\documentclass[11pt,reqno]{amsart}
\numberwithin{equation}{section}
\newtheorem{theorem}{Theorem}[section]
\newtheorem{lemma}[theorem]{Lemma}
\newtheorem{proposition}[theorem]{Proposition}
\newtheorem{corollary}[theorem]{Corollary}

\newtheorem{conjecture}[theorem]{Conjecture}
\theoremstyle{definition}
\newtheorem{example}[theorem]{Example}
\newtheorem{num_definition}[theorem]{Definition}
\newtheorem*{definition}{Definition}
\newtheorem{remark}[theorem]{Remark}

\newtheorem*{ack}{Acknowledgement}

\newenvironment{romenumerate}{\begin{enumerate}
 }{\end{enumerate}}

\newenvironment{thmxenumerate}{\begin{enumerate}
\setlength{\leftmargin}{0pt}
\setlength{\itemindent}{0pt}
 }{\end{enumerate}}

\newcounter{thmenumerate}
\newenvironment{thmenumerate}
{\setcounter{thmenumerate}{0}%
 \def\item{\par
 \refstepcounter{thmenumerate}\textup{(\roman{thmenumerate})\enspace}}
}
{}

\newcommand\pfitemx[1]{\par#1}
\newcommand\pfitem[1]{\pfitemx{(#1)}}
\newcommand{\refT}[1]{Theorem~\ref{#1}}
\newcommand{\refC}[1]{Corollary~\ref{#1}}
\newcommand{\refL}[1]{Lemma~\ref{#1}}
\newcommand{\refR}[1]{Remark~\ref{#1}}
\newcommand{\refS}[1]{Section~\ref{#1}}
\newcommand{\refSS}[1]{Subsection~\ref{#1}}
\newcommand{\refD}[1]{Definition~\ref{#1}}
\newcommand{\refSA}[1]{Appendix~\ref{#1}}
\newcommand{\refP}[1]{Proposition~\ref{#1}}
\newcommand{\refE}[1]{Example~\ref{#1}}
\newcommand{\refand}[2]{\ref{#1} and~\ref{#2}}
\newcommand\nopf{\qed}   

\newcommand\webcite[1]{\hfil\penalty0\texttt{\def~{\~{}}#1}~.\hfill\hfill}
\newcommand\arxiv[1]{\webcite{http://arXiv.org/#1}}

\newcommand\set[1]{\ensuremath{\{#1\}}}
\newcommand\bigset[1]{\ensuremath{\bigl\{#1\bigr\}}}
\newcommand\Bigset[1]{\ensuremath{\Bigl\{#1\Bigr\}}}
\newcommand\xpar[1]{(#1)}
\newcommand\bigpar[1]{\bigl(#1\bigr)}
\newcommand\Bigpar[1]{\Bigl(#1\Bigr)}
\newcommand\biggpar[1]{\biggl(#1\biggr)}

\newcommand\xfrac[2]{#1/#2}

\newcommand\parfrac[2]{\Bigpar{\frac{#1}{#2}}}
\newcommand\norm[1]{\ensuremath{\|#1\|}}
\newcommand\normpx[2]{\ensuremath{\|#1\|_{#2}}}
\newcommand\normp[1]{\normpx{#1}{p}}
\newcommand\normll[1]{\normpx{#1}{2}}
\newcommand\normoo[1]{\normpx{#1}{\infty}}
\newcommand\normHS[1]{\ensuremath{\norm{#1}_{HS}}}

\newcommand\ceil[1]{\lceil#1\rceil}
\newcommand\floor[1]{\lfloor#1\rfloor}
\newcommand\ntoo{\ensuremath{{n\to\infty}}}
\newcommand\mtoo{\ensuremath{{m\to\infty}}}
\newcommand\ktoo{\ensuremath{{k\to\infty}}}
\newcommand\ttoo{\ensuremath{{t\to\infty}}}
\newcommand\limn{\lim_{n\to\infty}}
\newcommand\limm{\lim_{m\to\infty}}
\newcommand\setdiff{\bigtriangleup}
\newcommand\iid{i.i.d.\spacefactor=1000}     
\newcommand\ie{i.e.\spacefactor=1000}
\newcommand\eg{e.g.\spacefactor=1000}
\newcommand\viz{viz.\spacefactor=1000}
\newcommand\cf{cf.\spacefactor=1000}
\newcommand{\as}{a.s.\spacefactor=1000}
\newcommand{\aex}{a.e.\spacefactor=1000}
\newcommand\whp{whp}
\newcommand\dto{\overset{\mathrm{d}}{\to}}
\newcommand\pto{\overset{\mathrm{p}}{\to}}
\newcommand\asto{\overset{\mathrm{a.s.}}{\to}}

\renewcommand\={:=}
\newcommand\bbR{\mathbb R}
\newcommand\bbC{\mathbb C}

\newcommand\bbZ{\mathbb Z}
\newcommand\bbT{\mathbb T}

\newcommand\E{\operatorname{\mathbb E{}}}
\renewcommand\P{\operatorname{\mathbb P{}}}
\newcommand\Var{\operatorname{Var}}

\newcommand\Exp{\operatorname{Exp}}
\newcommand\Po{\operatorname{Po}}
\newcommand\Bi{\operatorname{Bi}}
\newcommand\Be{\operatorname{Be}}

\newcommand\diam{\operatorname{diam}}
\newcommand\sign{\operatorname{sign}}
\newcommand\dd{\,d}
\newcommand\ddd{\partial}

\newcommand\intooo{\int_{-\infty}^\infty}

\newcommand\intoi{\int_0^1}
\newcommand\ints{\int_\sss}
\newcommand\sumin{\sum_{i=1}^n}
\newcommand\sumir{\sum_{i=1}^r}
\newcommand\ga{\alpha}
\newcommand\gb{\beta}
\newcommand\gd{\delta}

\newcommand\gam{\gamma}

\newcommand\gk{\kappa}
\newcommand\kk{\kappa}
\newcommand\kkb{\bar\kappa}
\newcommand\gl{\lambda}
\newcommand\gL{\Lambda}
\newcommand\go{\omega}

\newcommand\gs{\sigma}

\newcommand\eps{\varepsilon}

\newcommand\cA{\mathcal A}
\newcommand\cC{\mathcal C}

\newcommand\cE{\mathcal E}

\newcommand\cG{\mathcal G}
\newcommand\cL{{\mathcal L}}
\newcommand\cM{\mathcal M}

\newcommand\cP{\mathcal P}
\newcommand\cS{{\mathcal S}}

\newcommand\hmu{\widehat \mu}
\newcommand\hmup{\widehat \mu'}
\newcommand\hkkp{\widehat \kk'}
\newcommand\hvxs{\widehat \vxs}

\newcommand\tkkmux{{\tilde\kk_\mux}}
\newcommand\tkko{\tilde\kk_1}
\newcommand\tkkn{\tilde\kk_n}
\newcommand\kkn{\kk_n}
\newcommand\kkm[1]{\kk_m^{#1}}
\newcommand\kkmx[1]{\kk_{#1}^-}
\newcommand\tkkm{\hat\kk_m^{-}}
\newcommand\kkt{\kk_t}
\newcommand\kktx{\kk_t^{*}}
\newcommand\kkxM{{\kk}\bmin M} 
\newcommand\Pm{\cP_m}
\newcommand\downto{\searrow}
\newcommand\upto{\nearrow}
\newcommand\ett{\mathbf 1}
\newcommand\IN[1]{\ett[#1]}
\newcommand\bmin{\wedge}
\newcommand\bmax{\vee}
\newcommand\innprod[1]{\langle #1\rangle}
\newcommand\qi{^{-1}}
\newcommand\qh{^{1/2}}
\newcommand\qhi{^{-1/2}}
\newcommand\Gcjn{\ensuremath{{G^{1/j}_n(c)}}}
\newcommand\sss{\cS}
\newcommand\bfx{\mathbf{x}}
\newcommand\bfy{\mathbf{y}}
\newcommand\bfd{\mathbf{d}}
\newcommand\bfw{\mathbf{w}}
\newcommand\xs{\ensuremath{{\bfx}_n}}
\newcommand\ys{\ensuremath{{\bfy}_n}}
\newcommand\xss{\ensuremath{(\xs)_{n\ge 1}}}
\newcommand\nxss{\ensuremath{(\xs)}}
\newcommand\xssd{\ensuremath{(\xs')_{n\ge 1}}}
\newcommand\yss{\ensuremath{(\ys)_{n\ge 1}}}
\newcommand\vxs{\ensuremath{\mathcal V}}
\newcommand\pij{p_{ij}}
\newcommand\ppij{(\pij)}
\newcommand\gtx[3]{\ensuremath{G^{#1}(#2,#3)}}
\newcommand\gxx[2]{\ensuremath{G(#1,#2)}}

\newcommand\gnx[1]{\ensuremath{G(n,#1)}}
\newcommand\gnxx[1]{\ensuremath{G^\vxs(n,#1)}}
\newcommand\gnxxp[1]{\ensuremath{G^{\vxs'}(n,#1)}}
\newcommand\gnxxh[1]{\ensuremath{G^{\hvxs}(n,#1)}}

\newcommand\gnp{\gnx{p}}
\newcommand\gnk{\gnx{\kk}}
\newcommand\gnkx{\gnxx{\kk}}

\newcommand\gnkn{\gnx{\kk_n}}
\newcommand\gnkxn{\gnxx{\kk_n}}
\newcommand\gnkp{G^{\langle p\rangle}(n,\kk)}
\newcommand\gnkpp{G^{[p]}(n,\kk)}
\newcommand\gnc{\gnx{c/n}}
\newcommand\bpk{\ensuremath{\mathfrak X_{\kk}}}

\newcommand\bpkx{\ensuremath{\mathfrak X_{\kk}(x)}}
\newcommand\bpxx[1]{\ensuremath{\mathfrak X_{#1}(x)}}

\newcommand\bpkxp{\ensuremath{\widehat{\mathfrak X}_{\kk}(x)}}
\newcommand\bpkmod{\ensuremath{{\overline{\mathfrak X}}_{\kk}}}

\newcommand\Ga{\Gamma}
\newcommand\bpx[1]{\ensuremath{\mathfrak X_{#1}}}
\newcommand\rhox[1]{\ensuremath{\rho_{#1}}}
\newcommand\rhok{\rhox{k}}
\newcommand\rhogek{\rhox{\ge k}}
\newcommand\sgek{\sigma_{\ge k}}
\newcommand\extrho{\rho^+}
\newcommand\rhoe{\rhox{1+\eps}}
\newcommand\rhoo{\rhox{1}}
\newcommand\rhoex{\rhox{1+\eps}^*}
\newcommand\xphi{\phi}  
\newcommand\nv[1]{v_{#1}}
\newcommand\nvn{\nv{n}}
\newcommand\nvl{\nv{\gl}}
\newcommand\nx[1]{\ensuremath{N_{#1}}}
\newcommand\nk{\nx{k}}
\newcommand\ngek{\nx{\ge k}}
\newcommand\ngex[1]{\nx{\ge #1}}
\newcommand\nuni{\nu_n^{1}}
\newcommand\nunmi{\nu_{nm}^{1}}
\newcommand\Gtc{{G^2}}
\newcommand\Gtcm{{G^{2-}}}
\newcommand\Wc{{\overline W}}

\newcommand\EC{{EC}}
\newcommand\psimple{p_{\mathrm{simple}}}
\newcommand\pcut{p_{\mathrm cut}}
\newcommand\tk{T_{\kk}}

\newcommand\tkn{T_{\kk_n}}
\newcommand\tkx[1]{T_{\kk_{#1}}}
\newcommand\tkq{T_{\kk}'}

\newcommand\htk{\widehat T_{\kk}}
\newcommand\htkp{T_{\hkkp}}
\newcommand\Tx[1]{T_{#1}}
\newcommand\ttk{\widetilde T_{\kk}}
\newcommand\Phik{\Phi_{\kk}}
\newcommand\Phikn{\Phi_{\kk_n}}

\newcommand\Phix[1]{\Phi_{#1}}
\newcommand\trhokk{\tilde\rho_{\kk}}
\newcommand\rhokk{\rho_{\kk}}
\newcommand\rhokkx{\rho(\kk;x)}
\newcommand\rhokn{\rho_{\kk_n}}
\newcommand\rhokm{\rho_{\kk_m}}
\newcommand\rhokx[1]{\rho_{\kk_{#1}}}

\newcommand\rhock{\rhox{c\kk}}
\newcommand\muk{\mu_{\kk}}
\newcommand\tmu{\tilde\mu}
\newcommand\iik{\iint\kk}
\newcommand\hiik{\tfrac12\iik}
\newcommand\iikxy{\iint_{\sss^2}\kk(x,y)\dd\mu(x)\dd\mu(y)}
\newcommand\hiikxy{\tfrac12\iikxy}

\newcommand\rhokq{\rhokk^*}
\newcommand{\mucs}{$\mu$-continuity set}
\newcommand{\rfin}{regular finitary}

\newcommand\CA{\sss\setminus A}

\newcommand\infkappa{\inf\kappa}
\newcommand\bnk{\bar n_k}
\newcommand\step[1]{\medskip\par\emph{#1.}}
\newcommand\mux{\gd}
\newcommand\tgl{\widetilde G_\gl}
\newcommand\tgn{\widetilde G_n}

\newcommand\wnm{W_n^{(m)}}
\newcommand\wm{w_m}
\newcommand\qir{quasi-irreducible}
\newcommand\qirity{quasi-irreducibility}
\newcommand\glcrit{\gl_{\textrm{cr}}}
\newcommand\CS{Cauchy--Schwarz}
\newcommand\CSineq{\CS{} inequality}
\newcommand{\Holder}{H\"older}
\newcommand\iikk{\iint \kk^2}
\newcommand\hath{\widehat h}
\newcommand\Tr{\operatorname{Tr}}
\newcommand\qm{^{(m)}}
\newcommand\Xnf{neighbourhood function}
\newcommand\Lnf{$L$-\Xnf}
\newcommand\Onf{$1$-\Xnf}
\newcommand\Dnf{$D$-\Xnf}

\newcommand\edgeno{\zeta}
\newcommand\edgenokk{\edgeno(\kk)}
\newcommand\mui{\mu(\{i\})}
\newcommand\muj{\mu(\{j\})}
\newcommand\Gale[1]{\Ga_{\le #1}}
\newcommand\supxyn{\sup_{x,y,n}\gk_n(x,y)}
\newcommand\vnu{\boldsymbol{\nu}}

\newcommand\vNt{\mathbf{N}_t}

\newcommand\lntk{\log\norm{\tk}}
\newcommand\lntki{\log\norm{\tk}^{-1}}
\newcommand\lntkdi{\log\norm{T_{\kkd}}^{-1}}
\newcommand\ntkd{\norm{T_{\kkd}}}
\newcommand\kkd{{\hat\kk}}
\newcommand\mud{{\hat\mu}}
\newcommand\mudk{{\mud\{k\}}}
\newcommand\bpkd{\bpx{\kkd}}

\begin{document}

\title
{The phase transition in inhomogeneous random graphs}
\date{April 29, 2005; revised June 20, 2006}

\begingroup
\makeatletter \gdef\FNX#1{\@fnsymbol{#1}}
\endgroup
\newcommand\FN[1]{\textsuperscript{\FNX{#1}}}
\author{B\'ela Bollob\'as\FN1\FN4\FN6}
\author{Svante Janson\FN2\FN5}
\author{Oliver Riordan\FN3\FN4}
\thanks{\FN1Department of Mathematical Sciences, University of Memphis,
  Memphis TN 38152, USA}
\thanks{\FN2Department of Mathematics, Uppsala University,
 PO Box 480, SE-751 06 Uppsala, Sweden}
\thanks{\FN3Royal Society Research Fellow, Department of Pure Mathematics
and Mathematical Statistics, University of Cambridge, UK}
\thanks{\FN4Trinity College, Cambridge CB2 1TQ, UK}
\thanks{\FN5Churchill College, Cambridge CB3 0DS, UK}
\thanks{\FN6Research supported by NSF grants CCR 0225610 and DMS 0505550}


\subjclass[2000]{60C05; 05C80}

\begin{abstract}
\leftskip=-6pt \rightskip=\leftskip
The `classical' random graph models, in particular $G(n,p)$, are
`homogeneous', in the sense that the degrees (for example) tend to be
concentrated around a typical value.  Many graphs arising in
the real world do not have this property, having, for example,
power-law degree distributions. Thus there has been a lot of recent
interest in defining and studying `inhomogeneous' random graph models.

One of the most studied properties of these new models is their
`robustness', or, equivalently, the `phase transition' as an edge
density parameter is varied. For $G(n,p)$, $p=c/n$, the phase transition
at $c=1$ has been a central topic in the study of random graphs
for well over 40 years.

Many of the new inhomogeneous models are rather complicated; although
there are exceptions, in most cases precise questions such as
determining exactly the critical point of the phase transition are
approachable only when there is independence between the
edges. Fortunately, some models studied have this property already, and others
can be approximated by models with independence.

Here we introduce a very general model of an inhomogeneous random graph
with (conditional)
independence between the edges, which scales so that the number
of edges is linear in the number of vertices.  This scaling
corresponds to the $p=c/n$ scaling for $G(n,p)$ used to study the
phase transition; also, it seems to be a property of many large
real-world graphs. Our model includes as special cases many models
previously studied.

We show that, under one very weak assumption (that the expected
number of edges is `what it should be'), many properties
of the model can be determined, in particular the critical
point of the phase transition, and the size of the giant component
above the transition. We do this by relating our random graphs
to branching processes, which are much easier to analyze.

We also consider other properties of the model, showing,
for example, that when there is a giant component, it is
`stable': for a typical random graph, no matter
how we add or delete $o(n)$ edges, the size of
the giant component does not change
by more than $o(n)$.
\end{abstract}

\maketitle
\renewcommand\FN[1]{\relax}

\tableofcontents

\section{Introduction}\label{S1}

The theory of random graphs was founded in the late 1950s and early
1960s by Erd\H os and R\'enyi \cite{ER60,ER61},
who started the systematic study of the
space $\cG(n,M)$ of all graphs with $n$ labelled vertices and $M=M(n)$
edges, with all graphs equiprobable. (Usually, one writes
$G(n,M)$ for a random element of $\cG(n,M)$.)
At about the same time, Gilbert \cite{G59} introduced the closely
related model $\cG(n,p)$ of random graphs on $n$ labelled vertices:
a random $G(n,p)\in \cG(n,p)$ is obtained by selecting edges independently,
each with probability $p=p(n)$.
For many questions,
such as those
considered in this paper, the models are essentially equivalent
(if $p=M/\binom{n}{2}$, say). As
Erd\H os and R\'enyi are the founders of the theory
of random graphs, it is not surprising that both \gnp\ and
$G(n,M)$ are now known as Erd\H os--R\'enyi random graphs.

In addition to these two `classical' models, much attention
has been paid to the space of random $r$-regular graphs,
and to the space $\cG(k-\mathrm{out})$
of random directed graphs where each vertex has out-degree $k$,
and the undirected graphs underlying these.
All these random graph models are
`homogeneous' in the sense that all vertices are exactly
equivalent in the definition of the model. Furthermore,
in a typical realization, most vertices are in some
sense similar to most others. For example, the vertex degrees 
in $G(n,p)$ or in $G(n,M)$ do not vary very much: their distribution
is close to a Poisson distribution.

In contrast, many
large real-world graphs are highly inhomogeneous. One
reason is that the
vertices may have been `born' at different times, with
old and new vertices having very different properties.
Experimentally, the spread of degrees is often very large.
In particular, in many examples the degree distribution follows
a power law. In the last few years, this has led to the introduction
and analysis of many new random graph models designed to incorporate
or explain these features.
Recent work in this area perhaps
started from the observations of Faloutsos, Faloutsos and Faloutsos~\cite{FFF}
concerning scaling in real-life networks, in particular the power-law
distribution of degrees in the `internet graph',
and similar observations concerning the `web graph' made by
Kleinberg, Kumar, Raghavan, Rajagopalan and Tomkins~\cite{KKRRT}, 
and by Barab\'asi and Albert~\cite{BAsc}, who also looked at several other
real world graphs.
The latter two groups introduced two of the first models to explain
these observations, using the ideas of `copying' and of
`growth with preferential attachment', respectively.
Observations
of and proposed models for other networks followed, including
protein interaction networks, telephone call graphs,
scientific collaboration graphs and many others. Extensive
surveys of the mostly experimental or heuristic work
have been written by Barab\'asi and Albert~\cite{BAsurv} and
Dorogovtsev and Mendes~\cite{DMsurv}.

Some of the first rigorous mathematical results 
concerning (precisely defined variants of) these new models are
those of Bollob\'as and Riordan~\cite{diam},
Bollob\'as, Riordan, Spencer and Tusn\'ady~\cite{BRST},
Buckley and Osthus~\cite{BO}
and Cooper and Frieze~\cite{CF}.
For a partial survey of the rapidly growing 
body of rigorous work see Bollob\'as and
Riordan~\cite{BRsurv}. 
Needless to say, surveys in an active area such as this quickly become dated, 
and there are already too many rigorous results in the field to list here.

Perhaps the most striking and important result of Erd\H os and R\'enyi
concerns the sudden emergence of the `giant component', the phase transition
in $G(n,p)$ that occurs at $p=1/n$:
if $c>0$ is a constant, then the largest component of \gnx{c/n}
has order $O(\log n)$ \whp{} if $c<1$, and order $\Theta(n)$ \whp{} if $c>1$ 
(see \refS{Sresults} for the notation).
In particular, a
giant component of order $\Theta(n)$ exists (with high probability
as \ntoo) if and only if $c>1$.
Over twenty years later, Bollob\'as~\cite{BBtr} and \L uczak \cite{Lucz1} 
proved considerably sharper results about the exact nature of this phase
transition: in particular, they determined the exact size of the `window'
in which the transition takes place. Further, very detailed results
were proved by Janson, Knuth, \L uczak and Pittel~\cite{JKLP}; see
Bollob\'as~\cite{BB} and Janson,
{\L}uczak and Ruci\'nski~\cite{JLR} for numerous
results and references.

\medskip
Our main purpose in this paper is to lay the foundations
of a very general theory of inhomogeneous sparse random graphs. To this end
\begin{itemize}
\item
we shall define a general model that is sufficiently flexible to include
exactly many of the specific spaces of inhomogeneous random graphs
that have been studied in recent years,
\item
we shall establish a close connection between the component structure
of a random graph in this model, the survival probability of a related
branching process, and the norm of a certain operator,
\item
we shall use these connections to study the phase transition in our model,
examining especially the numbers of vertices and edges in the giant component,
and
\item
we shall prove results concerning the stability of the giant component
under the addition and deletion of edges.
\end{itemize}

In addition, we shall study various other properties
of our model, including the degree distribution, the numbers of paths
and cycles, and the typical distance between pairs of vertices
in the giant component.
Furthermore, we shall spell out what our results say about many specific
models that have been studied previously.

\medskip
Although we shall give many examples throughout the paper, to
motivate the definitions it may help to bear in mind one
particular example of the general class of models we shall study.
This example is the \emph{uniformly grown random graph}, or
\emph{$c/j$-graph}, $\Gcjn$. Here $c>0$ is a parameter that will
be kept constant as $n$ varies, and the graph $\Gcjn$ is the graph
on $[n]=\{1,2,\ldots,n\}$ in which edges are present
independently, and the probability that for $i\ne j$ the edge $ij$
is present is $\pij=\min\{c/\max\{i,j\},1\}$, or simply $c/j$ if
$i<j$ and $c \le 2$.

A natural generalization of \gnp\ that includes the example
above as a very special case is obtained by replacing the single parameter
$p$ by a symmetric $n\times n$ matrix $\ppij$ with $0\le\pij\le1$. We
write \gnx{\ppij} for the random graph with vertex set $[n]$ where
$i$ and $j$ are connected by an edge with probability $\pij$, and
these events are independent for all pairs $(i,j)$ with $i<j$;
see \cite[p.\ 35]{BB}.
We are interested in asymptotics as \ntoo, usually with
$\ppij=\bigpar{\pij(n)}$ depending on $n$.
It seems difficult to obtain substantial asymptotic results for \gnx{\ppij}
without further restrictions; the model is too general.
(However, for connectedness Alon \cite{Alon} proved a number of results.)

Here we are
mainly interested in random graphs where the average degree is
$\Theta(1)$; one of the main cases treated in this paper is
$\pij=\gk(i/n,j/n)/n$ for a suitable function $\gk$ on $(0,1]^2$.
Taking $\gk(x,y)=c/\max\{x,y\}$, we obtain $\Gcjn$.
Many other graphs studied earlier by different authors can also
be obtained by choosing $\gk(x,y)$ suitably; see \refS{Sapp},
and the forthcoming papers \cite{kerapp,quantum}. 
A precise definition of the random graphs treated here will be given in
\refS{Smodel}, and some simple examples in \refS{Sex}.

\medskip
The rest of the paper is organized as follows.
In \refS{Smodel} we define the model $\gnk$ we shall study,
along with the branching process $\bpk$ and integral operator $\tk$
to which we shall relate its component structure.

In \refSS{SSgc} we present our main results on the giant component
of $\gnk$: under certain weak
assumptions we obtain necessary and sufficient conditions for $\gnk$
to have a giant 
component, show that when the giant component exists it is unique,
and find its normalized size and number of edges.
Further results are presented in the following subsections, on
the `stability' of the giant component in \refSS{SSstable}, on
small components in \refSS{SSsmall}, on the degree sequence in
\refSS{SSdegrees}, 
and on the typical
distance between vertices of the giant component in \refSS{SSdist}.
In \refSS{SStrans} we turn to the phase-transition in \gnk; more precisely,
we examine the
growth rate of the giant component as it emerges.

Since our model is very general, and the definition rather lengthy, special
cases of the model play an important role in the paper; we have described one,
$\Gcjn$ already. In \refS{Sex} we give several further simple
examples, to illustrate  
the definitions and results of the previous sections. Towards the end
of the paper, 
we shall discuss several other special cases more extensively, in
particular describing 
the relationship to other models studied earlier; we consider these to
be applications 
rather than illustrations of the model, and so present them after the proofs.

The next several sections are devoted to the proofs of the main
results; the reader 
interested primarily in the applications may wish to skip straight to
\refS{Sapp}. 
We start by analyzing the branching process \bpk, proving results
about this process that will help us relate \gnk{} to \bpk.
The study of \bpk{} itself is not one of our main aims.
In \refS{S:branch1} we prove various lemmas needed in \refS{S:branch2} to
prove the results about \bpk{} that we shall use throughout the paper.

Next, we turn to preparatory results concerning \gnk{} itself, 
starting with simple approximation lemmas in \refS{Sapprox}; basic
results on the number of edges of \gnk\ are given in \refS{Sedges}.

Our main results about the existence and size of the giant component
are proved in
\refS{Smainpf}, using material from the previous sections;
the reader who is interested only in the
derivations of these results should read Sections \ref{Smodel} to
\ref{Smainpf}.
The number of edges in the giant component is determined in \refS{Sgiantedges}.

Sections \ref{Sstable} to \ref{Stransitionpf} are devoted to the proofs of the results
in Subsections~\ref{SSfirst_after_gc} to \ref{SSlast}:
broadly speaking, these proofs rely on the results up to
\refS{Smainpf}, but not
on each other, so the reader may safely omit any subset of these
sections.
The stability result is proved in \refS{Sstable},
the results on small components in \refS{S2ndpf},
the vertex degrees are studied in \refS{Sdegrees},
the distance
between vertices in \refS{Sdist}, and
the phase transition in \refS{Stransitionpf}. The latter results
may be viewed purely as statements about a branching process,
in which case the proofs need only the results of
Sections~\refand{S:branch1}{S:branch2}. 

In \refS{Sapp} we apply our general model to deduce results about
several specific models, in particular ones that have
been studied in recent years, and discuss
the relationship of our results to earlier work.

In \refS{Scycles} we give some simple results about paths
and cycles in \gnk, inspired by the work of Turova on a special case 
of the model (described in \refS{Sapp}), and show that a conjecture of
hers holds under mild conditions. 
In \refS{Sconcl} we discuss several
(at least superficially) related models as well as possible
future work. Finally, in the appendix we give some basic results
on random measures used throughout the paper.

\section{The model}\label{Smodel}

In this section we define the random graph model that we shall 
study throughout the paper, as well as a branching process and an
integral operator that will be key to
characterizing the component
structure of this random graph. This section also includes various
remarks on the definitions, including descriptions of several minor
variants; the 
formal definitions may be understood without reference to these
remarks. To make sense of the definitions, the reader may wish to keep
in mind the model $\Gcjn$ defined in the introduction.

Our model is an extension of one defined by S\"oderberg \cite{Sod1}.
Let $\sss$ be a separable metric space 
equipped with a Borel probability measure $\mu$. We shall often suppress
the measure $\mu$ in our notation, writing, for example, $\normp{\;}$
for the norm in $L^p(\sss)=L^p(\sss,\mu)$, and `a.e.\ on $\sss$' for
$\mu$-a.e.\ on $\sss$.
Much of the time (for example, when studying $\Gcjn$), we shall take
$\sss=(0,1]$ with $\mu$ Lebesgue measure.
Throughout the paper, the `kernel' $\kappa$ will be a symmetric
non-negative function on 
$\sss\times\sss$.
Further conditions on $\gk$ will be given in Definitions
\refand{Dg1}{Dg0}.

\medskip
For each $n$ we have a deterministic or random 
sequence $\xs=(x_1,\dots,x_n)$ of points in $\sss$.
Formally, we should write $\xs=(x_1^{(n)},\ldots,x_n^{(n)})$, say, as we
assume no relationship between the $i$th elements of $\xs$ and of ${\bfx}_{n'}$.
However, this notation would be rather cumbersome, and it will always
be clear which $\xs$ an $x_i$ is an element of.
Writing $\gd_x$ for the measure consisting of a point mass 
of weight $1$ at $x$, and
\begin{equation}
  \label{a2}
\nu_n\=\frac1n\sumin \gd_{x_i}
\end{equation}
for the empirical distribution of $\xs$, we shall assume that
$\nu_n$ converges in probability to $\mu$ as $n\to\infty$, with
convergence in the usual space of probability measures on $\sss$
(see, e.g., \cite{Bill}). This condition has a simple down-to-earth description 
in terms of the number of $x_i$ in certain sets $A$:
a set $A\subseteq\sss$ is a \emph{\mucs}
if $A$ is (Borel) measurable and $\mu(\ddd A)=0$, where
$\ddd A$ is the boundary of $A$. The convergence condition
$\nu_n\pto\mu$ means exactly that for every \mucs\ $A$,
\begin{equation}
  \label{a2a}
\nu_n(A)\=\#\set{i:x_i\in A}/n \pto \mu(A);
\end{equation}
see \refSA{Smeasure} for technical details.

One example where \eqref{a2a} holds is 
the random case, where the $x_i$ are independent and uniformly
distributed on $\sss$ with distribution $\mu$ (as in S\"oderberg~\cite{Sod1});
then \eqref{a2a} holds by the law of large numbers. 

We shall often
consider $\sss=(0,1]$ with the Lebesgue measure $\mu$; in this
case, condition \eqref{a2a} has to be verified only for intervals
(see \refR{RR}). For this pair  $(\sss, \mu)$ we shall have two
standard choices for the $(x_i)$: the deterministic case
$x_i=i/n$, and the random case where the $x_i$ are independent and
uniformly distributed on $(0,1]$. To express $\Gcjn$ as a special case 
of our model, we shall take $x_i=i/n$.

For later formal statements, we gather the preceding assumptions
into the following definitions.
\begin{definition}
A \emph{ground space} is a pair $(\sss,\mu)$, where
$\sss$ is a separable metric space and $\mu$ is a Borel probability
measure on $\sss$.
\end{definition}
\begin{definition}
A \emph{vertex space} $\vxs$ is a triple $(\sss,\mu,\xss)$, where
$(\sss,\mu)$ is a ground space and, for each $n\ge 1$, $\xs$
is a random sequence $(x_1,x_2,\ldots,x_n)$ of $n$ points of $\sss$,
such that \eqref{a2a} holds.
\end{definition}
Of course, we do not need $\xss$ to be defined for every $n$, but only
for an infinite set of integers $n$.
\begin{definition}
A \emph{kernel} $\kappa$ on a ground space $(\sss,\mu)$
is a symmetric non-negative
(Borel) measurable
function on
$\sss\times \sss$.
By a kernel on a vertex space $(\sss,\mu,\xss)$ we mean a kernel on
$(\sss,\mu)$.
\end{definition}

{}From now on, unless otherwise stated, we shall always write a vertex
space $\vxs$ as $(\sss,\mu,\xss)$, and $\xs$ as $(x_1,x_2,\ldots,x_n)$.
As noted above, the (distributions of) the individual $x_i$ depend on $n$;
in the notation we suppress this dependence as it will always be
clear which $\xs$ an $x_i$ is a member of.

Let $\kappa$ be a kernel on the vertex space $\vxs$. Given
the (random) sequence $(x_1,\dots,x_n)$, we let $\gnkx$ be the random
graph \gnxx{\ppij}
with
\begin{equation}
  \label{pij}
\pij\=\min\bigset{\gk(x_i,x_j)/n,1}.
\end{equation}
In other words, \gnkx{} has $n$ vertices \set{1,\dots,n} and, given
$x_1,\dots,x_n$, an edge $ij$ (with $i\neq j$) exists
with probability $\pij$, independently of all other
(unordered) pairs $ij$. Often, we shall suppress the dependence on
$\vxs$, writing $\gnk$
for $\gnkx$.
We have described one example already: if we take $\kk(x,y)=c/\max\{x,y\}$, 
with $\sss=(0,1]$, $\mu$ Lebesgue measure, and $x_i=x_i^{(n)}=i/n$,
then \eqref{pij} gives $\pij=\min\{c/\max\{i,j\},1\}$, so
$\gnkx$ is exactly the uniformly grown random graph $\Gcjn$ described
in the introduction. We shall discuss several other examples in
Sections \refand{Sex}{Sapp}.

\begin{remark}
The random graph \gnk{}=\gnkx{} depends not only on $\kk$ but also 
on the choice of $x_1,\dots,x_n$. Much of the time, our notation
will not indicate how the points $x_i$ are chosen, since this
choice is irrelevant for our results as long as \eqref{a2a} holds
and certain pathologies are excluded (see \eqref{t1b}, \refL{LE}
and \refE{Ebad}). The freedom of choice of $x_1,\dots,x_n$ gives
our model flexibility,
as shown by \refP{Pfinite}, \refT{Tdual}
and the examples in Sections
\refand{Sex}{Sapp}, but does not affect the asymptotic
behaviour. Of course, this asymptotic behaviour does depend very
much on $\sss$ and $\mu$.
\end{remark}

In order to make our results easy to apply, it will be convenient to
extend the definitions above in two ways, by allowing $\mu(\sss)$ to
take any value in $(0,\infty)$, and by allowing the number of vertices
of $\gnkx$ to be random, rather than exactly $n$.  As we shall see
later, this makes no essential difference, and we shall almost always
work with the $n$ vertex model in our arguments.  We shall consider the
`generalized' model only for the convenience of a reader wishing to
apply the results in the next section, obviating the need for a
separate reduction to the $n$ vertex model in each case. All other
readers may safely ignore the `generalized' model, including the
formal definitions that we now state.

\begin{definition}
A \emph{generalized ground space} is a pair $(\sss,\mu)$, where
$\sss$ is a separable metric space and $\mu$ is a Borel measure on
$\sss$ with $0<\mu(\sss)<\infty$. 
\end{definition}
Let $I\subset(0,\infty)$ be any unbounded set, the \emph{index set}
parametrizing our model. 
Usually, $I$ is the positive integers, or the positive reals. For
compatibility with 
our earlier definitions, we write $n$ for an element of $I$, even
though this need not be an integer. 
\begin{definition}
A \emph{generalized vertex space} $\vxs$ is a triple
$(\sss,\mu,(\xs)_{n\in I})$, where 
$(\sss,\mu)$ is a generalized ground space and, for each $n\in I$, $\xs$
is a random sequence $(x_1,x_2,\ldots,x_{\nvn})$ of points of $\sss$
of random length $\nvn\ge 0$,
such that \eqref{a2a} holds, i.e., such that 
\begin{equation}\label{nunA2}
 \nu_n(A)\=\#\set{i:x_i\in A}/n \pto \mu(A)
\end{equation}
as $n\in I$ tends to infinity
for every \mucs\ $A$;
equivalently, $\nu_n\pto\mu$.
\end{definition}

The definition of a \emph{kernel} $\kk$ on a generalized ground space
is exactly as before. 
Finally, given a kernel $\kk$ on a generalized vertex space, for $n\in
I$ we let 
$\gnkx$ be the random graph on $\{1,2,\ldots,\nvn\}$ in which, given
$\xs=(x_1,\ldots,x_{\nvn})$, 
each possible edge $ij$, $1\le i<j\le \nvn$, is present with probability
\begin{equation}\label{pij2}
 \pij\=\min\bigset{\gk(x_i,x_j)/n,1},
\end{equation}
and the events that different edges are present are independent.

Note that if $\vxs$ is a generalized ground space, then, applying
\eqref{nunA2} with $A=\sss$, we see that $\gnkx$ has $\mu(\sss)n+o_p(n)$ vertices.
In both \eqref{nunA2} and \eqref{pij2} we divide by $n$, rather than by the actual number
of vertices, or by $\mu(\sss) n$; this turns out to be most convenient normalization.
Roughly speaking, by conditioning on the sequences $(\xs)$, or by
adding $o_p(n)$ isolated vertices, we may assume without loss of generality
that the number of vertices is deterministic. Furthermore, multiplying $\kk$ and the index
variable $n$ by some constant factor, and dividing $\mu(\sss)$ by the same factor,
leaves the edge probabilities $\pij$ unchanged. As the condition \eqref{nunA2}
is also unaffected by this transformation, the only effect
on the model is to rescale the parameter $n$, and we may assume
without loss of generality 
that we have a vertex space rather than a generalized vertex space;
see \refSS{SSgvs}.

\begin{remark}\label{Rrandom}
We regard our random graphs as indexed by $n$, and consider what happens
as \ntoo. This is for notational convenience only; we
could consider graphs indexed by some other (possibly continuous)
parameter, $m$, say, such that the number of vertices $\nv{m}$ of the
graph with parameter $m$ tends to infinity. This superficial modification
is covered by the definitions above: instead of considering
graphs on $2^n$ vertices, say, one can always consider graphs on $n$
vertices with $n$ restricted to an `index set' $I$ consisting of the
powers of $2$. 

The generalized vertex space setting also allows the number of vertices
to be random. 
In other words, we may let \set{x_i} be a point process on
$\sss$, for example, a Poisson process of intensity $n$;
see Examples \refand{EPo}{Evp} and \refSS{SSTurova}.
\end{remark}

\begin{remark}\label{Rmeas0}
Changing $\kk$ on a set of measure zero may have a significant
effect on the graph $\gnkx$; see \refE{Ebad}, for instance.
Indeed, if the $x_i$ are deterministic, then $\gnkx$ depends only
on the values of $\kk$ on a discrete set. This means that in our
proofs we cannot just ignore measure zero sets
in the usual way. Later we shall impose very
weak conditions to control such effects; see \refR{Rpath}.
\end{remark}

Before turning to the key definitions, giving conditions under which we
can prove substantial results about $\gnkx$, let us make some remarks about
some minor variants of the model.

\begin{remark}\label{Rexp}
As an alternative to \eqref{pij} (or \eqref{pij2}),
we could use $\gk$ to define the
intensities of Poisson processes of edges, and ignore multiple
edges, so the probability $\pij$ that there is an edge
between $i$ and $j$ would be given by
\begin{equation}  \label{pij'}
\pij\=1-\exp\bigpar{-\gk(x_i,x_j)/n},
\end{equation}
rather than by \eqref{pij}. The results below are valid for this
version too; this can be shown either by checking that all
arguments hold with only trivial changes, or by defining
$\gk_n(x,y)\=n\bigpar{1-\exp\xpar{-\gk(x,y)/n}}$ and using the
setting in \refD{Dg2}.

Another alternative,
studied by Britton, Deijfen and Martin-L\"of \cite{BrittonDML}
in a special case (see \refSS{SSrank1}),
is to let
$\pij/(1-\pij)=\gk(x_i,x_j)/n$, \ie{}, to take 
\begin{equation}\label{pij''}
\pij\=\gk(x_i,x_j)/(n+\gk(x_i,x_j)).
\end{equation}
Again, the results below remain valid; we now define
$\gk_n(x,y)\=\gk(x,y)/(1+\gk(x,y)/n)$.
\end{remark}
\begin{remark}
  \label{Rmulti}
In this paper we treat only simple graphs. One natural variation
that yields a multigraph is to let the number of edges between $i$
and $j$ have a Poisson distribution with mean $\kk(x_i,x_j)/n$.
Under suitable conditions (\eg{} that $\kk$ is bounded), it is
easy to see that \whp{} there are no triple edges, and that the
number of double edges is $O_p(1)$; more precisely, it has an
asymptotic Poisson distribution with mean $\frac14\iint\kk^2$, see
\refS{Scycles}. The underlying simple graph is just the graph
defined in \refR{Rexp}.

Another variation (which can be combined with the previous one) is
to permit loops by allowing $i=j$ in the definition above. These
variations do not affect our results on component sizes.
\end{remark}

\begin{remark}\label{Rproc} 
Our model can be extended to a random graph process on a fixed vertex set 
describing
an inhomogeneously growing random graph: Start without any edges and,
given $\xs$, add edges at random times given by independent Poisson
processes with
intensities $\kk(x_i,x_j)/n$. (Ignore multiple edges.)
At time $t$, we obtain the version of
the random graph $\gnx{t\kk}$ given by \eqref{pij'};
\cf{} \refR{Rexp}.
Alternatively, we may add edges sequentially, with each new edge
chosen at random with probabilities proportional to $\kk(x_i,x_j)$;
this gives the same process except for a (random) change of time
scale.
\end{remark}

Without further restrictions, the model $\gnkx$ we have defined
is too general for us to prove meaningful results. Indeed,
the entire graph may be determined by the behaviour
of $\kk$ on a measure zero set.
Usually, $\kk$ is continuous, 
so this problem does not arise. However, there are natural examples with $\kk$
discontinuous, so we shall assume that $\kk$ is continuous
a.e.\ rather than everywhere. 
With this weaker condition, to
relate the behaviour of $\gnkx$ to that of $\kk$ we
shall need some extra assumptions. The behaviour of the total number
of edges turns out to be the key to the elimination of pathologies.

As usual, we write $e(G)$ for the number of edges in a graph $G$. Note
that 
$\E e\bigpar{\gnx{\ppij}}=\sum_{i<j}\pij$, so we have
\begin{equation}\label{Ee}
  \E e\bigpar{\gnkx}=\E\sum_{i<j} \min\bigset{\gk(x_i,x_j)/n,1}.
\end{equation}
In well behaved cases, this expectation is asymptotically $n\hiik$;
see, for example, 
\refL{LE}.

\begin{num_definition}\label{Dg1}
A kernel $\kk$ is \emph{graphical}
on a (generalized) 
vertex space $\vxs=(\sss,\mu,\nxss)$
if the following conditions hold:
\begin{romenumerate}
\item \label{T1a}
  $\gk$ is continuous \aex{} on $\sss\times\sss$;
\item \label{T1c}
$  \kk\in L^1(\sss\times\sss,\mu\times\mu)$;
\item  \label{T1d}
\begin{equation}
\label{t1b}
\frac1n  \E e\bigpar{\gnkx} \to \hiikxy.
\end{equation}
\end{romenumerate}
\end{num_definition}

Note that whether $\kk$ is graphical on $\vxs$ depends on the
sequences $\xs$. Also, as we shall see in \refR{Rckgraphical}, 
if $\kk$ is graphical on $\vxs$, then so is $c\kk$ for any constant $c>0$.
(This statement would be trivial without the $\min\{\cdot,1\}$ operation
in the right-hand side of \eqref{Ee}. With this, it is still not hard
to check.) 

\begin{remark}\label{Rpath}
Conditions \ref{T1a} and \ref{T1c} are natural technical conditions; at first
sight, condition \ref{T1d} is perhaps unexpected.
As we shall see, some extra condition is needed to exclude various pathologies;
see \refE{Ebad}, for example. Condition \ref{T1d} is in fact extremely weak:
the natural interpretation of $\gk$ is that it measures the
density of edges, so the integral should be the expected number of edges,
suitably normalized. Thus condition \ref{T1d} says that
$\gnkx$ has about the right number of edges, so
if \ref{T1d} does not hold, $\kk$ has failed to capture
even the most basic property of the graph. What is surprising, is that
this condition is enough: we shall show that the assumptions above are enough
for $\kk$ to capture many properties of the graph.

In fact, in many circumstances, condition \ref{T1d} is automatically satisfied.
Indeed, 
one of the two inequalities implicit in this definition, namely
\[
\liminf \frac1n \E e\bigpar{\gnkx}
\ge \hiikxy,
\]
always holds; see \refL{LE}. This lemma also shows that \ref{T1d}
holds whenever $\kk$ is bounded and $\vxs$ is a vertex space.
It 
also holds whenever $\vxs$ is a vertex space in which 
the $x_i$ are (pairwise) independent and
distributed according to $\mu$. Moreover, condition \ref{T1d} is likely to
hold, and to be easy to verify, for any particular model that is of
interest. \refP{PE} shows that when \ref{T1d} does hold, so the normalized
number of edges converges in expectation, then it also converges in
probability. 
Note also that \eqref{t1b} holds if and only 
if the corresponding relation holds for the variants of $\gnkx$ discussed
in \refR{Rexp}; see \refR{Rckgraphical}.

In conjunction with condition \ref{T1d},
condition \ref{T1c} 
says that the expected
number of edges is $O(n)$, so the (expected) average degree is $O(1)$.
There are interesting cases with more
edges, but they will not be treated here;
\cf{} \refS{Sconcl}.
\end{remark}

We can be somewhat more general and allow minor deviations in
\eqref{pij} by letting $\kk$ depend on $n$. This will ensure that our
results apply directly to the various variations on the model
discussed above. The conditions we shall need on a sequence
of kernels are contained in the next definition.

\begin{num_definition}\label{Dg2}
Let $\vxs=(\sss,\mu,\nxss)$ be a (generalized) vertex space and let $\kk$ be a
kernel on $\vxs$. A sequence $(\kkn)$ of kernels on
$(\sss,\mu)$ is \emph{graphical on $\vxs$ with limit $\kk$}
if, for \aex{}
$(y,z)\in\sss^2$,
\begin{equation}
  \label{t2a}
\text{$y_n\to y$ and $z_n\to z$ imply that }
\kk_n(y_n,z_n)\to\kk(y,z),
\end{equation}
$\kk$ satisfies conditions \ref{T1a} and \ref{T1c} of \refD{Dg1},
and
\begin{equation}
\label{t2b}
\frac1n  \E e\xpar{\gnxx{\kk_n}} \to \hiikxy.
\end{equation}
\end{num_definition}

Note that if $\kk$ is a graphical kernel on $\vxs$, then the sequence
$\kkn$ with $\kkn=\kk$ for every $n$ is a graphical sequence on $\vxs$
with limit $\kk$.
Much of the time, members of a graphical sequence of kernels on $\vxs$ are
themselves graphical on $\vxs$.

\medskip
Much of the time, the conditions in \refD{Dg2} will be all we shall 
need to prove results about $\gnxx{\gk_n}$. However, when we come
to the size of the giant component, one additional condition will be needed.

\begin{num_definition}\label{Dg0}
A kernel $\kk$ on a (generalized) ground space $(\sss,\mu)$
is \emph{reducible} if
\begin{equation*}
  \text{$\exists A\subset \sss$ with
$0<\mu(A)<\mu(\sss)$ such that $\kappa=0$ a.e.\ on $A\times(\sss\setminus
A)$};
\end{equation*}
otherwise $\gk$ is \emph{irreducible}. Thus $\gk$ is irreducible
if
\begin{equation}
\label{t1ay}
  \text{$A\subseteq \sss$ and $\kappa=0$ a.e.\ on $A\times(\sss\setminus A)$
implies $\mu(A)=0$ or $\mu(\sss\setminus A)=0$}.
\end{equation}
\end{num_definition}

Roughly speaking, $\kk$ is reducible if the vertex set of $\gnkx$ can
be split into two parts so that the probability of an edge from one
part to the other is zero, and irreducible otherwise.
For technical reasons, we consider a slight weakening of irreducibility.

\begin{num_definition}\label{Dg0a}
A kernel $\kk$ on a (generalized) ground space $(\sss,\mu)$ is \emph{\qir} if
there is a \mucs{} $\sss'\subseteq\sss$ 
with $\mu(\sss')>0$ such that the restriction
of $\kk$ to $\sss'\times\sss'$
is irreducible, and $\kk(x,y)=0$ if $x\notin \sss'$ or $y\notin
\sss'$.
\end{num_definition}

\begin{remark}
  \label{Rquasi}
Given a \qir\ kernel $\kk$ and the associated graph $G_n=\gnkx$,
we can consider the irreducible restriction $\kk'$ of $\kk$ to
$\sss'\times\sss'$, and the corresponding graph
$G_n'=G^{\vxs'}(n,\kk')$ obtained from $G_n$ by deleting the
vertices with types in $\sss\setminus\sss'$; these vertices are
isolated in $G_n$. This graph is an instance of our model with 
a generalized vertex space $\vxs'$; note that the number $N'$ of
vertices of $G_n'$ 
is random. Thus, we may reduce suitable questions about \qir\ kernels
to the irreducible case. In our main results, the reader will lose
nothing by reading irreducible instead of \qir. We state some of
the results for the \qir\ case because this is all we need in the
proofs (even without removing isolated vertices as above), and we
sometimes need 
the \qir\ case of one result to prove the irreducible case of
another.
\end{remark}

\subsection{A branching process}\label{SS:branch}
Let $\kk$ be a kernel on a (generalized) ground space $(\sss,\mu)$.
To study the component structure of \gnk, we shall use the multi-type
Galton--Watson branching process with type space $\sss$, where a
particle of type $x\in \sss$ is replaced in the next generation by
a set of particles distributed as a Poisson process on $\sss$ with
intensity $\kk(x,y)\dd\mu(y)$. (Thus, the number of children with
types in a subset $A\subseteq\sss$ has a Poisson distribution with
mean $\int_A \kk(x,y)\,d\mu(y)$, and these numbers are independent
for disjoint sets $A$ and for different particles; see, \eg{},
Kallenberg \cite{Kall}.) We denote this branching process, started with a
single particle of type $x$, by \bpkx. When $\mu(\sss)=1$, so $\mu$ is
a probability measure, 
we write $\bpk$ for the
same process with the type of the initial particle random,
distributed according to $\mu$.

Let $\rhok(\kk;x)$ be the probability that the branching process
\bpkx\ has a total population of exactly $k$ particles, and let
$\rhogek(\kk;x)$ be the probability that the total population is
at least $k$. Furthermore, let $\rho(\kk;x)$ be the
probability that the branching process survives for eternity. If
the probability that a particle has infinitely many children is 0, then
$\rho(\kk;x)$ is equal to $\rho_{\infty}(\kk;x)$, the probability that
the total population is infinite; see \refR{R:bb1}.

Set
\begin{align} \label{rho}
\rhok(\kk)\=\int_\sss\rhok(\kk;x)\dd\mu(x),
&&&
\rho(\kk)\=\int_\sss\rho(\kk;x)\dd\mu(x),
\end{align}
and define $\rhogek(\kk)$ analogously.
Thus, if $\mu(\sss)=1$, then $\rho(\kk)$ is the survival probability
of the branching process $\bpk$.
Note that multiplying $\kk$ by a constant factor $c$ and dividing $\mu$ by the same factor
leaves the branching process $\bpk(x)$, and hence $\rho(\kk;x)$ and $\rhok(\kk;x)$,
unchanged. However, $\rho(\kk)$, for example, is decreased by a factor of $c$.

\begin{remark}\label{R:bp}
As we shall see later, the branching process \bpkx{} arises naturally 
when exploring a component of $\gnk$ starting at a vertex of type~$x$;
this is directly analogous to the use of the single-type Poisson branching
process in the analysis of the Erd\H os-R\'enyi graph $G(n,c/n)$. In models with a fixed degree sequence,
a related `size-biased' branching process arises, as it matters how we reach
a vertex. Here, due to the independence between edges, there is no size-biasing.
\end{remark}

\subsection{An integral operator}
Given a kernel $\kk$ on a (generalized) ground space $(\sss,\mu)$,
let $\tk$ be the integral operator on $(\sss,\mu)$ with kernel $\kk$,
defined by
\begin{equation}
  \label{tk}
(\tk f)(x)=\int_\sss \kk(x,y)f(y)\dd\mu(y),
\end{equation}
for any (measurable) function $f$ such that this integral is defined (finite or
$+\infty$) for \aex{} $x$.
As usual, we need never consider non-measurable functions; in future,
we shall assume without comment that all functions considered are measurable.
Note that $\tk f$ is defined for every $f\ge0$, with $0\le\tk
f\le\infty$. If $\kk\in L^1(\sss\times\sss)$, as we shall assume throughout,
then $\tk f$ is also defined for every bounded $f$; in this case  $\tk f\in
L^1(\sss)$
and thus $\tk f$ is finite a.e.

We define
\begin{equation}\label{tnorm}
  \norm{\tk}\=\sup\bigset{\normll{\tk f}: f\ge0,\,\normll{f}\le1}
\le\infty.
\end{equation}
When finite, $\norm{\tk}$ is the norm of $\tk$ as an operator in
$L^2(\sss,\mu)$;
it is infinite  if $\tk$ does not define a bounded operator
in $L^2$.
Trivially,
$\norm{\tk}$ is at most the Hilbert--Schmidt norm of $\tk$:
\begin{equation}\label{HSdef}
 \norm{\tk} \le \normHS{\tk}\=\norm{\kk}_{L^2(\sss\times\sss)}
 = \left(\iint_{\sss^2} \kk(x,y)^2\dd\mu(x)\dd\mu(y)\right)^{1/2}.
\end{equation}

We also define the non-linear operator $\Phik$ by
\begin{equation}
  \label{Phik}
\Phik f \= 1-e^{-\tk f}
\end{equation}
for $f\ge0$.
Note that for such $f$ we have $0\le \tk f\le\infty$, and thus
$0\le\Phik f\le 1$. We shall
characterize the survival probability $\rho(\kk;x)$, and thus
$\rho(\kk)$, in terms of the 
non-linear operator $\Phik$, showing essentially that the function
$x\mapsto \rho(\kk;x)$  
is the maximal fixed point of the non-linear operator $\Phik$;
see \refT{Trho}.

\section{Main results}\label{Sresults}

In this section we present our main results describing 
various properties of the general model $\gnxx{\kkn}$;
some further general results will be given in the
later sections devoted to individual properties. In 
\refS{Sapp}, we shall present results for special
cases of the model, including several that have been
studied previously.

All our results are asymptotic, and all
unspecified limits are taken as \ntoo. We use the following
standard notation: for (deterministic) functions $f=f(n)$ and $g=g(n)$,
we write $f=O(g)$ if $f/g$ is bounded, $f=\Omega(g)$ if $f/g$ is bounded 
away from zero, i.e., if $g=O(f)$,
and $f=\Theta(g)$ if $f=O(g)$ and $g=O(f)$. We write $f=o(g)$ if $f/g\to 0$.

Turning to sequences of events and random variables,
we say that an event holds \emph{with high probability} (\whp), if
it holds with probability tending to 1 as $n\to\infty$.
(Formally, it is a {\em sequence} $E_n$ of events that may hold \whp, but
the $n$ is often suppressed in the notation.)
We write $\pto$ for convergence in probability. Thus, for example, 
if $a\in\bbR$, then $X_n\pto a$ if and only if, for every $\varepsilon >0$, the
relations $X_n>a-\eps$ and $X_n<a+\eps$ hold \whp{}.

We shall
use $O_p$, $o_p$ and $\Theta_p$ in the standard way (see \eg{}
Janson, {\L}uczak and Ruci\'nski~\cite{JLR}); for example,
if $(X_n)$ is a sequence of random variables, then 
$X_n=O_p(1)$ means ``$X_n$ is bounded in probability'' and
$X_n=o_p(1)$ means that $X_n \pto0$.  Given a function $f(n)>0$,
we shall write $X_n=O(f(n))$ \whp{} if there
exists a constant $C<\infty$ such that $|X_n| \le Cf(n)$ \whp.
(This is written $X_n=O_C(f(n))$ in \cite{JLR}.) Note that this is stronger
than $X_n=O_p(f(n))$; the two statements can be written as $\exists C
\forall \eps \limsup_n \P(|X_n|>Cf(n)) <\eps$ and $\forall \eps
\exists C \limsup_n \P(|X_n|>Cf(n)) <\eps$, respectively. We shall
use $X_n=\Theta(f(n))$ \whp\ similarly.

We denote the orders of the components of a graph $G$ by
$C_1(G)\ge C_2(G) \ge\dots$, with $C_j(G)=0$ if $G$ has fewer than
$j$ components. We let $\nk(G)$ denote the total number  of 
vertices in components of order $k$, and write $\ngek(G)$ for
$\sum_{j\ge k}\nx{j}(G)$, the number of vertices in components of
order at least $k$. 

We shall write
$a\bmin b$ and $a\bmax b$ 
for $\min\{a,b\}$ and $\max\{a,b\}$, and use the same notation for the
pointwise minimum or maximum of two functions.

\medskip
As noted in the previous section, a reader who wishes to understand 
the following results, rather than apply them to a specific model,
may safely ignore all references to generalized vertex spaces.

\subsection{Existence, size and uniqueness of the giant component}\label{SSgc}

Our first result gives a necessary and sufficient condition
for the existence of a giant component in our model.

\begin{theorem}
   \label{T2}
Let $(\kkn)$ be a graphical sequence of kernels on a (generalized) vertex space
$\vxs$ with limit $\kk$.
\begin{thmxenumerate}
\item
If\/ $\norm{\tk}\le1$, then
$C_1\bigpar{\gnxx{\kk_n}}=o_p(n)$, while if\/
$\norm{\tk}>1$, then
$C_1\bigpar{\gnxx{\kk_n}}=\Theta(n)$ \whp.
\item
For any $\eps>0$, \whp\ we have 
\begin{equation}\label{t2upper}
  \frac1n C_1\bigpar{\gnxx{\kk_n}} \le \rho(\kk)+\eps.
\end{equation}
\item\label{T2_3}
If $\gk$ is \qir, then
\begin{equation}\label{t2bound}
\frac1n  C_1\xpar{\gnxx{\kk_n}} \pto \rho(\kk).
\end{equation}
\end{thmxenumerate}
In all cases $\rho(\kk)<\mu(\sss)$; furthermore, $\rho(\kk)>0$ if
and only if $\norm{\tk}>1$.
\end{theorem}
This result will be proved in \refS{Smainpf}, along with an additional result,
\refT{T1A}, concerning the distribution of the types of the vertices
in the giant component. 
We have included the final statement about $\rho(\kk)$ for ease of
future reference, even though it is purely a statement about the
branching process $\bpk$.
As remarked above, $\rho(\kk)$ can be found from
the solutions of the non-linear equation $f=\Phik(f)$;
see \refT{Trho}.

\refT{T2} has several immediate consequences. As customary,
we say that a sequence of random graphs $G_n$ (with $\Theta(n)$ 
vertices
in $G_n$) has a \emph{giant component} (\whp) if
$C_1(G_n)=\Theta(n)$ \whp. For simplicity we state these
results in the form where the kernel $\kkn$ is independent of $n$.

\begin{corollary}
  \label{C1}
Let $\kk$ be a graphical kernel on a (generalized) vertex space $\vxs$,
and consider the random graphs $\gnxx{c\kk}$ where $c>0$ is a constant.
Then the threshold for the existence of a giant component is
$c=\norm{\tk}\qi$.
More precisely, if $c\le\norm{\tk}\qi$,
then $C_1\bigpar{\gnxx{c\kk}}=o_p(n)$,
while if $c>\norm{\tk}\qi$ and $\kk$ is irreducible,
then $C_1\bigpar{\gnxx{c\kk}}=\rho(c\kk)n+o_p(n)=\Theta_p(n)$.
\end{corollary}

\begin{corollary}
  \label{C2}
Let $\kk$ be a graphical kernel on a (generalized) vertex space
$\vxs$. Then the 
property that
$\gnxx{c\kk}$ has \whp{} a giant component holds for every $c>0$ if
and only if
$\norm{\tk}=\infty$. Otherwise it has a finite threshold $c_0>0$.
\end{corollary}

The corollaries above are immediate from \refT{T2}, the observation 
that $\norm{\Tx{c\kk}}=c\norm{\tk}$, and the fact
that $\kk$ graphical on $\vxs$ implies $c\kk$ graphical on $\vxs$ (see
\refR{Rckgraphical}).
In the light of the results above, we say that a kernel $\gk$ is
\emph{subcritical} if $\norm{\tk}<1$, \emph{critical} if
$\norm{\tk}=1$, and \emph{supercritical} if $\norm{\tk}>1$. We use
the same expressions for a random graph $\gnk$ and a
branching process $\bpk$.

The next result shows that the number of edges in the graph at the
point where the giant component emerges is maximal
in the classical Erd\H os--R\'enyi case, or the slightly more general
`homogeneous case' described in \refE{Ehomo}; see \refS{Stransitionpf} for
the proof. (In this result we do need $\mu(\sss)=1$ as a normalization.)

\begin{proposition}\label{PC6}
Let $\kkn$ be a graphical sequence of kernels on a vertex space 
$\vxs$ with limit $\kk$, and assume that $\kk$ is critical, \ie{}
$\norm{\tk}=1$.
Then $\frac1n e(\gnxx{\kkn})\pto\hiik\le 1/2$, with equality in the uniform
case $\kk=1$;
more precisely, equality holds if and only if 
$\int_\sss \kk(x,y)\dd\mu(y)=1$ for a.e.\ $x$.
\end{proposition}

We can also determine the asymptotic number of edges in the giant
component. As this is not always uniquely defined,
for any graph $G$, let $\cC_1(G)$ be the {\em largest component} of $G$,
i.e., the component with most vertices, chosen according to any fixed
rule if there 
is a tie. In order to state the next result concisely,
let
\begin{equation}\label{Zdef}
 \edgenokk \=
 \frac12
 \iint_{\sss^2}\kk(x,y)\bigpar{\rho(\kk;x)+\rho(\kk;y)-\rho(\kk;x)\rho(\kk;y)}
 \dd\mu(x)\dd\mu(y).
\end{equation}
Note that the bracket above is the probability that, given 
independent branching processes $\bpkx$ and $\bpk(y)$, at least one survives.
Intuitively, given that a certain edge is present in $\gnxx{\kkn}$,
this edge is in the giant component if, when exploring the rest of the graph from its end-vertices,
there is at least one from which we can reach many vertices.

\begin{theorem}\label{Tedges}
Let $(\kkn)$ be a graphical sequence of kernels on a (generalized) vertex space
$\vxs$ with
\qir\ limit $\kk$.
Then
\begin{equation}\label{eedge}
\frac1n  e\bigpar{\cC_1\xpar{\gnxx{\kk_n}}} \pto \edgenokk.
\end{equation}
\end{theorem}

This result will be proved in \refS{Sgiantedges},
together with some properties of $\edgenokk$. 

\medskip
Under our assumptions, the giant component is \whp{} unique when it exists;
the second largest component is much smaller. Indeed,
as we shall show in \refS{Smainpf},
only $o_p(n)$ vertices
are in `large' components other than the largest.

\begin{theorem}
   \label{T2b}
Let $(\kkn)$ be a graphical sequence of kernels on a (generalized) vertex space
$\vxs$ with \qir\ limit $\kk$,
and let $G_n=\gnxx{\kkn}$.
If\/ $\go(n)\to\infty$ and $\go(n)=o(n)$, then
\begin{equation} \label{t2b1}
\sum_{j\ge2:\; C_j(G_n)\ge\go(n)}  C_j(G_n) =o_p(n).
\end{equation}
In particular,
\begin{equation}
 C_2(G_n)=o_p(n).
\label{c2small}
\end{equation}
\end{theorem}

\begin{remark}\label{Rred}
If $\gk$ is reducible and the $x_i$ are (absolutely)
continuous random
variables, then $\gnkx$ decomposes into two (or more)
disjoint parts that can be regarded as $\gtx{\vxs_i}{n_i}{\kk_i}$, 
for suitable $\vxs_i$, $n_i$ and $\kk_i$. By considering the parts
separately, many of our results for the irreducible case can be
extended to the reducible case; note, however, that each of
the parts may contain a giant component, so it is possible to have
$C_2=\Theta_p(n)$. The restriction to the case where the $x_i$ are
continuous, which includes the Poisson case of \refE{EPo},
is necessary
unless we impose a further restriction on $\kk$:
in general there may be
a subset $A\subset\sss$ of measure zero which always contains some $x_i$,
and this can link the subgraphs $\gtx{\vxs_i}{n_i}{\kk_i}$.
Worse still, such an $A$ may contain an $x_i$ with probability bounded
away from $0$ and $1$, so $\frac{1}{n}C_1(\gnkx)$ need not converge in
probability.
\end{remark}

\begin{remark}\label{Rmono}
If $\kk$ and $\kk'$ are two kernels on the same vertex space with
$\kk\le \kk'$,  
then $\gnkx$ and $\gnxx{\kk'}$ are random graphs on the same vertex
set, and there is a natural 
coupling between them in which
$\gnkx$ is always a subgraph of $\gnxx{\kk'}$, i.e., a coupling with
$\gnk\subseteq\gnx{\kk'}$. 
Similarly, one can couple the corresponding branching processes so
that every particle 
present in one is present in the other, i.e., so that
$\bpk\subseteq \bpx{\kk'}$. Thus $\rho(\kk)\le\rho(\kk')$.
If $\kk$ is irreducible 
and $\rho(\kk)>0$, then it follows from \refT{Trho} and \refL{L1c}
that $\rho(\kk')>\rho(\kk)$ unless $\kk'=\kk$ a.e.
Similarly, the threshold $c_0(\kk')\=\norm{T_{\kk'}}\qi$ is at
most $c_0(\kk)\=\norm{T_{\kk}}\qi$. Here, however, somewhat
surprisingly, we may have $c_0(\kk')=c_0(\kk)$
even if $\kk'>\kk$; see \refSS{SSchkns}. On the other hand, it is
easily seen that if $\tk$ is compact,
$\kk$ is irreducible,
and $\kk'>\kk$ on a set of
positive measure, then $\norm{T_{\kk'}}
>\norm{T_{\kk}}$ and thus $c_0(\kk')<c_0(\kk)$.
\end{remark}

\subsection{Stability}  \label{SSstable}\label{SSfirst_after_gc}

The giant component of $G_n=\gnxx{\kkn}$ is stable in
the sense that its size does not change much if we add or delete a
few edges; this is made precise in the following theorem. Note
that the edges added or deleted do not have to be random or
independent of the existing graph; they can be chosen by an
adversary after inspecting the whole of $G_n$. Also,
we may delete vertices instead of (or as well as) edges.

\begin{theorem}\label{TrobustNEW}
Let $(\kkn)$ be a graphical sequence of kernels on a (generalized) vertex space
$\vxs$ with irreducible
limit $\kk$, and let $G_n=\gnxx{\kk_n}$. For every
$\eps>0$ there is a $\gd>0$ (depending on $\kk$) such that, \whp{},
\begin{equation}\label{estableNEW} 
 (\rho(\kk)-\eps)n \le C_1(G_n') \le (\rho(\kk)+\eps)n
\end{equation}
for every graph $G_n'$ that may be obtained from $G_n$ by deleting at most
$\gd n$ vertices and their incident edges, and then
adding or deleting at most $\gd n$ edges.
\end{theorem}

In particular, if $G_n'$ is a graph on $V(G_n)$ with
$e(G_n'\setdiff G_n)=o_p(n)$ then
\begin{equation*}
 C_1(G_n')=C_1(G_n)+o_p(n)=\rho(\kk)n+o_p(n).
\end{equation*}
\refT{TrobustNEW} is proved in \refS{Sstable}.
Clearly, in proving the first inequality in \eqref{estableNEW}, we
may assume that $G_n'\subseteq G_n$, and in proving the second
that $G_n\subseteq G_n'$. The latter case will be easy
to deal with using \refT{T2b}.
Proving the first inequality amounts to showing that, \whp, the giant
component of $G_n$ cannot be cut into two pieces of size
at least $\Theta(n)$ by deleting $o(n)$ vertices and then $o(n)$ edges.
For edge deletion,
Luczak and McDiarmid~\cite{LMcD} gave a very simple
proof of this result in
the Erd\H os-R\'enyi case, which adapts easily to the finite
type case and hence (using our general results) to the full
generality of \refT{TrobustNEW}. This proof is presented
in \refS{Sstable}.

Another approach to proving \refT{TrobustNEW} involves reducing this
statement to an equivalent statement about the two-core, which is very
easy to prove in the uniform case. This reduction involves relating
the two-core to the branching process, using results that we believe
are of interest in their own right, presented in \refS{Sstable}.
Unfortunately, while the general
case of the two-core result can be proved by branching process
methods, the proof is very complicated, so we shall not give it.

\begin{remark}
\refT{TrobustNEW} may be viewed as a statement about the vulnerability
of large-scale networks to attack by an adversary who knows the
detailed structure of the network, and attempts to disconnect the
network into small pieces by deleting a small fraction of the vertices
or edges.  The vulnerability of `scale-free' networks to such attacks
has been considered by many people; see, for example,
\cite{AJBattack,CNSWfragile,CohenAttack, BRcoupling}; it turns out
that such networks are much more resilient to random failures than
$G(n,c/n)$, but also more vulnerable to attack. In general, the
flexibility available to the attacker makes rigorous analysis
difficult, although a result for the Barab\'asi-Albert model was given
in~\cite{BRcoupling}. \refT{TrobustNEW} shows in particular that, for
\gnkx, the network is at most a constant factor more vulnerable than a
homogeneous network: a constant fraction of the vertices or edges must be
deleted to destroy (or significantly shrink) the giant component.
\end{remark}

\begin{remark}\label{Rstable}
As pointed out by Britton and Martin-L\"of~\cite{BML}, in
the case of vertex deletion \refT{TrobustNEW} also has the
following interpretation: suppose that $G_n$ represents the network
of contacts that may allow the spread of an infectious disease from
person to person, and that we wish to eliminate
the possibility of an epidemic by vaccinating some of
the population. Even if the entire network of contacts is known,
if the source of the infection is not known, a significant (constant,
as $n\to\infty$) proportion of the population must be vaccinated:
otherwise, there is still a giant component in the graph on the
unvaccinated people, and if the infection starts at one of its vertices,
it spreads to $\Theta(n)$ people.
\end{remark}

\subsection{Bounds on the small components}\label{SSsmall}

For the classical random graph $\gnx{c/n}$ it is well-known that
in the subcritical ($c<1$) case, $C_1=O(\log n)$ \whp{}, 
and that in the supercritical ($c>1$) case,
$C_2=O(\log n)$ \whp{}; 
see \cite{BB,JLR}, for example.
These bounds do not always hold in the general framework we are considering
here, but if we add some
conditions, then we can improve the estimates
$o_p(n)$
in \refT{T2} and \eqref{c2small} to $O(\log n)$ \whp.
As before, we write $G_n$ for $\gnxx{\kkn}$.

\begin{theorem}
   \label{T4}
Let $(\kkn)$ be a graphical sequence of kernels on a (generalized) vertex space
$\vxs$ with limit $\kk$.
\begin{romenumerate}
\item
If\/ $\gk$ is subcritical,
\ie{}, $\norm{\tk}<1$, and $\sup_{x,y,n}\gk_n(x,y)<\infty$,
then
$C_1\xpar{G_n}=O(\log n)$ \whp.
\item \label{T4p2}
If\/ $\gk$ is supercritical,
\ie{}, $\norm{\tk}>1$, $\gk$ is irreducible,
and either $\inf_{x,y,n}\gk_n(x,y)>0$
or $\sup_{x,y,n}\gk_n(x,y)<\infty$,
then
$C_2\xpar{G_n}=O(\log n)$ \whp.
\end{romenumerate}
\end{theorem}

\refT{T4} is proved in \refS{S2ndpf}.
Note that in part (ii) we draw the same conclusion from
the very different assumptions $\inf_{x,y,n}\kk_n(x,y)>0$ and
$\sup_{x,y,n}\kk_n(x,y)<\infty$. There is no similar result for
the subcritical case (part (i))
assuming only that
$\inf_{x,y,n}\kk_n(x,y)>0$: \cite[Theorems 1 and  2]{BJR} show that the
random graph $\Gcjn$ with $0<c<1/4$ is subcritical  and satisfies
$C_1(\Gcjn)=n^{\Theta(1)}$ \whp.

\subsection{Degree sequence}\label{SSdegrees}

We next turn to the degrees of the vertices of
$G_n=\gnkxn$, where $\kkn\to\kk$. As we shall see, the degree
of a vertex of a given type $x$ is asymptotically Poisson
with a mean
\begin{equation}\label{gldef}
 \gl(x)\=\ints\kk(x,y)\dd\mu(y)
\end{equation}
that depends on $x$. This leads to a mixed Poisson
distribution for the degree $D$ of a (uniformly chosen) 
random vertex of $G_n$.
We write $Z_k$ 
for the number of vertices of $G_n$ with degree $k$.
\begin{theorem}
    \label{TD}
Let $(\kkn)$ be a graphical sequence of kernels on a (generalized) vertex space
$\vxs$ with limit $\kk$, and let $G_n=\gnkxn$.
For any fixed $k\ge0$, 
\begin{equation*}
  Z_k/n \pto \ints\frac{\gl(x)^k}{k!}e^{-\gl(x)}\dd\mu(x),
\end{equation*}
where $\gl(x)$ is defined by \eqref{gldef}.
Equivalently,
\begin{equation*}
  Z_k/|V(G_n)| \pto \P(\Xi=k),
\end{equation*}
where 
$\Xi$ has the mixed Poisson distribution
$\ints\Po(\gl(x))\dd\mu(x)/\mu(\sss)$. 
\end{theorem} 

In other words, if $D$ is the degree of a random vertex of
$G_n=\gnkxn$, and we normalize so that $\mu(\sss)=1$, then 
\begin{equation*}
  \cL(D\mid G_n)\pto\cL(\Xi)= \ints\Po(\gl(x))\dd\mu(x).
\end{equation*}

As we shall show in \refS{Sdegrees} and Subsections \refand{SSrtij}{SSrank1},
our model includes natural examples of `scale-free' random graphs,
where the degree distribution has a power-law tail.
We believe that when it comes to modelling 
real-world graphs with, for example, observed power laws for vertex degrees,
our model provides an interesting and flexible alternative to existing
models based on generating graphs with a given degree sequence (e.g.,
Molloy and Reed \cite{MR1,MR2}), or given expected degrees (e.g.,
Aiello, Chung and Lu \cite{AielloCL}).

\subsection{Distances between vertices}\label{SSdist}

Next, we consider the distances between vertices of $G_n=\gnxx{\kkn}$ where,
as usual, $\kkn$ is a graphical sequence of kernels on $\vxs$ with limit $\kk$.
Let us write $d(v,w)$ for the graph distance between two vertices of $G_n$,
which we take to be infinite if they lie in different components.
Note that
\begin{equation}\label{nconn}
 \bigl|\bigl\{ \{v,w\}: d(v,w)<\infty \bigr\}\bigr| 
= \sum_i \binom{C_i(G_n)}{2},
\end{equation}
where $\{v,w\}$ denotes an unordered pair of distinct vertices of $G_n$.

Under certain conditions, we can give upper and lower bounds on
$d(v,w)$ for almost 
all pairs with $d(v,w)<\infty$.

\begin{theorem}\label{Tdist}
Let $\kkn$ be a graphical sequence of kernels on a (generalized)
vertex space \vxs\ with 
limit $\kk$, with $\norm{\tk}>1$.
Let $G_n=\gnxx{\kkn}$, and let $\eps>0$ be fixed.
\begin{romenumerate}
\item\label{Tdist0}
If $\kk$ is irreducible, then
\[
 \bigl|\bigl\{ \{v,w\}: d(v,w)<\infty \bigr\}\bigr|
 = \frac{C_1(G_n)^2}2 + o_p(n^2) = \frac{\rho(\kk)^2n^2}2 + o_p(n^2).
\]
\item\label{Tdistu}
If $\supxyn<\infty$, then
\[
 \bigl|\bigl\{ \{v,w\}: d(v,w)\le (1-\eps)\log n/\log\norm{\tk} \bigr\}\bigr| 
= o_p(n^2).
\]
\item\label{Tdistl}
If $\kk$ is irreducible and $\norm{\tk}<\infty$, then
\[
\bigl|\bigl\{ \{v,w\}: d(v,w)\le (1+\eps)\log n/\log\norm{\tk} \bigr\}\bigr| 
= \rho(\kk)^2n^2/2+o_p(n^2).
\]
\item\label{Tdistl2}
If $\kk$ is irreducible and $\norm{\tk}=\infty$, then
there is a function $f(n)=o(\log n)$ such that
\begin{equation}\label{pi}
 \bigl|\bigl\{ \{v,w\}: d(v,w)\le f(n) \bigr\}\bigr| 
= \rho(\kk)^2n^2/2+o_p(n^2).
\end{equation}
\end{romenumerate}
\end{theorem}

Note that part \ref{Tdist0} is immediate from \eqref{nconn} and
Theorems \refand{T2}{T2b}. Related earlier results  
are discussed briefly in \refS{Sdist}.

\medskip
In the finite-type non-critical case, we can give
an asymptotic formula for the `diameter' of $G_n$, i.e., for
\[
 \diam(G_n) \= \max\{d(v,w)\: :\: v,w\in V(G),\:d(v,w)<\infty\},
\]
the maximum of the diameters of the components of $G_n$.
This turns out to depend not only on the norm of $\tk$, but also 
on the norm of the operator associated to the `dual kernel' $\kkd$.

\begin{num_definition}\label{Ddual}
Let $\kk$ be a supercritical kernel on a (generalized) ground space
$(\sss,\mu)$.  
The {\em dual kernel} is the kernel $\kkd$ on the generalized ground
space $(\sss,\mud)$  
defined by $\kkd(x,y)=\kk(x,y)$, with $\dd\mud(x) = (1-\rho(\kk;x))\dd\mu(x)$.
\end{num_definition}

Note that $\kkd$ and $\kk$ are identical as functions on
$\sss\times\sss$. However, they are defined 
on different generalized ground spaces. Hence, the operators $\tk$ and
$T_{\kkd}$ have 
(in general) different norms. If we wish to consider only ground spaces,
we may renormalize, defining $\kkd'$ on $(\sss,\mud')$ by
$\kkd'(x,y)=(1-\rho(\kk))\kk(x,y)$ and $\dd\mud'(x) =
(1-\rho(\kk;x))/(1-\rho(\kk))\dd\mu(x)$. 
The choice of normalization does not affect the norm of the operator:
$\norm{T_{\kkd}}=\norm{T_{\kkd'}}$. 

The relevance of the dual kernel is that it describes the `small'
components of $\gnkx$; see 
\refS{S2ndpf}. The distribution of these small components is
essentially the same 
as the distribution of trees hanging off the two-core of the giant
component, which affects 
the diameter of $\gnkx$.

\begin{theorem}\label{th_diam}
Let $\kk$ be a kernel on a (generalized) vertex space $\vxs=(\sss,\mu,(\xs))$, 
with $\sss=\{1,2,\ldots,r\}$ finite and 
$\mui>0$ for each $i$, and let $G_n=\gnxx{\kk}$.   
If\/ $0<\norm\tk<1$, then
\[
 \frac{\diam(G_n)}{\log n} \pto \frac{1}{\lntki}
\]
as $n\to\infty$.
If\/ $\norm\tk>1$ and $\kk$ is irreducible, then
\[
 \frac{\diam(G_n)}{\log n} \pto   \frac{2}{\lntkdi} +\frac{1}{\lntk},
\]
where $\kkd$ is the dual kernel to $\kk$.
\end{theorem}

Note that we do not require $\kk$ to be graphical on $\vxs$: if $\vxs$
is a vertex space,  
then, as $\sss$ is finite, any kernel $\kk$ on $\vxs$ is graphical;
see \refR{Rrfin}. 
If $\vxs$ is a generalized vertex space, then $\kk$ need not be graphical.
However, by conditioning on the sequences $(\xs)$, we can reduce to
the vertex space 
case; see \refSS{SSgvs}.

The assumptions of \refT{th_diam} are much more restrictive than those 
of our other results: we require the type space to be finite.
Note, however, that even the single type case of this result,
concerning the classical random graph $G(n,c/n)$, is non-trivial;
it answers in the negative a question of Chung and Lu~\cite{ChungLu:diam}.
This special case of \refT{th_diam} was proved independently by
Fernholz and Ramachandran~\cite{FR:diam}, again as a special case of a result
for a more general model. The nature of their model makes their
proof much more difficult than that of \refT{th_diam}; see \refSS{SSdiam}.

\subsection{The phase transition}\label{SStrans}\label{SSlast}

Finally, we turn to the phase transition in $\gnkx$, where the giant 
component first emerges. As usual, to study the transition, we should
vary a single density parameter. Here, it is most natural
to fix a graphical kernel $\kk$ on a vertex space $\vxs$,
and to study $\gnxx{c\kk}$ for a real parameter $c>0$,
as in \refC{C1}. Instead, we could consider random subgraphs of $\gnkx$
obtained by keeping each edge, or edge vertex, independently with
probability $p$, and use $p$ as the parameter; as we shall see
in Examples \refand{Eep}{Evp}, all three approaches are equivalent,
so we shall use the first.

By \refT{T2}, the size of the largest component of \gnxx{c\kk} is
described by the function $\rho(c\kk)$, which is 0 for
$c\le c_0\=\norm{\tk}\qi$ and strictly positive for larger $c$.
With $\vxs$ and $\kk$ fixed, let us denote this function by
$\rho(c)$, $c>0$.
We shall see (from \refT{TappB}) that $\rho(c)$ is continuous on $(0,\infty)$.

Since $\rho(c)=0$ for $c\le c_0$ but not for larger $c$, the function
$\rho$ is not analytic at $c_0$; in physical terminology, {\em there is a
phase transition at} $c_0$.

For the classical Erd\H os--R\'enyi random graph \gnx{c/n} (obtained
with $\kk=1$), it is well-known that
$\rho$ is continuous but
the first derivative has a jump at $c_0=1$; more precisely, $\rho'$
jumps from 0 to $\rho'_+(c_0)=2$.
For finite $d$, we shall say that the phase transition in $\gnkx$ 
has \emph{exponent} $k$ if $\rho(c_0+\eps)=\Theta(\eps^k)$ as
$\eps\downto 0$.
As we have just noted,
in \gnx{c/n} the phase transition has exponent $1$.
If $\rho(c_0+\eps)=o(\eps^k)$ for all $k$, we say that the phase
transition has \emph{infinite exponent}.
We are deliberately avoiding
the physical term `order',
as it is not used in a consistent way in this context.
In other contexts, discontinuous phase transitions are possible; see,
for example, 
Aizenman, Chayes, Chayes and Newman~\cite{ACCN}.

It was shown in \cite{BJR} (see also
Dorogovtsev, Mendes and Samukhin~\cite{DMS-anomalous} and Durrett
\cite{Durrett}) that in the case $\sss=(0,1]$ and $\kk(x,y)=1/(x\bmax
y)$, the phase transition `is of infinite order', i.e., has
infinite exponent (see \refSS{SSDubins} for more details). We shall see
in \refSS{SSrank1} that it is also possible to have a phase
transition with any finite exponent larger than 1 (including
non-integer values).

The next theorem shows that the phase transition has exponent 1
for a wide class of kernels $\kk$, including all bounded $\kk$.
We also prove that for this class there is no other phase transition:
as $\rho=0$ on $(0,c_0)$, it is trivially analytic there,
and we shall prove that $\rho$ is analytic on $(c_0,\infty)$.
As $\rho$ is defined in terms of the branching process, rather than
a graph, we do not need a vertex space for the statement of the next result;
to deduce conclusions for graphs of the type we consider, we should let
$\kk$ be an irreducible graphical kernel on a vertex space $\vxs$,
satisfying the
additional condition \eqref{t5a} below. (Also, there is no need to consider 
generalized ground spaces, as we may trivially normalize so that $\mu(\sss)=1$
by multiplying $\kk$ by $\mu(\sss)$ and dividing $\mu$ by the same
factor -- this 
leaves the branching process unchanged.) 
When we say that a function $f$ defined on the reals is \emph{analytic}
at a point $x$, we mean that 
there is a neighbourhood of $x$ in which $f$
is given by the sum of a convergent power series;
equivalently, 
$f$ extends to a complex analytic
function in a complex neighbourhood of $x$.

\begin{theorem}\label{T5}
Let $\kk$ be a kernel on a ground
space $(\sss,\mu)$. Suppose
that $\kk$ is irreducible, and that
  \begin{equation}\label{t5a}
\sup_x \int_\sss \kk(x,y)^2\dd\mu(y) <\infty.
  \end{equation}
\begin{thmxenumerate}
\item The function $c\mapsto\rho(c)\=\rho(c\kk)$ is analytic
except at
  $c_0\=\norm{\tk}\qi$.
\item
The linear operator $\tk$ has an eigenfunction $\psi$ of eigenvalue
$\norm{\tk}<\infty$,
and every such eigenfunction is bounded and satisfies
\begin{equation}\label{t5}
  \rho(c_0+\eps)
=2c_0\qi\frac{\int_\sss\psi \int_\sss\psi^2}{\int_\sss\psi^3}
\eps+O(\eps^2),
\qquad \eps>0,
\end{equation}
so $\rho'_+(c_0)=
2c_0\qi\xfrac{\int_\sss\psi \int_\sss\psi^2}{\int_\sss\psi^3}>0$
and $\rho$ has a phase transition at $c_0$ with exponent $1$.
\end{thmxenumerate}
\end{theorem}

The proof is given in \refS{Stransitionpf}.
Note that \eqref{t5a} implies
that $\kk\in L^2\subseteq L^1$.
\refT{T5} has an easy consequence concerning the extremality of the
Erd\H os--R\'enyi random graphs, also proved in \refS{Stransitionpf}.

\begin{corollary}\label{C5}
Let $\kk$ be an irreducible kernel on a ground space $(\sss,\mu)$ such that 
\eqref{t5a} holds, and let
$c_0\=\norm{\tk}\qi>0$.
Then $c_0\rho'_+(c_0)\le 2$, with equality
in the classical Erd\H{o}s--R\'enyi case;
more precisely, equality holds if and only if 
$\int_\sss \kk(x,y)\dd\mu(y)=1$ for a.e.\ $x$.
\end{corollary}

Let $\kk$ be an irreducible graphical kernel on a vertex space $\vxs$; let us
assume \eqref{t5a} and, as a normalization, that $c_0=1$.
Letting $c$ increase from the threshold $c_0$, \refC{C5} says that
the giant component of $\gnxx{c\kk}$ has maximal growth-rate
in the Erd\H os--R\'enyi case, and, more generally, in the
`homogeneous' case treated in \refE{Ehomo} below.
In this example, the vertex
degrees are more or less the same, so there is no first-order
inhomogeneity 
in the graph; any inhomogeneity in vertex degrees leads
to a slower growth.

\begin{remark}\label{Retrans}
By \refT{Tedges}, the number of edges in the giant component
of $\gnx{c\kk}$ near the phase transition is asymptotically determined
by the behaviour of the function $\edgeno(c\kk)$ as $c\downto c_0$.
As we shall show in \refP{Pedges}, if $\norm{\tk}<\infty$, then
$\edgeno(c\kk)/\rho(c\kk)\to1$ as $c\downto c_0\=\norm{T_{\kk}}\qi$.
In particular, under the conditions of \refT{T5}, there is a phase
transition of exponent 1 in $\edgeno$ too.
(In addition, the proof of \refT{T5} will show that
$c\mapsto\edgeno(c\kk)$ is also analytic except at $c_0$.)
In the case $\norm{\tk}=\infty$, when $c_0=0$,
it is not always true that $\edgeno(c\kk)\sim\rho(c\kk)$ as $c\downto c_0$:
this will be shown by \refE{E2large2}. An important case
when this does hold is described in \refSS{SSrank1}.
\end{remark}

\section{Examples}\label{Sex}

In this section we give several simple examples of the random
graph model we study; these examples are chosen to illustrate the
definitions and the scope of the model, as well as various
pathologies that may occur. In subsequent sections we shall refer
to several of these examples;
in particular, many of our proofs will be based
on the `finite-type' case.
Further examples of interest in
their own right are discussed at length in \refS{Sapp}, as
applications of our results. We often suppress the dependence on
$\vxs$, writing $\gnk$ for $\gnkx$.

\begin{example}
{\bf The Erd\H os-R\'enyi random graph.}
  \label{E1}
If $\gk=c$ is constant, then the edge probabilities $\pij$ given 
by \eqref{pij} are all equal to $c/n$ (for $n>c$).
Thus any choice of vertex space gives the
classical Erd\H os--R\'enyi random graph \gnx{c/n}.
The simplest choice is to let $\sss$ consist of a single point. Then
the operator 
$\tk$ is simply multiplication by $c$, so $\norm{\tk}=c$ and \refC{C1} yields
the classical result that there is a phase
transition at $c=1$. Furthermore,
the function $\rho(c;x)$ reduces to the single value $\rho(c)$,
and the survival probability $\rho(c)$ of the branching
process $\bpx{c}$ is given by the formula
\begin{equation}\label{e1}
  \rho(c)=1-e^{-c\rho(c)},
\qquad\text{with }
\rho(c)>0 \text{ if } c>1;
\end{equation}
this classical branching process result is the simplest case of
\refT{Trho} below. 
Returning to the graph, in this case \refT{T2} reduces to the
classical result of Erd\H os and R\'enyi~\cite{ER60}. 
\end{example}

\begin{example}\label{Ebib}
{\bf The homogeneous bipartite random graph.}
Set $\sss=\{1,2\}$, $\mu\{1\}=\mu\{2\}=1$, and let
$\vxs=(\sss,\mu,(\xs))$ be a generalized  
vertex space in which $\xs$ consists of $n$ vertices of type $1$ and
$n$ vertices of type $2$. 
Let $\kk$ be defined by $\kk(1,1)=\kk(2,2)=0$ and $\kk(1,2)=\kk(2,1)=c$.
Then $\gnkx$ is the random bipartite graph $G(n,n;c/n)$ with $n$
vertices in each class, 
where each possible edge between classes is present with probability
$c/n$, independently 
of the other edges. While it is natural to use a generalized vertex
space to describe 
this example, it is not necessary: the same graph can be written as $\gnkx$
in another way: 
take $\mu\{1\}=\mu\{2\}=1/2$, and let $\vxs$ be a vertex space where
$\xs$ is defined 
only for $n$ even, and then consists of $m=n/2$ vertices
of each type. Let $\kk(1,1)=\kk(2,2)=0$ as before, and $\kk(1,2)=\kk(2,1)=2c$,
so the edge probabilities are $2c/n=c/m$.
\end{example}

\begin{example}
{\bf The finite-type case.}
  \label{Efin}
Let $\sss=\set{s_1,\dots,s_r}$ be finite. Then $\kk$ is an $r\times r$
matrix. In this case, $\gnk$ has vertices of $r$ different types (or
colours), say
$n_i$ vertices of type $i$, with two vertices of types $i$ and $j$
joined by an edge with probability $n\qi\kk(i,j)$ (for $n\ge\max\kk$).
The condition \eqref{a2a} means that $n_i/n\to\mu_i$ for each $i$ (in
probability if the $n_i$ are random), where $\mu_i\=\mu\set{i}\ge0$.

This case has been studied by S\"oderberg \cite{Sod1,Sod2,Sod3,Sod4},
who noted our \refT{T2} in this case (with $\kkn=\kk$ for all $n$).
\end{example}

Most of our proofs will be based on a disguised form of this case,
described by the following definition.

\begin{num_definition}\label{Drfin}
A kernel $\kk$ on a (generalized) ground space $(\sss,\mu)$ is {\em \rfin}
if $\sss$ has a finite partition into sets $S_1,\dots,S_r$ such that $\kk$ is
constant on each $S_i\times S_j$, where each $S_i$ is a \mucs, \ie{},
is measurable and has $\mu(\partial S_i)=0$.
\end{num_definition}

Clearly, if $\kk$ is \rfin{} on $(\sss,\mu)$ then 
the random graph $\gnkx$ has the same distribution
as a finite-type graph $\gnxxp{\kk'}$,
$\vxs'=(\sss',\mu',(\ys))$: take
$\sss'=\set{1,\dots,r}$, let $y_k=i$ whenever $x_k\in S_i$,
and define $\mu'\set{i}$ and
$\kk'(i,j)$ in the obvious way. Let $n_i=\#\set{l:x_l\in
S_i}=n\nu_n(S_i)$, where $\nu_n$ is as in \eqref{a2}. The numbers
$n_i$ may be random, but since each $S_i$ is a
\mucs, \eqref{a2a} yields $n_i/n=\nu_n(S_i)\pto\mu(S_i)$,
so $\vxs'$ is a (generalized) vertex space.

\begin{remark}\label{Rrfin}  
Let us note for later
that a finite-type or \rfin{} kernel $\kk$ on a vertex space \vxs{}
is automatically graphical on \vxs; conditions \refand{T1a}{T1c}
of \refD{Dg1} are trivial in this case, while condition \ref{T1d}
holds in the much more general case of $\kk$ bounded; see \refL{LE}.
This observation does not extend to generalized vertex spaces: 
there may be a very large number of vertices with some small probability,
so the expectation in \eqref{t1b} need not converge, or even exist;
see \refR{Rgenexp}. 
\end{remark}

\begin{example}
{\bf The homogeneous case.}
  \label{Ehomo}
Generalizing the Erd\H os-R\'enyi and homogeneous bipartite cases above, 
let $(\sss,\mu)$ be an arbitrary (generalized) ground space, and
let $\kk$ be such that
$\int_\sss\kk(x,y)\dd\mu(y)$ is essentially independent of $x\in\sss$,
i.e., that
\begin{equation}\label{homo}
 \int_\sss \kk(x,y)\dd\mu(y)=c
\qquad\text{for \aex{} $x$},
\end{equation}
for some constant $c$.
(This says roughly that, asymptotically, all vertices have the same
average degree.) Then $\tk 1=c$ \aex, so the constant function 1 is a positive
eigenfunction with eigenvalue $c$, and thus $\norm{\tk}=c$, and by
\refT{T2} there is a giant component (and
$\rho(\kk)>0$) if and only if $c>1$.

Normalizing (if necessary) so that $\mu(\sss)=1$, 
in the branching process, apart from particles
with types in a measure zero set, which arise in $\bpk$ with probability 0,
the number of children of each particle
has a $\Po(c)$ distribution. Hence, ignoring the types of the particles,
the distributions of the process $\bpk$ and the single-type process
$\bpx{c}$ are the same. In particular, $\rho(\kk)=\rho(c)$, so
$\rho(\kk)=\rho(c)$ is given by \eqref{e1} in this case too.
If $\kk$ is irreducible, the global behaviour of $\gnk$ is
thus exactly the same as that of $\gnc$, at least in terms of the
size of the giant component. The local behaviour can be quite
different, though. For example, $\gnk$ may
have many more triangles or other small cycles than $\gnc$; see
\refE{EChomo}. On the other hand, by \refT{TD},
the vertex degrees have an
asymptotic $\Po(c)$ distribution just as in \gnc.

A natural example of such a homogeneous $\kk$ is given by taking
$\sss$ as $(0,1]$ (now 
better regarded as the circle $\bbT$),
$\mu$ as Lebesgue measure, and $\kk(x,y)=h(x-y)$ for an even
function $h\ge0$ of period 1. For example, $h$ can be constant on
a small interval $(-\gd,\gd)$ and vanish outside it;
this gives a modification of $\gnc$ where only ``short'' edges are
allowed.

More generally, $\sss$ can be any compact homogeneous space, for
example a sphere, with Haar measure $\mu$ and an invariant metric
$d$, and $\kk(x,y)$ a function of the distance $d(x,y)$.
\end{example}

\begin{example}\label{Etri}
Take $\sss=(0,1]$ with $\mu$ the Lebesgue measure, 
and let $x_i=i/n$.
 Set $\gk(x,y)=\ett[x+y\le 1]$ and consider the kernel
$c\kk$, so that
\begin{equation*}
\pij=
\begin{cases}
  c/n, & i+j\le n;
\\
 0, & i+j > n.
\end{cases}
\end{equation*}
Thus $\gnx{c\kk}$ can be obtained from the random graph \gnx{c/n}
by deleting all edges $ij$ with $i+j>n$.

The operator $\tk$ is compact, and it easy to see that it has
eigenvalues $(-1)^k\go_k\qi$ and eigenfunctions
$\cos(\go_k x)$, with $\go_k=(k+1/2)\pi$, $k=0,1,\dots$. Hence
$\norm{\tk}=2/\pi$ and the critical value is $c_0=\pi/2$.
\refT{T5} shows that at the critical value we have
$c_0\rho'_+(c_0)=3/2$.
\end{example}

\begin{example}
{\bf I.i.d. vertices.}
  \label{EU}
For any ground space $(\sss,\mu)$, we can obtain a vertex space by 
taking $x_1,\dots,x_n$ to be
\iid{} random points in $\sss$ with
distribution $\mu$. (This has been proposed by S\"oderberg
\cite{Sod1}.) 
In this case
\begin{equation*}
  \begin{split}
 \E e\bigpar{\gnk}
&=
\frac{n(n-1)}2 \iint_{\sss^2} \frac{\kk(x,y)\bmin n}{n} \dd\mu(x)\dd\mu(y)
\\&
< \frac n2 \iikxy.
  \end{split}
\end{equation*}
Hence, by \refL{LE} below, \eqref{t1b} always holds, and to verify
that a kernel $\kk$ is graphical,
we only have to check conditions \ref{T1a} and \ref{T1c} in \refD{Dg1}.
Similarly, for a sequence of kernels,
\eqref{t2b} holds provided $\iint\kk_n\to\iint\kk$.
\end{example}

\begin{example}
{\bf Poisson process graph.}
  \label{EPo}
For any generalized ground space $(\sss,\mu)$ and any $\gl>0$, let 
${\bfx}_\gl=(x_1,\ldots,x_{\nvl})$ be the points of a Poisson
process on $\sss$ with intensity measure $\gl\mu$.  In other words,
$\nvl$ has a Poisson distribution $\Po(\gl\mu(\sss))$, and, given
$\nvl$, the points $x_i$ are \iid{} as in \refE{EU}. Then
$(\sss,\mu,({\bfx}_{\gl}))$ is a generalized vertex space. 
Here, it is
natural to write $\gl$ rather than $n$ for an element of the index set
$I=(0,\infty)$. 
Note that \eqref{nunA2} holds because
$\gl\nu_\gl(A)\sim\Po(\gl\mu(A))$.
This is the
canonical example of a generalized vertex space, and one of the main
reasons for allowing a random number of vertices.

Let $\kk$ be a kernel on $(\sss,\mu)$, so, given ${\bfx}_\gl$, the
edge probabilities 
in the graph $G^\vxs(\gl,\kk)$ are given by
\[
 \pij =\min\bigset{\gk(x_i,x_j)/\gl,1},
\]
for $1\le i<j\le \nvl$.
As in \refE{EU}, \eqref{t1b} always holds. To see this, let
$\vxs'=(\sss,\mu',\yss)$, where $\mu'=\mu/\mu(\sss)$ 
is the normalized version of $\mu$, and $\ys$ consists of $n$ \iid
points of $\sss$ 
chosen with distribution $\mu'$.
Given that $\nvl=n$, the distribution of $G^\vxs(\gl,\kk)$
is exactly that of $G^{\vxs'}(n,(n/\gl)\kk)$. In particular, as
\eqref{t1b} holds for the latter graph, 
\begin{align*}
  \E\bigpar{e(G^\vxs(\gl,\gk))\mid \nvl=n}
 &\sim \frac{n}{2} \iint_{\sss^2}\frac{n\kk}{\gl} \dd\mu'(x)\dd\mu'(y) \\
 &= \frac{n^2}{2\gl\mu(\sss)^2} \iikxy.
\end{align*}
As $\gl\to\infty$ we have $\E(\nvl^2)\sim(\gl\mu(\sss))^2$, and
\eqref{t1b} follows. 
Hence, as in \refE{EU}, a kernel $\kk$ on $(\sss,\mu)$ is graphical
on $\vxs$ if and only if conditions \ref{T1a} and \ref{T1c} of \refD{Dg1} hold.

In this Poisson process example, it is easy to see that allowing a random 
number of vertices makes the model only superficially more
general. Indeed, renormalizing so that $\mu(\sss)=1$, since
$\nvl\sim\Po(\gl)$, we can regard $\nvl$ as a random function of $\gl$,
which is increasing (a Poisson process), and then $\nvl/\gl\asto1$ as
$\gl\to\infty$. It follows that if we condition on the process
$\nvl(\gl)$, then \refT{T2} applies a.s.\ to the corresponding
graphs $G^{\vxs'}(\nvl,(\nvl/\gl)\kk)$ on the (ungeneralized) vertex
space $\vxs'$. Thus, conditioning on $\nvl(\gl)$,
\begin{equation}
  \label{epo}
\gl\qi C_1(G^\vxs(\gl,\kk)) \pto \rho(\kk)
\qquad\text{as } \gl\to\infty
\end{equation}
holds a.s.
It follows that \eqref{epo} holds unconditionally too.
Other properties can be treated similarly. We shall see later, in
\refSS{SSgvs}, 
that all our results for generalized
vertex spaces can be reduced to the vertex space case.
\end{example}

\begin{example}
{\bf Edge percolation.}
  \label{Eep}
Let $\kk$ be an irreducible graphical kernel on a (generalized) vertex
space 
$\vxs$ with $\norm{\tk}>1$, and let $0< p\le 1$. Independently of
everything else, keep each edge in \gnk{} with probability $p$ and
delete it with probability $1-p$. Denote the resulting
graph by $\gnkp$.

This random graph $\gnkp$ is nothing but $\gnx{\tkkn}$, where
\begin{equation*}
\tkkn(x,y)\=p\bigpar{\gk(x,y)\bmin n}.
\end{equation*}
Clearly, $x_n\to x$ and $y_n\to y$ imply $\tkkn(x_n,y_n)\to
p\kk(x,y)$, provided $(x,y)$ is a point of continuity of $\kk$.
Furthermore, $\frac1n\E e\bigpar{\gnkp}=\frac pn\E e\bigpar{\gnk}\to p\hiik$.
Hence,
$(\tkkn)$ is a graphical sequence with limit $p\kk$, so
\refT{T2} applies with $\kk$ replaced by $p\kk$, and
\begin{equation*}
  n\qi C_1\bigpar{\gnkp}
\pto
\rho(p\kk).
\end{equation*}
In particular, $\gnkp$ has \whp{} a component of order $\Theta(n)$
if and only if $\norm{\Tx{p\kk}}>1$, \ie{}, if $p>\norm{\tk}\qi$.
Thus, as expected, we obtain the same threshold for
edge percolation in $\gnk$ (meaning that there remains a giant
component) as for the existence of a giant component in
$\gnx{p\kk}$; see \refC{C1}.

Of course, the same conclusions follow if we start with the more general
setting of \refD{Dg2}.
\end{example}

\begin{example}
{\bf Vertex percolation.}
  \label{Evp}
Again, let $\kk$ be an irreducible graphical kernel on a vertex space 
$\vxs=(\sss,\mu,(\xs))$ with
$\norm{\tk}>1$, and let $0< p\le 1$.
Independently of everything else,
keep each vertex  in \gnkx{} with
probability $p$ and delete it with probability $1-p$.
Denote the resulting graph by $\gnkpp$.
This graph is again an instance of our model with a generalized vertex space.
Indeed, writing $\ys$ for the subsequence of $\xs$ corresponding
to the vertices that were not deleted,
$\vxs'=(\sss,p\mu,(\ys))$ is a generalized vertex space, and
$\gnkpp$ has exactly the distribution of $G^{\vxs'}(n,\kk)$. 
Since the kernel $\kk$ is graphical on $\vxs$, and
$\E \bigpar{e\bigpar{\gnkpp}} = p^2\E\bigpar{e\bigpar{\gnkx}}$,
the kernel $\kk$ is also graphical on $\vxs'$, so our results
apply to $G^{\vxs'}(n,\kk)$ and hence to $\gnkpp$.

Here, one must be a little careful with the normalization: the norm
of $\tk$ defined with respect to $(\sss,p\mu)$ is $p$ times $\norm{\tk}$,
the norm defined with respect to $(\sss,\mu)$.
In particular, \refT{T2} tells us that $\gnkpp$ has \whp{} a component of order
$\Theta(n)$ if and only if $p\norm{\Tx{\kk}}>1$, \ie{} if
$p>\norm{\tk}\qi$. We thus obtain the same threshold for vertex
percolation in $\gnk$ as for edge percolation in
\refE{Eep}.

Once again, we could have started with the setting of \refD{Dg2};
we could also have started with a generalized vertex space.
\end{example}

Note that we can obtain the Poisson graph
$\tgl(\kk)$ in \refE{EPo} as a limit of the vertex percolation
model $\gnkpp$ in \refE{Evp} if we
take $p=\gl/n$ and let \ntoo.

\medskip
Our next example shows that even in the supercritical, irreducible case, the
second largest component may be rather large -- certainly much larger than
$O(\log n)$ as in the Erd\H os--R\'enyi case.

\begin{example} \label{E2large+}
{\bf Large second component.}
Let $\sss=\set{1,2,3,\dots}$ with $\mu\set{k}=2^{-k}$, and let
$x_1,\dots,x_n$ be \iid{} random points in $\sss$ with
distribution $\mu$. Let $(\eps_k)_1^{\infty}$ be a
sequence of positive numbers tending to zero, to be
chosen below. Set $\kk(k,k)=2^{k+1}$ for $k\ge1$,
$\kk(1,k)=\kk(k,1)=\eps_k$ for $k\ge2$, and $\kk(i,j)=0$
otherwise.
Note that $\kk\in L^1(\sss\times\sss,\mu\times\mu)$;
as noted in \refE{EU}, from
our choice of $x_i$ it follows that $\kk$ is graphical
on $(\sss,\mu,\xss)$.

For each $k\ge1$, the graph $\gnk$ contains $n_k\sim\Bi(n,2^{-k})$
vertices of type $k$, forming a random subgraph $H_k$ which has
the distribution of the Erd\H os--R\'enyi graph $G(n_k,2^{k+1}/n)$.
Each potential edge between $H_1$
and $H_k$ is present with probability $\eps_k/n$.
Note that $n_k=\bnk+O_p(\bnk\qh)$, where $\bnk=\E n_k=n/2^k$, and thus each
$H_k$ is (\whp) supercritical.
In particular, \whp{} $C_1(\gnk)\ge C_1(H_1)\ge cn$ for some $c>0$, so
$\gnk$ is supercritical.

Let $k_n\to\infty$ with $\log_2 n-k_n\to\infty$, so that $\bar
n_{k_n}\to\infty$. Let us choose the $\eps_k$ so that
$\eps_{k_n}\le n^{-2}$. Then the expected number of edges between
$H_1$ and $H_{k_n}$ is $\E(n_1n_{k_n})\eps_{k_n}/n \le
n\eps_{k_n}\to0$, so \whp{} $H_{k_n}$ is isolated in $\gnk$. As
$n_{k_n}=\bar n_{k_n}+O_p(\bar n_{k_n}\qh)$, we may couple the
$\gnk$ for different $n$ so that
\begin{equation}\label{nov}
 n_{k_n} = \bar n_{k_n} +O(\bar n_{k_n}\qh)
\end{equation}
holds a.s. (Here the implicit constant is random.) We may then condition
on $n_{k_n}$, assuming that $n_{k_n}$ is deterministic, and that \eqref{nov}
holds.

Clearly, $H_{k_n}$ is a uniform Erd\H os--R\'enyi random graph
$G(n_{k_n},2^{k_n+1}/n)$.
As $2^{k_n+1}/n\sim 2/n_{k_n}$, this graph is supercritical (for large $n$),
and has a largest component of order $(c+o_p(1))n_{k_n}$ for some constant $c$.
Thus,
$$
C_1\bigpar{H_{k_n}}
=(c+o_p(1)) n_{k_n}
=(c+o_p(1))\bar  n_{k_n}
=(c+o_p(1))n/2^{k_n}.
$$
Given any function $\go(n)$ with $\go(n)=o(n)$, we can choose $k_n$
so that $2^{k_n}\go(n)/n\to0$;
it follows that \whp{}
$C_2(\gnk)\ge C_1(H_{k_n})
>(c/2)n/2^{k_n}
>\go(n)$.
Thus, the $o_p(n)$ bound in \refT{T2b} is best possible.
\end{example}

The final example in this section shows that when $\norm{\tk}=\infty$,
the ratio of the number of edges to the number of vertices in the
giant component of $G(n,c\kk)$ need not tend to 1 as $c\to 0$. In fact,
it may even tend to $\infty$.

\begin{example} \label{E2large2}
{\bf Dense giant component.}
Let $\sss$, $\mu$ and $\xs$ be as in \refE{E2large+}, and let
$\kk(1,k)=\kk(k,1)=1$ for $k\ge1$, $\kk(k,k)=4^k/k^2$ for $k\ge2$,
and $\kk(i,j)=0$ otherwise.
Again $\kk\in L^1$, so $\kk$ is graphical.
Let $c>0$ be small but fixed and consider
$\gnxx{c\kk}$. Let $k_0$ be the smallest integer such that
$2^{k_0}/k_0^2>1/c$; taking $c$ small enough, we may assume that
$k_0\ge10$.

Using the notation of \refE{E2large+},
if $k\ge 2$, then $H_k$ forms a random subgraph of the type
$G(n_k,c4^k/(k^2n))$. Since $n_kc4^k/(k^2n)\approx c2^k/k^2$,
this
subgraph is a supercritical Erd\H os--R\'enyi graph if $k\ge k_0$, and
if $k\ge k_0+1$, classical results show that \whp{} $H_k$ contains a
component of order $\Theta(n_k)=\Theta(2^{-k}n)$ with
$\Theta(n_k^2c4^k/(k^2 n) )=\Theta(n c/k^2)$ edges; throughout
this example the implicit constants in $\Theta(\cdot)$
and $O(\cdot)$ notation
do not depend on $c$.
Each of these components is \whp{} of order $n$, so they are
subsets of the giant component of \gnkx. Summing over
$k=k_0+1,\dots,2k_0$, the giant component thus has at least
$\Theta(n c/k_0)$ edges,
so $\edgeno(c\kk)=\Omega(c/k_0)$; see \refT{Tedges}.

To bound the number of vertices in the giant component, condition on
$x_1,\dots,x_n$ and say that a
vertex of type $k$ is \emph{light} if $k\le k_0-3$, and \emph{heavy} otherwise.
The total number of heavy vertices is $O(n2^{-k_0})$ \whp.
Furthermore, it is easy to check that if $c$ is small enough,
then the expected number of edges to light vertices 
{}from each heavy vertex is at most $1/2$, as is
the expected degree of each light vertex.
Each light vertex in the giant component has to
be connected to some heavy vertex by a path whose other vertices all
are light. As the expected number of such
paths starting at a given heavy vertex is at most $\sum_{l\ge1} (1/2)^l=1$,
the expected number of light vertices in the giant component is
$O(n2^{-k_0})$ too. Hence, the number of vertices in the
giant component is $O_p(n2^{-k_0})$,
so $\rho(c\kk)=O(2^{-k_0})=O(c/k_0^2)$.
Consequently, $\edgeno(c\kk)/\rho(c\kk)=\Omega(k_0)=\Omega(\log(1/c))$,
as $c\to0$. In particular,
$\edgeno(c\kk)/\rho(c\kk)\to\infty$ as $c\downto c_0=0$;
see \refR{Retrans}.
\end{example}

\section{Branching process lemmas}\label{S:branch1}

In this section and the next
we study the Poisson branching processes $\bpkx$
and $\bpk$ defined in \refSS{SS:branch}, and their survival
probabilities. These turn out to be given by the solutions
to a certain non-linear
functional equation \eqref{eq}. Let us briefly recall some
definitions.

Let $(\mu,\sss)$ be a (generalized) ground space. 
The branching process $\bpkx$ is a multi-type
Galton--Watson branching processes with type space $\sss$:
a particle of type $y\in \sss$ is replaced in the next generation by
its `children', a set of
particles whose types are distributed as a Poisson process on $\sss$
with intensity
$\kk(y,z)\dd\mu(z)$. The zeroth generation of $\bpkx$ consists of a
single particle of type $x$. Note that the distribution of $\bpkx$ is unaffected
if we multiply $\kk$ by a constant and divide $\mu$ by the same
constant; thus, we may assume without loss of generality that $\mu(\sss)=1$.
We shall make this assumption throughout this section.
In this normalized case, the branching process \bpk{} is just
the process \bpkx{} started with a single particle whose (random) type
is distributed according to the probability measure $\mu$.

Here, we have no need for the metric or topological structure of $\sss$; in
this section and the next,
$\sss$ can be any measurable space equipped with a
probability measure $\mu$. We assume, as before, that the kernel $\kk$ is a
measurable symmetric non-negative function on $\sss^2$.
We shall also assume that $\kk\in L^1(\sss\times\sss,\mu\times\mu)$,
\ie{}, that $\iint\kk<\infty$.

Let us recall our notation for the survival probabilities of
particles in $\bpkx$. We write $\rhok(\kk;x)$ for the probability
that the total population consists of exactly $k$ particles, and
$\rhogek(\kk;x)$ for the probability that the total population
contains at least $k$ particles. Furthermore,
$\rho(\kk;x)$ is the probability that the branching
process survives for eternity.

We write  $\rhok(\kk)$, $\rhogek(\kk)$ and $\rho(\kk)$ for the
corresponding probabilities for $\bpk$, so that, e.g.,
$\rhok(\kk)=\int_\sss\rhok(\kk;x)\dd\mu(x)$.

We start with a trivial observation that will
enable us to eliminate certain pathologies.

\begin{lemma}
  \label{Lae}
If\/ $\kk=\kk'$ a.e., then
$\rho(\kk;x)=\rho(\kk';x)$
and $\rhogek(\kk;x)=\rhogek(\kk';x)$ hold for \aex{} $x$;
hence
$\rho(\kk)=\rho(\kk')$ and $\rhogek(\kk)=\rhogek(\kk')$.
\end{lemma}

\begin{proof}
There is a measure zero set
$N\subset\sss$ such that if $x\notin N$,
then $\kk(x,y)=\kk'(x,y)$ for \aex{} $y$. It follows that if we start the
processes $\bpk$ and $\bpx{\kk'}$
at the same $x\notin N$, the processes will be identical in
distribution.
Hence $\rho(\kk;x)=\rho(\kk';x)$ and
$\rhogek(\kk;x)=\rhogek(\kk';x)$ for all $x\notin N$,
and the result
follows from \eqref{rho}.
\end{proof}

For the sake of convenience, in this section we impose one more
assumption on $\kk$, namely that
\begin{equation}
  \label{bb1}
\int_{\sss} \kk(x,y)\dd\mu(y)<\infty
\end{equation}
for \emph{every} $x\in \sss$. This assumption loses no generality,
as \eqref{bb1} holds for \aex{} $x$, since $\iint\kk<\infty$.
Writing $N$ for the measure zero set of $x$ such that \eqref{bb1}
does not hold, define $\kkb$ by setting $\kkb(x,y)=0$ if $x\in N$
or $y\in N$, and $\kkb(x,y)=\kk(x,y)$ otherwise. Then $\kk=\kkb$
\aex{}, so by \refL{Lae} we have $\rho(\kkb)=\rho(\kk)$ and so on.

All the assumptions above apply to all the kernels considered
below, denoted $\kk_1$, $\kk'$, etc. In this section, unless
explicitly stated, we do \emph{not} assume that $\kk$ is
irreducible.

\begin{remark}\label{R:bb1}
Condition \eqref{bb1}
means that a particle of type $x$ has a finite number of children
in the branching process.
As we are assuming \eqref{bb1} for all $x$, a particle survives for
eternity (has descendants in all future generations) if and only if
it has infinitely many descendants.
In other words,
$\rho(\kk;x)=\rhox{\infty}(\kk;x)$.
\end{remark}

\begin{remark}
Our process is a very special branching process
since we assume that the children of a particle are distributed
according to a Poisson process; in particular, the number of
children has a Poisson distribution. Other branching processes,
and other functional equations, appear when studying random graphs
with dependencies between edges, as in
\cite{robust,giant2,Rsmall}, but will not be considered here.

Note also that even with the Poisson assumption, our processes are
special.
For the branching process, there is no reason to assume $\kk$ to be
symmetric; 
moreover, $\mu$ may be any $\gs$-finite measure, and the
hypothesis $\kk\in L^1$ could be weakened to \eqref{bb1} for
\aex{} $x$ (or perhaps removed completely). We shall, however,
consider only the special case just defined; this will be useful
in the proofs. We have not yet investigated to what
extent the results generalize and remark only that in
non-symmetric situations, the norm $\norm{\tk}$ should be replaced
by the spectral radius.
\end{remark}

There is an abundant literature on branching processes with
different types; see, for example, the book by Mode \cite{Mode}.
However, we have not found the results we need in the generality
required here, so for the sake of completeness we give full
proofs, although the results are only minor
extensions of known results; see, for example, \cite[Chapter
6]{Mode}.

We start with the connection between our branching process and the
operator $\Phik$ defined in \eqref{Phik}.

\begin{lemma}\label{Lpoi}
Consider the random offspring of a single particle of type $x$; let
$N$ be the number of children, and
denote their types by $(\xi_i)_{i=1}^N$.
If\/ $g$ is a measurable function on $\sss$ with $0\le g\le 1$, then
\begin{equation}\label{poi}
  \E\prod_{i=1}^N\bigpar{1-g(\xi_i)}
=e^{-(\tk g)(x)}
=
1-(\Phik g)(x).
\end{equation}
\end{lemma}

\begin{proof}
This is a standard formula for Poisson processes; see, for
example, Kallenberg~\cite[Lemma 12.2(i)]{Kall}, taking
$f=-\ln(1-g)$. 
For completeness, we include the simple proof.
Indeed, let $\nu=\nu_x$ be the measure
defined by $\dd\nu(y)=\kk(x,y)\dd\mu(y)$. Then $N\sim\Po(\nu(\sss))$
and, given $N$, the types $\xi_i$ of the children are \iid{} with
the renormalized distribution $\nu'=\nu/\nu(\sss)$.
Hence, given $N$, the conditional expectation of
$\prod_{i=1}^N\bigpar{1-g(\xi_i)}$ 
is just
\begin{multline*}
  \prod_{i=1}^N\E\bigpar{1-g(\xi_i)} = \bigpar{\E(1-g(\xi_1))}^N \\
 = \left(1-\int_\sss g(y)\dd\nu'(y)\right)^N
 = \bigpar{1-(\tk g)(x)/\nu(\sss)}^N.
\end{multline*}
Using $\P(N=n)=e^{-\nu(\sss)} \nu(\sss)^n/n!$ and taking the expectation, the result follows.
\end{proof}

Our next aim is to study the fixed points of $\Phik$, i.e., the
solutions of the equation
\begin{equation}
  \label{eq}
f=\Phik f\=1-e^{-\tk f},
\end{equation}
where $f$ is a non-negative function on $\sss$.

\begin{remark}
  \label{Rae}
If $f=g$ \aex, then $\Phik f=\Phik g$. In particular, if $f=\Phik
f$ \aex, then $\Phik f=\Phik(\Phik f)$; thus, if $f$ satisfies
\eqref{eq} \aex, then there is a solution $\bar f$ to \eqref{eq}
(\viz{} $\Phik f$) such that $f=\bar f$ a.e. This shows that it
makes no essential difference if we require \eqref{eq} to hold
only \aex{} (which might be natural from an $L^2$ perspective). We
shall, however, find it convenient to interpret \eqref{eq} and
similar relations as holding everywhere unless we explicitly state
otherwise. Similarly, if $\kk=\kkb$ a.e., then for any solution
$f$ to $f=\Phik f$ there is a unique $\bar f$ with the properties
that $\bar f=\Phix{\kkb}\bar f$ and $\bar f=f$ a.e.
\end{remark}

Note that $\Phik$ is monotone: if $0\le f\le g$  \aex{} then $\tk
f\le \tk g$ and thus $\Phik f\le \Phik g$.

In the lemma below, 1 denotes the function with
constant value 1.

\begin{lemma}\label{L1b}
  \begin{thmenumerate}
\item
For $m\ge0$ the probability that a particle of type $x$ has
descendants in at least $m$ further generations is
$(\Phik^m1)(x)$.
\item
As \mtoo, $(\Phik^m 1)(x)\downto\rho(\kk;x)$.
\item
The function $\rho_\kk=\rho_\kk (x)=\rho (\kk; x)$ is a solution
of \eqref{eq}, i.e., satisfies $\Phik\rho_\kk=\rho_\kk$.
\item
The function $\rho_\kk$ is the maximum solution of \eqref{eq}: if
$f$ is any other solution,  then $\rho_\kk(x)\ge f(x)$ for every
$x$.
  \end{thmenumerate}
\end{lemma}

\begin{proof}
  \pfitem{i}
Let
$g_m(x)$ be this probability. Then, with $g=g_m$,
the left-hand side of
\eqref{poi} is the probability that none of the
children of $x$
has descendants in
at least $m$ generations,
\ie{}, the probability $1-g_{m+1}(x)$ that $x$ does not have descendants
in $m+1$ generations. Thus $g_{m+1}=\Phik g_m$, and
the result follows by induction, since $g_0(x)=1$.

  \pfitem{ii}
An immediate consequence of (i).

  \pfitem{iii}
This follows by the same argument as (i) (and is also a
consequence of (ii) and dominated convergence).

   \pfitem{iv}
Suppose that $f$ is a solution of \eqref{eq}. Then
$f=\Phik f\le1$, and thus $f=\Phik^m f\le \Phik^m1$ for every
$m$. Hence, $f\le\rho_\kk$ follows from (ii).
\end{proof}

\begin{remark}
  If we do not impose \eqref{bb1} for all $x$, then
 (iii), i.e., $(\Phik\rho_\kk)(x)=\rho_\kk(x)$,
 could fail for $x$ in the measure zero set for which \eqref{bb1} does not
  hold. This is because a particle of type
$x$ has infinitely many children, which may have finite but unbounded
  lines of descendants; for an example, take $\sss=(0,1]$ and $\kk=1$
    except that $\kk(x,1)=\kk(1,x)=1/x$.
\end{remark}

We continue to study the functional equation \eqref{eq}.

\begin{lemma}\label{LPhi}
Suppose that $f\ge0$ with $f=\Phik f$. Then
\begin{romenumerate}
  \item \label{lphia}
$0\le f<1$;
\item \label{lphid}
$\tk f\ge f$, with strict inequality when $f(x)>0$;
\item \label{lphib}
$\tk f\le f/(1-f)$, with strict inequality when $f(x)>0$;
\item \label{lphic}
if $\kk$ is irreducible, then
either $f=0$ everywhere or $f>0$ \aex
\end{romenumerate}
\end{lemma}

\begin{proof}
\pfitem{i} We have $f(x)=1-e^{-(\tk f)(x)}\le1$. Hence,
$(\tk f)(x) \le (\tk 1)(x)=\int_{\sss} \kk(x,y)\dd\mu(y)
<\infty$, where the second inequality is just our assumption
\eqref{bb1}. Therefore, $f(x)=1-e^{-(\tk f)(x)}<1$.
\pfitem{ii} This is immediate from $f={1}-e^{-\tk f}\le \tk f$, with
equality only when $\tk f=0$.
\pfitem{iii} We have $e^{-\tk f}=1-f$, and thus, as $f<1$,
\begin{equation*}
\tk f \le e^{\tk f}-1
=\frac1{1-f} -1
=\frac {f}{1-f},
\end{equation*}
with equality only when $\tk f=0$.
\pfitem{iv}
Let $A\=\set{x\in \sss:f(x)=0}$.
For $x\in A$, $(\Phik f)(x)=f(x)=0$, and thus $(\tk f)(x)=0$.
Hence $\gk(x,y)=0$ for \aex{} $y\notin A$.
Consequently, $\gk=0$ \aex{} on $A\times(\sss\setminus A)$, which by
\eqref{t1ay} implies $\mu(A)=0$ or $\mu(A)=1$.
In the latter case, $f=0$ \aex, and thus $f=\Phik f=0$.
\end{proof}

In the next two lemmas we consider irreducible $\kk$.

\begin{lemma}\label{L2}
Suppose that $\kk$ is irreducible.
Suppose further that $f=\Phik f$ and $g=\Phik g$ with $0\le f\le g$.
Then either $f=0$ or $f=g$.
\end{lemma}

\begin{proof}
By \refL{LPhi}\ref{lphic} we may assume that $f>0$ a.e.

Let $h=(g-f)/2\ge0$; thus $f+h=(f+g)/2$. The function $t\mapsto
1-e^{-t}$ is strictly concave;
in particular, 
$1-e^{-(t+u)/2}\ge \frac{1}{2}\left((1-e^{-t})+(1-e^{-u})\right)$.
Hence,
\begin{equation}\label{q1}
  \begin{split}
 \Phik\parfrac{f+g}{2}
& =1-e^{-\tk ((f+g)/2)} =1-e^{-(\tk f+\tk g)/2}
\\&
\ge \tfrac12\left((1-e^{ -\tk f}) + \bigpar{1-e^{-\tk g}}\right)
=\tfrac12(f+g) =f+h,
  \end{split}
\end{equation}
with strict inequality at every point where $f<g$ and thus
$\Phik f<\Phik g$ and $\tk f<\tk g$.
On the other hand,
\begin{equation}\label{q2}
\begin{split}
 1-\Phik\parfrac{f+g}{2}
&
= e^{-\tk (f+h)}
= e^{-\tk f} e^{-\tk h}
=(1-f)e^{-\tk h}
\\&
\ge
(1-f)(1-\tk h).
  \end{split}
\end{equation}
Combining \eqref{q1} and \eqref{q2}, we find that
\begin{equation*}
  (1-f)(1-\tk h)
\le 1-(f+h)=1-f-h
\end{equation*}
and thus
\begin{equation}\label{q3}
  (1-f)\tk h
\ge h
\end{equation}
with strict inequality when $g>f$.

Suppose now that $g>f$ on a set of positive measure. Then,
inequality \eqref{q3}, the fact that $f>0$ \aex{}, and
\refL{LPhi}\ref{lphib} imply that
\begin{equation}\label{q4}
  \int_\sss f \,\tk h \dd\mu
>
  \int_\sss f \frac{h}{1-f} \dd\mu
=
  \int_\sss h \frac{f}{1-f} \dd\mu
\ge
  \int_\sss h \,\tk f \dd\mu.
\end{equation}
Note that the integrals above are finite because $\kk\in L^1$ and
$f,h\le1$. However, as $\gk$ is symmetric, $\tk$ is a symmetric
operator, and so $\int_\sss f \,\tk h \dd\mu = \int_\sss h \,\tk f
\dd\mu$, contradicting \eqref{q4}. This shows that $g=f$ \aex{}
and thus $f=\Phik f=\Phik g=g$.
\end{proof}

\begin{lemma}\label{L3}
Suppose that $\kk$ is irreducible.
Then $f=0$ and $f=\rho_\kk$ are the only solutions to \eqref{eq};
these solutions may coincide.
\end{lemma}

\begin{proof}
By parts (iii) and (iv) of \refL{L1b}, the function $\rho_\kk$ is
a solution of \eqref{eq}, and  $0\le f\le \rho_\kk$ for every
solution $f$ of \eqref{eq}. The result follows by \refL{L2}.
\end{proof}

It remains to decide whether $\rho_\kk=0$ or not.
Recall that $\norm{\tk }$ is defined in \eqref{tnorm}.
We shall show that $\rho_\kk=0$ if and only if $\norm{\tk}\le 1$.

\begin{lemma}  \label{L4}
If\/ $\norm{\tk }\le1$, then $\rhokk=0$.
\end{lemma}

\begin{proof}
Suppose that $f$ is a solution of \eqref{eq}, and that
we do not have
$f=0$ a.e.  \refL{LPhi}\ref{lphid} implies that $\tk f\ge f$,
with $\tk f>f$ on a set of positive measure, and hence that
$\normll{\tk f}>\normll{f}$, contradicting $\norm{\tk }\le1$.
Consequently, if $f$ is a solution of \eqref{eq}, then $f=0$ \aex,
and thus $f=\Phik f=0$, so the only solution is $f=0$. In
particular, $\rho_\kk=0$ since $\rho_\kk$ is a solution by
\refL{L1b}.
\end{proof}

It remains to show that if $\norm{\tk}>1$, then $\rho_\kk$ is not
identically zero. We proceed
in several steps.

\begin{lemma}
  \label{L1c}
If\/ $f\ge0$ and $\Phik f\ge f$, then $\Phik^m f \upto g$ as
\mtoo, for some $g\ge f\ge0$ with $\Phik g=g$.
\end{lemma}

\begin{proof}
By induction, $f\le \Phik f\le \Phik^2 f\le\dots$. Since
$0\le \Phik^m f\le1$, the limit $g(x)\=\limm (\Phik^m f)(x)$ exists for
every $x$, and
$g\ge0$.
Monotone convergence yields
\begin{equation*}
  (\tk g)(x)=\limm \int_\sss\gk(x,y)(\Phik^m f)(y)\dd\mu(y)
= \limm (\tk (\Phik^m f))(x)
\end{equation*}
and thus
\begin{equation*}
  (\Phik g)(x)= \limm (\Phik(\Phik^m f))(x)=g(x).
\end{equation*}
\vskip-\baselineskip
\end{proof}

\begin{lemma}
  \label{L5}
If there is a bounded function $f\ge0$, not \aex{} $0$, such that
$\tk f\ge(1+\gd)f$ for some $\gd>0$, then
$\rho_\kk>0$ on a set of positive measure.
\end{lemma}

\begin{proof}
Let $M=\sup f<\infty$. Fix $\eps>0$ with $(1-M\eps)(1+\gd)\ge1$.
Since $-\log(1-x)\le x/(1-x)$ we have
\begin{equation*}
  -\ln(1-\eps f)
\le  \frac1{1-\eps M} \eps f
\le (1+\gd)\eps f
\le \eps \tk f,
\end{equation*}
and thus
\begin{equation*}
  \Phik(\eps f) = 1-e^{-\eps \tk f} \ge 1-(1-\eps f)=\eps f.
\end{equation*}
By \refL{L1c}, there exists a solution $g$ to $\Phik g=g$ with
$g\ge \eps f$, and thus $g$ not \aex{} $0$. By part (iv) of
\refL{L1b}, we have $\rho_\kk\ge g$.
\end{proof}

\begin{remark}
The proof of \refL{L5} shows that
$\rho_\kk \ge \frac{\gd}{1+\gd}\frac{f}{\sup f}$.
In particular, this immediately implies Theorem 10 of \cite{BJR}.
\end{remark}

We should like to find an eigenfunction of $\tk$ with eigenvalue
greater than~1, so that we can apply \refL{L5}. If the
Hilbert--Schmidt norm of $\tk$ (see \eqref{HSdef}) is finite, then
a standard result gives us such an eigenfunction.

\begin{lemma}\label{Lcomp}
If\/ $\normHS{\tk} <\infty$,
then $\tk$ is compact and has an eigenfunction $\psi\in L^2(\sss)$,
$\psi\ge0$, with
eigenvalue $\norm{\tk}$.

If, in addition, $\kk$ is irreducible, then
$\psi>0$ a.e.,
and every eigenfunction with
eigenvalue $\norm{\tk}$ is a multiple of $\psi$.
\end{lemma}

\begin{proof}
Suppose that $\normHS{\tk}<\infty$. It is well-known (see \eg{}
\cite[XIV.6, p. 202]{BBLinAnal}) that $\tk$ is then compact,
and so has an eigenfunction $g\in L^2$ with
eigenvalue of modulus $\lambda:=\norm{\tk}$. Then
\begin{equation*}
\tk |g| \ge
|\tk g|
=\gl|g|
\qquad\text{a.e.,}
\end{equation*}
and since $\norm{\tk}=\gl$ we must have $\tk|g|=\gl|g|$ a.e. Hence
$\psi\=|g|$ is an eigenfunction with eigenvalue $\gl=\norm{\tk}$.

Now suppose that $\kk$ is irreducible, with $\normHS{\tk}<\infty$,
and let $h$ be any (real)
function in $L^2$ with $\tk h=\gl h$ a.e. By the argument above,
$\tk|h|=\gl|h|$ a.e.\ holds as well.  Let
$A\=\set{|h|=0}$. Then $\tk|h|=\gl|h|=0$ \aex{} on $A$, so
  $\kk=0$ \aex{} on $A\times(\sss\setminus A)$ and \eqref{t1ay} yields
$\mu(A)=0$ or 1.
Hence either $h=0$ \aex{} or
$h\neq 0$ a.e.
In particular, taking $h=g$ we see that $\psi>0$ a.e.

Returning to a general $h$ satisfying $\tk h=\gl h$ a.e.,
as $\tk(|h|+h)=\gl(|h|+h)$ \aex\ by linearity,
we can apply the argument above to
$|h|+h$, deducing that either $h>0$ \aex{} or $h\le 0$ a.e.
Finally, applying this to $h-a\psi$, with $a$ chosen such that
$\int(h-a\psi)\dd\mu=0$, we see that $h-a\psi=0$ a.e.
\end{proof}
The second part of
\refL{Lcomp} will be needed only in \refS{Stransitionpf}.

After this preparation, it is easy to show that if $\norm{\tk}>1$
then $\rho_\kk>0$ on a set of positive measure.

\begin{lemma}\label{L6}
  If\/ $1<\norm{\tk}\le\infty$, then
$\rho_\kk>0$ on a set of positive measure. Thus
\eqref{eq} has at least one non-zero solution.
\end{lemma}

\begin{proof}
Since $\norm{\tk}>1$, there is function $f\in L^2$ with
$\normll{f}=1$ and $\normll{\tk f}>1$. As $\tk |f|\ge |\tk f|$,
we may assume that $f\ge0$. Let $T_N$ be the
integral operator on $\sss$ with the truncated kernel
$\gk_N(x,y)\=\gk(x,y)\bmin N$, $N\ge1$. By monotone convergence,
$T_N f\upto \tk f$ as $N\to\infty$, and thus $\normll{T_N f}\upto
\normll{\tk f}>1$. We can thus choose an $N$ such that
$\normll{T_Nf}>1=\normll{f}$, and thus $\norm{T_N}>1$. Set
$\gd=\norm{T_N}-1>0$.

Since the kernel $\gk_N$ is bounded and $\mu$ is a finite measure,
by \refL{Lcomp} $T_N$ has an eigenfunction $\psi\in L^2(\sss)$ with
$\psi\ge0$ and  
\begin{equation}
  \label{qg}
T_N \psi =\norm{T_N} \psi=(1+\gd)\psi.
\end{equation}
Since the kernel $\gk_N$ is bounded, it  follows that $T_N \psi$
is a bounded function. Indeed, $(T_N\psi)(x)\le N\int_S\psi\dd\mu
=N\normpx{\psi}{1}\le N\normll{\psi}<\infty$. 
From \eqref{qg} it follows
that $\psi$ is bounded.

Since $\gk\ge\gk_N\ge0$, we have, using \eqref{qg} again,
\begin{equation*}
\tk \psi \ge
T_N \psi
=(1+\gd)\psi,
\end{equation*}
and the result follows by \refL{L5}.
\end{proof}

The final lemma of this section will enable us to reduce 
the reducible case to the irreducible one.

\begin{lemma}
  \label{Lred}
Let $\kk$ be a symmetric measurable function on $\sss\times\sss$.
Then there exists a partition
$\sss=\bigcup_{i=0}^N \sss_i$
with $0\le N\le \infty$ such that each $\sss_i$ is measurable,
$\mu(\sss_i)>0$ for $i\ge 1$,
the restriction of
$\kk$ to $\sss_i\times\sss_i$ is irreducible for each $i\ge1$, and
$\kk=0$ \aex{} on
$(\sss\times\sss)\setminus\bigcup_{i=1}^N (\sss_i\times\sss_i)$.
\end{lemma}

Note that $\kk=0$ \aex{} on $\sss_0\times\sss_0$.

\begin{proof}
Let $\cG$ be the family of all measurable subsets $A\subseteq\sss$
such that $\kk=0$ \aex{} on $A\times(\sss\setminus A)$. It is
easily verified that $\cG$ is a $\gs$-field; thus $(\sss,\cG,\mu)$
is a finite measure space. Hence there exists a partition
$\sss=\bigcup_{i=0}^N \sss_i$ with $0\le N \le\infty$ and each
$\sss_i\in\cG$ such that for each $i\ge1$ the set $\sss_i$ is an atom in
$(\sss,\cG,\mu)$ with positive measure,
while $\sss_0$ is non-atomic, \ie, contains no atoms
with non-zero measure. (We allow $\sss_0=\emptyset$.)
Here `$\sss_i$ is an atom' means that if
$A\subseteq\sss_i$ with $A\in\cG$, then $\mu(A)=0$ or
$\mu(A)=\mu(\sss_i)$; this is equivalent to \eqref{t1ay}, so
$\kk$ is irreducible on $\sss_i\times\sss_i$ for each $i\ge1$.

Finally, since $\sss_0$ is non-atomic, for every positive integer $M$ there
exists a partition $\sss_0=\bigcup_{j=1}^M T_j$ with $T_j\in\cG$ and
$\mu(T_j)=\mu(\sss_0)/M$. Then $\kk=0$ \aex{} on $T_i\times T_j$ when
$i\neq j$, and thus
\begin{equation*}
  (\mu\times\mu)\bigset{(x,y)\in\sss_0\times\sss_0: \kk(x,y)\neq0}
\le
  \left(\mu\times\mu\right)\Bigpar{\bigcup_{j=1}^M (T_j\times T_j)}
=M\parfrac{\mu(\sss_0)}{M}^2.
\end{equation*}
Letting $M\to\infty$, we see that $\kk=0$ \aex{} on $\sss_0\times\sss_0$.
\end{proof}

\begin{remark}\label{Rallsolns}
One application of \refL{Lred}
is a generalization of Lemma \ref{L3} to arbitrary $\kk$.
With $\sss_i$ as in \refL{Lred}, let $J$ be the set of indices $i$
such that the restriction of the operator $\tk$ to $L^2(\sss_i)$ has
norm strictly greater than $1$.
Then there are $2^{|J|}$ solutions of \eqref{eq}, where
$0\le|J|\le\infty$: for every subset $J'\subseteq J$, there is exactly
one solution that equals $\rhokk$ \aex{} on $\bigcup_{i\in J'}\sss_i$ and
vanishes \aex{} elsewhere. (This is easily seen using the argument in
the proof of \refT{TappB} below.)
\end{remark}

\section{Branching process results}\label{S:branch2}

In this section we collect the branching process results we shall
use. These are all simple consequences of the lemmas in the
previous section.
In this section, $\kk$ will always be a kernel
on a measure space $(\sss,\mu)$, \ie, a symmetric non-negative
measurable function on $\sss\times\sss$.
Unless explicitly stated otherwise, $\mu$ will 
be a probability measure, i.e., $\mu(\sss)=1$.
We shall assume that $\kk\in L^1$; as noted in the previous section,
it follows that \eqref{bb1} holds \aex{} $x$. We do not assume
that \eqref{bb1} holds for every $x$ except when explicitly stated.

\begin{theorem}\label{T:GW}
Suppose that $\kk$ is a kernel on the space $(\sss,\mu)$,
that $\kk\in L^1$, and that \eqref{bb1} holds for every $x$. Then
the function $\rho_\kk$ defined by $\rho_\kk(x)=\rho(\kk;x)$ is
the maximum solution of \eqref{eq}. Furthermore:
\begin{romenumerate}
\item
If\/ $\|\tk \|\le1$, then $\rho(\kk;x)=0$ for every $x$, and
\eqref{eq} has only the zero solution.
\item
If\/ $1<\|\tk \|\le\infty$, then $\rho(\kk;x)>0$ on a set of positive
measure.
If, in addition, $\kk$ is irreducible, then $\rho(\kk;x)>0$ for \aex{} $x$,
and $\rho(\kk;x)$ is the only non-zero solution of \eqref{eq}.
\end{romenumerate}
In particular, $\rho(\kk)>0$ if and only if\/ $\|\tk \|>1$.
\end{theorem}

\begin{proof}
The first statement is just part (iv) of \refL{L1b}. The remaining
statements
follow directly from Lemmas \ref{L4}, \ref{L6},
\refL{LPhi}\ref{lphic} and \ref{L3}.
\end{proof}

The next result is essentially a restatement of \refT{T:GW}, in the setting 
of the results in \refS{Sresults}. Thus, $\mu$ will not necessarily
be a probability measure, and we 
shall not require that \eqref{bb1} holds; this makes very little
difference. This result gives the promised characterization of
$\rho(\kk;x)$ and $\rho(\kk)$ in terms of a functional equation,
in the full generality of the setting of \refT{T2}.

Recall that \eqref{Phik} defines
$\Phik$ only for non-negative functions; we thus consider only
non-negative solutions to \eqref{t1e} below.

\begin{theorem}\label{Trho}
Let $\kk$ be a kernel on a (generalized) ground space $(\sss,\mu)$, 
with $\kk\in L^1(\sss\times\sss,\mu\times\mu)$.
There is a (necessarily unique) maximum solution $\trhokk$ to
\begin{equation}
  \label{t1e}
\Phik(\trhokk)=\trhokk,
\end{equation}
i.e., a solution that pointwise dominates all other solutions.
Furthermore, $\rho(\kk;x)=\trhokk(x)$ for \aex\ $x$,
and
\begin{equation} \label{Phirho}
 \Phik(\rho_\kk)=\rho_\kk\quad a.e.,
\end{equation}
where the function $\rho_\kk$ is defined by $\rho_\kk(x)\=\rho(\kk;x)$.

If\/ $\norm{\tk}\le 1$, then $\trhokk$ is identically zero, and this
is thus the only solution to \eqref{t1e}.
If\/ $\norm{\tk}>1$, then $\trhokk$ is positive
on a set of positive measure. Thus $\rho(\kk)>0$ if and only
if $\norm{\tk}>1$.

If\/ $\norm{\tk}>1$ and $\kk$ is irreducible, then $\trhokk$ is the
unique non-zero solution to \eqref{t1e}, and $\trhokk=\rho_\kk>0$
a.e.
\end{theorem}

\refT{Trho} follows almost immediately from \refT{T:GW} and
\refL{Lae}.

\begin{proof}
Multiplying $\kk$ by a constant factor and dividing $\mu$ by the 
same constant factor does not affect the definition of the branching process
$\bpk(x)$. Hence, the function $\rho_\kk$ is not affected by this rescaling.
As the operators $\tk$ and $\Phik$ are also unchanged, we may assume without loss
of generality that $\mu(\sss)=1$.
As noted in \refS{S:branch1}, since $\kk\in L^1$ there is a
kernel $\kkb$ with $\kkb=\kk$ \aex{},
such that \eqref{bb1} holds for $\kkb$ for every $x$.
Applying \refT{T:GW} to the kernel $\kkb$, the result follows by
\refL{Lae} and \refR{Rae}.
\end{proof}

We now study monotonicity and continuity properties of
$\rho(\kk;x)$ and $\rho(\kk)$ when $\kk$ is varied.
For the rest of the section, we assume that $\mu(\sss)=1$. 
As usual, we
say that a sequence of functions $f_n$ \emph{increases (\aex)} to
a function $f$ if for every $x$ (\aex{} $x$) the sequence $f_n(x)$
is monotone increasing and converges to $f(x)$. As before, we
write $\rho_\kk$ for the function given by
$\rho_\kk(x)\=\rho(\kk;x)$. We start with a trivial lemma.

\begin{lemma}\label{L12}
  If\/ $\kk_1\le\kk_2$, then $\rhokx1\le\rhokx2$.
\end{lemma}

\begin{proof}
Immediate by coupling the branching processes.
\end{proof}

\begin{theorem}\label{TappB}
  \begin{thmenumerate}
\item
 Let $(\kappa_n)_1^{\infty}$ be a sequence of kernels on $(\sss,\mu)$
 increasing
 \aex{} to $\kappa$.
Then
$\rhokn\upto \rhokk$ for \aex{} $x$ and
$\rho(\kk_n)\upto\rho(\kk)$.
\item
 Let $(\kappa_n)_1^{\infty}$ be a sequence of kernels on $(\sss,\mu)$
 decreasing 
 \aex{} to $\kappa$.
Then
$\rhokn\downto \rhokk$ for \aex{} $x$ and
$\rho(\kk_n)\downto\rho(\kk)$.
 \end{thmenumerate}
\end{theorem}

\begin{proof}
  On the measure zero set where $\kk_n\not\to\kk$, redefine all $\kk_n$ and
  $\kk$ to be 0.
By \refL{Lae}, this does not affect the conclusions, so we may assume
$\kk_n\upto\kk$ or $\kk_n\downto\kk$ everywhere.
Similarly, we may assume that \eqref{bb1} holds for every $x$, for
each $\kk_n$ and for $\kk$.
 It suffices to prove
the conclusions for $\rhokk$:
the conclusions for $\rho(\kk)$
follow from \eqref{rho} and dominated
  convergence.

\pfitem{i}
We choose a partition $\sss=\bigcup_{i=0}^N \sss_i$ as in \refL{Lred},
and redefine $\kk$ and all $\kk_n$ to be 0 on
$(\sss\times\sss)\setminus\bigcup_{i=1}^N (\sss_i\times\sss_i)$;
this only changes the kernels on a set of measure zero,
so we may again apply \refL{Lae}.
Now $\rho_{\kk_n}=\rho_\kk=0$ on $\sss_0$.
We may consider each $\sss_i$, $i\ge1$, separately, and
we may thus assume without loss of generality that $\kk$ is
irreducible.
The only problem is that the restriction
of $\mu$ to $\sss_i$ does not have total mass $1$, but this is not a
real problem,
and can be handled
by renormalizing, \ie,
dividing the measure by $\mu(\sss_i)$ and
multiplying all kernels by the same factor; 
as remarked earlier,
this operation does
not affect the branching process.

We have shown that we may assume that $\kk$ is irreducible;
let us do so.
By \refL{L12}, if $m\le n$, then $\rhokm\le\rhokn$.
Thus $(\rhokn)$ 
is an increasing sequence of functions, all bounded
by 1, so the limit $\rhokq(x)\=\limn\rhokn(x)$ exists everywhere.
By monotone convergence,
\begin{equation*}
  \begin{split}
  (\tk  \rhokq)(x)
&
=\int_\sss\gk(x,y)\rhokq(y)\dd\mu(y)
=\limn \int_\sss\gk_n(x,y)\rhokn(y)\dd\mu(y)
\\&
=\limn   (\tkn \rhokn)(x),
  \end{split}
\end{equation*}
so $\Phik\rhokq=\limn \Phikn\rhokn=\limn\rhokn=\rhokq$.
Hence, by \refL{L3}, either $\rhokq=\rhokk$, and we are done, or
$\rhokq=0$.
In the latter case, each $\rhokn=0$, and thus, by \refL{L6},
$\norm{\tkn}\le1$.

Hence, if $f\in L^2$ with $f\ge 0$ and $\normll f\le1$, then
$\normll{\tkn f}\le 1$.
Monotone convergence shows that,  as \ntoo,
$\tkn f\upto \tk f$
and
$\normll{\tkn f}\upto\normll{ \tk f}$.
Consequently, $\normll{\tk f}\le1$ for each such $f$, and thus
$\norm{\tk}\le1$.
By \refT{T:GW}, $\rhokk=0$ in this case, so $\rhokk=\rhokq$ in this case
too.

\pfitem{ii}
This is similar.
Now $(\rhokn)$ 
is a decreasing sequence
of functions, and $\rhokq(x)\=\limn\rhokn(x)$ still exists everywhere.
By dominated convergence,
$  (\tk  \rhokq)(x)=\limn   (\tkn \rhokn)(x)$,
so $\Phik\rhokq=\limn \Phikn\rhokn=\limn\rhokn=\rhokq$.
In other words, $\rhokq$ satisfies
\eqref{eq}.
Furthermore, by \refL{L12} again, $\rhokn\ge\rhokk$, so $\rhokq\ge\rhokk$.
Since $\rhokk$ is a maximal
solution to \eqref{eq} by \refL{L1b},
$\rhokq =\rhokk$.
\end{proof}

\begin{theorem}\label{TappC}
  \begin{thmenumerate}
\item
 Let $(\kappa_n)_1^{\infty}$ be a sequence of kernels on $(\sss,\mu)$
 increasing 
 \aex{} to $\kappa$.
Then, for every $k\ge 1$,
$\rhogek(\kk_n;x)\upto\rhogek(\kk;x)$ for \aex{} $x$ and
$\rhogek(\kk_n)\upto\rhogek(\kk)$.
\item
 Let $(\kappa_n)_1^{\infty}$ be a sequence of kernels on $(\sss,\mu)$
 decreasing
 \aex{}
to  $\kappa$.
Then, for every $k\ge1$,
$\rhogek(\kk_n;x)\downto\rhogek(\kk;x)$ for \aex{} $x$ and
$\rhogek(\kk_n)\downto\rhogek(\kk)$.
 \end{thmenumerate}
\end{theorem}

\begin{proof}
As in the proof of \refT{TappB}, we may assume that
$\kk_n\upto\kk$ or $\kk_n\downto\kk$ everywhere, and that \eqref{bb1}
always holds.

\pfitem{i}
Let $\kk_0\=0$ and $\Delta\kk_n\=\kk_n-\kk_{n-1}$, $n\ge1$. 
The children of a particle of type $x$ are given by a Poisson
process with intensity
$\kk(x,y)\dd\mu(y)= \sum_n \Delta\kk_n(x,y)\dd\mu(y)$, which can be
represented as the sum of independent Poisson processes with
intensities $\Delta\kk_n(x,y)\dd\mu(y)$. We label the children in the
$n$th of these processes  by $n$, and give the initial `root' vertex label 0.
This gives a labelling of all
particles in the branching process $\bpk(x)$ (which starts with a single
particle of type $x$) such that the
subset of all particles that, together with all their ancestors, have
labels at most $n$ gives the branching process
$\bpx{\kk_n}(x)$.
Consequently (using this coupling of the processes), the family tree of
the initial particle in $\bpx{\kk_n}(x)$
will grow to its family tree in $\bpkx$ as \ntoo.
Hence $\rhogek(\kk_n;x)\upto\rhogek(\kk;x)$ and
$\rhogek(\kk_n)\upto\rhogek(\kk)$.

\pfitem{ii}
We may similarly label all particles in $\bpx{\kk_1}(x)$ with labels
\set{1,2,\dots,\infty} such that
$\bpx{\kk_n}(x)$ [$\bpkx$] consists of all particles that, together with their
ancestors, have labels at least $n$ [$\infty$].
By \refR{R:bb1},
a particle always has a finite number of children, so a particle
survives for eternity
if and only if it has infinitely many descendants.
By \refT{TappB} we have
\begin{equation}\label{pconv}
 \rho(\kk_m;x)\downto\rho(\kk;x)
\end{equation}
for \aex{} $x$.
Fix any $x$ for which \eqref{pconv} holds.
Writing $|\bpxx{\kk_n}|$ for the total population of the branching
process $\bpxx{\kk_n}$,
whenever $|\bpxx{\kk_m}|<\infty$ for some $m$,
we have $|\bpxx{\kk_n}|\downto
|\bpkx|$ as $n\to\infty$;
indeed, for large $n$ the entire processes $\bpx{\kk_n}(x)$ and
$\bpkx$ coincide.
{}From \eqref{pconv}, with probability $1$ either
$|\bpkx|=\infty$, in which
case $|\bpxx{\kk_n}|\ge |\bpkx|=\infty$ for all $n$, or there is an $m$
with $|\bpxx{\kk_m}|<\infty$, in which case
$|\bpxx{\kk_n}|=|\bpkx|$ for all large enough $n$.
Thus, the events $|\bpxx{\kk_n}|\ge k$ converge \aex{} to
$|\bpkx|\ge k$,
and  $\rhogek(\kk_n;x)\downto\rhogek(\kk;x)$.
\end{proof}

Suppose that $\kk$ is supercritical (\ie, that $\norm{\tk}>1$),
and assume for simplicity
that \eqref{bb1} holds for every $x$. Consider the branching
process \bpkx{} starting with a particle of type $x$, and classify
its children in the first generation according to whether they
have infinitely many descendants or not.
By the 
properties of Poisson processes, this exhibits the
children as the union of two independent Poisson processes with
intensities $\kk(x,y)\rhokk(y)\dd\mu(y)$ and $\kk(x,y)(1-\rhokk(y))\dd\mu(y)$
respectively, where the first litter consists of the children
with infinitely many descendants, or, equivalently,
those whose descendants live for ever.

The process $\bpkx$ eventually becomes extinct if and only if the
first litter is empty. It follows that if $\bpkxp$ denotes the
branching process \bpkx{} conditioned on extinction, then \bpkxp{} is
itself a multi-type Galton--Watson
branching process, where the set of children of a particle of type $z$ is
given by a Poisson process with intensity
$\kk(z,y)(1-\rhokk(y))\dd\mu(y)$.
This is another instance of the situation studied here, with $\mu$
replaced by $\hmu$ defined by $\dd\hmu(y)\=(1-\rhokk(y))\dd\mu(y)$,
except that $\hmu$ is not a probability measure -- this is
unimportant since we can normalize and consider
$\hkkp\=(1-\rho(\kk))\kk$ and
$\hmup\=(1-\rho(\kk))\qi\hmu$; see \refD{Ddual} and the discussion 
following.

The process $\bpkxp$ dies out by construction, and is thus subcritical
or critical. \refE{E2large1} shows that it can be critical (even when
$\kk$ is irreducible). In many cases, however, $\bpkxp$ is
subcritical; we give one simple criterion.

\begin{lemma}
  \label{Ldual}
Suppose that $\kk$ is irreducible and that $\norm{\tk}>1$.
If $g\ge0$ is integrable and such that $\tk\bigpar{(1-\rhokk)g}\ge g$
\aex, then $g=0$ a.e.
\end{lemma}

\begin{proof}
We may assume that \eqref{bb1} holds for every $x$.
By \refT{T:GW} and \refL{LPhi}\ref{lphib},
  $(1-\rhokk)\tk\rhokk<\rhokk$ a.e.
If $g>0$ on a set of positive measure, then
\begin{equation*}
  \int_\sss g\rhokk\dd\mu
>\int_\sss g(1-\rhokk)\tk\rhokk\dd\mu
=\int_\sss \rhokk\tk\bigpar{g(1-\rhokk)}\dd\mu
\ge \int_\sss \rhokk g\dd\mu,
\end{equation*}
a contradiction.
\end{proof}

\begin{theorem} \label{Tdualbp}
Suppose that $\kk$ is a \qir\ kernel on $(\sss,\mu)$, 
and that $\norm{\tk}>1$.
Let $\hmu$ be the measure defined by
$\dd\hmu(y)=(1-\rhokk(y))\dd\mu(y)$, and let $\htk$ be the
corresponding integral operator
$$
\htk g\= \int_\sss \kk(x,y)g(y)\dd\hmu(y)=\tk\bigpar{(1-\rhokk)g}.
$$
Then $\norm{\htk}_{L^2(\hmu)}\le1$.

If, in addition, $\iint_{\sss^2} \kk(x,y)^2\dd\mu(x)\dd\mu(y)<\infty$,
then
$\norm{\htk}_{L^2(\hmu)}<1$.
\end{theorem}

Note that with $\hkkp=(1-\rho(\kk))\kk$ and $\hmup=(1-\rho(\kk))\qi\hmu$ as above, 
we have $\htk=\htkp$, where $\htkp$ is defined by
$\htkp g(x)\= \int_{\sss}\hkkp(x,y)g(y)\dd\hmup(y)$.
Thus
$\norm{\htk}_{L^2(\hmu)}=\norm{\htkp}_{L^2(\hmup)}$.

\begin{proof}
We may assume \eqref{bb1}
and that $\kk$ is irreducible.
The discussion above and \refT{T:GW} show that
$\norm{\htkp}_{L^2(\hmup)}\le1$, as $\bpkxp$ dies out by construction.
(An analytic proof is easily given
too, using a truncation of $\kk$ and the argument below for the second part.)

For the second part, the additional assumption implies that
\linebreak
$\iint_{\sss^2} \hkkp(x,y)^2\dd\hmup(x)\dd\hmup(y)<\infty$, so
$\normHS{\htkp}<\infty$. \refL{Lcomp} shows that $\htk=\htkp$ has an
eigenfunction $g\ge0$ with eigenvalue $\norm{\htk}$, and thus
$$
\norm{\htk} g
=\htk g
=\tk\bigpar{(1-\rhokk)g}
\qquad\text{a.e.}
$$
If $\norm{\htk}\ge1$, this contradicts \refL{Ldual}.
\end{proof}

With a few exceptions, in the rest of the paper we shall not refer 
directly to the lemmas in \refS{S:branch1}; the results in this section
describe the properties of the branching process we shall use.

\section{Approximation}\label{Sapprox}

In this section we introduce certain upper
and lower approximations to a kernel $\gk$ on a (generalized) ground space 
$(\sss,\mu)$, in preparation for the study of the random
graph $\gnxx{\kkn}$.
Recall that $\sss$ is a separable metric space, and
that $\mu$ is a Borel measure on $\sss$ with $0<\mu(\sss)<\infty$. We
usually assume  
that $\mu(\sss)=1$; in this section, this makes no difference.
Here the metric and topological structure of $\sss$ will be important.

Given a sequence of finite partitions
$\Pm=\set{A_{m1},\dots,A_{mM_m}}$, $m\ge1$, of $\sss$
and an $x\in\sss$, we define $i_m(x)$ by
\begin{equation} \label{imdef}
 x\in A_{m,i_m(x)}.
\end{equation}
As usual, for $A\subset \sss$ we write $\diam(A)$ for 
$\sup\{d(x,y): x,y\in A\}$, where $d$ is the metric on
our metric space $\sss$.

\begin{lemma}
  \label{LP}
Let $(\sss,\mu)$ be a (generalized) ground space. 
There exists a sequence of finite partitions
$\Pm=\set{A_{m1},\dots,A_{mM_m}}$, $m\ge1$, of $\sss$ such that
\begin{romenumerate}
\item\label{LPa}
each $A_{mi}$ is measurable and $\mu(\partial A_{mi})=0$;
\item \label{LPb}
for each $m$, $\cP_{m+1}$ refines $\Pm$, \ie{}, each $A_{mi}$ is a
union $\bigcup_{j\in J_{mi}}A_{m+1,j}$ for some set $J_{mi}$;
\item \label{LPc}
for \aex{} $x\in\sss$, $\diam(A_{m,i_m(x)})\to0$ as \mtoo,
where $i_m(x)$ is defined by \eqref{imdef}.
\end{romenumerate}
\end{lemma}

\begin{proof}
If $\sss=(0,1]$ and $\mu$ is continuous, \eg{}, $\mu$ is the
Lebesgue measure, we can take $\Pm$ as the dyadic partition into
intervals of length $2^{-m}$. If $\sss=(0,1]$ and $\mu$ is
arbitrary, we can do almost the same; we
  only shift the endpoints of the intervals a little when necessary to
  avoid point masses of $\mu$.

In general, we can proceed as follows.
Let $z_1,z_2,\dots$ be a
dense sequence of points in $\sss$. For any  $z_i$,
the balls $B(z_i,r)$, $r>0$, have disjoint boundaries, and thus
all except at most a countable number of them are \mucs s.
Consequently, for every $m\ge1$ we may choose balls
$B_{mi}=B(z_i,r_{mi})$ that are \mucs s and have radii satisfying
$1/m<r_{mi}<2/m$. Then, $\bigcup_i B_{mi}=\sss$, and if we define
$B'_{mi}\=B_{mi}\setminus\bigcup_{j<i}B_{mj}$, we obtain for each
$m$ an infinite partition $\set{B'_{mi}}_1^\infty$ of $\sss$ into
\mucs{s}, each with diameter at most $4/m$. 
To get a finite partition, we choose $N_m$ large enough
to ensure that, with $B'_{m0}\=\bigcup_{i>N_m}
B'_{mi}$, we have $\mu(B'_{m0})< 2^{-m}$; then
$\set{B'_{mi}}_{i=0}^{N_m}$ is a partition of $\sss$ for each $m$,
with $\diam(B'_{mi})\le 4/m$ for $i\ge 1$.

Finally, we let $\cP_m$ consist of all intersections 
$\bigcap_{l=1}^m B'_{l i_l}$ with $0\le i_l\le N_l$; then conditions
\ref{LPa} and \ref{LPb} are satisfied. Condition \ref{LPc} follows
from the Borel--Cantelli Lemma: as $\sum_m\mu(B'_{m0})$ is finite,
a.e.\ $x$ is in finitely many of the sets $B'_{m0}$. For any such $x$,
if $m$ is large enough then $x\in B'_{mi}$ for some $i\ge 1$,
so the part of $\cP_m$ containing $x$ has diameter at most
$\diam(B'_{mi})\le 4/m$. 
\end{proof}

Recall that a kernel $\gk$ on $(\sss,\mu)$
is a symmetric measurable function on
$\sss\times\sss$. Fixing a sequence of partitions
with the properties described in \refL{LP},
we can define sequences of lower and upper
approximations to $\kk$ by
\begin{align}
\kkm-(x,y)&\=\inf\set{\kk(x',y'):x'\in A_{m,i_m(x)},\;y'\in A_{m,i_m(y)}},
\label{k-}
\\
  \kkm+(x,y)&\=\sup\set{\kk(x',y'):x'\in A_{m,i_m(x)},\;y'\in A_{m,i_m(y)}}.
\label{k+}
\end{align}
We thus replace $\kk$ by its infimum or supremum on each
$A_{mi}\times A_{mj}$. As $\kkm+$ might be $+\infty$, we shall
use it only for bounded $\kk$.

By \refL{LP}\ref{LPb},
\begin{align*}
  \kkm-\le\kk_{m+1}^-
\quad\text{and}\quad
  \kkm+\ge\kk_{m+1}^+.
\end{align*}
Furthermore, if $\kk$ is continuous \aex\ then, by \refL{LP}\ref{LPc},
\begin{align}\label{kk}
  \kkm-(x,y)\to\kk(x,y)
\text{ and }
  \kkm+(x,y)\to\kk(x,y)
\text{ for \aex{} $(x,y)\in\sss^2$}.
\end{align}
Since $\kkm-\le\kk$, we can obviously construct our random graphs so that
$\gnx{\kkm-}\subseteq\gnk$;
in the sequel we shall assume this. Similarly, we shall
assume that $\gnx{\kkm+}\supseteq\gnk$ when $\kk$ is bounded.

If $(\kkn)$ is a graphical sequence of kernels with limit $\kk$,
we define instead
\begin{equation} \label{kn-}
\kkm-(x,y)\=\inf\set{(\kk\bmin\kk_n)(x',y'):
  x'\in A_{m,i_m(x)},\;y'\in A_{m,i_m(y)}, \; n\ge m}.
\end{equation}

By \refL{LP}\ref{LPb}, we have $\kkm-\le\kk_{m+1}^-$, and from
\refL{LP}\ref{LPc} and \eqref{t2a} we see that
\begin{equation}\label{kkk}
  \kkm-(x,y)\upto\kk(x,y) \text{\quad as \mtoo, for \aex{} $(x,y)\in\sss^2$}.
\end{equation}
Moreover, when $n\ge m$ we have
\begin{equation}\label{kknm}
 \kk_n\ge\kkm-,
\end{equation} and we may assume that
$\gnx{\kkm-}\subseteq\gnx{\kk_n}$.

For a uniformly bounded graphical sequence
$(\kkn)$ of kernels with limit $\kk$, we similarly define
\begin{equation}\label{kn+}
\kkm+(x,y)\=\sup\set{(\kk\bmax\kk_n)(x',y'):
  x'\in A_{m,i_m(x)},\;y'\in A_{m,i_m(y)}, \; n\ge m} <\infty.
\end{equation}

Relations corresponding to \eqref{kkk} and \eqref{kknm} hold for
$\kkm+$; we collect  
these and an additional result in the following lemma.

\begin{lemma}\label{L:uapprox}
Let $(\kkn)_{n\in I}$ be a graphical sequence of kernels on a
(generalized) vertex space $\vxs$ 
with limit $\kk$, and suppose that $\sup_{x,y,n}\kkn(x,y)<\infty$.
Then there is a sequence $\kkm+$, $m=1,2,\ldots$, 
of \rfin\ kernels on $\vxs$
with the following properties.
\begin{romenumerate}
\item\label{L:uapprox1}
We have $\kkm+(x,y)\downto\kk(x,y)$  as $\mtoo$ for \aex{} $(x,y)\in\sss^2$.
\item\label{L:uapprox2}
Whenever $n\ge m$ we have $\kkm+(x,y)\ge \kkn(x,y)$ for \emph{every}
$(x,y)\in\sss^2$.
\item\label{L:uapprox3}
$\norm{\Tx{\kkm+}}\downto \norm{\tk}$ as $\mtoo$.
\end{romenumerate}
\end{lemma}

\begin{proof}
Let $\Pm=\set{A_{m1},\dots,A_{mM_m}}$, $m\ge1$, be a sequence of partitions
with the properties described in \refL{LP}, and define  
$\kkm+(x,y)$ by \eqref{kn+}.
(If $\kk_n=\kk$ for all $n$, this is just \eqref{k+}.)
Then \ref{L:uapprox2} holds trivially.
By \refL{LP}\ref{LPc} and \eqref{t2a},
$\kkm+\downto\kk$ a.e., proving \ref{L:uapprox1}.
Finally, by dominated convergence,
$\normHS{\Tx{\kkm+}-\tk}\to0$.
Hence,
\begin{equation*}\label{tk+}
\norm{\tk} \le \norm{\Tx{\kkm+}}
\le
\norm{\tk}+\norm{\Tx{\kkm+}-\tk}
\le \norm{\tk} +\normHS{\Tx{\kkm+}-\tk}
\downto \norm{\tk},
\end{equation*}
proving \ref{L:uapprox3}.
\end{proof}

We finish this section with a result for lower approximations 
corresponding to \refL{L:uapprox}, but with one additional
ingredient: for lower approximations to be
useful we shall often need them to be \qir.

\begin{lemma}\label{L:approx}
If\/ $(\kkn)_{n\in I}$ is a graphical sequence of kernels on a
(generalized) vertex space 
$\vxs$ with limit $\kk$,
there is a sequence $\tkkm$, $m=1,2,\ldots$, 
of regular finitary kernels on
$\vxs$ with the
following properties.
\begin{romenumerate}
\item\label{L:approxa}
If $\kk$ is \qir{}, then so is $\tkkm$ for all large $m$.
\item\label{L:approxb}
We have $\tkkm(x,y)\upto\kk(x,y)$  as $\mtoo$ for \aex{} $(x,y)\in\sss^2$.
\item\label{L:approxc}
Whenever $n\ge m$ we have $\tkkm(x,y)\le \kkn(x,y)$ for \emph{every}
$(x,y)\in\sss^2$.
\end{romenumerate}
\end{lemma}

Before turning to the proof, let us note that the conclusions of the
lemma are obvious for suitably
`nice' kernels $\kk$ (or sequences $\kkn\to\kk$);
for example if $\kk$ is continuous, $\sss$ is compact and $\kk>0$.
Indeed, if we
partition $\sss$ into finitely many pieces $\sss_i$
in a suitable way, we may then set
$\tkkm(x,y)=\inf\{\kk(x',y')\ :\ x'\in \sss_i,\ y'\in \sss_j\}$
whenever $x\in \sss_i$
and $y\in\sss_j$.
Note also that in the application we shall need condition
\ref{L:approxc} for \emph{every} $(x,y)\in\sss^2$:
while changes in a kernel $\kk$ on a set of measure zero do not affect
the branching process $\bpk$,
they can affect the graph $\gnkx$.

\begin{proof}[Proof of \refL{L:approx}]
We may assume that $\kk>0$ on a set of positive measure, as otherwise
we may take $\tkkm=0$ for every $m$
and there is nothing to prove.
We shall construct the sequence $\tkkm$ in two stages.

Let $\Pm=\set{A_{m1},\dots,A_{mM_m}}$, $m\ge1$, be a sequence of partitions
with the properties described in \refL{LP}.
If $\kkn=\kk$ for all $n$,
we start with $\kkm-$ defined in
\eqref{k-}.
In general, with a sequence $\kk_n$, we use instead
the definition \eqref{kn-}.

Each $\kkm-$ is of the \rfin{} type treated above, and the $\kkm-$ have
two of the properties required for the $\tkkm$, namely (ii) and (iii),
by \eqref{kkk} and \eqref{kknm}, respectively. However, (i) may fail,
as some $\kkm-$ may be reducible. From now on we shall assume
that $\kk$ is \qir, as otherwise we may take $\tkkm=\kkm-$.
In fact, without loss of generality we may assume that $\kk$ is
irreducible. Indeed, it suffices to prove this case as, given a \qir{}
$\kk$, we may then apply the result to the irreducible restriction
to $\sss'\times\sss'$, and extend the approximating $\tkkm$
obtained to $\sss\times\sss$ by taking them to be zero off
$\sss'\times\sss'$. We shall thus assume that $\kk$ is irreducible.

If $\kkm-=0$ \aex{} for every $m$, then $\kk=0$ \aex{} by \eqref{kkk},
contradicting our
assumption.
We may thus assume that there exists an $m_0$ such that
$\kkmx{m_0}>0$ on a set of positive measure.
We consider only $m\ge m_0$, and assume for notational convenience
that $m_0=1$.
Thus there exist $i$ and $j$ (possibly
equal) with $\mu(A_{1i}),\mu(A_{1j})>0$ and $\kkmx{1}>0$ on
$A_{1i}\times A_{1j}$. From now on we fix such a pair $i$ and $j$.

For $m\ge1$, let
$E_m\=\bigcup\set{A_{mi}:\mu(A_{mi})=0}$, noting that $\mu(E_m)=0$, and let
$B_m$
be the set of all $x\in\sss$ such that for some $k\ge1$
there exists a sequence $x_0,\dots,x_k$
with  $x_0=x$, $x_k\in A_{1i}$, $\kkm-(x_{l-1},x_l)>0$,
and $x_l\notin E_m$ for
$l=1,\dots,k$. (Note that $x=x_0$ may belong to $E_m$.)
Since $\kkm-$ is constant on each $A_{mp}\times A_{mq}$,
$B_m$ is a union of some of the sets $A_{mp}$.
It is easily seen that
$B_m\subseteq B_{m+1}$, that
$B_m\supseteq B_1 \supseteq A_{1j}$,
that the restriction of $\kkm-$ to $B_m$ is irreducible and
that $\kkm-=0$ on $(B_m\setminus E_m)\times(\sss\setminus B_m)$ and
thus \aex{} on
$B_m\times(\sss\setminus B_m)$.

Let $B\=\bigcup_1^\infty B_m$.
If $n\ge m$, then $B_m\subseteq B_n$ and thus $\kkmx{n}=0$ \aex{} on
$B_m\times(\sss\setminus B)\subseteq B_n\times(\sss\setminus B_n)$.
Letting $n\to\infty$, \eqref{kkk} shows that $\kk=0$ \aex{} on
$B_m\times(\sss\setminus B)$. Letting now $m\to\infty$ (taking
the union) yields $\kk=0$ \aex{} on $B\times(\sss\setminus B)$. Since
$\kk$ is irreducible, it follows by \eqref{t1ay} that $\mu(B)=0$ or
$\mu(\sss\setminus B)=0$.  
As
$B\supseteq B_{1}\supseteq A_{1j}$, we have $\mu(B)>0$, so
$\mu(\sss\setminus B)=0$. 
In other words, \aex{} $x\in B=\bigcup_m B_m$.

Now define
\begin{equation*}
  \tkkm(x,y)=\kkm-(x,y)\ett[x\in B_m]\ett[y\in B_m].
\end{equation*}
Thus $\tkkm$ is 0 off $B_m\times B_m$, and the restriction to $B_m$ is
by construction irreducible and of the \rfin{} type, so condition
\ref{L:approxa} of the lemma is satisfied.
Furthermore, by \eqref{kkk} and the fact that $B_m\upto B$ with
$\mu(\sss\setminus B)=0$, 
we have
$  \tkkm(x,y)\upto\kk(x,y)$  as $\mtoo$ for \aex{} $(x,y)\in\sss^2$,
so \ref{L:approxb} holds.

Finally, if $n\ge m$, then $\tkkm\le\kkm-\le\kk_n$, so \ref{L:approxc} holds.
\end{proof}

\section{The number of edges}\label{Sedges}

In this section we consider circumstances in which the condition
\eqref{t1b} or \eqref{t2b} 
on the convergence of the number of edges in $\gnkx$ does, or does not, hold.
In doing so, we shall make frequent use of the approximating kernels
$\kkm-$ and $\kkm+$ defined for a single kernel $\kk$ by \eqref{k-}
and \eqref{k+}, 
and for a sequence by \eqref{kn-} and \eqref{kn+}.
As before, we shall always write the (generalized)
vertex space $\vxs$ under
consideration as $(\sss,\mu,\nxss)$,
unless otherwise specified.

\begin{lemma}
  \label{LE}
Let $\kk$ be an \aex{} continuous kernel on a (generalized) vertex
space $\vxs$. Then 
  \begin{equation}\label{le1}
\liminf \frac1n \E e\bigpar{\gnkx}
\ge \hiikxy.
  \end{equation}
If\/ $\kk$ is a bounded \aex{} continuous kernel on a vertex space $\vxs$, then
  \begin{equation}\label{le2}
\lim \frac1n \E e\bigpar{\gnkx}
= \hiikxy.
  \end{equation}
\end{lemma}

\begin{proof}
Write $\gnk$ for $\gnkx$. Consider first the \rfin{} case defined in
\refD{Drfin}. For $n\ge\max\kk$, conditioning on $n_1,\dots,n_r$ we
have
\begin{multline}  \label{le3}
\frac1n \E e\bigpar{\gnk\mid n_1,\dots,n_r}
=\frac1{2n} \sum_{i,j=1}^r (n_in_j-n_i\gd_{ij}) \frac1n \kk(i,j)
\\
\pto\tfrac1{2} \sum_{i,j=1}^r  \kk(i,j)\mu(S_i)\mu(S_j) =\hiikxy.
\end{multline}
Taking expectations and applying Fatou's Lemma, it follows that 
\eqref{le1} holds in this case.

In general, to prove \eqref{le1}
we use \refL{LP} and the approximation \eqref{k-}. For
every $m$, by the case just treated,
\begin{equation*}
\liminf_{\ntoo} \frac1n \E e\bigpar{\gnk}
\ge
\liminf_{\ntoo} \frac1n \E e\bigpar{\gnx{\kkm-}}
=\tfrac12\iint_{\sss^2} \kkm-. 
\end{equation*}
As \mtoo,  the monotone convergence theorem implies that$\iint
\kkm-\to \iint \kk$, and  \eqref{le1} follows.

If $\vxs$ is a vertex space and $\kk$ is \rfin, then the left-hand 
side of \eqref{le3} is bounded by $\max\kk/2$, so by the dominated convergence
theorem we have
$\frac1n \E e\bigpar{\gnk}\to \hiik$.
In general, if $\vxs$ is a vertex space and $\kk$ is bounded,
we can use $\kkm+$ in place of $\kkm-$ to show that
$\limsup _{\ntoo} \frac1n \E e\bigpar{\gnk}\allowbreak
\le \hiik$.
\end{proof}

\begin{remark}\label{Rgenexp}
Condition \eqref{le2} may fail for a generalized vertex space $\vxs$,
even if $\kk$ 
is constant. The problem is that the definition of a generalized
vertex space only 
imposes `\whp{} conditions' on the number of vertices, giving no control on the
distribution in the $o(1)$ probability case that these conditions
fail, and hence 
giving no control on expectations. In particular, with $\kk$ identically $1$,
the expected number of edges is essentially $\frac{1}{2n}\E(|V(\gnkx)|^2)$,
and we have no control over this expectation -- it can even be infinite.

When the number $\nvn$ of vertices is sufficiently concentrated (for example
Poisson), this problem does not arise. 
Indeed, \eqref{le2} holds whenever $\kk$ is bounded and $\Var(\nvn/n)\to0$;
since $\nvn/n\pto\mu(\sss)$ by assumption,
the variance condition is easily shown to be
equivalent to $\E(\nvn/n)^2\to\mu(\sss)^2$, and to imply
uniform integrability of $(\nvn/n)^2$;
see \eg{} \cite[Proposition 4.12]{Kall}.
(If the parameter $n$ is not restricted to integers,
we may have to consider a sequence of indices $n$.)  
Since the left-hand 
side of \eqref{le3} is bounded by $\max\kk\,(\nvn/n)^2$, which is
also uniformly integrable, we may take the expectation in \eqref{le3}
and obtain \eqref{le2}.

Our main results concern statements that hold whp, and
convergence in probability 
of various quantities. For such statements, a small chance of a very large
number of vertices is not a problem.
\end{remark}

The following lemma shows that the condition to be graphical is
essentially equivalent to a statement about approximations with
bounded kernels. 
\begin{lemma}
  \label{LEM}
Let $\kk$ be a bounded \aex{} continuous kernel on a (generalized) vertex
space $\vxs$.

If $\vxs$ is a vertex space, then 
$\kk$ is graphical if and only if
\begin{romenumerate}
\item \label{lemeps}
for every $\eps>0$ there exists an $M<\infty$ such that
  \begin{equation*}
\limsup \frac1n \E e\bigpar{\gnkx\setminus\gnxx{\kkxM}}
\le \eps.
  \end{equation*}
\end{romenumerate}

In general,
$\kk$ is graphical if and only if
\ref{lemeps} holds together with
\begin{romenumerate}\setcounter{enumi}{1}
  \item \label{lemM}
for every $M<\infty$,
  \begin{equation*}
\limsup \frac1n \E e\bigpar{\gnxx{\kkxM}}
\le \hiikxy .
  \end{equation*}
\end{romenumerate}
\end{lemma}
\begin{proof}
  It is obvious that (i) and (ii) imply that
$\limsup \frac1n \E e\bigpar{\gnxx{\kk}}
\le \hiikxy$, which
together with \refL{LE} shows that $\kk$ is graphical.

Conversely, suppose that $\kk$ is graphical on $\vxs$.
Then, from the definition of graphicality (see \eqref{t1b}), 
$\E e\bigpar{\gnkx}\to\frac12\iint\kk$. 
Applying \refL{LE} to $\kkxM$, it follows that
  \begin{multline*}
\limsup \frac1n \E e\bigpar{\gnkx\setminus\gnxx{\kkxM}}
\\
= \hiik-
\liminf \frac1n \E e\bigpar{\gnxx{\kkxM}}
\le \tfrac12\iint(\kk-\kkxM)<\eps
  \end{multline*}
if $M$ is large enough.
\end{proof}
Note that \ref{lemM} almost always holds by \refL{LE} and \refR{Rgenexp}.
Arguing as in the proof of \refL{LEM}, one can show that \ref{lemM}
can be replaced by 
the condition that each $\kkxM$ be graphical; we omit the details.

\begin{remark}\label{Rckgraphical}
\refL{LEM} implies that,
if $\kk$ is a graphical kernel on a (generalized) 
vertex space $\vxs$ and $0<c<\infty$, then $c\kk$ is also graphical on
$\vxs$. Indeed, it suffices to check condition \ref{T1d} of
\refD{Dg1}, namely that $\frac1n \E e\bigpar{\gnx{c\kk}} \to
\tfrac12\iint c\kk(x,y)$. Without the $\min\{\cdot,1\}$ in the formula
\eqref{Ee}, 
this would be immediate from the same condition for $\kk$; indeed,
the claim that $c\kk$ is graphical is equivalent to the claim that replacing
this $1$ with $1/c$ does not affect $\E e\bigpar{\gnk}$ by more than $o(n)$.

Since $c\kk\bmin cM=c(\kkxM)$, it is obvious that
condition \ref{lemM} of \refL{LEM}
holds for $c\kk$ if and only if it holds for $\kk$.
Moreover, if $n\ge M$, 
\begin{equation*}
\tfrac1n\E\bigpar{e(\gnkx\setminus\gnxx{\kkxM})\mid\xs}
=n^{-2}\sum_{i<j} \bigpar{\kk(x_i,x_j)\bmin n - \kk(x_i,x_j)\bmin M}.
\end{equation*}
It is clear that if we replace $\kk$ by $c\kk$ and $M$ by $cM$, and
assume $n\ge 2(1\bmax c)M$, then this sum changes by at most a
constant factor. Hence condition \ref{lemeps} of \refL{LEM} also
holds for $c\kk$ if and only if it holds for $\kk$.

Lemmas \refand{LE}{LEM} hold also for the variants of $\gnkx$ defined
in \refR{Rexp}, by the same proofs. Moreover, it is easily seen that
conditions (i) and (ii) of \refL{LEM} hold for one of these versions if and only if
they hold for $\gnkx$. Hence $\kk$ is graphical if and only if the
analogue of \eqref{t1b} for one of these variants holds.

The results above can be extended to sequences $(\kk_n)$ satisfying
\eqref{t2a}, using the approximations $\kkm-$ defined by 
\eqref{kn-}. In particular, a similar argument shows 
that if $(\kkn)$ is a graphical sequence of kernels on a (generalized)
vertex space 
$\vxs$ with limit $\kk$, and $(c_n)$ is a sequence of positive reals
with $c_n\to c>0$, 
then $(c_n\kkn)$ is graphical on $\vxs$ with limit $c\kk$.
\end{remark}

Let us emphasize that relation \eqref{le2}, i.e., condition \eqref{t1b} 
from the definition of graphicality, often holds for unbounded $\kk$ too,
and for generalized vertex spaces $\vxs$. 
One example is when the $x_i$ are random as in \refE{EU}; another is the
Poisson process 
case in \refE{EPo}. A rather different example is the
following.

\begin{example}
 \label{E42}
Suppose that $\sss=(0,1]$, $\mu$ is the Lebesgue measure and
$x_i=i/n$; this vertex space was considered in
\refE{Etri} and will be used in several further examples in
\refS{Sapp}, with several different kernels. Suppose that
  $\kk$ is decreasing in each variable, so $\kk(x,y)\ge \kk(x',y')$ when
  $x\le x'$ and $y\le y'$.
Then
\begin{equation*}
  \frac1n\E e\bigpar{\gnk}
=\frac{1}{2n^2}\sum_{i\neq j} \kk(i/n,j/n)\bmin n
\le \frac12\int_0^1\int_0^1\kk(x,y)\dd x\dd y.
\end{equation*}
Hence \eqref{le1} implies that \eqref{le2} holds in this case.
Note that this includes both 
\eqref{b3a} and \eqref{b3b}.
\end{example}

We next give a simple example where \eqref{t1b}, i.e., \eqref{le2}, fails.

\begin{example}
  \label{Ebad}
Take again $\sss=(0,1]$, let $\mu$ be the Lebesgue
measure and set $x_i=i/n$. Let $0<\gd<1$ be constant and define $\kk$
by
\begin{equation*}
  \kk(x,y)=
  \begin{cases}
m & \text{if $x \bmin y=1/m$ and $m\ge1$;}
\\
\gd & \text{otherwise}.
  \end{cases}
\end{equation*}
Note that $\kk=\gd$ a.e., and hence $\rho(\kk)=\rho(\gd)=0$ as for
$\gnx{\gd/n}$ in \refE{E1}; furthermore, $\kk$ is continuous a.e.

Now, $\kk(1/n,j/n)=n$ for every $j\le n$. Hence, $\gnk$ contains
the star consisting of all edges $1j$, $1<j\le n$, so $\gnk$ is
connected and $C_1(\gnk)=n$, although, as remarked above,
$\rho(\kk)=0$. Consequently, \eqref{t2bound} fails in 
this case. Note that all assumptions of \refT{T2}
are satisfied except \eqref{t1b}; indeed, 
$e(\gnk)\ge n-1$
while
$ \iikxy=\gd$.

We can modify this example to make $\kk$ continuous on $(0,1]^2$:
for $0<\eps<1/4$, let
  $\kk_\eps(x,y)=\xphi(x\bmin y)$ with $\xphi(1/m)=m$,
$\xphi(1/m\pm\eps m^{-4})=\gd$, and $\xphi$ linear in between. If
$\eps$ is small enough, then $\norm{\tkx{\eps}}<1$ (because the
Hilbert--Schmidt norm satisfies
$\normHS{\tkx{\eps}}\to\gd$ as $\eps\to0$ by dominated convergence);
thus
$\rho(\kk_\eps)=0$, although $C_1(\gnx{\kk_\eps})=n$.
\end{example}

We next give a result on the number of edges conditioned on
$\xs$; 
this time we consider a sequence $(\kk_n)$ of kernels.

\begin{lemma}
  \label{LEcond}
Let $(\kkn)$ be a graphical sequence of kernels on a (generalized)
vertex space 
$\vxs$ with limit $\kk$. Then
\begin{equation*}
\frac1n  \E\Bigpar{e\bigpar{\gnxx{\kk_n}}\Bigm|\xs}
\pto \hiikxy.
\end{equation*}
\end{lemma}

\begin{proof}
Let $W_n\= \E\bigpar{e\bigpar{\gnxx{\kk_n}}\bigm|\xs}/n$ and
$w\=\hiikxy$. By our assumption \eqref{t2b}, we have
$\E W_n\to w$.

Define $\kkm-$ by \eqref{kn-}. By \eqref{le3}, applied to $\kkm-$,
\begin{equation*}
\wnm\= \E\bigpar{e\bigpar{\gnxx{\kkm-}}\bigm|\xs}/n
\pto \wm:=\tfrac12\iint_{\sss^2}\kkm-(x,y)\dd\mu(x)\dd\mu(y).
\end{equation*}
Let $\eps>0$ be given.
By \eqref{kkk} and monotone
convergence, $\wm\to w$ as \mtoo, so we
may choose $m$ such that $w_m>w-\eps$.
For $n\ge m$ we have $W_n\ge\wnm$, and hence 
\begin{equation*}
  \P\bigpar{W_n<w-2\eps}
\le
\P\bigpar{\wnm<\wm-\eps}
\to0
\qquad
\text{as }\ntoo.
\end{equation*}
Hence, writing $f_-$ for $-(f\bmin 0)$, we have
$(W_n-w)_-\pto0$ and, by dominated convergence,
$\E(W_n-w)_-\to0$.
Consequently, $\E|W_n-w|=2\E(W_n-w)_- +\E(W_n-w)\to0$.
\end{proof}

\begin{remark}
  \label{Rcond}
Recalling \eqref{a2a} or \eqref{nunA2},
the convergence condition for the empirical 
distribution $\nu_n$ of the types of the vertices in a 
(generalized) vertex space,
we have
$\nu_n\pto\mu$ and $W_n\pto w$ (in the notation of the proof above),
where $\nu_n$ and $W_n$ are functions of $\xs$.
Coupling the $\xs$
for different $n$ appropriately
(a simple application of
the Skorohod coupling theorem \cite[Theorem 4.30]{Kall}),
or considering appropriate subsequences,
we may assume that
$\nu_n\to\mu$ and $W_n\to w$ a.s.
Consequently, we may condition on $\xs$ and assume that \eqref{a2a}
and \eqref{t2b} still hold.
In other words, after conditioning on $\xs$, $\vxs$ is still a
(generalized) vertex 
space and $(\kk_n)$ is still graphical with limit $\kk$.
By conditioning in this way we may thus assume that $\xs$ is
deterministic; see \refSS{SSgvs}.
\end{remark}

Our next result shows that the number of edges is concentrated,
so that the actual number converges as well as its mean.

\begin{proposition}
    \label{PE}
Let $(\kkn)$ be a graphical sequence of kernels on a (generalized)
vertex space  
$\vxs$ with limit $\kk$. Then
  \begin{equation*}
 \frac1n e\bigpar{\gnkxn}
\pto \hiikxy.
  \end{equation*}
\end{proposition}

\begin{proof}
Let $G_n=\gnkxn$ and, as above, $W_n=\E(e(G_n)\mid\xs)/n$.
Conditioned on $\xs$, the number
$e(G_n)$ of edges is a sum of independent $\Be(\pij)$ variables,
and thus $\Var(e(G_n)\mid\xs) \le \E(e(G_n)\mid\xs)$.
Hence, using \eqref{t2b},
\begin{equation*}
  \E\bigpar{e(G_n)/n-W_n}^2
=n^{-2}\E\bigpar{\Var(e(G_n)\mid\xs)}
\le n^{-2}\E\bigpar{e(G_n)}
\to0.
\end{equation*}
Consequently, $e(G_n)/n-W_n\pto0$, and the result follows by \refL{LEcond}.
\end{proof}

Finally, we note that small sets of vertices do not connect to too
many edges.
For this we need a simple lemma.

\begin{lemma}
  \label{Ldeg2}
Let $\kk$ be a bounded kernel on a (generalized) vertex space $\vxs$,
and let $G_n=\gnkx$.
Then $\sum_i d^2_{G_n}\!(i)=O(n)$ \whp.
\end{lemma}

\begin{proof}
  Let $a\=2\mu(\sss)$ and $b\=\sup\kk(x,y)<\infty$.
Then \whp{} $v_n\le an$, and thus, in the natural coupling,
$\gnkx\subseteq G(\floor{an},b/n)$ \whp.
Consequently, it suffices to prove the result for $G(\floor{an},b/n)$
or, changing the notation slightly, for $G(n,c/n)$ for every fixed
$c>0$. However, for any graph $G$, $\sum_i d_G(i)=2 e(G)$ and 
$\sum_i d_G(i)(d_G(i)-1)=2 P_2(G)$, twice the number of paths of
length 2 (\cf{} \refS{Scycles}). It is well-known,
and easy to prove, that $e(G(n,c/n))/n\pto \ga_1$ and 
$P_2(G(n,c/n))/n\pto \ga_2$ for some constants $\ga_1,\ga_2$ 
(depending on $c$), see \eg{} \cite[Chapter 3 and Theorem 6.5]{JLR}.
Consequently, with $C\=2\ga_1+2\ga_2+1$, 
$\sum_i d^2_{G(n,c/n)}(i) <Cn$ whp.
(Alternatively, we may use \refT{TPC}.)
\end{proof}

\begin{proposition}\label{Ponon}
Let $(\kkn)$ be a graphical sequence of kernels on a (generalized)
vertex space  
$\vxs$ with limit $\kk$.
Given $\eps>0$, there is a $\delta>0$ so that \whp\
the sum of the degrees 
of any set of at most
$\gd n$ vertices of $G_n=\gnkxn$
is at most $\eps n$.
In particular, any set of $o_p(n)$ vertices of $G_n$ has $o_p(n)$ neighbours.
\end{proposition}

\begin{proof}
Let $d_G(v)$ denote the degree of vertex $v$ in a graph $G$, 
and  let $\tkkm$ be as in \refL{L:approx}. 
Since $\tkkm\upto\kk$ \aex, 
$\iint\tkkm \to \iik$, and thus we can choose $m$ such that
$\iint\tkkm > \iik-\eps$. 
Let $G_n'\=\gnxx{\tkkm}$.
For $n\ge m$ we have $\tkkm\le\kkn$, and we may as usual assume that
$G_n'\subseteq G_n$. Moreover, \refP{PE} applies to both $G_n$ and
$G'_n$, so 
\begin{eqnarray*}
 \frac{1}{2n} \sum_{i\in V(G_n)} \bigpar{d_{G_n}(i)-d_{G'_n}(i)}
&=&
\frac{1}{n}e\bigpar{\gnkxn\setminus\gnxx{\tkkm}} \\
&\pto& \tfrac12\iint_{\sss^2}\kk -\tfrac12\iint_{\sss^2}\tkkm
<\eps.
\end{eqnarray*}
Hence, whp
\begin{equation}\label{oden}
  \sum_{i\in V(G_n)} \bigpar{d_{G_n}(i)-d_{G'_n}(i)}<2\eps n.
\end{equation}

By \refL{Ldeg2}, applied to $G_n'$, there
is a constant $C<\infty$ such that \whp{} 
$\sum_{i} d^2_{G_n'}(i) < Cn$. Hence, if $\gd=\eps^2/C$, the \CSineq{}
shows that \whp{} for every set $A\subseteq V(G_n)$ with $|A|\le\gd n$,
\begin{equation*}
  \sum_{i\in A} d_{G'_n}(i)
\le
\Bigpar{|A|  \sum_{i\in A} d^2_{G'_n}(i)}\qh
< \bigpar{|A|Cn}\qh \le \eps n.
\end{equation*}
Combining this with \eqref{oden}, we obtain
$\sum_{i\in A} d_{G_n}(i) <3\eps n$, \whp{} for all such $A$, 
and the result follows by replacing $\eps$ by $\eps/3$.
\end{proof}

\subsection{Generalized vertex spaces}\label{SSgvs}

Our main results concern graphical sequences of kernels on generalized
vertex spaces, expressing properties of the graphs $\gnxx{\kkn}$ in
terms of the limiting kernel $\kk$. As noted earlier, it is intuitively
clear that we lose no generality by restricting our attention to vertex spaces.
Furthermore, as noted in \refR{Rcond}, we may assume that the vertex types $\xs$
are deterministic.
As we shall now see, \refL{LEcond} and a simple probabilistic lemma given
in the appendix imply precise forms of these assertions.
We start by showing
that we may take the sequences $\xs$ to be deterministic.
 
Let $\kkn$ be a graphical sequence of kernels on a generalized vertex space $\vxs$
with limit $\kk$. As noted in 
\refR{Rcond}, by coupling appropriately we may assume that, after conditioning
on $(\xs)$, the triple $(\sss,\mu,(\xs))$, in which the sequences $\xs$ are now
deterministic, is (\as) still a generalized vertex space, and that
$\kkn$ is (\as) 
graphical on this space with limit $\kk$.
Almost all our results assert that (given some $\eps>0$) a certain
event $\cE_n$ 
holds \whp; recall that statements like $X_n\pto a$ and $X_n=o_p(a_n)$
can be expressed in this form. The $\xs$ deterministic case of such a result
then implies that (\as) the conditional probabilities
$\P(\cE_n\mid\xs)$ tend to 1. 
Taking expectation with respect to the random sequences $\xs$, it follows
by dominated convergence that $\P(\cE_n)\to1$, \ie, the result
holds also for random $\xs$.

Some of our results are of the form $X_n=O(a_n)$ \whp. Again, it suffices to
prove such a result for deterministic $\xs$; the general case then follows by
\refL{Loc2}, with $Y_n=\nu_n$, and $\cM_n$ the set of all measures of the form $n\qi\sum_1^N\gd_{x_i}$,
a subset of the metric space $\cM$ of all
finite Borel measures on $\sss$. The key point is that $\nu_n$
determines $\xs$ up to relabelling the vertices, and that the conditional
distribution of the unlabelled graph $\gnxx{\kkn}$ given $\xs$ does not
depend on the labelling, or on $\bfx_{n'}$, $n'\ne n$.

We now turn to the simple reduction from generalized vertex spaces
to vertex spaces. Although the arguments apply to all our main results,
for definiteness, we shall illustrate them with one particular
example: we shall show that statement \ref{T2_3} of \refT{T2}, namely
\begin{equation}\label{t2bound2}
\frac1n  C_1\xpar{\gnxx{\kk_n}} \pto \rho(\kk),
\end{equation}
follows
from the same statement restricted to the case that $\vxs$ is a vertex
space.

Let $\vxs=(\sss,\mu,(\xs)_{n\in I})$ be a generalized
vertex space, and let $\kkn$, $n\in I$, be a graphical sequence of
kernels on $\vxs$ with limit $\kk$.  As noted in \refS{Smodel},
purely formal manipulations show
that taking $\mu(\sss)=1$ loses no generality, although one must
be a little careful with the introduction of normalizing factors.
To spell this out pedantically, let
$I'=\mu(\sss)I=\{\mu(\sss)n: n\in I\}$, let $\mu'=\mu/\mu(\sss)$ be
the normalized version of the measure $\mu$, and let $\vxs'$ be the
generalized vertex space $(\sss,\mu',({\bfy}_m)_{m\in I'})$ defined by
${\bfy}_m={\bfx}_{m/\mu(\sss)}$, so the sequences $(\xs)$ and 
$({\bfy}_m)$ are 
identical except for our rescaling of the index set. Writing $\kk'$
for $\mu(\sss)\kk$ and $\kk_m'$ for $\mu(\sss)\kk_{m/\mu(\sss)}$, for
$n\in I$ the graphs $\gnxx{\kkn}$ and $G^{\vxs'}(m,\kk_m')$, $m=\mu(\sss)n$,
have exactly the same distribution. Also (as a consequence), the sequence
$\kkn'$ is graphical 
on $\vxs'$ with limit $\kk'$, so our main results,
in particular \refT{T2},
apply to the model $G^{\vxs'}(m,\kk_m')$.

Multiplying $\kk$ by the constant factor $\mu(\sss)$ and dividing $\mu$
by the same factor leaves the branching process $\Phik$, and hence
the survival probability $\rho(\kk;x)$, unchanged,
and so divides $\rho(\kk)$ by a factor $\mu(\sss)$. Multiplying
the index variable $n$ by $\mu(\sss)$ divides the left-hand side of
\eqref{t2bound2} 
by the same factor, so this relation for $\gnxx{\kkn}$
follows from the same relation for the model $G^{\vxs'}(m,\kk_m')$.

Apart from the rather trivial normalization above, there are two further
differences between vertex spaces and generalized vertex spaces.
One is that, in the former, the index set is discrete, indeed a subset
of the integers. 
This makes very little difference: for any result of the form
$f(G_n)\pto a$, $n\in I$, it suffices to consider `thin' index sets
$I$, say discrete sets $\{i_1,i_2,\ldots,\}$ with $i_1\ge 100$ and
$i_{t+1}\ge 2i_t$. 
Indeed, if $f(G_n)\pto a$ fails, there is an $\eps>0$ and
an unbounded set $I'\subset I$
with $\P(|f(G_n)-a|\ge \eps)\ge \eps$ for every $n\in I'$, and then
$f(G_n)\pto a$ fails along any subsequence of $I'$, and hence along
at least one thin sequence. Thus, in all our main
results we need only consider `thin' index sets.

The final extension allowed by generalized vertex spaces is a that
the number of vertices in $\xs$ may be random, rather than exactly $n$.
As noted at the start of the section, we may assume that each $\xs$
is deterministic, and in particular that the number $\nvn$ of vertices
is deterministic. This does not quite give a vertex space, as
we need not have $\nvn=n$:
instead, taking $A=\sss$ in \eqref{nunA2}, we have $\nvn/n\to 1$.
Rescaling the indexing parameter as above, replacing 
$n$ by $\nvn$ (after taking a subsequence if necessary) and
multiplying $\kkn$ by a factor $\nvn/n$,
does not affect the distribution of the graph, so the
resulting kernels 
are still graphical with limit $\kk$. Hence, our results for vertex
spaces apply. 
In particular, using \eqref{t2bound2} for vertex spaces, we find that
\[
 \frac1{\nvn}  C_1\xpar{\gnxx{\kk_n}} \pto \rho(\kk).
\]
As $n\sim v_n$, this implies \eqref{t2bound2}.

We have shown that it suffices to 
prove \eqref{t2bound2}, i.e., part \ref{T2_3} of \refT{T2}, for vertex
spaces in which the sequences $\xs$ are deterministic; this was our 
aim in this subsection.
Similar comments apply to all our results.

\section{The giant component}\label{Smainpf}

In this section we prove our main results, Theorems~\refand{T2}{T2b}
of \refSS{SSgc},
concerning the existence, size and uniqueness of the giant component
in the random graph $\gnxx{\kk_n}$. The basic strategy
will be to relate the neighbourhoods of a vertex of $\gnxx{\kkn}$
to the branching process, by exploring these neighbourhoods step by step.
In the context of random graphs,
this step-by-step exploration and comparison with a branching process,
which now is standard, was 
perhaps first used by Karp~\cite{Karp}, who
applied it to study the size of the giant component in
random directed graphs; similar ideas were used earlier in other
contexts, for example by Kendall \cite{Kendall} in the study of epidemics.

Let us first recall some notation.
We shall work with the branching process \bpk\ defined in \refSS{SS:branch}
and studied in Sections \refand{S:branch1}{S:branch2}.
As before, when the branching process is started with a single
particle of type $x$ we
denote it $\bpkx$.
Unless explicitly stated otherwise, $\kk$ will be a kernel on a vertex
space $\vxs=(\sss,\mu,\xss)$;  
most of the time we shall not consider generalized vertex spaces.
We shall assume that
$\kk\in L^1$, i.e., that $\iint\kk<\infty$. Any additional assumptions on $\kk$
(such as irreducibility) will be stated.

Recall that $\rhogek(\kk;x)$ is the probability that $\bpk(x)$
contains at least $k$ particles in total (in all generations taken
together),
and $\rhok(\kk;x)$ is the probability that $\bpk(x)$ 
contains exactly $k$ particles in total, while $\rho(\kk;x)$ is the
probability that $\bpk(x)$ 
survives for eternity, i.e., for infinitely many generations.
Starting the process with a particle of
random type with distribution $\mu$, the corresponding
probabilities for $\bpk$ are $\rhogek(\kk)$, $\rhok(\kk)$, and $\rho(\kk)$.

A key step in our proofs will be an additional result,
relating the fixed-size components
of $\gnxx{\kkn}$ to the branching process $\bpk$.
As before, we write $N_k(G)$ for the number of vertices
of a graph $G$ in components of order $k$, and
$\ngek(G)$ for $\sum_{j\ge k}N_j(G)$, the number
of vertices in components of order at least $k$.

\begin{theorem}
   \label{T3}
Let $(\kkn)$ be a graphical sequence of kernels on a vertex space
$\vxs$ with limit $\kk$.
If\/ $k\ge1$ is fixed, then
$\nk\bigpar{\gnxx{\kkn}}/n\pto\rhok(\gk)$.
\end{theorem}

\begin{remark}\label{Rcompfinite}
In \cite{giant2,Rsmall}, results similar to \refT{T2} were proved
(for special $\kk$ but with more complicated dependencies) using a
careful coupling of the discovery process of the random graph and
the limiting  branching process; here we shall do this coupling
only in the simple case of finitely many types (\refE{Efin}); the
general case will then follow by approximation and monotonicity
arguments. In particular, we shall show that any $G_n=\gnkx$
contains a $G_n'=\gnxx{\kk'}$, where $\kk'$ may be regarded as a
kernel defined on a finite set $\sss$, such that
$C_1(G_n')$ is no more than $o_p(n)$ smaller than $C_1(G_n)$; a
formal statement is given below. This reduces many questions
concerning the very general model $\gnkx$ to the much simpler
`finite-type' case.
\end{remark}

\begin{proposition}\label{Pfinite}
Let $(\kkn)$ be a graphical sequence of kernels on a vertex space
$\vxs$ with \qir\ limit $\kk$.
Given any $\eps>0$, there is a vertex space $\vxs'=(\sss',\mu',\yss)$
with $\sss'$ finite and a \qir\ kernel $\kk'$ on $\sss'\times\sss'$ with
the following properties: $\rho(\kk')\ge \rho(\kk)-\eps$,
the graphs $G_n=\gnxx{\kk_n}$ and $G_n'=G^{\vxs'}(n,\kk')$
can be coupled
so that $G_n'\subseteq G_n$ for sufficiently large $n$,
and $C_1(G_n')/n\pto\rho(\kk')$.
\end{proposition}

The assertion concerning $C_1(G_n')$ will follow from the other assertions
and \refT{T2}. However, we shall prove \refP{Pfinite}
as a step towards the proof of \refT{T2}.
This is an example where \qirity\ is forced on us: if we assume $\kk$
is irreducible, we still cannot insist that $\kk'$ is irreducible.

\medskip
We now turn to the proofs.
We start by
giving two elementary results that will be useful below. The first
concerns
$\ngek(G)$,  the number of vertices of a graph $G$ that are in
components of order at least $k$. Note that for any graph $G$ and
any $k\ge1$,
\begin{equation}
 C_1(G) \le \max\bigset{k, \ngek(G)},
 \label{eL0}
\end{equation}
since if $C_1(G)\ge k$ then $\ngek(G)\ge C_1(G)$.

\begin{lemma}  \label{LA}
If\/ $k\ge 2$ and  $G$, $G'$ are two graphs with $G\subseteq G'$, then
\begin{equation*}
  \ngek(G)\le \ngek(G') \le \ngek(G)+2k\bigpar{e(G')-e(G)}.
\end{equation*}
\end{lemma}

\begin{proof}
  If we add a single edge to $G$, the set of vertices belonging to
  components of orders $\ge k$ will either remain the same or increase
  by the inclusion of one or two smaller components; hence
  $\ngek(G)$ will increase by at most $2(k-1)$. The result follows by
  iterating $e(G')-e(G)$ times.
\end{proof}

\begin{lemma}\label{L00}
As $\ktoo$,
$\rhogek(\kk;x)\downto\rho(\kk;x)$ a.e.\ $x$,
and
$\rhogek(\kk)\downto\rho(\kk)$.
\end{lemma}

\begin{proof}
As $\kk\in L^1$, \eqref{bb1} holds \aex{} $x$. By \refL{Lae}, we
may assume that \eqref{bb1} holds for every $x$. Then every particle
in the branching process $\bpk$ has a finite number of children,
so a particle survives for eternity if and only if
it has infinitely many descendants, and the result follows.
\end{proof}

Now we turn to the main part of this section, which concerns the
connection between the order of the giant component of $\gnkx$ and
the survival probability $\rho(\kk)$.

We begin by studying the case when $\sss$ is finite. It will turn
out that this case gives essentially everything, using our
monotonicity results and \refL{LA}. We use the notation in
\refE{Efin}.
We shall assume that we have a fixed $\kk$, as in
\refD{Dg1},
rather than a convergent sequence $\kkn$ as in \refD{Dg2} and
\refT{T2}. In addition, we shall assume that the
matrix $\kk$ is irreducible and that $\mui>0$ for every $i$.
As observed by S\"oderberg \cite{Sod1}, we then can adapt the
standard branching process argument for the classical random graph
\gnx{c/n}, see, \eg{}, \cite[Section 5.2]{JLR}. The details are as
follows.

\begin{lemma}\label{L:step1}
Let $\kappa$ be a kernel on the vertex space
$\vxs=(\sss,\mu,\xss)$,
where $\sss=\{1,2,\ldots,r\}$, and
suppose that $\mui>0$ for every $i$.
Writing $G_n$ for $\gnkx$, if $\kk$ is irreducible we have
\begin{equation}
  \label{nc}
C_1(G_n)/n\pto\rho(\kk).
\end{equation}
Whether or not $\kk$ is irreducible, for any fixed $k$ we have
\begin{equation}\label{n7}
\ngek\bigpar{G_n}/n\pto\rhogek(\gk).
\end{equation}
\end{lemma}

\begin{proof}
Recall that we have $n_i$ vertices of type $i$, $i=1,\dots,r$, and that
$n_i/n\pto \mu_i=\mu\set{i}$.
Coupling the graphs (or just the $\xs$)
for different $n$ appropriately, we may of course
assume that $n_i/n\to \mu_i$ a.s.
{}From now on we condition on $n_1,\dots,n_r$; we may thus assume
that $n_1,\dots,n_r$ are deterministic with $n_i/n\to\mu_i$.

Let $\go(n)$ be any function such that $\go(n)\to\infty$ and $\go(n)/n\to0$.
(Although it might seem more natural to fix $\go(n)=\log\log n$, say,
we shall need this flexibility
in the choice of $\go(n)$ later.)
We call a component of $G_n\=\gnkx$ \emph{big} if it has at least
$\go(n)$ vertices. Let $B$ be the union of the big components,
so $|B|=\ngex{\go(n)}(G_n)$.

Fix $\eps>0$. We may assume that $n$ is so large that
$\go(n)/n<\eps\mu_i$ and $|n_i/n-\mu_i|<\eps\mu_i$ for every $i$; thus
$(1-\eps)\mu_in<n_i<(1+\eps)\mu_in$. We may also assume that $n>\max\kk$,
as $\kk$ is a function on the finite set $\sss\times \sss$.

Select a vertex and explore its component in the usual way, one vertex
at a time. We first reveal all edges from the initial vertex, and put
all neighbours that we find in a list of unexplored vertices; we then
choose one of these and reveal its entire neighbourhood, and so on. Stop when
we have found at least $\go(n)$ vertices (so $x\in B$), or when
there are no unexplored vertices left (so we have found the entire
component and $x\notin B$).

Consider one step in this exploration, and assume that we
are about to reveal the neighbourhood of a vertex $x$ of type $i$. Let
us write $n_j'$ for the number of unused vertices of type $j$ remaining.
Note that
$n_j\ge n_j'\ge n_j-\go(n)$, so
\begin{equation}
  \label{n1}
(1-2\eps)\mu_j < n_j'/n < (1+\eps)\mu_j.
\end{equation}
The number of new neighbours of $x$ of type $j$ has a binomial
$\Bi\bigpar{n_j',\kk(i,j)/n}$ distribution, and the numbers for
different $j$ are independent. The total variation distance
between a binomial $\Bi(n,p)$ distribution and the Poisson
distribution with the same mean is at most $p$, see, \eg{}, the
first inequality in Barbour, Holst and Janson~\cite[(1.23)]{SJI}.
Hence the total variation distance between the binomial
distribution above and the Poisson distribution
$\Po\bigpar{\kk(i,j)n_j'/n}$ is at most $\kk(i,j)/n=O(1/n)$. Also,
by \eqref{n1},
\begin{equation*}
 (1-2\eps)\kk(i,j)\mu_j \le \kk(i,j)n_j'/n \le (1+\eps)\kk(i,j)\mu_j.
\end{equation*}
Since we perform at most $\go(n)$ steps in the exploration, we
may, with an error probability of $O(\go(n)/n)=o(1)$, couple the
exploration with two multi-type branching processes
$\bpx{(1-2\eps)\kk}$ and $\bpx{(1+\eps)\kk}$ such that the first
process always finds at most as many new vertices of each type as
the exploration, and the second process finds at least as many.
Consequently, for a vertex $x$ of type $i$,
\begin{equation}
  \label{n2}
\rhox{\ge\go(n)}\bigpar{(1-2\eps)\kk;i}+o(1)
\le \P(x\in B)
\le
\rhox{\ge\go(n)}\bigpar{(1+\eps)\kk;i}+o(1).
\end{equation}
Note for later (after \eqref{ndiff}) that,
as for any constant $C$
the Poisson distribution with mean $C$ has probability $o(1/n)$ of
exceeding $\log n$, the probability that we find more than $\log
n$ new neighbours in one step is $O(1/n)$. It follows that the
probability that we reach more than $\go(n)+\log n$ vertices
during the exploration is $o(1)$. (Informally, we cannot
`overshoot' by more than $\log n$.)

Since $\go(n)\to\infty$, by \refL{L00} we have
$\rhox{\ge\go(n)}(\kk';i)\to\rho(\kk';i)$ for every kernel $\kk'\in L^1$, so
we can rewrite \eqref{n2} as
\begin{equation*}
\rho\bigpar{(1-2\eps)\kk;i}+o(1)
\le \P(x\in B)
\le
\rho\bigpar{(1+\eps)\kk;i}+o(1).
\end{equation*}
Letting $\eps\to0$ we find, using  \refT{TappB},
that if $x$ is of type $i$, then the probability that the component
containing $x$
is big satisfies
\begin{equation}
  \label{n4}
\P(x\in B)
\to
\rho(\kk;i).
\end{equation}
(Recall that we are conditioning on the types of the vertices, treating
the numbers $n_i$ of vertices of type $i$ as deterministic, and assuming
that $n_i/n\to \mu_i$.)
Summing over all vertices $x$ we find
\begin{equation}
  \label{n5}
\begin{split}
\frac1n \E |B|
&
=\frac1n \sum_x \P(x\in B)
=\frac1n \sumir n_i \P(x\in B\mid \text{$x$ is of type $i$})
\\&
\to \sumir \mu_i\rho(\kk;i)
=\rho(\kk).
\end{split}
\end{equation}
Note that this limit is independent of the choice of $\go(n)$ in the
definition of $B$. Hence, if we define $B'$ using another such
function $\go'(n)$, it follows from \eqref{n5} (considering
$\go\bmin\go'$ and $\go\bmax \go'$) that
\begin{equation}
\label{ndiff}
\E|B\setdiff B'|/n\to0.
\end{equation}

Next, start with two distinct vertices $x$ and $y$, of types $i$
and $j$, say, and explore their components as above, again
stopping each exploration if we find $\go(n)$ vertices.
Assume for the moment that $\go(n)$ is small, say
$\go(n)=\log n$. The probability that during the truncated
exploration we find a connection between the two components is
$O(\go(n)^2/n)+o(1)=o(1)$. (Here we use the fact noted after
\eqref{nc}, that we are not likely to overshoot: with probability
$1-o(1)$, at every stage, even after stopping the exploration of
one component because it has become too large, the explored parts
of the components contain at most $\go(n)+\log n$ vertices.) As
before, fix $\eps>0$. For $n$ large enough, ignoring the
possibility of joining the truncated components of $x$ and $y$, we
can couple the two explorations as above with independent
branching processes (with $(1-3\eps)\kk$ for the lower bound) to
obtain
\begin{multline*}
\rhox{\ge\go(n)}\bigpar{(1-3\eps)\kk;i}
\rhox{\ge\go(n)}\bigpar{(1-3\eps)\kk;j}+o(1)
\\
\le
\P(x,y\in B)
\le
\rhox{\ge\go(n)}\bigpar{(1+\eps)\kk;i}
\rhox{\ge\go(n)}\bigpar{(1+\eps)\kk;j}+o(1).
\end{multline*}
Letting $\eps\to 0$, it follows, as above, that
\begin{equation}\label{rkij}
\P(x,y\in B) \to \rho(\kk;i)\rho(\kk;j);
\end{equation}
therefore, summing over all pairs of vertices $x,y$, we find that
\begin{equation*}
\begin{split}
\frac1{n^2} \E |B|^2
&
=\frac1{n^2} \sum_{x\neq y} \P(x,y\in B) + \frac1{n^2}\E|B|
\to \sum_{i,j=1}^r \mu_i\mu_j\rho(\kk;i)\rho(\kk;j)
=\rho(\kk)^2.
\end{split}
\end{equation*}
Combining this and \eqref{n5}, we see that $\Var(|B|/n)\to0$, and
thus that
\begin{equation}
\label{n6}
|B|/n\pto\rho(\kk).
\end{equation}
So far, we have assumed that $\go(n)$ was small. However, by
\eqref{ndiff}, having proved \eqref{n6} for one choice of $\go(n)$
it follows that \eqref{n6} holds for every choice of $\go(n)$
satisfying $\go(n)\to\infty$ and $\go(n)=o(n)$.

For any choice of $\go(n)$ with $\go(n)\to\infty$ and $\go(n)=o(n)$,
equation \eqref{n6}
gives the upper bound on the size $C_1(G_n)$ of
the largest component claimed in \eqref{nc}, since
$C_1(G_n)\le \max\{\go(n),|B|\}$ by \eqref{eL0}. In other words,
for any $\eps>0$,
\begin{equation}\label{ncub}
 |C_1(G_n))|/n\le \rho(\kk)+\eps
\end{equation}
holds \whp.

To obtain the matching lower bound, it remains to show that all
but $o_p(n)$ vertices in $B$ belong to a single
component. (We note that this is the only place where the
irreducibility of $\kk$ is needed.)
We first consider the simpler case where $\kk(i,j)>0$ for every
$i$ and $j$; we shall return to the general case afterwards. We
shall reveal the edges in $G_n$ in two rounds: given $0<\eps<1$,
we may take independent graphs $G_{n,0}$ and $G_{n,1}$ on the same
vertex set, with the distributions of $\gnx{(1-\eps)\kk}$ and
$\gnx{\eps\kk}$ respectively, so that $G_{n,0}\cup G_{n,1}\subseteq
G_n$. We shall think of $G_{n,0}$ as containing almost all the
edges of $G_n$, and $G_{n,1}$ as containing a few edges we
initially keep in reserve.

Recalling that $|\sss|=r$, set $\go(n)=rn^{2/3}$, and let $B_0$ be the
union of the big
components in $G_{n,0}$. From \eqref{n6}, applied with $(1-\eps)\kk$
in place of $\kk$, \whp\ we have
\begin{equation}\label{b0}
 |B_0|/n\ge \rho((1-\eps)\kk)-\eps.
\end{equation}
We claim that \whp\ all vertices of $B_0$ lie in a single component in $G_n$.
To see this, we condition on $G_{n,0}$
and use the random graph $G_{n,1}$: let $x,y\in B_0$ be vertices in distinct
components $C_x$, $C_y$ of $G_{n,0}$. As $\go=rn^{2/3}$, there are
$1\le i,j\le r$
such that $C_x$ contains a set $C_x^i$ of at least $n^{2/3}$ vertices
of type $i$, and $C_y$ a set $C_y^j$ of at least $n^{2/3}$ vertices
of type $j$. Now the probability that $G_{n,1}$ does not
contain a $C_x^i-C_y^j$ edge is
$(1-\kappa(i,j)/n)^{|C_x^i||C_y^j|}=\exp(-\Omega(n^{1/3}))=o(n^{-2})$.
As there are at most $n^2$ pairs to consider, it follows that \whp\
all vertices of $B_0$
lie in a single component of $G_n$, and hence, from \eqref{b0}, that
\begin{equation}\label{c1b0}
 |C_1(G_n)|/n\ge \rho((1-\eps)\kk)-\eps
\end{equation}
holds \whp.

The case when some $\kappa(i,j)$ may be zero is only slightly more
complicated. This time,
we replace $G_{n,1}$ by $r$ independent graphs $G_{n,l}$ with the distribution
of $\gnx{\eps\kappa/r}$. Given $C_x^i$ and $C_y^j$ as above,
the irreducibility of $\kappa$ implies that there is a sequence of types,
$i=i_1,i_2,\ldots,i_t=j$, such that $\kappa(i_l,i_{l+1})>0$
for all $l$.
As there are only $r$ types, we may suppose that $t\le r+1$ (note that
we may have $i=j$). Let $A_1=C_x^i$,
and, for $2\le l\le t-1$, let $A_l$ be the set of vertices of type $i_l$
adjacent to $A_{l-1}$ in $G_{n,l-1}$.
As $|A_1|=\Omega(n^{2/3})$ and $\kappa(i_1,i_2)>0$, the expected size
of $A_2$ is $\Omega(n^{2/3})$;
furthermore, from a standard Chernoff bound, with probability
$1-\exp(-\Omega(n^{2/3}))=1-o(n^{-2})$
we have
$|A_2|\ge \E|A_2|/2$, say.
Iterating, we see that
for some $c>0$ we have $|A_{t-1}|\ge cn^{2/3}$
with probability $1-o(n^{-2})$.
Finally, we find an edge in $G(n,t-1)$ from $A_{t-1}$ to $C_y^j$ with
very high probability, as above,
establishing \eqref{c1b0} in this case as well.

Letting $\eps\to0$ and using \refT{TappB}, the right-hand side of
\eqref{c1b0} tends
to $\rho(\kk)$, so \eqref{c1b0} proves the lower bound on $C_1(G_n)$
claimed in \eqref{nc}.
Combining this with the upper bound \eqref{ncub}, equation \eqref{nc} follows.

To prove \eqref{n7}, observe that
if we replace $\go(n)$ by a fixed number $k$ in the
argument leading to \eqref{n6} above,
and use \refT{TappC} instead of \refT{TappB},
we obtain \eqref{n7} instead of \eqref{n6}.
Note that this argument has not made use of the irreducibility of $\kk$ either.
\end{proof}

Note that the first part of \refL{L:step1} and \refT{Trho}
imply \refT{T2} in the case when
$\sss$ is finite, $\mu\set{i}> 0$ for every $i\in\sss$,
$\kkn=\kk$ for every $n$,
and $\kk$ is irreducible.

\medskip
We next consider the \rfin{} case in \refD{Drfin}; let us recall the
definition. 
A kernel $\kk$ on a vertex space $\vxs$ is \emph{regular finitary}
if $\sss$ may be partitioned into a finite number $r$ of \mucs s
$S_1,\ldots,S_r$
so that $\kk$ is constant on each $S_i\times S_j$.
A \mucs\ is a measurable set $A\subseteq\sss$ with $\mu(\partial A)=0$.
We next prove an extension of \refL{L:step1} to this \rfin{} case.

\begin{lemma}\label{L:step2}
Let $\kk$ be a \rfin{} kernel on a vertex space $\vxs$, and let
$G_n=\gnkx$. Then \eqref{n7} holds.
If $\kk$ is irreducible, then \eqref{nc} holds.
\end{lemma}

\begin{proof}
As noted in \refE{Efin}, the regular finitary case differs only in
notation from the finite case, so it suffices to prove that the conclusions
of \refL{L:step1} hold without the assumption that each $\mu\{i\}>0$.
Due to the generality of our model, we cannot just ignore sets
of measure zero; see \refR{Rmeas0}.

Using the notation of \refL{L:step1},
let us say that a type $i$ is \emph{bad} if $\mu_i=0$, and let
$\sss'\=\set{i\in\sss:\text{$i$ is not bad}}$.
Conditioning on the sequences $n_i$ as in the proof of \refL{L:step1},
if $i$ is a bad type then $n_i/n\to0$. Hence,
if we eliminate all vertices of bad type, we are left with a random
graph $G_n'=\gxx{n'}{(n'/n)\kk'}$,
where $\kk'$ is the restriction of $\kk$ to
$\sss'\times\sss'$ and $n'/n\to1$. It is easily seen that
$\rho(\kk')=\rho(\kk)$, and
$\rhogek(\kk')=\rhogek(\kk)$.
The expected degree of any vertex is at most $\max\kk<\infty$,
so the expected number of edges with at least one bad endpoint is
$o(n)$. Hence, \refL{LA} shows that for each fixed $k$,
$\E\bigpar{\ngek(G_n)-\ngek(G_n')}=o(n)$. Consequently,
\eqref{n7} holds for $G_n$ because it holds for $G_n'$.

Similarly, applying \eqref{nc} to $G_n'$, we see that if $\eps>0$
then
$C_1(G_n)/n \ge C_1(G_n')/n \allowbreak > \rho(\kk)-\eps$ \whp.
In the opposite direction, \eqref{n7} yields that
for every $\eps>0$ and $k\ge1$, \whp{}
$\ngek(G_n)/n \le \rhogek(\kk)+\eps$, and \eqref{eL0} implies
$C_1(G_n)/n \le \rhogek(\kk)+\eps$ \whp.
Further, by \refL{L00}, we have
$\rhogek(\kk)\downto\rho(\kk)$ as $\ktoo$.
Taking $k$ large enough, we find that $C_1(G_n)/n \le \rho(\kk)+2\eps$ \whp,
and \eqref{nc} follows.
\end{proof}

For technical reasons, we prove a slight extension of \refL{L:step2} to the
\qir{} \rfin{} case; \cf{} \refR{Rquasi}.

\begin{lemma}\label{L:step3}
Let $\kk$ be a regular finitary kernel on a vertex space $\vxs$.
Suppose that $\kk$ is \qir, i.e., that there is
a \mucs{} $\sss'\subseteq\sss$ such that $\kk$ restricted
to $\sss'$ is irreducible and $\kk=0$ off $\sss'\times\sss'$.
Then \eqref{nc} 
holds for $G_n=\gnkx$.
\end{lemma}

\begin{proof}
We may ignore all vertices with types not in $\sss'$, since they
will be isolated, and consider the restriction of our model to
$\sss'$. Note that we now have $n'$ vertices, with
$n'/n\pto\mu(\sss')$. The case $\mu(\sss')=0$ is trivial, and
otherwise we can consider the normalized measure $\mu/\mu(\sss')$ on
$\sss'$ and the kernel $\kk'=\mu(\sss')\kk$ on $\sss'\times\sss'$.
It is
easily checked that \refL{L:step2}
implies that
\eqref{nc} 
holds for $G_n$ in this case as well.
\end{proof}

It turns out that most of the work is behind us; roughly speaking,
to prove \refT{T2} we shall approximate with the \rfin{} case and use
\refL{L:step2}. There are some complications, as we must ensure
irreducibility of the approximations, but these
have already been dealt with in \refS{Sapprox}: when
$\kk_n$ is a graphical sequence of kernels with \qir\ limit $\kk$,
\refL{L:approx} gives us a sequence of \qir\ \rfin\ kernels
$\tkkm$ approaching $\kk$ from below. Furthermore, $\tkkm\le \kk_n$
when $n\ge m$, so we may and shall assume that
\begin{equation}\label{samuel}
 \gnx{\tkkm}\subseteq\gnx{\kk_n}
\end{equation}
for $n\ge m$.
This will allow us to apply \refL{LA}.

We are now in a position to prove our main results.
We start with the approximation result \refP{Pfinite}, which shows
that for many purposes we need only consider the finite-type case.

\begin{proof}[Proof of \refP{Pfinite}]
We use the kernels $\tkkm$ constructed in \refL{L:approx}.
{}From
\refL{L:approx}\ref{L:approxa}\ref{L:approxb}
and \refT{TappB}, if $m$ is large enough then
$\tkkm$ is \qir\ and
$\rho(\tkkm)\ge \rho(\kk)-\eps$.
Fix such an $m$.
We may regard the regular finitary kernel
$\tkkm$ as a kernel on a finite set
$\sss'$, so the graph
$\gnxx{\tkkm}$ has the required distribution for $G_n'$.
{}From  \refL{L:approx}\ref{L:approxc} and \eqref{samuel} we
can couple $G_n'$ and $G_n$
so that $G_n'\subseteq G_n$
whenever $n\ge m$.
Finally, from
\refL{L:step3}, we have $C_1(G_n')/n\pto\rho(\kk')$ as required.
\end{proof}

Next, it will be convenient to prove a restatement of \refT{T3}.

\begin{lemma}
   \label{LT3}
Let $(\kkn)$ be a graphical sequence of kernels on a vertex space
$\vxs$ with limit $\kk$,
and let $k\ge1$ be fixed. Then
\begin{equation}\label{E:LT3}
 \ngek\bigpar{\gnxx{\kk_n}}/n\pto\rhogek(\gk).
\end{equation}
\end{lemma}

Note that \refL{LT3} immediately implies \refT{T3},
as $\nk=\ngek-\ngex{k+1}$.

\begin{proof}
As before, to avoid clutter we suppress the dependence on $\vxs$,
writing $\gnx{\cdot}$ for $\gnxx{\cdot}$. We shall also
write $G_n$ for $\gnx{\kk_n}=\gnxx{\kk_n}$.
We may assume that $k\ge2$, since the case $k=1$ is trivial.
We use the $\tkkm$ constructed in \refL{L:approx}.

For each $m$, by \refL{L:step2} we have
\begin{equation}\label{emma}
\ngek\bigpar{\gnx{\tkkm}}/n\pto\rhogek(\tkkm).
\end{equation}
Let $\eps>0$.
Since $\rhogek(\tkkm)\to\rhogek(\kk)$ as $\mtoo$ by \refT{TappC}, we
can choose $m$ such that
$\rhogek(\tkkm)>\rhogek(\kk)-\eps$.
Using \eqref{samuel}, it follows from \eqref{emma} that \whp
\begin{equation}\label{jesper}
\ngek\bigpar{\gnx{\kk_n}}/n \ge \ngek\bigpar{\gnx{\tkkm}}/n
>\rhogek(\kk)-\eps,
\end{equation}
proving the lower bound claimed in \eqref{E:LT3}.

To prove the upper bound, consider any $\eta>0$.
By monotone convergence,
\begin{equation*}
\iint_{\sss^2}\tkkm(x,y)\dd\mu(x)\dd\mu(y)
\to
\iint_{\sss^2}\kk(x,y)\dd\mu(x)\dd\mu(y)
\end{equation*}
as \mtoo. Hence we may choose $m$ such that $\iint\kk-\iint\tkkm <
\eta/k$. We now fix this $m$.

By \eqref{t2b} and \refL{LE} (applied to the bounded kernel $\tkkm$),
we have
\begin{multline*}
\E \bigpar{e(G_n)-e\bigpar{\gnx{\tkkm}}}/n
\\
\to
\tfrac12\iint_{\sss^2}\kk(x,y)\dd\mu(x)\dd\mu(y)
-
\tfrac12\iint_{\sss^2}\tkkm(x,y)\dd\mu(x)\dd\mu(y)
<\eta/2k.
\end{multline*}
Hence, for $n$ large,
\begin{equation}\label{fewextra}
 \E \bigpar{e(G_n)-e\bigpar{\gnx{\tkkm}}}/n<\eta/2k,
\end{equation} and by
\eqref{samuel} and \refL{LA},
$\E \bigpar{\ngek(G_n)-\ngek\bigpar{\gnx{\tkkm}}}/n<\eta$.
Hence, for large $n$, using also \eqref{emma},
\begin{equation*}
  \begin{split}
\P\bigpar{\ngek(G_n)/n>\rhogek(\kk)&+2\eps}
\le
\P\bigpar{\ngek(G_n)/n>\rhogek(\tkkm)+2\eps}
\\&
\le
\P\bigpar{\ngek(\gnx{\tkkm})/n>\rhogek(\tkkm)+\eps}
\\&\qquad
+\P\bigpar{\bigpar{\ngek(G_n)-\ngek\bigpar{\gnx{\tkkm}}}/n>\eps}
\\&
<\eta+\eta/\eps.
  \end{split}
\end{equation*}
Letting $\eta\to0$, we find
$\ngek(G_n)\le \rhogek(\kk)+2\eps$ \whp, which together with
\eqref{jesper}
completes the
proof of the lemma.
\end{proof}

As noted above, \refT{T3} is just a reformulation of \refL{LT3}.
We are now ready to prove \refT{T2}.

\begin{proof}[Proof of \refT{T2}]
As noted in \refSS{SSgvs}, without loss of generality we may assume that 
$\vxs$ is a vertex space, rather than a generalized vertex space.
As above we write $G_n$ for $\gnxx{\kk_n}$, and consider the
approximating kernels 
$\tkkm$ constructed in \refL{L:approx}.

We first observe that
\eqref{rho}, \refT{Trho} and \refL{LPhi}\ref{lphia} 
imply that $\rho(\kk)<1$,
and that $\rho(\kk)>0$ if and only if $\norm{\tk}>1$.

Next we prove the upper bound \eqref{t2upper} on the size of the giant
component
of $G_n$. Fix $\eps>0$.
By \eqref{eL0} and \refL{LT3},
for every fixed $k\ge1$, \whp
\begin{equation}\label{david}
  C_1(G_n)/n \le k/n + \ngek(G_n)/n <\eps + \rhogek(\kk)+\eps.
\end{equation}
By \refL{L00},
as $\ktoo$,
$\rhogek(\kk)\downto\rho(\kk)$.
Hence we may choose $k$ so large that $\rhogek(\kk)<\rho(\kk)+\eps$,
and
\eqref{david} yields
$C_1(G_n)/n < \rhogek(\kk)+3\eps$ \whp, proving \eqref{t2upper}.

For \qir\ $\kk$, the lower bound on the size of the giant
component claimed in \eqref{t2bound} follows from \refP{Pfinite}.
Alternatively, we may argue as in the proof of \refL{LT3}.
Fix $\eps>0$.
By \refT{TappB},
$\rho(\tkkm)\to\rho(\kk)$ as $\mtoo$, so we
can choose $m$ such that
$\rho(\tkkm)>\rho(\kk)-\eps$, and then, by \refL{L:step3}
and \eqref{samuel}, \whp
\[
C_1\bigpar{\gnx{\kk_n}}/n \ge C_1\bigpar{\gnx{\tkkm}}/n
>\rho(\kk)-\eps.
\]
Together with \eqref{t2upper}, this proves
the convergence claimed in \eqref{t2bound}.

It remains to prove part (i) of \refT{T2}.
When $\norm{\tk}\le1$ we have $\rho(\kk)=0$, so \eqref{t2upper} yields
$C_1(G_n)=o_p(n)$, as required.
Suppose that $\norm{\tk}>1$. Recall that $\tkkm\upto\gk$ a.e.,
by \refL{L:approx}\ref{L:approxb}. It follows as in the proof of
\refL{L6} that $\norm{T_{\tkkm}}>1$ if $m$ is large enough. 
Let us fix such an $m$. As $\tkkm$ is of the regular finitary type,
there is a finite partition $\sss=\bigcup_{i=0}^r \sss_i$ of $\sss$
into \mucs s such that the restriction $\kk'_i$
 of $\tkkm$ to $\sss_i\times\sss_i$
is irreducible for $1\le i\le r$, and $\tkkm$ is zero \aex{}
off $\bigcup_{i=1}^r\sss_i\times\sss_i$. (This can be regarded
as an application of \refL{Lred} with $\sss$ finite. However,
the lemma is trivial in this case.)
As $T_{\tkkm}$ operates separately on each $\sss_i$, we have
$\norm{T_{\tkkm}}=\max_i\norm{T_{\kk'_i}}$, where 
$\norm{T_{\kk'_i}}$ is defined either on the generalized
ground space $(\sss_i,\mu|_{\sss_i})$, or,
equivalently, by extending $\kk'_i$ to $\sss\times\sss$ 
by setting $\kk'_i(x,y)=0$ if $x\notin \sss_i$ or $y\notin \sss_i$.
In particular, there is an $i$
with $\norm{T_{\kk'_i}}>1$.
Extending $\kk_i'$ to $\sss\times\sss$ as above,
$\kk'_i$ is a supercritical \qir\ kernel on $\sss$
of the regular finitary type, with
$\kk_n\ge\tkkm\ge\kk'_i$ for large $n$. Hence,
by \eqref{samuel} and \refL{L:step3}, we have
\[
 C_1(G_n)/n \ge C_1(\gnx{\kk'_i})/n \pto \rho(\kk'_i) >0,
\]
and $C_1(G_n)=\Theta(n)$ \whp{} follows, completing the proof.
\end{proof}

We next prove \refT{T2b}, showing
that the second largest component has size $o_p(n)$.

\begin{proof}[Proof of \refT{T2b}]
Let $(\kkn)$ be a graphical sequence of kernels on a (generalized) vertex space
$\vxs$ with \qir\ limit $\kk$, and $\go(n)$ a function
satisfying  $\go(n)\to\infty$ and $\go(n)=o(n)$. Our task is to show that
\begin{equation*}
\ngex{\go(n)}(G_n) \=
\sum_{j\ge1:\;C_j(G_n)\ge\go(n)}  C_j(G_n)
 =n\rho(\gk)+o_p(n).
\end{equation*}
Then \eqref{t2b1} follows by \refT{T2}.
In turn, \eqref{c2small} follows immediately, taking $\omega(n)=\log n$, say.
As before, we may assume that $\vxs$ is a vertex space. 

Let $\eps>0$.
For an upper bound on $\ngex{\go(n)}(G_n)$, fix a large $k$ such that
$\rhogek(\kk)<\rho(\kk)+\eps$.
For large $n$ we have $\go(n)>k$ and thus by \refL{LT3} \whp{}
$\ngex{\go(n)}(G_n)/n\le \ngek(G_n)/n < \rho(\kk)+2\eps$.

For a lower bound, assume that $\rho(\kk)>0$. Then, by \refT{T2},
\whp{}
$C_1(G_n)>\tfrac12\rho(\kk)n>\go(n)$, so by \refT{T2} again, \whp{}
$\ngex{\go(n)}(G_n)/n\ge C_1(G_n)/n>\rho(\kk)-\eps$. This is trivially
true if $\rho(\kk)=0$ too.

Since $\eps$ was arbitrary, the proof is complete.
\end{proof}

We now turn to a result giving the distribution of the
types of the vertices
making up the giant component; to state this, we need some more definitions.

Let $\cC_1(G_n)$ be the largest component
of $G_n=\gnxx{\kk_n}$, i.e., the component with most vertices,
chosen by any rule if there is a tie. (Thus, if a sequence $(G_n)$
has a unique giant component, then $\cC_1(G_n)$ is this
giant component.) Let $\nuni\=\frac1n
\sum_{i\in\cC_1(G_n)}\gd_{x_i}$ be the random measure with total
mass $C_1(G_n)/n$ that describes the distribution of the points
$x_i$ corresponding to the vertices in the largest component. We
equip the space of finite positive Borel measures on $\sss$ with
the weak topology; see \refSA{Smeasure}. 
As before, we write $\rhokk$ for the
function defined by $\rhokk(x)=\rho(\kk;x)$.

\begin{theorem}
   \label{T1A}
Let $(\kkn)$ be a graphical sequence of kernels on a (generalized)
vertex space 
$\vxs$ with \qir{}
limit $\kk$.
Then $\nuni\pto\mu_\kk$ in the space of finite measures on $\sss$
with the weak topology, where
$\mu_\kk$ is the measure on $\sss$ defined by
$\dd\mu_\kk=\rhokk\dd\mu$.
In other words, for every \mucs\ $A$,
\begin{equation}\label{t1a}
\nuni(A)=\frac1n  \#\set{i\in \cC_1(G_n):x_i\in A}
\pto \mu_\kk(A)=\int_A\rho(\kk;x)\dd\mu(x),
\end{equation}
where $G_n=\gnxx{\kkn}$. 
Furthermore, if $f:\sss\to \bbR$ is continuous $\mu$-\aex{} and
satisfies
\begin{equation}\label{emmsan}
  \frac1n \sum_{i\in V(G_n)} |f(x_i)|
\pto \int_\sss |f|\dd\mu
<\infty,
\end{equation}
then 
\begin{equation}\label{t1f}
  \frac1n \sum_{i\in \cC_1(G_n)} f(x_i)
\pto \int_\sss f\dd\mu_\kk=\int_\sss f(x)\rho(\kk;x)\dd\mu(x).
\end{equation}
In particular, \eqref{t1f} holds for every bounded and $\mu$-\aex{} continuous 
$f:\sss\to \bbR$.
\end{theorem}

Condition \eqref{emmsan} is very natural and often easy to verify;
for example, if $\vxs$ is a vertex space in which the $x_i$ are \iid,
as in \refE{EU}, 
or a generalized vertex space in which $\xs$ is a Poisson process, as
in \refE{EPo}, then 
\eqref{emmsan} holds for every integrable $f$ by the law of large
numbers.
Similarly, if $\sss=(0,1]$ and $x_i=i/n$, then \eqref{emmsan}
holds for every decreasing integrable positive $f$.
Note that some restriction on $f$ is needed for \eqref{t1f};
it is not hard to construct 
an example where \eqref{t1f}  fails, and so does \eqref{emmsan}.

\begin{proof}
We begin by proving the first statement.
We proceed in several steps, as before. Arguing as in \refSS{SSgvs},
we may assume 
without loss of generality that $\vxs$ is a vertex space.

First we assume that the conditions of \refL{L:step1} are satisfied:
$\sss$ is finite, $\kk$ is fixed and irreducible, and $\mu_i=\mui>0$
for every $i$.
We use the notation of the proof of \refL{L:step1}; in particular,
$\go$ is some
function with $\go(n)\to\infty$ and $\go(n)=o(n)$, and $B$
is the set of vertices of $G_n=\gnkx$ in `big' components, i.e.,
components of order
at least $\go(n)$.

Let $V_i$ be the set of vertices of type $i$. The arguments leading to
\eqref{n6}
in the proof of \refL{L:step1} yield also
$|B\cap V_i|/n\pto\rho(\kk;i)\mu_i$; see \eqref{n4} and \eqref{rkij}.

If $\rho(\kk)>0$, then the conclusion \eqref{nc} of \refL{L:step1}
implies that \whp{}
$\cC_1(G_n)\subseteq B$, and thus \eqref{nc} and \eqref{n6} imply that
$|B \setdiff \cC_1(G_n)|/n\pto0$. This is clearly true when $\rho(\kk)=0$
too, and implies that
$|\cC_1(G_n)\cap V_i|/n\pto\rho(\kk;i)\mu_i$ for every $i$, which is
exactly \eqref{t1a}.

The result extends to the case when some $\mu_i=0$ as before.
Thus \eqref{t1a} holds in the irreducible \rfin{} case considered
in \refL{L:step2},
provided
$A$ is one of the sets $\sss_i$ in the partition or a union of such
sets. In fact, $A$ may be any \mucs, since we may replace the
partition $\set{\sss_i}$ by $\set{\sss_i\cap A,\,\sss_i\setminus A}$,
noting that all parts are \mucs s.
Similarly, the extension to the \qir\ case is immediate, as in
\refL{L:step3}.

We now turn to the general case.
Note that $\muk(\sss)=\rho(\kk)$.
Assume that $\rho(\kk)>0$; otherwise the result is trivial (with
$\muk=0$) by \refT{T2}.

Fix a \mucs{} $A$.
Use a sequence of partitions
$\Pm$ as in \refL{LP},
and consider the finitary approximation
$\tkkm$ given by \refL{L:approx} for some $m$.
Let $\nunmi$ be the random measure $\nuni$ defined for $\gnx{\tkkm}$.
(As before, we suppress the dependence on $\vxs$.)
By the finitary
case completed above,
\begin{equation}\label{cecilia}
\nunmi(A)\=\frac1n  \#\bigset{i\in \cC_1\bigpar{\gnx{\tkkm}}:x_i\in A}
\pto \mu_{\tkkm}(A)
\end{equation}
for every fixed $m$.

Let $\eps>0$ and choose $m$ so large that $\rho(\tkkm)>\rho(\kk)-\eps$
and
$\rho(\tkkm)>0$ (see \refT{TappB}).
Then, applying Theorem \ref{T2} to $\tkkm$
and \eqref{c2small} of Theorem \ref{T2b} to $\kk$,
\whp{} $C_1(\gnx{\tkkm}) > \tfrac12\rho(\tkkm) n >
C_2(G_n)$. Recalling the coupling
\eqref{samuel}, it follows that
the largest component of $\gnx{\tkkm}$ is contained in the largest
component of $G_n$,
\ie{}, $\cC_1(\gnx{\tkkm})\subseteq \cC_1(G_n)$, and thus $\nunmi\le
\nuni$.
Consequently, from \eqref{cecilia}, \whp
\begin{equation}
  \label{sofie}
\nuni(A)
\ge \nunmi(A)
\ge
\mu_{\tkkm}(A)-\eps
\ge
\mu_\kk(A)-2\eps,
\end{equation}
because
$\mu_\kk(A)-\mu_{\tkkm}(A) \le \mu_\kk(\sss)-\mu_{\tkkm}(\sss)
=\rho(\kk)-\rho(\tkkm)<\eps$.

Since $\CA$ also is a \mucs, we may replace $A$ by $\CA$ in
\eqref{sofie} and obtain that \whp
\begin{equation*}
\nuni(\CA)
\ge
\mu_\kk(\CA)-2\eps.
\end{equation*}
Since $\nuni(\sss)=C_1(G_n)/n$ and $\muk(\sss)=\rho(\kk)$, this and
\refT{T2} show that \whp
\begin{equation*}
\nuni(A)=C_1(G_n)/n-\nuni(\CA)
\le \rho(\kk)+\eps -\mu_\kk(\CA)+2\eps
=\muk(A)+3\eps.
\end{equation*}
This and \eqref{sofie} yield $\nuni(A)\pto\muk(A)$, so we have shown
that \eqref{t1a} holds
for this $A$. 
We have shown that \eqref{t1a} holds for an arbitrary \mucs{} $A$,
which yields
$\nuni\pto\mu_\kk$ by \refL{Lmeas}.

Turning to the second part of the lemma, 
note that the left-hand sides of \eqref{emmsan} and \eqref{t1f} are
equal to $\int |f|\dd\nu_n$ and $\int f\dd\nuni$, respectively.
If $f$ is bounded and $\mu$-\aex{} continuous, \refL{Lmeas} thus shows that 
these relations follow from \eqref{nunA2} and \eqref{t1a}, respectively.

To deduce \eqref{t1f} from \eqref{emmsan} for unbounded $f$,
we use the truncations $f_M\=(|f|\bmin M)\sign(f)$.
Let $\eps>0$.
By monotone convergence, $\int|f_M|\dd\mu\to\int|f|\dd\mu$ as
$M\to\infty$. Thus, we can choose $M$ such that 
$\int|f_M|\dd\mu>\int|f|\dd\mu-\eps$.
Since \eqref{emmsan} holds for bounded $\mu$-\aex{} continuous
functions, it holds for 
$f_M$, so
\begin{equation*}
\frac1n \sum_{i\in V(G_n)}\bigpar{ |f(x_i)|-|f_M(x_i)|}
\pto \int_\sss |f|\dd\mu - \int_\sss |f_M|\dd\mu <\eps.
\end{equation*}
Hence the left-hand side is at most $\eps$ \whp.
Consequently, 
\whp
	\begin{multline*}
\left|\frac1n \sum_{i\in \cC_1(G_n)}f(x_i)
-\frac1n \sum_{i\in \cC_1(G_n)}f_M(x_i)\right|
\le	
\frac1n \sum_{i\in \cC_1(G_n)}\bigl|f(x_i)-f_M(x_i)\bigr|
\\
\le	
\frac1n \sum_{i\in V(G_n)}\bigl|f(x_i)-f_M(x_i)\bigr|
=
\frac1n \sum_{i\in V(G_n)}\bigpar{ |f(x_i)|-|f_M(x_i)|}
\le\eps,
	\end{multline*}
with the first two inequalities holding unconditionally.
Note that
\begin{equation*}
\left|\int_\sss f\dd\muk - \int_\sss f_M\dd\muk\right| 
\le
\int_\sss |f-f_M|\dd\mu 
=
\int_\sss |f|\dd\mu - \int_\sss |f_M|\dd\mu 
<\eps.
\end{equation*}
Since $\eps>0$ is arbitrary, and \eqref{t1f} holds for each $f_M$,
relation \eqref{t1f} for $f$ follows by a standard
$3\eps$-argument. 
\end{proof}

\section{Edges in the giant component}\label{Sgiantedges}

The main aim of this section is to prove \refT{Tedges}, which claims
that if $\kkn$ is a graphical sequence of kernels on a (generalized) vertex
space $\vxs$ with \qir\ limit $\kk$,
then $\frac{1}{n}e(\cC_1(\gnxx{\kkn}))\pto \edgeno(\kk)$,
where $\edgeno(\kk)$ is defined in \eqref{Zdef} as
\begin{equation}\label{Zdef2}
 \edgenokk \=
 \frac12
 \iint_{\sss^2}\kk(x,y)\bigpar{\rho(\kk;x)+\rho(\kk;y)-\rho(\kk;x)\rho(\kk;y)}
 \dd\mu(x)\dd\mu(y).
\end{equation}

Before turning to the proof, we briefly examine the behaviour of $\edgenokk$,
giving two alternative formulae for $\edgenokk$,
together with upper and lower bounds in terms of $\rhokk$.

{}From the symmetry of $\kk$ and the definition \eqref{tk} of $\tk$,
\eqref{Zdef2} 
is equivalent to
\[
 \edgenokk
 = \ints (1-\rho_\kk/2) T_\kk\rho_\kk \dd\mu,
\]
where, as usual,
$\rho_\kk$ is the function defined by $\rho_\kk(x)\=\rho(\kk;x)$.
By relation \eqref{Phirho} of \refT{Trho} and the definition
of $\Phik$ in \eqref{Phik}, it follows that
\begin{equation}\label{Zform2}
 \edgenokk
 = \ints\Bigpar{1-\frac{\rho(\kk;x)}2}\ln\parfrac{1}{1-\rho(\kk;x)}\dd\mu(x).
\end{equation}

Writing $G_n$ for $\gnxx{\kk_n}$, note that the assumptions
of \refT{Tedges} include convergence of the expectation of $e(G_n)/n$. 
As shown in \refP{PE}, an easy consequence of these assumptions is that
\begin{equation}\label{ecnv}
 e(G_n)/n \pto \frac12\iint_{\sss^2}\kk.
\end{equation}
In the light of \eqref{ecnv}, relation \eqref{eedge}
is equivalent to the assertion that
number of edges not in the giant component is
\[
\frac n2\iint_{\sss^2} (1-\rho(\kk;x))\kk(x,y)(1-\rho(\kk;y))
  \dd\mu(x)\dd\mu(y)+o_p(n).
\]

In any connected graph,
the number of edges is at least the number of
vertices minus 1;
hence $\edgenokk\ge\rho(\kk)$. In fact,
\refT{Trho} has
the following simple consequence.
\begin{proposition}\label{Pedges}
Let $\kk$ be a kernel on a (generalized) ground space $(\sss,\mu)$. Then 
\begin{equation*}
  \rho(\kk)
\le
\edgenokk
\le\tfrac12(\norm{\tk}+1)\rho(\kk)
\le\norm{\tk}\rho(\kk).
\end{equation*}
Furthermore, the first two inequalities are strict
when $\rho(\kk)>0$.
\end{proposition}

\begin{proof}
If $0<s<1$, then $s<(1-s/2)\ln(1/(1-s))$, as is easily verified by
computing the Taylor series.
Thus,
\[
 \rhokkx \le  \Bigpar{1-\frac{\rho(\kk;x)}2}\ln\parfrac{1}{1-\rho(\kk;x)},
\]
with strict inequality when $\rhokkx>0$. Integrating with
respect to $\mu$, the left-hand side
becomes $\rho(\kk)$, while, from \eqref{Zform2}, the right-hand side
becomes $\edgenokk$.  Thus $\rho(\kk)\le\edgenokk$, with strict
inequality if $\rho(\kk)>0$.

In the other direction,
if $0<s<1$, then $(1-s)\ln(1/(1-s))<s-s^2/2$, as can again
be verified by computing the Taylor series.
Hence,
\[
 \Bigpar{1-\frac s2}\ln\parfrac{1}{1-s} < s-\frac{s^2}2
  + \frac s2\ln\parfrac{1}{1-s}.\]
Substituting $s=\rhokkx$ and integrating, it follows that
\[
 \edgenokk
 \le \ints\left(
  \rhokkx-\frac12\rhokkx^2+\frac12\rhokkx\ln\parfrac{1}{1-\rhokkx}
 \right)\dd\mu(x),
\]
with strict inequality when $\rhokkx>0$.
Writing $\rhokk$ for the function
defined by $\rhokk(x)\=\rhokkx$, from \eqref{Phirho}
and the definition \eqref{Phik} of $\Phik$, we have
\[
 \ln\parfrac{1}{1-\rhokkx} = (\tk\rhokk)(x).
\]
It follows that
\begin{align*}
  \edgenokk
&\le
\rho(\kk)-\tfrac12\innprod{\rhokk,\rhokk}+\tfrac12\innprod{\rhokk,\tk\rhokk}
\\&
\le
\rho(\kk)+\tfrac12\bigpar{\norm{\tk}-1}\ints\rhokk^2\dd\mu(x)
\le
\tfrac12\bigpar{\norm{\tk}+1}\rho(\kk),
\end{align*}
with strict inequality unless $\edgenokk=\rho(\kk)=0$.
\end{proof}

Our proof of \refT{Tedges}
will be very similar to that of \refT{T2}, except that
we need to consider
certain branching process expectations $\sigma(\kk)$ and $\sgek(\kk)$ in place
of $\rho(\kk)$ and $\rhogek(\kk)$.
In preparation for the proof, we shall relate $\edgenokk$ to the branching
process $\bpk$ via $\sigma(\kk)$.
As before, we assume that $\kk$ is
a kernel on $(\sss,\mu)$
with $\kk\in L^1$; in particular, it is convenient here to normalize
so that $\mu(\sss)=1$.

Let $A$ be a Poisson process on $\sss$,
with intensity given by
a finite measure
$\lambda$, so that $A$ is a random
multi-set on $\sss$.
If $g$ is a bounded measurable function on multi-sets on $\sss$, it is
easy to see that
\begin{equation}\label{eX}
 \E\left(|A|g(A)\right) = \int_{\sss} \E g(A\cup\{y\}) \dd\lambda(y).
\end{equation}
(This is a simple consequence of the well-known fact that the Palm
  distribution equals the distribution of   $A\cup\set y$.
To show \eqref{eX} directly, note that we may construct $A$ as
  follows: first
decide the total number $N$ of points in $A$, according to a Poisson
$\Po(c)$ distribution with mean $c=\gl(\sss)$. Then let
$(a_i)_{i=1}^N$ be a sequence of \iid{} random points of $\sss$,
each distributed according to the normalized form $\gl/c$
of $\gl$, and take $A=\{a_1,\ldots,a_N\}$. Let $\nu$ be the measure
(on finite sequences of points in $\sss$) associated to $(a_i)_{i=1}^N$,
and let $\nu'$ be the measure with density $N\dd\nu$.
Recalling that
if $Z$ has a $\Po(c)$ distribution, then $k\P(Z=k)=c\P(Z-1=k)$,
we find that
$\nu'/c$
may be constructed by taking $N-1$ to have a $\Po(c)$ distribution,
and then taking the $a_i$ \iid{} as before, or, equivalently,
by constructing a sequence
according to $\nu$ and appending a new random point with the distribution
$\gl/c$.
Neglecting the order of the points, \eqref{eX} follows.)

Let $X(x)$ denote the first generation of the branching process $\bpkx$.
Thus $X(x)$ is given by a Poisson process on $\sss$ with intensity
$\kk(x,y)\dd\mu(y)$.
Suppose that \eqref{bb1} holds, so $X(x)$ is finite.
Let $\sigma(\kk;x)$ denote the expectation of $|X(x)|\ett[|\bpkx|=\infty]$,
recalling that under the assumption \eqref{bb1}, the branching process
$\bpkx$ dies out if and only if $|\bpkx|<\infty$.
Then
\begin{align*}
 \int_\sss \kk(x,y)\dd\mu(y) &- \sigma(\kk;x)
 = \E\big(|X(x)|\ett[|\bpkx|<\infty]\big) \\
 &= \E\left( |X(x)| \prod_{z\in X(x)} (1-\rho(\kk;z)) \right) \\
 &= \int_\sss \kk(x,y) (1-\rho(\kk;y))
    \E\left(\prod_{z\in X(x)} (1-\rho(\kk;z)) \right) \dd\mu(y)\\
 &= \int_\sss \kk(x,y) (1-\rho(\kk;y)) (1-\rho(\kk;x)) \dd\mu(y).
\end{align*}
Here the penultimate step is from \eqref{eX}; the last step uses the
fact that the branching process dies
out if and only if none of the children of the initial particle survives.
Writing $X$ for the first generation of $\bpk$, let
\[
 \sigma(\kk) \= \E\left( |X|\ett[|\bpk|=\infty] \right)
=\int_{\sss} \gs(\kk;x)\dd\mu(x).
\]
Then, integrating over $x$ and subtracting from $\iint\kk(x,y)$, we obtain
\begin{equation}\label{sform}
 \sigma(\kk) =
 \iint_{\sss^2}\kk(x,y)\big(1-(1-\rho(\kk;x))(1-\rho(\kk;y))\big)
  \dd\mu(y)\dd\mu(x),
\end{equation}
i.e., $\sigma(\kk)=2\edgenokk$, where $\edgenokk$ is
defined in \eqref{Zdef}.

\begin{lemma}\label{el1}
Let $\kk$ be a \qir\ kernel on a ground space $(\sss,\mu)$,
with $\kk\in L^1$.
If\/ $(\kkn)_1^\infty$ 
is a sequence of kernels that increase to $\kk$ a.e., then
$\sigma(\kkn)\to \sigma(\kk)<\infty$.
\end{lemma}

\begin{proof}
This is immediate from \refT{TappB}(i), \eqref{sform}, the fact that
$\iint\kk<\infty$, and dominated convergence.
\end{proof}

As we shall see next, $\sigma(\kk)$ is the limit of the expectations
\[
 \sgek(\kk) \= \E\left( |X|\ett[|\bpk|\ge k] \right).
\]

\begin{lemma}\label{el2}
With $\kk\in L^1$ fixed,
\begin{equation}\label{sup}
 \sgek(\kk) \downto \sigma(\kk)\quad \text{as } k\to\infty.
\end{equation}
\end{lemma}

\begin{proof}
We have $|X| \ge |X|\ett[|\bpk|\ge k] \downto
|X|\ett[|\bpk|=\infty]$. As $\E|X|=\iint\kk(x,y)<\infty$,
the result follows by dominated convergence.
\end{proof}

Using the above lemmas we can prove \refT{Tedges}. As the argument is
very similar to that for
\refT{T2}, we give only an outline.

\begin{proof}[Proof of \refT{Tedges}]
As usual, we may assume without loss of generality that $\vxs$ is a
vertex space.  
Let $M_{\ge k}(G)$ denote the number of edges of a graph $G$ that lie
in components of order
at least $k$.

We start with the case when $\sss$ is finite
and $\kk$ is irreducible,
writing $G_n$ for $\gnxx{\kkn}$.
Let $d(x)$ denote the degree of a vertex $x$ of $G_n$. Using the local
coupling of the neighbourhood
of a random vertex $x$ to the branching process $\bpk$ described in
the proof of \refL{L:step1},
considering $\E(d(x)\ett[x\in B])$ in place of $\P(x\in B)$, the proof
of \refL{L:step1} yields
the relations
\begin{equation}\label{enc}
 2e(\cC_1(G_n))/n\pto \sigma(\kk)
\end{equation}
and
\begin{equation}\label{enck}
 2M_{\ge k}(G_n)/n \pto \sgek(\kk),
\end{equation}
corresponding to \eqref{nc} and \eqref{n7}.
As before, the same formulae in the \qir{}
regular finitary setting of \refL{L:step3} follow.

To complete the proof, we consider the approximating kernels $\tkkm$
constructed in \refL{L:approx}.
By \refL{L:approx}\ref{L:approxb} and \refL{el1} we have
$\sigma(\tkkm)\to\sigma(\kk)$.
Applying \eqref{enc} to $\tkkm$ and using the coupling
$\gnxx{\tkkm}\subseteq\gnxx{\kkn}$, $n\ge m$,
it follows that
for any $\eps>0$,
\begin{equation}\label{ebd1}
 e(\cC_1(G_n))/n \ge \sigma(\kk)/2 -\eps
\end{equation}
holds \whp. This is exactly the lower bound claimed in \eqref{eedge}.

For the upper bound, we claim first that, for each fixed $k$,
\[
 2M_{\ge k}(\gnxx{\kkn})/n \pto \sgek(\kk).
\]
The argument is exactly as for \eqref{E:LT3}, except that in place of
\eqref{fewextra} we show that there is an $m$ for which
$\E(e(G_n)-e(G(n,\tkkm)))/n<\eta/(2k^2)$, and in place of \refL{LA} we
use the fact
that, for $k\ge 1$,
adding an edge to a graph cannot change $M_{\ge k}$ by more than
$2\binom{k-1}{2}+1\le k^2$.
The rest of the proof is as for \refT{T2}, using\
\[
 e(\cC_1(G_n))/n \le k^2/n + M_{\ge k}(G_n)/n
\]
in place of \eqref{david} and \refL{el2} instead of \refL{L00}.
\end{proof}

\section{Stability}\label{Sstable}

This section is devoted to the proof of
the `stability' result, \refT{TrobustNEW},
which states that deleting a few vertices and their incident
edges, and then adding or deleting a few edges, does not change
the size  of the giant component of $G_n=\gnxx{\kkn}$ significantly.
As usual, 
without loss of generality we may restrict our attention to the case
where $\vxs$ is a vertex space; we shall return to this later.
For the moment, we shall ignore vertex deletion; our aim is thus to
prove the following special case of \refT{TrobustNEW}. 

\begin{theorem}
  \label{TrobustOLD}
Let $(\kkn)$ be a graphical sequence of kernels on a vertex space 
$\vxs$ with irreducible 
limit $\kk$, and let $G_n=\gnxx{\kk_n}$. For every
$\eps>0$ there is a $\gd>0$ (depending on $\kk$) such that, \whp{},
\begin{equation}\label{Estab}
 (\rho(\kk)-\eps)n \le C_1(G_n') \le (\rho(\kk)+\eps)n
\end{equation}
for every graph $G_n'$ on $V(G_n)=[n]$ 
with 
$e(G_n'\setdiff G_n)\le \gd n$.
\end{theorem}

We shall see later (at the end of \refSS{SSLMcD}) that
\refT{TrobustNEW} follows.  
As noted in \refSS{SSstable}, to prove \refT{TrobustOLD}
it suffices to consider separately the cases where edges are added and
where edges are deleted.
More precisely, as $G_n'\cap G_n\subseteq G_n'\subseteq G_n'\cup G_n$, it
suffices to prove the upper bound
in \eqref{Estab} for $G_n'\supseteq G_n$, and the lower bound for
$G_n'\subseteq G_n$. 

\medskip
The upper bound is easy. Indeed, by \refL{L00},
$\rhogek(\kk)\downto\rho(\kk)$ as \ktoo. Thus, given $\eps>0$,
we may choose $k$ such that
$\rhogek(\kk)\le \rho(\kk)+\eps/3$. By \refL{LT3}, \whp{}
$\ngek(G_n) \le (\rho(\kk)+\eps/2)n$.
Taking $\gd=\eps/4k$, it follows by \refL{LA} that \whp{}
$\ngek(G_n') \le (\rho(\kk)+\eps)n$, which implies
the upper bound in \eqref{Estab}.

\medskip
For the lower bound, our aim is to show that \whp{}
\begin{equation}\label{esl}
 C_1(G_n-E) \ge (\rho(\kk)-\eps)n
\end{equation}
for every $E\subseteq E(G_n)$ with $|E|\le\gd n$.

We may assume that $\rho(\kk)>0$, as otherwise there is nothing to prove.
As in the proof of \refT{T2}, it suffices to consider the
regular finitary case; in fact, given
$\eps>0$, by \refP{Pfinite} there is a vertex space $\vxs'$ with finite
type space and a \qir\ kernel $\kk'$ on $\vxs'$ with
$\rho(\kk')>\rho(\kk)-\eps/2$ such that we may consider $\gnxxp{\kk'}$
as a subgraph of $G_n$.
It suffices to prove that
there is a $\gd>0$ such that removing at most $\gd n$ edges from
$\gnxxp{\kk'}$ leaves \whp\ a graph
with a component of order at least $(\rho(\kk')-\eps/2)n$. Replacing
$\eps$ by $2\eps$, this is exactly \eqref{esl},
but with $\gnxx{\kkn}$ replaced by $\gnxxp{\kk'}$.
Thus we may assume that
$G_n=\gnkx$, where $\kk$ is a \qir{} kernel
on a finite set $\sss=\{1,2,\ldots,r\}$.
In fact, by rescaling,
as in the proof of \refL{L:step3},
we may assume that $\kk$ is irreducible.
Finally, as in the proof of \refL{L:step2}, we may assume that
$\mu\{i\}>0$ for every $i$,
as there are $o_p(n)$ edges incident with types $i$ with
$\mu\{i\}=0$. In other words,
we may assume the setting of \refL{L:step1}. We shall do so for the
rest of this section; thus $G_n=\gnkx$, where $\vxs=(\sss,\mu,\xss)$
is a vertex space, and
\begin{equation}\label{sassump}
 \sss=\{1,2,\ldots,r\}, \ \mu\{i\}>0\ \forall i,\ \hbox{$\kk$ is
 irreducible, and $\norm{\tk}>1$}.
\end{equation}

In a paper studying the bisection width of sparse random graphs,
Luczak and McDiarmid~\cite{LMcD} proved \eqref{esl} for the
Erd\H os-R\'enyi case, where $|\sss|=1$ or $\kk$ is constant.
Their proof adapts easily to the finite-type
case, from which, as shown above, \refT{TrobustOLD} follows. We present
this proof in \refSS{SSLMcD}.

A different, perhaps more natural, approach to proving \eqref{esl} is
to work with the branching process $\bpk$,
using the coupling of vertex neighbourhoods in $G_n$
with $\bpk$ to reduce \eqref{esl}
to an equivalent statement for the two-core, \refL{l2core} below.
The latter statement
has a very simple proof in the uniform case. We present this approach
here, in \refSS{SS2core} below,
because the intermediate results, relating properties of the two-core
to the branching process, are likely to be of interest in their own right.
Unfortunately, while \refL{l2core} can be proved
in the general case by branching process methods, our proof is rather
complicated. 
As the result follows from \refT{TrobustOLD}, which can be proved
more simply by the method of Luczak and McDiarmid, we omit the proof.
A reader interested only in the proof of \refT{TrobustOLD} can safely
omit \refSS{SS2core}.

\subsection{Counting cuts in the giant component}\label{SSLMcD}

In this subsection we prove \refT{TrobustOLD}, and then deduce
\refT{TrobustNEW}. 
Apart from the straightforward adaptations to non-constant $\kk$,
the argument for \refT{TrobustOLD} is that of  Luczak and
McDiarmid~\cite{LMcD}. 
We start with a deterministic lemma whose statement and proof are
taken verbatim from \cite{LMcD}.

\begin{lemma}\label{LLMcD}
For any $\eps>0$, there exist $\eta_0=\eta_0(\eps)>0$ and $n_0$
such that the following holds. For all $n\ge n_0$, and for
all connected graphs $G$ with $n$ vertices, there are at most
$(1+\eps)^n$ bipartitions of $G$ with at most $\eta_0n$ cross
edges.
\end{lemma}
\begin{proof}
Let $T$ be an arbitrary spanning tree of $G$.
Any $2$-partition $S$, $\overline S$ of $T$ is determined uniquely
by the corresponding set of cross edges, together with the
specification for each cross edge of which of its endpoints is in $S$.
For as $T$ is connected, the cross edges specify a nonempty
subset $S^*$ of $S$, and then $S$ is the set of vertices $v$
such that there is a path from $v$ to one of the vertices in $S^*$
where this path does not use any of the cross edges.
(If $v\in \overline S$ then no path from $v$ to $S^*$
can avoid the cross edges, and if $v\in S$ then any shortest
path from $v$ to $S^*$ avoids the cross edges.)
Hence, since $T$ has $n-1$ edges, the number of $2$-partitions
of $T$ (and hence also of $G$) with at most $\eta n$ cross edges
is no more than
\[
 \sum_{j\le \eta n} 2^j\binom{n}{j} =
 O\left( n2^{\eta n}\left( (1-\eta)^{1-\eta}\eta^\eta \right) ^{-n} \right),
\]
assuming $\eta\le 1/2$.
Now let $\eps>0$. As $\eta\to 0$, $2^\eta/((1-\eta)^{1-\eta}\eta^\eta)\to 1$.
Hence, for $\eta$ sufficiently small and $n$ sufficiently large, there are
at most $(1+\eps)^n$ partitions with at most $\eta n$ cross edges.
\end{proof}

Recall our assumptions \eqref{sassump}, that $\sss=\{1,2,\ldots,r\}$,
$\mu\{i\}>0$ for every $i$, $\kk$ is irreducible and
$\norm{\tk}>1$. As usual, we condition on $\xs$, so we may assume that
$\xs$ is deterministic for every $n$, so there are $n_i$ vertices of
type $i$, with $n_i/n\to\mu\{i\}>0$ as $\ntoo$.

The main additional ingredient needed to adapt the proof of \cite{LMcD}
to non-constant kernels is the following simple lemma.

\newcommand\grklet{\theta}
\begin{lemma}\label{lexp}
Suppose that the assumptions \eqref{sassump} hold.
For any $\eps>0$ there is a $\grklet=\grklet(\kk,\eps)>0$ 
with the following property.
If $n$ is large enough then,
whenever $V_1$, $V_2$ are disjoint sets of at least $\eps n$
vertices of $G_n=\gnk$ such that $V_1\cup V_2$
contains at least $\eps n$ vertices of each type,
the expected number of edges
from $V_1$ to $V_2$ in $G_n$ is at least $\grklet n$.
\end{lemma}
\begin{proof}
We assume that $n\ge \max\kk$.
Let $\eps'=\min\{\eps/r,\eps/2\}$,
and let
\[
 \grklet=(\eps')^2\min\{\kk(i,j):\kk(i,j)>0\}>0. 
\]
There are types $i$ and $j$
such that $V_1$ contains at least $\eps'n$ vertices of type $i$, and $V_2$
at least $\eps'n$ vertices of type $j$.
As $\kk$ is irreducible, there is a sequence $i=i_0,i_1,\ldots,i_t=j$
such that $\kk(i_s,i_{s+1})>0$ for each $s$.
For each $i_s$, from our condition on $V_1\cup V_2$,
one or both of $V_1$ and $V_2$ must contain at least $\eps'n$ vertices
of type $i_s$. It follows that for some $s$, $V_1$ contains
at least $\eps'n$ vertices of type $i_s$, and $V_2$ contains
at least $\eps'n$ vertices of type $i_{s+1}$.
But then the expected number of edges from $V_1$ to $V_2$ is at least
$(\eps' n)^2 \kk(i_s,i_{s+1})/n \ge \grklet n$, as required.
\end{proof}

Using \refL{lexp}, the proof of Lemma 2 in \cite{LMcD} adapts immediately
to our setting. Note that we use different notation (in particular,
Greek letters) from \cite{LMcD}, for consistency with the rest
of the present paper.

\begin{proof}[Proof of \refT{TrobustOLD}]
As noted at the start of the section, 
it suffices to prove \eqref{esl}, assuming
that \eqref{sassump} holds.

Given $\gd,\eps>0$, by an \emph{$(\eps,\gd)$-cut} in a graph $G$
we shall mean a partition $(W,\Wc)$ of
the vertex set of $G$ with $|W|$, $|\Wc|\ge\eps|G|$, such that $G$
contains at most $\gd|G|$ edges from $W$ to $\Wc$.
We know from \refT{T2} that $\frac1n C_1(G_n)\pto \rho(\kk)>0$, so proving
\eqref{esl} is equivalent to showing that 
for any $\eps>0$ there is a
$\gd=\gd(\eps)>0$ such that \whp\ the giant component of $G_n$ has no
$(\eps,\gd)$-cut. 

Given $0<\gamma<1$, let $G_1$, $G_2$ be independent graphs with the
distributions 
of $\gnx{(1-\gamma)\kk}$ and $\gnx{\gamma\kk}$, respectively. We may and shall
couple the pair $(G_1,G_2)$ with $G_n\sim \gnk$ so that $G_1\cup G_2\subseteq G_n$.
(The union has almost the distribution of $G_n$; the only difference
arises from the possibility of $G_1$ and $G_2$ sharing edges.)

Fix $\eps>0$.
Recall that $\kk$ is supercritical, so $\rho(\kk)>0$.
Furthermore, $\kk$ is irreducible, so $\rho(\kk;i)>0$ for each $i$.
By \refT{TappB}, for each $i$ we have
$\rho((1-\gamma)\kk;i)\upto \rho(\kk;i)$ as $\gamma\to 0$.
Let us fix
a $\gamma$ such that 
\[
 \rho((1-\gamma)\kk;i)\ge (1-\eps/3)\rho(\kk;i)
\]
holds for every $i$. Thus, $\rho((1-\gamma)\kk)\ge (1-\eps/3)\rho(\kk)$.

Following (in this respect) the notation of \cite{LMcD}, let $U$
and $U_1$ denote the largest components of $G_n$ and $G_1$ respectively,
chosen according to any rule if there is a tie.
Then, by Theorems \refand{T2}{T2b}, the events
\[
 A_1 \= \{ |U_1| \ge (1-\eps/2)\rho(\kk) n \}
\]
and
\[
 A_2 \= \{ U_1\subseteq U \hbox{ and }|U_1|\ge (1-\eps/2)|U| \}
\]
hold \whp; for the condition $U_1\subseteq U$, note that $U_1$ must
be contained in some component of $G_n$, and \whp\ only $U$
is large enough.

Let $\eps_1=\min\{ \rho(\kk;i)\mu\set{i} : i\in \sss\}/2>0$.
By \refT{T1A}, the event 
\[
A_1'\=\set{\text{$U_1$ contains at least $\eps_1n$ 
vertices of each type $i$}}
\]
holds \whp. Without loss of generality, we may assume
that $\eps<\eps_1$. Let $\nu=\gamma\grklet(\kk,\eps\rho(\kk)/2)$,
where $\grklet$ is the function appearing in \refL{lexp}.
If $n$ is large enough then, by \refL{lexp},
whenever $A_1'$ holds,
if we partition the vertex set of $U_1$ into two parts $V_1$, $V_2$ each of
size at least $\eps\rho(\kk)n/2<\eps_1n$, then the expected number
of edges in $G_2$ from $V_1$ to $V_2$ is at least $\nu n$.

Continuing exactly as in \cite{LMcD}, but keeping our notation for the
relevant constants, 
let $\eta>0$ satisfy
\[
 1+2\eta \le \exp(\nu/8),
\]
and let $\delta>0$ be the minimum of $\nu/4$ and
$\frac{1}{2}\eta_0(\eta)$ (from 
\refL{LLMcD}).
Let
\[
 A_3 \= \{U \hbox{ has an $(\eps,\delta)$-cut in }G_n\},
\]
and
\[
 A_4 \= \{ U_1 \hbox{ has an $(\eps/2,2\delta)$-cut in }G_n\}.
\]
We claim that $A_2\cap A_3\subseteq A_4$. Indeed, suppose
that $A_2$ holds and that $U$ has an $(\eps,\delta)$-cut into
$B\cup C$. Let $B_1=B\cap U_1$ and $C_1=C\cap U_1$. Then $U_1$
has a partition into $B_1\cup C_1$, both $|B_1|$ and $|C_1|$
are at least
\[
 \eps |U| - (|U|-|U_1|) \ge \eps |U_1|/2,
\]
and the number of cross edges is at most $\delta|U|\le 2\delta|U_1|$,
so $A_4$ holds, proving the claim.
As $A_2$ holds \whp, and our aim is to show that $\P(A_3)\to 0$,
it thus suffices to show that $\P(A_4)\to 0$.

Let us condition on $G_1$, assuming that $A=A_1\cap A_1'$ holds.
By \refL{LLMcD} and our choice of $\delta$, there are at most
$(1+\eta)^n$ $(\eps/2,2\delta)$-cuts of $U_1$ in $G_1$.
Consider any one such cut, partitioning $U_1$ into $B\cup C$, say.
Let $X_2$ be the number of edges of $G_2$ from $B$ to $C$.
Recalling that $G_1$ and $G_2$ are independent,
as noted above, $\E(X_2)\ge \nu n \ge 4\delta n$.
As $X_2$ has a binomial distribution, a standard Chernoff estimate
implies that
\[
 \P(X_2\le 2\delta n)\le \P(X_2\le \E(X_2)/2) \le \exp(-\E(X_2)/8)
 \le \exp(-\nu n/8).
\]
As $G_2\subseteq G_n$, the probability that $B\cup C$
is an $(\eps/2,2\delta)$-cut of $U_1$ in $G_n$
is at most $\exp(-\nu n/8)$.
Hence, conditional on $G_1$ and assuming that $A$ holds,
\[
 \P(A_4\mid G_1) \le (1+\eta)^n \exp(-\nu n/8) \le (1+\eta)^n(1+2\eta)^{-n} 
= o(1).
\]
As the estimate above holds uniformly for all $G_1$ such that $A$ holds,
it follows that
$\P(A_4\mid A)=o(1)$. As $A$ holds \whp, this shows
that $\P(A_4)\to 0$, as required.
\end{proof}

As noted earlier, it is easy to deduce \refT{TrobustNEW} from 
\refT{TrobustOLD}. Recall that the only differences between these results are
that in \refT{TrobustNEW} we allow $\vxs$ to be a generalized vertex space,
and we allow the deletion of vertices as well as the addition and
deletion of edges. 
\begin{proof}[Proof of \refT{TrobustNEW}]
We first show that, as usual, we lose no generality by assuming that 
$\vxs$ is a vertex space.  Although this is not obvious at first
sight, the general arguments in \refSS{SSgvs} apply.  Indeed, the only
potential problem arises when we condition on the sequences $(\xs)$,
since $\gd$ might depend on $(\xs)$. However, fixing $\eps$ and
defining $X_n$ as the smallest number of changes (edge/vertex
deletions or edge additions) that can be made to $G_n$ to obtain a
graph $G_n'$ for which \eqref{estableNEW} fails, then
Theorem~\ref{TrobustNEW} states exactly that, for any $\eps>0$, we
have $n/X_n = O(1)$ \whp. As noted in \refSS{SSgvs}, in proving that
any function of $G_n$ is $O(1)$ \whp, we may assume that the sequences
$(\xs)$ are deterministic, by conditioning and applying \refL{Loc2}.

From now on we assume that $\vxs$ is a vertex space.
Turning to vertex deletion, given an $\eps>0$, let $\gd>0$ be such that
the conclusion of \refT{TrobustOLD} holds.
By \refP{Ponon}, there is a $\gd'>0$ such that
the event $\cE$ that any $\gd'n$ vertices of $G_n$ are incident with
at most $\gd n/2$ edges 
holds \whp.
Set $\gd''=\min\{\gd', \gd/2\}$.

Let $G_n'$ be any graph obtained from $G_n$ by deleting at most $\gd''
n\le \gd' n$ vertices, 
and then adding and deleting at most $\gd'' n\le \gd n/2$ edges. 
If $\cE$ holds,
then replacing the deleted vertices as isolated vertices to obtain a graph
$G_n''$ on $V(G_n)$, 
we have
\[
 |E(G_n'')\setdiff E(G_n)| \le \gd n/2+\gd n/2 = \gd n.
\]
Hence, by \refT{TrobustOLD}, \whp\ every such $G_n'$ satisfies
\eqref{Estab}, which is exactly \eqref{estableNEW}. This  
completes the proof of \refT{TrobustNEW}.
\end{proof}

The above proof of \refT{TrobustOLD} is much simpler than any proof we have
been able to find based directly on branching process methods.
However, the branching process approach does give additional insight
into the relationship between the giant component and two-core of $G_n$
and the branching process $\bpk$.

\subsection{Branching process analysis of the two-core}\label{SS2core}

Throughout this subsection we work with a kernel $\kk$ on a vertex space 
$\vxs$ satisfying the assumptions
\eqref{sassump}. As usual, we assume without loss
of generality that the number $n_i$ of vertices
of each type $i$ is deterministic, with $n_i/n\to \mui$ as $n\to\infty$.
The cornerstone of the branching process approach is
the following form of the coupling between the
neighbourhood exploration process in
$G_n$ and the branching process $\bpk$.

\begin{lemma}\label{l:fincouple}
There is a function $L_0=L_0(n)\to \infty$ such that
we may couple the neighbourhood exploration
process of a random vertex $v$ of $G_n=\gnkx$ with the branching
process $\bpk$ so that \whp\ they agree
for the first $L_0$ generations.
\end{lemma}
The sense of agreement is that there is a bijection between the
vertices of $G_n$ at distance at most
$L_0$ from $v$ and the first $L_0$ generations of $\bpk$ mapping $v$ to
the initial particle and preserving type
and adjacency, where
particles in the branching process are adjacent if one is a child of the other.
\begin{proof}
The argument is the same as the proof of \eqref{n2},
except for the error bounds. Note that it suffices to consider the case $L_0$
fixed. With $L_0$ fixed, the total number of vertices
encountered has bounded expectation, so we may abandon
the coupling if we reach more than $\log n$ vertices, say, in
the neighbourhood exploration.
At every step, the number of unused vertices of type $j$
is $\mu\{j\}n+o(n)$. Using this estimate
in place of \eqref{n1}, we may couple the number of new neighbours
of each type found with a corresponding $\Po(\kk(i,j)\mu\{j\})$
random variable so as to agree with probability $1-o(1)$.
 As the expected total
number of steps is $O(1)$, the total error probability is $o(1)$.
\end{proof}

If we have $\mu\{j\}n+O(1)$ vertices of type $j$, \refL{l:fincouple}
holds for any $L_0=o(\log n)$.

As in the proof of \refL{L:step1}, the coupling easily extends to
the $L_0$-neighbourhoods of two vertices. Given $G_n$, let $v$ and $w$
be chosen independently 
and uniformly at random from the vertices of $G_n$.

\begin{lemma}\label{coup2}
There is an $L_0(n)\to\infty$ such that
we may couple $(G_n,v,w)$ with two independent copies $\bpk$, $\bpk'$
of the branching 
process $\bpk$ so that \whp\ the first $L_0$ neighbourhoods of $v$ and
of $w$ agree with the 
first $L_0$ generations of $\bpk$ and of $\bpk'$, respectively.
\end{lemma}

We omit the proof, noting only that for $L_0$ fixed, the probability
that $v$ and $w$ 
are within graph distance $2L_0$ is $o(1)$.

The next step is to find a way of applying the coupling results above
to expectations of functions of the
neighbourhoods. This will require some care, due to the possible large
contribution to an expectation from the low probability event that the
coupling fails.

We consider functions $f(v,G)$ defined on a pair $(v,G)$, where $G$ is
a graph in which 
each vertex has a type from $\sss=\{1,2,\ldots,r\}$, and $v$
is a distinguished vertex of $G$, the {\em root}.
We call such a function an
{\em \Lnf} if it is invariant under type preserving
rooted-graph isomorphisms and depends only
on the subgraph of $G$ induced by vertices within a fixed distance $L$
of $v$.  
We define $f(\bpk)$ by evaluating $f$ on the branching process
in the natural way: form a graph from the branching process as above,
and take the initial 
particle as the root. Thus \refL{l:fincouple} implies that we can
couple $(G_n,v)$ with $\bpk$ 
so that $f(v,G_n)=f(\bpk)$ \whp\
for every \Lnf\ $f$.

Given an \Lnf\ $f$, let
\[
 S_n \= \frac1n \sum_{v\in V(G_n)} f(v,G_n).
\]
Also, for $v$ and $w$ independent random
vertices of $G_n$, let $X_n=f(v,G_n)$ and $Y_n=f(w,G_n)$.
Note that $\E(S_n)=\E(X_n)=\E(Y_n)$.

\begin{theorem}\label{t:cexp}
Let $\vxs$ be a vertex space with finite type space $\sss$, and let
$\kk$ be a kernel on $\vxs$. 
If $f$ is an \Lnf\ such that $\sup_n \E(X_n^4)<\infty$, then $S_n\pto
\E(f(\bpk))$. 
\end{theorem}

\begin{proof}
Let $X=f(\bpk)$.
By \refL{l:fincouple} we may couple $X_n$ and $\bpk$ so that
$\P(X_n\ne X)\to 0$,
and hence $X_n\pto X$.
Since $\sup_n \E(X_n^4)<\infty$ implies that the variables $X_n$ are
uniformly integrable, it follows that
\begin{equation}\label{eSnx}
 \E(S_n)=\E(X_n)\to \E(X);
\end{equation}
see \cite[Lemma 4.11]{Kall}, for example. 

Let $Y$ be an independent copy of $X$. From \refL{coup2} we may couple
$(X_n,Y_n)$ with $(X,Y)$ so that 
$\P((X_n,Y_n)\ne (X,Y))\to 0$. In particular, $X_nY_n\pto XY$.
As
\[
 \E((X_nY_n)^2)=\E(X_n^2Y_n^2) \le \sqrt{\E(X_n^4)\E(Y_n^4)} =\E(X_n^4) \le C,
\]
for some $C<\infty$, the variables $X_nY_n$ are also uniformly integrable,
so
$\E(X_nY_n)\to \E(XY)$. But $\E(S_n^2)=\E(X_nY_n)$ by linearity of expectation,
while $X$ and $Y$ are independent and have the same distribution. Thus
$\E(S_n^2)\to \E(XY)=\E(X)^2$. 
Together with \eqref{eSnx}, this proves the result.
\end{proof}

\begin{remark}\label{R:cexp}
\refT{t:cexp} can be applied to any \Lnf\ $f$ bounded by
a polynomial of
the number of vertices within distance $L$ of $v$.
Indeed, the number of vertices at
distance $t$ from $v$ in $G_n$ is stochastically dominated
by the number $N_t$ of particles in generation $t$ of a Galton-Watson branching process
in which the number of children of each particle has a $\Bi(n,\max\kk/n)$ distribution.
As a $\Bi(1,p)$ distribution is stochastically dominated by a $\Po(1,-\log(1-p))$
distribution, if $n$ is large enough then $N_t$ is dominated by $N_t'$, the number
of particles in generation $t$ of the single-type Poisson branching process
$\bpx{2\max\kk}$. The probability generating function of $N_t'$
is obtained by iterating that of the Poisson distribution $t$ times. As all moments
of a Poisson distribution are finite, it follows that all moments of $N_t'$ are finite,
so the fourth moment of any power of $N_t'$ is finite.
\end{remark}

\refP{Ponon} states that, given $\eps>0$, there is a $\gd>0$ such that 
\whp{} any set of at most $\gd n$ vertices of $G_n$ is incident with
at most $\eps n$ edges. A very special case of \refT{t:cexp} gives an
alternative proof of this result (under the more restrictive
assumptions of the present section). Indeed, writing $X$ for the
number of particles in the first generation of \bpk, since $\E X$ is
finite, we have $\E(X\IN{X>M})\to 0$ as $M\to\infty$.  Given $\eps>0$
there is thus an $M$ for which $\E(X\IN{X>M})<\eps/3$.  Writing
$d_G(v)$ for the degree of a vertex $v$ in the graph $G$, let
$f(v,G)=d_G(v)\IN{d_G(v)>M}$; clearly, $f$ is a \Onf.  By
\refR{R:cexp}, \refT{t:cexp} applies to $f$, so
\[
 S_n = \frac1n \sum_{v\in V(G_n)\,:\, d(v)>M} d(v) \pto \E(X\IN{X>M})<\eps/3.
\]
Hence $S_n\le \eps/2$ \whp.
Set $\gd=\eps/(2M)$, and let $W$ be a set of at most $\gd n$ vertices of $G_n$.
Then
\begin{eqnarray*}
 \sum_{w\in W} d(w) &\le& \sum_{w\in W\,:\, d(w)>M} d(w) + M|W| \\
  &\le& nS_n + \eps n/2 \le \eps n
\end{eqnarray*}
whenever $S_n\le \eps/2$ holds, so \whp{} any set of at most $\gb n$ vertices
of $G_n$ are incident with at most $\eps n$ edges.

\smallskip
Our next result is a simple observation concerning short cycles.
As before, we assume throughout that \eqref{sassump} holds.
\begin{lemma}\label{lshort}
Let $L=L(n)=o(\log n)$. The probability
that a random vertex $v$ of
$G_n$ is  within distance $L$ of a cycle of length at most $L$ is $o(1)$.
\end{lemma}
\begin{proof}
As all edge probabilities are bounded by $p=\max\kk/n$,
the expected number of vertices $v$ at distance $d\ge 0$ from a cycle
of length $l\ge 3$ is at most
$n^{l+d}p^{l+d}\le (\max\kk)^{l+d}$. Summing over $l,d\le L$, the
expectation is $o(n)$.
\end{proof}

The \emph{two-core} $C^2(G)$ of a graph $G$ is the maximal subgraph of
$G$ with minimum degree at least 2.
Equivalently, $C^2(G)$ consists of those vertices and edges of $G$
that lie in some cycle in $G$,
or on a path joining two vertex-disjoint cycles.
We shall work with the two-core $\Gtc\=C^2(G_n)$ of $G_n$.
To do so, we need to relate certain properties of $\Gtc$
to the branching process $\bpk$. In the light of
\refL{l:fincouple}, it will be useful to have
a reasonably accurate ($o(1)$ error probability) `local'
characterization of when a vertex $v$ is in the two-core.
We shall need similar results for vertices not in the two-core, but
connected to it by short paths.
Note that \refT{T2b} gives us a corresponding characterization for the
giant component: for a suitable $L(n)\to \infty$,
up to an error probability of $o(1)$, a vertex $v$ is in the giant
component if and only if it
is in a component of size at least $L$, 
and using \refL{l:fincouple}, it is easy to check
that \whp\ when this condition holds the
{\em $L$-distance set}, the set of vertices
at graph distance exactly $L$ from $v$,
is non-empty, so there is a path of length $L$
starting at $v$. (We omit the details
as we use this statement only to motivate what follows, not
in the proof.) For the two-core, we need two vertex-disjoint paths.

Let $L=L(n)$ be a function tending to infinity slowly, to be chosen below.
For a vertex $v$ of $G_n$ and an integer $d\ge 0$, let $TC_d(v)$ be
the event that $v$ is at graph
distance at most $d$ from the two-core $\Gtc$ of $G_n$. Thus $TC_0(v)$ is the
event $v\in \Gtc$.
Let $LTC_d(v)$ be the `localized' event that there is a vertex $w$ at
distance $d'\le d$ from $v$ joined by
two vertex-disjoint paths of length $L$ to vertices at distance $d'+L$ from
$v$. Thus, as we explore the neighbourhoods of $v$ successively, 
$LTC_d(v)$ is the event that after $d'\le d$ steps we reach
a vertex $w$ (which we expect to be the closest vertex of the two-core to $v$)
with two neighbours in the next generation each of which has neighbours
for at least $d-1$ further generations.

\begin{lemma}\label{ltc}
Let $d\ge 0$ be fixed, and let $v$ be a random vertex of $G_n$.
Provided $L(n)$ tends to infinity sufficiently slowly, the event
$TC_d(v)\setdiff LTC_d(v)$,
i.e., the event
that one of $TC_d(v)$ and $LTC_d(v)$ holds but not the other, has
probability $o(1)$.
\end{lemma}
\begin{proof}
Assume, as we may, that $L=o(\log n)$.
We start with the case $d=0$.
Let us say that a cycle is \emph{short} if it has length at most $2L$.
By \refL{lshort}, the probability that $v$ is within distance $L$ of a
short cycle is $o(1)$.
If $TC_0(v)$ holds, i.e., $v$ is in the two-core, then $v$ is in a
cycle, or on a path joining two
vertex-disjoint cycles.
Assuming that $v$ is not close to a short cycle, in either case we can
find two vertex-disjoint paths of length $L$ starting from $v$, so
$LTC_0(v)$ holds. Hence $\P(TC_0(v)\setminus LTC_0(v))=o(1)$.

The reverse bound is more difficult, as what we need is an equivalent
for the two-core of \refT{T2b}, which states
that almost all vertices in largish
components are in a single giant component. In fact, we can use
\refT{T2b}.  Suppose that $LTC_0(v)\setminus TC_0(v)$ holds. Note that
$v$ is not in a cycle by definition of the two-core.  Let $w_1$, $w_2$
be two neighbours of $v$ joined by vertex-disjoint paths to vertices $x_1$,
$x_2$ at distance $L$ from $v$.  In $G_n-v$, there is no path from
$w_1$ to $w_2$;
otherwise, there would be a cycle in
$G_n$ containing $v$.  Hence, at least one of $w_1$ and $w_2$, let us
say $w_1$, is not in the giant component of $G_n-v$. (Here, by the
giant component we mean the largest component, chosen according to any
fixed rule if there is a tie.)  But $w_1$ is in a component of size at
least $L$, as witnessed by the path $w_1x_1$. In summary, if
$LTC_0(v)\setminus TC_0(v)$ holds, so does the event $E(v)$ that $v$
is adjacent in $G_n$ to a vertex $w$ in an intermediate component of
$G_n-v$, i.e., a component other than the largest having size at least
$L$.  As the random vertex $v$ is chosen independently of $G_n$, the
graph $G_n-v$ is an $(n-1)$-vertex graph to which \refT{T2b}
applies. Hence, taking $\omega(n)=L(n)$, by \refT{T2b} the number of
vertices $w$ of $G_n-v$ in intermediate components is
$o_p(n)$. Conditioning on $G_n-v$ tells us nothing about the
edges from $v$ to $G_n-v$. As $\kk$ is bounded, it follows that $E(v)$ has
probability $o(1)$.  Thus $\P(LTC_0(v)\setminus TC_0(v))=o(1)$,
completing the proof in the case $d=0$.

The general case follows using \refP{Ponon}. If $TC_d(v)\setminus
LTC_d(v)$ holds, then $v$ is within distance $d+L$ of a vertex on a
short cycle. Hence, by \refL{lshort}, $\P(TC_d(v)\setminus
LTC_d(v))=o(1)$.  If $LTC_d(v)\setminus TC_d(v)$ holds, then $v$ is within
distance $d$ of a vertex $v'$ for which
$LTC_0(v')\setminus TC_0(v')$ holds.
By the case $d=0$ above, $o_p(n)$ vertices $v'$ have
this property, and the result follows by applying \refP{Ponon} $d$
times.
\end{proof}

We now turn to the branching process equivalents of the events $TC_d$
and $LTC_d$. Considering the branching process $\bpk$ (started with a
single particle of random type), let $DS_d$ be the event that there is
a particle in some generation $d'\le d$ which has at least two children with
descendants in all future generations. Similarly, let $LDS_{d,L}$ be
the event that there is a particle $x$ in generation up to $d$,
say in generation $d'$, such that $x$ has two children each of which
has one or more descendants in generation $d'+L$, i.e., $L$ generations
after $x$.
Note that $LDS_{d,L}$ depends only on the first $d+L$
generations of the branching process.  Suppose that $L(n)$ grows
slowly enough that \refL{l:fincouple} applies with $2L$ in place of
$L$.  Then for any fixed $d$ we have $d+L\le 2L$ for large enough $n$,
and, with $v$ a random vertex of $G_n$ as before, from
\refL{l:fincouple} we have
\begin{equation}\label{ltclds}
 \P(LTC_d(v))=\P(LDS_{d,L(n)})+o(1) \quad\text{as \ntoo}.
\end{equation}
Note that for $d$ and $L$ fixed, the event $LDS_{d,L}$, which is defined in
terms of the branching process, does not depend on $n$, so
$\P(LDS_{d,L})$ is a constant.
For each $d$, as $L$ increases the events $LDS_{d,L}$ decrease
to the event $DS_d$.  Hence,
\begin{equation}\label{llim}
 \lim_{L\to\infty} \P(LDS_{d,L})=\P(DS_d).
\end{equation}
Suppose now that $L(n)$ tends to infinity sufficiently slowly
that \eqref{ltclds} and \refL{ltc} hold. Then, from \eqref{llim},
$\P(LDS_{d,L(n)})\to \P(DS_d)$ as $\ntoo$. Hence, from \eqref{ltclds},
$\P(LTC_d(v))=\P(DS_d)+o(1)$. Finally, using \refL{ltc} we obtain
\[
 \P(TC_d(v)) = \P(DS_d) +o(1).
\]

Considering two random vertices $v$, $w$ of $G_n$ and using
\refL{coup2} instead of \refL{l:fincouple},
we obtain $\P(TC_d(v)\cap TC_d(w))=\P(DS_d)^2+o(1)$ similarly.

Writing $TC_d$ for the set of vertices of $G_n$ for which $TC_d(v)$
holds, i.e., for the set of vertices
within distance $d$ of the two-core, it follows that
$\E |TC_d|/n \to \P(DS_d)$ and $\E (|TC_d|/n)^2 \to (\P(DS_d))^2$,
and thus
\begin{equation}\label{tcd}
 |TC_d|/n \pto \P(DS_d).
\end{equation}
As $\norm{\tk}>1$, the branching process \bpk\ is supercritical, so
$\P(DS_0)>0$. Hence, taking $d=0$ in
\eqref{tcd},
\begin{equation}\label{tclarge}
 |V(\Gtc)|/n \pto \P(DS_0) >0.
\end{equation}

The reason for considering the two-core $\Gtc$ of $G_n$ is that 
\refT{TrobustNEW} boils down to a statement about $\Gtc$. Roughly speaking,
the largest component of $G_n$ consists of the two-core with some
trees hanging off it, and it is easy to see what effect deleting
edges from the trees has on the size of the largest component. The question
is what happens when edges are deleted from the two-core.

\begin{lemma}\label{l2core}
Suppose that \eqref{sassump} holds,
i.e., $\sss=\{1,2,\ldots,r\}$, $\mu\{i\}>0$ for each $i$,
$\kk$ is irreducible, and $\norm{\tk}>1$.
Let $\Gtc$ be the two-core of $\gnkx$.
For any $\eps>0$ there is a $\gd>0$ such that the following statement
holds \whp:
for any set $W\subset V(\Gtc)$ with $|W|\ge \eps n$ and
$|V(\Gtc)\setminus W|\ge \eps n$ there
are more than $\delta n$ edges of $\Gtc$ joining $W$ to $V(\Gtc)\setminus W$.
\end{lemma}

In other words, under the assumptions of \refL{L:step1}, if $\norm{\tk}>1$ then
the two-core $\Gtc$ cannot be cut into two large (size $\Theta(n)$)
pieces by a small set of edges. Note that the two-core itself is large
by \eqref{tclarge}.
As the proof of \refL{l2core} is rather long, we first show
that it implies \refT{TrobustNEW}.

\begin{proof}[Deduction of \refT{TrobustNEW} from \refL{l2core}]
We have already shown (at the end of \refSS{SSLMcD}) that, using \refP{Ponon}, 
\refT{TrobustNEW} can be deduced from \refT{TrobustOLD}.
As noted at the start of the section, in proving \refT{TrobustOLD} we may assume
that \eqref{sassump} holds, and it suffices to prove \eqref{esl}.
{}From now on, let us fix the quantity
$\eps>0$ appearing in \eqref{esl}.

The events
$DS_d$ form an increasing sequence, and their union is contained in
the event $S$ 
that the branching process $\bpk$ survives (contains points in all
generations).  Also, $S\setminus \bigcup_d DS_d$ is the event that the
process survives, but with only a single infinite line of descent.
From basic properties of Poisson processes, starting from a particle
of type $x$, 
the types its surviving children,  
i.e., its children that have descendants in all later generations,
form a Poisson process on $\sss$ with intensity $\kk(x,y)\rho(\kk;y)\dd\mu(y)$.
In particular, the number of such children is Poisson with some mean
$\gl(x)>0$. 
It follows that, conditional on a particle surviving, the probability
that it has 
at least two surviving children is positive, and hence, as the type
space $\sss$ is finite, 
bounded way from zero. Hence $\P(S\setminus \bigcup_d DS_d)=0$, so
$\P(DS_d)\upto\P(S)=\rho(\kk)$,
and there
is a constant $D$ such that $\P(DS_D)\ge \rho(\kk)-\eps/3$.
{}From \eqref{tcd}, for any fixed $D$
the set $TC_D$ of vertices within
distance $D$ of the two-core has size $\P(DS_D)n+o_p(n)$, so
\whp\
\begin{equation}\label{tcd1}
 |TC_D| \ge (\rho(\kk)-\eps/2)n.
\end{equation}

Let $\eps'<\P(DS_0)/3$ be a small positive constant to be chosen later,
and let $\gd=\min\{\eps',\gd(\eps')\}$, where $\gd(\cdot)$ is the function
appearing in \refL{l2core}.
Let us delete an arbitrary set $E$ of at most $\gd n$ edges from $G_n$,
leaving a graph $G_n'$.  Let $\Gtcm\subseteq \Gtc$ be the largest
remaining connected part of $\Gtc$. We claim that
\[
 |V(\Gtc)\setminus V(\Gtcm)|\le \eps'n
\]
holds \whp.
Note that $|V(\Gtc)|\ge 3\eps'n$ \whp\ by \eqref{tclarge}.
If this inequality holds and every component of $\Gtc\setminus E$ has size
at most $|V(\Gtc)|-\eps'n$, then there is a union $H$ of components of
$\Gtc\setminus E$ with between $\eps'n$ and $|V(\Gtc)|-\eps'n$ vertices:
indeed, if the largest component has at least $\eps'n$ vertices, this
will do as $H$. Otherwise, every component has at most $\eps'n$ vertices,
and the smallest union $H$ with at least $\eps'n$ vertices will do.
The existence of an $H$ with the stated properties
has probability $o(1)$ by \refL{l2core}, proving the claim.

Let $X$ be the component of $G_n'$ containing $\Gtcm$.
If $v\in TC_D\setminus X$ then, as $v\in TC_D$, there is a path in $G_n$ of length
at most $D$ from $v$ to a vertex of $\Gtc$. Taking any such path,
as $v\notin X$, either the path ends in a vertex of $\Gtc\setminus \Gtcm$,
of which there are \whp\ at most $\eps'n$, or
it contains an edge of $E$, and hence contains an endvertex of such an edge;
there are at most $2\gd n\le 2\eps'n$ such endvertices.
In particular, \whp\ all $v\in TC_D\setminus X$ are within distance $D$ (in $G_n$) of some set
of at most $3\eps'n$ vertices.
Applying \refP{Ponon} $D$ times, it follows that if we choose $\eps'$ small enough,
then \whp{}
$ |TC_D\setminus X| \le \eps n/2.$
Using \eqref{tcd1} it follows that \whp{}
\[
 C_1(G_n')\ge |X|\ge (\rho(\kk)-\eps)n,
\]
completing the proof of \refT{TrobustNEW}.
\end{proof}

It remains only to prove \refL{l2core}.
The uniform case ($\kk$ constant) has a simple proof, presented below.

\begin{proof}[Proof of \refL{l2core}, uniform case]
In a moment we shall restrict to the uniform case; for now, we assume
\eqref{sassump}.

As $\norm{\tk}>1$, the branching process $\bpk$ is supercritical, so
$\rho(\kk)>0$.  From irreducibility, it follows that $\rho(\kk;i)>0$
for every $i$.  Consider the event $TS$
that the initial particle has
exactly three children that survive.  As the initial particle has
positive probability of having exactly three children, $\P(TS)>0$.
Arguing as for \eqref{tclarge}, one can show that
the number of vertices of degree exactly 3 in $\Gtc$ is $\P(TS)n+o_p(n)$;
in fact, both statements are special cases of \refL{l:gtcexp} below.

We shall
condition on the vertex set and (labelled) degree
sequence of $\Gtc$. In other words, we shall condition on the sequence 
${\bfd}=(d(1),\ldots,d(n))$, where $d(i)$
is the degree in $\Gtc$ of the
vertex $i$ and $d(i)=0$ if $i\notin V(\Gtc)$. Let us write $n_2$ for
$|V(\Gtc)|=|\{i:d(i)>0\}|$, $m_2$ for $e(\Gtc)=\frac12\sum d(i)$, and
$n_{\ge3}$ for the number of $i$ for which $d(i)\ge 3$.
Note that \whp\ $n_{\ge3}\ge \P(TS)n/2>0$. Also,
$e(G_n)=O(n)$ \whp{} (for example, by
\refP{PE}), and $G_n$ has maximum degree
$O(\log n)$ \whp. Thus there are positive constants $\eps_1$ and $C$,
depending only on $\kk$, such that
\begin{equation}\label{typ}
 n_{\ge3}\ge \eps_1 n, \quad m_2\le Cn_2,
\quad \text{and }\max_i d(i)\le C\log n,
\end{equation}
hold \whp.

{}From now on we consider the uniform case, where $\kk$ is constant, or,
equivalently, $|\sss|=1$. This is just the usual Erd\H os--R\'enyi
random graph $G_n=\gnx{c/n}$, with $c>1$.
We condition on ${\bfd}$, assuming, as we may, that the conditions
\eqref{typ} hold. 
In the uniform case it is easy to
see that (given ${\bfd}$) the graph $\Gtc$ is uniformly distributed
among all graphs with degree sequence ${\bfd}$. This is because any
graph can be decomposed into its two-core and a collection of
vertex-disjoint trees, each sharing at most one vertex with the two-core.
Hence, any two graphs $H_1$, $H_2$ with degree sequence ${\bfd}$ can
be extended in exactly the same ways to graphs $H_1'$, $H_2'$
with vertex set $[n]$
so that
$H_i'$ has two-core $H_i$.  As corresponding graphs $H_1'$, $H_2'$
have the same number of edges, they are equally likely in the model
$\gnx{c/n}$. Summing over the possible extensions, $H_1$ and $H_2$ are
equally likely to arise as $\Gtc$.

Let $H$ be the random multigraph with degree sequence ${\bfd}$
generated by the configuration model of~\cite{B1}.  In other words,
for each vertex $i$ we take $d(i)$ `stubs', and we pair the $2m_2$
stubs randomly, with all $(2m_2-1)!!$ pairings equally likely. For
every pair in this pairing, we take an edge between the corresponding
vertices. This generates a multigraph $H$ with degree sequence ${\bf
d}$, where $H$ may contain loops and multiple edges. Let $\psimple$ be
the probability that $H$ is simple.  From our assumptions on ${\bfd}$
it is easy to check that $\psimple=\exp(-o(n))$.
This very crude lower bound is all that we shall need; the much
stronger bound $\exp(-O((\log n)^2))$ follows from
\eqref{typ} and the general results in McKay~\cite{McKay}.
In fact, using the fact that $\sum_i d(i)^2=O(n)$ \whp, one can show
that $\psimple=\Theta(1)$.
Given that $H$ is simple, it is uniformly
distributed among all simple graphs with degree sequence ${\bfd}$,
i.e., $H$ has the distribution of $\Gtc$. Hence, to show that $\Gtc$
has a certain property \whp\ it suffices to show that $H$ has the
property with probability $1-o(\psimple)$.

For $W\subseteq V(\Gtc)$ let $\Wc=V(\Gtc)\setminus W$. Let $S(W)$ denote
the set of stubs associated to $W$, so $|S(W)|=\sum_{i\in W}d(i)$. For
a set $S$ of stubs, let $\overline S=S(V(\Gtc))\setminus S$, and let
$\pcut(S,\overline S)$ denote the probability that, in the random
pairing, every stub in $S$ is paired with another stub in $S$. If
$V(\Gtc)$ has a partition $W$, $\Wc$ with at most $\gd n$ edges
between $W$ and $\Wc$, then there is an $S\subseteq S(W)$ with
$|S(W)\setminus S|\le\gd n$ so that every stub in $S$ is paired with
another stub in $S$. Hence, the expected number of such partitions with
$W$, $\Wc$ large is at most
\begin{equation}\label{ecdef}
 \EC \= \sum_{W\subset V(\Gtc): |W|,|\Wc|\ge \eps n}
 \ \sum_{S\subseteq S(W): |S(W)\setminus S|\le\gd n} \pcut(S,\overline S),
\end{equation}
and it suffices to show that $\EC$ is $o(\psimple)$ if we choose $\gd$
small enough.  Now $\pcut(S,\overline S)=0$ if $|S|$ is odd, and
otherwise
\begin{equation}\label{pcutu}
 \pcut(S,\overline S)
= \frac{(|S|-1)!! \,(|\overline S|-1)!!}{(2m_2-1)!!} \le \binom{m_2}{|S|/2}^{-1}.
\end{equation}
Every vertex of $\Gtc$ has degree at least 2, so for any $W$ we have
$|S(W)|/2\ge |W|$.  As, from \eqref{typ}, there are at least $\eps_1
n$ vertices of degree $3$ in $\Gtc$, either $|S(W)|/2\ge |W|+\eps_1
n/4$ or $|S(\Wc)|/2 \ge |\Wc|+\eps_1 n/4$.  It follows that whenever
$|W|$, $|\Wc|\ge \eps n$ we have
\[
 \binom{m_2}{|S(W)|/2} \ge \exp(4an)\binom{n_2}{|W|},
\]
where $a>0$ is a constant depending only on $\eps_1$ and $\eps$, and
hence only on $\kk$ and $\eps$.  Choosing $\gd$ small enough, it
follows that
\begin{equation}\label{pcut2}
 \binom{m_2}{|S|/2} \ge \exp(3an)\binom{n_2}{|W|},
\end{equation}
whenever $S\subseteq S(W)$ with $|S|\ge |S(W)|-\gd n$.  Given $w=|W|$,
there are at most $\binom{n_2}{w}$ choices for $W$. Given $W$, there
are, crudely, at most $m_2\binom{2m_2}{\gd n}$ choices for $S\subseteq
S(W)$ with $|S(W)\setminus S|\le \gd n$.
Hence, from \eqref{ecdef}, \eqref{pcutu} and \eqref{pcut2},
\[
 \EC \le
\sum_{\eps n\le w\le n_2-\eps n} \binom{n_2}{w} m_2\binom{2m_2}{\gd n}
\exp(-3an)\binom{n_2}{w}^{-1}.
\]
Choosing $\gd$ small enough, it follows that $\EC=O(n^2\exp(-2an))$,
and hence that $\EC\le \exp(-an)$ for $n$ large enough. As
$\psimple=\exp(-o(n))$, this completes the proof in the uniform
($|\sss|=1$) case.
\end{proof}

One might hope that any argument for the uniform case would adapt 
easily to the finite-type case.  However, we have been unable to find
a simple extension of the argument above. 
Our branching-process based argument for the general
case is somewhat involved and rather lengthy, and we shall not present it.
This is because adapting a proof due to Luczak and McDiarmid~\cite{LMcD}
gives the much simpler proof of \refT{TrobustOLD} in \refSS{SSLMcD}, 
which in turn immediately implies
the general case of \refL{l2core}.
We believe, however, that the results in this subsection are likely
to be useful for determining other properties of the two-core.

\medskip
We close this section with a final result, \refL{l:gtcexp} below,
relating any `local' property of the two-core to the branching process
$\bpk$. This will require a little introduction.

Recall that, by \refL{l:fincouple}, if $L\to\infty$ sufficiently
slowly then we may couple the $2L$-neighbourhood of a random vertex $v$ of
$G_n$ with the first $2L$ generations of the branching process $\bpk$
so that they agree with probability $1-o(1)$. By \refL{ltc}, for almost every
vertex $v$, $v$ is in the two-core if and only if there are two
disjoint paths of length $L$ starting at $v$. This allows us to adapt
the coupling, and hence \refT{t:cexp}, to the two-core.

For $v\in \Gtc$ and $t\ge 0$, let
$\Ga_t(v,\Gtc)$ be the set of vertices of $\Gtc$ at graph distance
$t$ from $v$, and set $\Ga_t(v,\Gtc)=\emptyset$ if $v\notin \Gtc$. Note
that any vertex on a path joining two vertices of $\Gtc$ is in $\Gtc$,
so $\Ga_t(v,\Gtc)$ is the $t$-distance set
of $v$ in the graph $\Gtc$.
Let $\bpkmod$ be
obtained from $\bpk$ in two steps: first, delete any particle that does
not have descendants in all future generations.  Then, if the initial
particle has only one remaining child, delete everything; we write
$\bpkmod=\emptyset$ in this case. We obtain a certain branching
process $\bpkmod$ having the following properties whenever
$\bpkmod\ne\emptyset$: the first particle has at least two children,
and every later particle at least one child.

For constant $D$, if $L=L(n)\to\infty$ then up to an $o(1)$ error probability, the first $D+L$
generations of $\bpk$ determine the first $D$ generations of $\bpkmod$:
consider surviving to generation $D+L$ instead of
surviving forever. Let $v$ be a random vertex of $G_n$.
If $L\to\infty$ sufficiently slowly then,
by \refL{ltc} and \refP{Ponon} (applied $D$ times),
the probability that $v$ is within distance
$D$ of a vertex $w$ for which one of $TC_0(w)$ and $LTC_0(w)$ holds but not the other is $o(1)$.
Using \refL{l:fincouple}, it follows that the $\Gtc$-neighbourhoods
$(\Ga_t(v,\Gtc))_{t=0}^D$ of a random vertex $v\in G_n$ can be coupled
with the first $D$ generations of $\bpkmod$ so as to agree with probability $1-o(1)$.
Similarly, \refL{coup2} implies its equivalent for $\Gtc$ and $\bpkmod$.
Using these two results, an analogue of \refT{t:cexp} follows.
In the result below
we take $f(\bpkmod)$ to be zero when $\bpkmod$ is empty. The proof follows
exactly that of \refT{t:cexp}, so we omit it. 

\begin{lemma}\label{l:gtcexp}
Let $D$ be fixed, and let 
$f=f(v,G)$ be a \Dnf\ bounded by
a polynomial of
the number of vertices
within distance $D$ of $v$. Then
\[
 S_n \= \frac1n \sum_{v\in\Gtc} f(v,\Gtc) \pto \E(f(\bpkmod)).
\]
\end{lemma}

Of course, the condition on $f$ could be replaced by a fourth-moment condition as in \refT{t:cexp}.
As an immediate consequence
of \refL{l:gtcexp}, we can describe, for example, the typed degree
sequence of $\Gtc$, taking $f_{i,d}(v,G)$ to be the
\Onf\ taking the value $1$ when $v$ has degree $d$ and type $i$, and $0$ otherwise.
\refL{l:gtcexp} can be used as the basis of a proof of \refL{l2core}, but as noted above, the
details are rather involved; see the first version of this paper, at
\arxiv{math.PR/0504589v1}

\begin{remark}
Recently, Riordan~\cite{Rkcore} proved
an analogue of \refL{l:gtcexp} for the $k$-core,
again using a direct coupling of the neighbourhood exploration process
with a branching process.
\end{remark}

\section{Bounds on the small components}\label{S2ndpf}

In this section we prove \refT{T4}, i.e., that the sizes of the small
components of $G_n=\gnkx$ are $O(\log n)$ \whp\ under certain assumptions.
As before, by the \emph{giant component} in a graph $G_n$ we mean the
unique largest component, provided it has $\Theta (n)$ vertices --
all other components are {\em small}.
Thus, in the supercritical case ($\norm{\tk}>1$), where there is a giant
component, a small component is any component other than
the largest, so our aim is to prove an upper bound on $C_2(G_n)$.
In the strictly subcritical case ($\norm{\tk}<1$), all components
are small, and our aim is to show that $C_1(G_n)=O(\log n)$ \whp.

We shall prove three results that together imply
\refT{T4}, namely Theorems \ref{T4b}, \ref{T4a} and \ref{T4c} below.
In the supercritical case, we shall also prove
a more general result,
\refT{Tdual} below,
describing the distribution of the
graph formed from $G_n$ by deleting the giant (more precisely,
largest) component, $\cC_1(G_n)$. This description is in terms
of another instance of our general model, involving the dual kernel
$\kkd$ defined in \refD{Ddual}.

\begin{theorem}
 \label{Tdual}
Let $(\kkn)$ be a graphical sequence of kernels on a (generalized) vertex space
$\vxs$ with \qir\ limit $\kk$,
with $\norm{\tk}>1$.
Let $G_n=\gnxx{\kkn}$, and let $G_n'$ be the graph obtained from $G_n$
by deleting
all vertices in the largest component
$\cC_1(G_n)$.
There is a generalized vertex space $\hvxs=(\sss,\hmu,(\ys))$ with
$\hmu$ given by $\dd\hmu(x)=(1-\rho(\kk;x))\dd\mu(x)$, 
such that $G_n'$ and $\gnxxh{\kkn}$ can be coupled to agree \whp.
Furthermore, the sequence $(\kkn)$ is graphical on $\hvxs$ with
\qir{} limit $\kk$.
\end{theorem}

If we wish, we can renormalize so that $(\sss,\hmu)$ becomes 
a ground space; see the comment after \refD{Ddual}. However,
the resulting graph still has a random number of vertices, so we
cannot insist that $\hvxs$ is a vertex space.

\medskip
Theorem~\ref{Tdual} is the natural
generalization to our context of the old `duality result' of
Bollob\'as~\cite{BBtr} for the Erd\H os--R\'enyi model $\gnx{c/n}$
that was the basis of the
study of the phase transition there (see also
{\L}uczak~\cite{Lucz1}, Janson, Knuth, {\L}uczak and
Pittel~\cite{JKLP}, and the books \cite{BB,JLR}).

\begin{remark}\label{RC2crit}
We know that the random graph $\gnxxh{\kkn}$ in \refT{Tdual} cannot be
supercritical, since otherwise $G_n$ would have a second giant
component.
It may be critical, see \refE{E2large1}, but is typically subcritical;
one sufficient condition for subcriticality
is given in the next result.
\end{remark}

\begin{theorem}\label{Tdualsub}
Under the assumptions of \refT{Tdual},
if, in addition,
$\iint_{\sss^2}\kk(x,y)^2\dd\mu(x)\dd\mu(y)<\infty$, then
$\gnxxh{\kkn}$ is subcritical.
\end{theorem}

\begin{proof}
As usual, we may normalize so that $\mu(\sss)=1$. The result 
then follows immediately from \refT{Tdualbp}.
\end{proof}

\begin{example} \label{E2large1}
As in \refE{E2large+}, let
$\sss=\set{1,2,3,\dots}$ with $\mu\set{k}=2^{-k}$, and let
$x_1,\dots,x_n$ be \iid{} random points in $\sss$ with
distribution $\mu$. Let $(\eps_k)_1^{\infty}$ be a
sequence of positive numbers tending to zero, to be
chosen below. Set $\kk(1,1)=4$,
$\kk(k,k)=2^{k}$ for $k\ge2$ (instead of $2^{k+1}$ in \refE{E2large+}),
$\kk(1,k)=\kk(k,1)=\eps_k$ for $k\ge2$, and $\kk(i,j)=0$
otherwise.

Furthermore, again as in \refE{E2large+}, let $H_k$ be the subgraph of $\gnk$
induced by the $n_k\sim\Bi(n,2^{-k})$ vertices of type $k$.
Then, conditional on $n_k$, each $H_k$ has the distribution
of the Erd\H os--R\'enyi graph $G(n_k,\kk(k,k)/n)$. As before,
$H_1$ is supercritical, which implies that $\gnk$ is supercritical.

Let $k_n\to\infty$ slowly, and choose $\eps_k$ such that
$\eps_{k_n}\le n^{-2}$; to be specific, set $k_n\=\ceil{\log_2\log
n}$ and $\eps_k\=\exp(-2^{k+2})$. Then \eqref{nov} implies that
$H_{k_n}$ is a critical Erd\H os--R\'enyi graph: the
edge probability is not exactly one over the number $n_{k_n}$ of
vertices, but is
$(1+o(n_{k_n}^{-1/3}))/n_{k_n}$, which is within the ``scaling
window''. It follows, as in \refE{E2large+}, that \whp\
$C_2(\gnk)\ge C_1(H_{k_n})=\Theta_p\bigpar{\bar  n_{k_n}^{2/3}}
>n^{2/3}/\log n$.

It is easy to see, analytically from \refT{Trho} or
probabilistically from \refT{T1A},
that $\rho(k)\to0$ as $k\to\infty$.
Consider the
graph $\gnxxh{\kk_n}$ in \refT{Tdual};
the norm of the corresponding integral operator $\htk$ on $L^2(\hmu)$
is at least the norm when restricted
to $\set{k}$, which is exactly $1-\rho(k)$, $k\ge2$.
Hence the norm is at least
1.
By \refR{RC2crit}, it follows that the norm of $\htk$ is exactly 1, i.e.,
that $\gnxxh{\kk_n}$ is critical.
\end{example}

We now turn to the proofs of Theorems \refand{T4}{Tdual},
starting with the subcritical case
of \refT{T4}, which we restate below.

\begin{theorem}
   \label{T4b}
Let $(\kkn)$ be a graphical sequence of kernels on a (generalized) vertex space 
$\vxs$ with limit $\kk$.
If\/ $\gk$ is subcritical,
\ie{}, $\norm{\tk}<1$, and $\sup_{x,y,n}\gk_n(x,y)<\infty$,
then
$C_1\xpar{G_n}=O(\log n)$ \whp.
\end{theorem}

\begin{proof}
As usual, we may assume that $\vxs$ is a vertex space; see \refSS{SSgvs}. 
Consider first the case when $\kkn=\kk$ for all $n$ and ${\sss}$ is finite, or,
equivalently, the \rfin\ case.
In this case, the result follows by comparing the neighbourhood 
exploration process of a vertex in the graph to a subcritical
branching process: this comparison is similar
to that made in the proof of \refL{L:step1}. This time,
setting $\go(n)=A\log n$, where $A$ is a (large) constant to be chosen below,
instead of the upper bound in \eqref{n2} we claim that
\begin{equation}\label{uexact}
 \P(x\in B)
\le
\rhox{\ge\go(n)}\bigpar{(1+2\eps)\kk;i},
\end{equation}
for all sufficiently large $n$.

This can be proved by the comparison argument used for \eqref{n2}, except
that in the final step, instead of using the total variation distance between
the binomial and Poisson distributions, we note that
if $\preceq$ denotes stochastic domination and 
$p'=-\log(1-p)$, so that $p'=p+O(p^2)$ as $p\to0$, then
$\Bi(1,p)\preceq\Po(p')$, and thus
$\Bi(m,p)\preceq\Po(mp')$ for every $m$.
Hence,
for $n$ large enough,
$\Bi(n_j',\kk(i,j)/n)\preceq \Bi(n_j,\kk(i,j)/n)
\preceq \Po\bigpar{(1+2\eps)\kk(i,j)\mu_j}$.
Note that in the argument
leading to \eqref{n2}, and hence in our proof of \eqref{uexact}, we do not
assume that $\kk$ is irreducible.

For $\eps$ small enough, the branching process $\bpx{(1+2\eps)\kk}$ is
subcritical.
Therefore,
there is an $a>0$ such that
\begin{equation}\label{ineq-GW}
\rhogek((1+2\eps)\kk;i)\le e^{-ak}
\end{equation}
holds for all $i$ and all $k\ge 1$.
Inequality \eqref{ineq-GW} is
undoubtedly  well known for finite-type Galton-Watson processes,
but for the sake of completeness we sketch a proof.
For $z>0$ let $g_z(x):=\E z^{|\bpkx|}$, where
$|\bpkx|$ denotes the total population of the branching process.
Using $\E w^{\Po(\lambda)}=e^{\lambda(w-1)}$ and independence of the
Poisson numbers of particles of each type in the first generation, we have
$g_z=ze^{\tk(g_z-1)}$. If $\kk$ is subcritical,
then by the implicit function theorem this functional equation has
a finite solution for $z$ in a
neighbourhood of 1, say for $|z-1|<\gd$, and it follows by an argument
similar to the proof of \refL{L1b} that indeed
$g_z(x)<\infty$ for all $x$ when $z<1+\gd$.

Taking $A=2/a$ for some $a$ for which inequality \eqref{ineq-GW} holds,
and recalling \eqref{uexact},  $\P(x\in B)\le
n^{-2}$ follows, and so \whp\ there are no vertices in $B$, i.e., there
is no component of order $A\log n$ or larger.

For the general case, it suffices to bound $\kappa_n$ for large
$n$ from above by a subcritical finitary $\kappa'$.
Recalling that $\sup_{x,y,n}\kkn(x,y)<\infty$,
let $\kkm+$ be defined as in \refL{L:uapprox}.
As $\norm{T_\kappa}<1$, by \refL{L:uapprox} we have
$\norm{\Tx{\kkm+}}<1$ for $m$ large enough,
and we can take
$\kk'=\kappa_m^+$.
\end{proof}

We now turn to the supercritical case, proving two results
that together imply part \ref{T4p2} of \refT{T4}.

\begin{theorem}
   \label{T4a}
Let $(\kkn)$ be a graphical sequence of kernels on a (generalized)
vertex space
$\vxs$ with irreducible limit $\kk$.
If\/ $\gk$ is supercritical,
\ie{}, $\norm{\tk}>1$, and $\inf_{x,y,n}\gk_n(x,y)>0$,
then
$C_2\xpar{G_n}=O(\log n)$ \whp.
\end{theorem}

\begin{proof}
As usual, we may assume that $\vxs$ is a vertex space. 
As the kernel $\kappa$ is supercritical, we have $\norm{T_\kappa}>1$,
so there is an integer $k$
such that $(1-1/k)\norm{T_\kappa}>1$. 
Fix such a $k$ throughout the proof.

Recall that $G_n=\gnxx{\kkn}$ is constructed by choosing in an
appropriate manner
a (deterministic or random) sequence
$\xs=(x_1,\ldots,x_n)$ giving the types of the vertices, and then
constructing
the edges using the kernel $\kkn$.
Independently of $G_n$, let us partition $V(G_n)=[n]$ 
into $k$ subsets $S_i$
in a random way, by independently assigning each vertex to a random
subset.
In other words, we construct $G_n$ (which has vertex set $[n]$)
and then partition its vertex set into $k$ classes
$S_i$.  Let $S_i^c\=[n]\setminus S_i$.

Let $A$ be a very large constant, to be chosen later.
We aim to prove that
the event $E$ that $G_n$ contains a component with more than $A \log
n$ and at most $n/A$
vertices has probability $o(1)$.
Every component $\cC$ of $G_n$ meets
some $S_i$ in at least $|\cC|/k$
vertices.
Let us say that a component $\cC$ of $G_n$ is \emph{bad} if it has at
least $A\log n$ and at most $n/A$
vertices, and meets $S_1$ in at least $A\log n/k$ vertices. Then, as
all $S_i$ are equivalent,
it suffices to prove that \whp\ $G_n$ has no bad component.

To this end, consider the subgraph $G_n'$ of $G_n$ induced by the 
vertices in $S_1^c$.
Let $\mu'=(1-1/k)\mu$, and let $\ys$ be the (random)
subsequence of $\xs$ corresponding to those vertices $i$ with $i\notin S_1$.
Note that $\vxs'=(\sss,\mu',\yss)$ is a generalized vertex space:
condition \eqref{a2a} for $\vxs'$
follows from the same condition for $\vxs$ and the random choice of
$S_1$.  Also, $G_n'$ has exactly the distribution of $G^{\vxs'}(n,\kkn)$.

Since $\E e(G_n')=(1-1/k)^2\E e(G_n)$,
and the sequence $(\kkn)$ is graphical on $\vxs$ with limit $\kk$, it
is graphical on $\vxs'$ with the same limit $\kk$.
As $\kk$ is irreducible on $(\sss,\mu)$, it is irreducible on $(\sss,\mu')$.
Let us write $\kk'$ for the kernel $\kk$ when viewed as a kernel on $(\sss,\mu')$.
Thus $T_{\kk'}$ is an operator on $L^2(\sss,\mu')$,
and, by assumption, $\norm{T_{\kappa'}}=(1-1/k)\norm{\tk}>1$.
Let $\cC_1$ be the largest component of $G_n'$ (chosen
according to any rule if there is a tie). Then, by \refT{T2},
\whp\
$\cC_1$ contains at least $\rho(\kappa')n/2$ vertices, and
$\rho(\kk')>0$.

{}From now on, we shall assume that $\cC_1$ has at least
$\rho(\kappa')n/2$ vertices, and choose $A$ so that
$1/A<\rho(\kappa')/2$.
If a component $\cC$ of $G_n$ is bad, then it cannot contain $\cC_1$, and
hence sends no edges to $\cC_1$.
Let $G_n''$ be the spanning subgraph of $G_n$
obtained from $G_n$ by deleting all edges between $\cC_1$ and vertices
in $S_1$.
If $\cC$ is a bad component of $G_n$, then all edges of $\cC$ are present
in $G_n''$, so $\cC$ is a component
of $G_n''$. Hence, the probability that $G_n$ has a bad component is
bounded by $o(1)$ (the probability that
$\cC_1$ is too small) plus the probability that some component $\cC\ne \cC_1$
of $G_n''$ containing at least
$A \log n/k$ vertices of $S_1$ sends no edges to $\cC_1$ in
$G_n$. Conditioning on $\xs$, $S_1$ and $G_n''$,
we have not tested any edges between $\cC_1$ and $S_1$, so these
edges are present independently,
each with its original probability. These individual probabilities are
all at least $(\infkappa)/n$, where $\infkappa\=\inf_{x,y,n}\gk_n(x,y)>0$
by assumption.

Thus, for any of the at most $n$ components $\cC\ne \cC_1$ of $G_n''$ with
at least $M=A \log n/k$ vertices in $S_1$,
the probability that $\cC$ sends no edges to $\cC_1$ is at most
\[
 (1-(\inf \kappa)/n)^{M|\cC_1|}.
\]
As $\infkappa>0$, $|\cC_1|\ge\rho(\kappa')n/2$ and $M=A \log n/k$, we can
make this probability
$o(n^{-100})$ by choosing $A$ large enough.
\end{proof}

Next, we turn to our general result
\refT{Tdual} on the distribution of the
small components
of a supercritical graph $G_n=\gnxx{\kkn}$. This will then be used to
prove the final statement in \refT{T4}, restated as \refT{T4c} below.

\begin{proof}[Proof of \refT{Tdual}]
The result is a simple consequence of Theorems \ref{T2}, \ref{T1A},
\refand{Tedges}{T2b}.
Indeed, given $G_n=\gnxx{\kkn}$, let $\cC_1$ be the largest component 
of $G_n$ (chosen according to any fixed rule if there is a tie), and,
for each $n$, let $\ys$ be the subsequence of $\xs$ consisting
of those $x_i$ for which $i\notin \cC_1$.
Then, by \refT{T1A}, $\hvxs=(\sss,\hmu,(\ys))$ is a generalized vertex space,
where $\hmu$ is defined by $\dd\hmu(x)=(1-\rho(\kk;x))\dd\mu(x)$;
indeed, the only non-trivial condition to verify is \eqref{nunA2}, which
is immediate from the same condition for $\vxs$ and \eqref{t1a}.

Next, we must show that the graphs $G_n'$ and $\gnxxh{\kkn}$ may be
coupled so that their edge-sets agree \whp. This is easy to see: as
usual, we condition throughout on $\xs$, treating $\xs$ as
deterministic. If we condition also on the vertex set $V$ of $\cC_1$,
the only information this gives about edges of $G_n$ inside
$V^c=V(G_n)\setminus V$ 
is that $G_n'=G_n[V^c]$ contains no component
larger than $|V|$ (and that certain components of order exactly $|V|$
are ruled out).  Without this condition, $G_n[V^c]$ would have exactly
the distribution of $\gnxxh{\kkn}$,
so it suffices to prove
that $C_1\bigpar{\gnxxh{\kkn}}<|V|$ holds \whp. In fact, as
$|V|\ge\rho(\kk)n/2$ \whp, it suffices to prove that
\[
 \eta\=\P\left(C_1\bigpar{\gnxxh{\kkn}} \ge \rho(\kk)n/2\right)
\]
tends to $0$ as $n\to\infty$.

Recall that the sequence $\xs$ is deterministic.  Pick a vertex $j$ of
$G_n$ at random, and explore its component in $G_n$ in the usual way,
by finding the neighbours of $j$, then the neighbours of the neighbours,
and so on.
Let $V_j$ be the vertex set of this component. The nature of the
exploration process ensures that, given $V_j$, the edges of
$G_n[V_j^c]$ are present independently with their unconditional
probabilities. In other words, writing $\ys'$ for the subsequence 
of $\xs$ consisting of those $x_i$ with $i\notin V_j$, and setting
$\hvxs'=(\sss,\hmu,(\ys'))$, the edge sets of $G_n[V_j^c]$ and
$G^{\hvxs'}(n,\kkn)$ have the same distribution.

Let us condition on the event $E$ that $|V_j|\ge \rho(\kk)n/2$.
As \refT{T2} can be applied to $G_n$, and $j$ was chosen at random,
we see that
$\P(E)$ is bounded
away from 0. Also, appealing to \refT{T2b}, we see that
\whp\ $G_n$ has a unique component
of size at least $\rho(\kk)n/2$. Thus, given $E$, we have $\ys=\ys'$ \whp.
Hence, with probability $\P(E)(o(1)+\eta)$, the graph $G_n$ contains
two components of order at least $\rho(\kk)n/2$. By \refT{T2b}, this
probability is $o(1)$, so we have $\eta=o(1)$, as required.
Finally, as $\kk$ is \qir\ on $(\sss,\mu)$, it is \qir\ on
$(\sss,\hmu)$. 

It remains only to show that the sequence
$(\kkn)$ is graphical on
$\hvxs$ with limit $\kk$.  Now $(\kkn)$ is graphical on $\vxs$ with
limit $\kk$, and
all conditions of \refD{Dg2} for $\hvxs$
apart from the last, \eqref{t2b},
follow immediately from the corresponding conditions for $\vxs$.
In other words, we must show that
\begin{align*}
 \E e\big(\gnxxh{\kkn}\big)/n
&\to
  \iint_{\sss^2}\kappa(x,y)\dd\hmu(x)\dd\hmu(y)
\\
 &=
 \iint_{\sss^2}\kappa(x,y)(1-\rho(\kk;x))(1-\rho(\kk;y))\dd\mu(x)\dd\mu(y).
\end{align*}
As we may couple $G_n'$ and $\gnxxh{\kkn}$ to agree \whp, and as the
number of edges
in either graph is bounded by that in $G_n$,
and thus, divided by $n$, is uniformly integrable,
it suffices to prove the
same limiting formula
for $\E e(G_n')/n$. This is immediate from the definition of $G_n'$,
condition \eqref{t2b}
for $G_n$, and \refT{Tedges}.
\end{proof}

\begin{theorem}
   \label{T4c}
Let $(\kkn)$ be a graphical sequence of kernels on a (generalized) vertex space
$\vxs$ with irreducible limit $\kk$.
If\/ $\gk$ is supercritical,
\ie{}, $\norm{\tk}>1$, and $\sup_{x,y,n}\gk_n(x,y)<\infty$,
then
$C_2\xpar{G_n}=O(\log n)$ \whp.
\end{theorem}

\begin{proof}
By \refT{Tdual},
\whp\ $C_2(G_n)=C_1(\gnxxh{\kk_n})$, and the result
  follows from Theorems \refand{Tdualsub}{T4b}.
\end{proof}

\begin{remark}
One might think that $C_2(G_n)=O(\log n)$ would hold \whp\ for
any supercritical kernel $\kk$,
i.e., that Theorems \refand{T4a}{T4c} would hold without the condition
$\inf_{x,y,n}\kk_n(x,y)>0$ or
$\sup_{x,y,n}\kk_n(x,y)<\infty$,
at least
for $\kkn=\kk$, say. However, it is easy to construct counterexamples,
similar to \refE{Ebad},
by taking
$x_i=i/n\in\sss=(0,1]$ and modifying a suitable kernel
to introduce a
largish star with centre $1$ not joined to the giant component.

Such pathologies are not the only counterexamples:
in the case where the $x_i$ are uniformly distributed,
for any function $\go(n)=o(n)$,
\refE{E2large+} provides an example of a supercritical random graph
with $C_2(G_n)>\go(n)$ \whp; thus
the $o_p(n)$ bound in \refT{T2b} is best possible.
\end{remark}

\section{Vertex degrees}\label{Sdegrees}

In this section we turn to the vertex degrees, proving \refT{TD}: 
if $\kkn$ is a graphical sequence of kernels
on a vertex space \vxs\ with limit $\kk$, and we define $\gl(x)$ by
\[
 \gl(x)\=\ints\kk(x,y)\dd\mu(y),
\]
then, writing $Z_k$ for the number of vertices of degree $k$ in $\gnxx{\kkn}$,
our aim is to show that
\begin{equation*}
  Z_k/n
\pto
\P(\Xi=k)
=
\ints\frac{\gl(x)^k}{k!}e^{-\gl(x)}\dd\mu(x),
\end{equation*}
where $\Xi$ has the mixed Poisson distribution $\ints\Po(\gl(x))\dd\mu(x)$.

In fact, \refT{TD} is stated for a generalized vertex space, and 
includes limiting results both for $Z_k/n$ and for
$Z_k/|V(G_n)|$. Since $|V(G_n)|/n\pto \mu(\sss)$, these results are
equivalent. As usual, the statement for generalized vertex spaces
reduces to that for vertex spaces (see \refSS{SSgvs}), and what we
must prove is exactly the statement above.

\begin{proof}[Proof of \refT{TD}]
Consider first the \rfin{} case in \refD{Drfin}.
Take a vertex $v$ of type $i$, let $D_v$ be its degree, and let
$D_{v,j}$ be the number of edges from $v$ to vertices of type $j$,
$j=1,\dots,r$; thus $D_v=\sum_jD_{v,j}$.
Assume that $n\ge\max\kk$ and condition on $n_1,\dots,n_r$.
Then the $D_{v,j}$ are independent for $j=1,\dots,r$, and
$D_{v,j}\sim\Bi\bigpar{n_j-\gd_{ij},\kk(i,j)/n}\dto\Po\bigpar{\mu_j\kk(i,j)}$;
hence
\begin{equation*}
  D_v\dto\Po\Bigpar{\sum_j\mu_j\kk(i,j)}
=\Po(\gl(i)),
\end{equation*}
as $\gl(i)=\int\kk(i,j)\dd\mu(j)=\sum_j\kk(i,j)\mu_j$.
Consequently,
\begin{equation*}
\P(D_v=k)
\to
\P\bigpar{\Po(\gl(i))=k}
= \frac{\gl(i)^k}{k!}e^{-\gl(i)}.
\end{equation*}
Let $Z_{k,i}$ be the number of vertices in \gnkx{} of type $i$ with degree $k$.
Then, still conditioning on $n_1,\dots,n_r$,
\begin{equation*}
\frac1n  \E Z_{k,i}
=\frac1n n_i\P(D_v=k)
\to \mu_i\P\bigpar{\Po(\gl(i))=k}.
\end{equation*}
It is easily checked that
$\Var(Z_{k,i}\mid n_1,\dots,n_r)=O(n)$.
Hence
\begin{equation*}
  \frac1n Z_{k,i} \pto \P\bigpar{\Po(\gl(i))=k}\mu_i,
\end{equation*}
and thus, summing over $i$,
\begin{equation*}
  \frac1n Z_{k}
=\sum_i  \frac1n Z_{k,i}
\pto \sum_i \P\bigpar{\Po(\gl(i))=k}\mu_i
=\P(\Xi=k).
\end{equation*}

This proves the theorem in the \rfin{} case.
In general, define $\kkm-$ by \eqref{kn-}.
Let $\eps>0$ be given. From \eqref{kkk} and
monotone convergence, there is an $m$ such that
\begin{equation}\label{ju}
\iint_{\sss^2}\kkm-(x,y)\dd\mu(x)\dd\mu(y)>
\iint_{\sss^2}\kk(x,y)\dd\mu(x)\dd\mu(y)-\eps.
\end{equation}
For $n\ge m$ we have $\kkm-\le \kk_n$ by \eqref{kknm},
so we may assume that
$\gnx{\kkm-}\subseteq\gnx{\kk_n}$. (Here, as usual,
we suppress the dependence on $\vxs$.)
Then, using \refP{PE} twice and \eqref{ju},
\begin{multline*}
\frac1ne\bigpar{\gnx{\kk_n}\setminus\gnx{\kkm-}}
=
\frac1ne\bigpar{\gnx{\kk_n}}-\frac1ne\bigpar{\gnx{\kkm-}}
\\
\pto
\tfrac12\iint_{\sss^2}\kk(x,y)\dd\mu(x)\dd\mu(y)
-
\tfrac12\iint_{\sss^2}\kkm-(x,y)\dd\mu(x)\dd\mu(y)
<\frac{\eps}2,
\end{multline*}
so \whp\
$e\bigpar{\gnx{\kk_n}\setminus\gnx{\kkm-}}<\eps n$.
Let us write $Z_k\qm$ for the number of vertices
of degree $k$ in $\gnx{\kkm-}$. It follows that \whp
\begin{equation}\label{ju1}
  |Z_k\qm-Z_k|<2\eps n.
\end{equation}
Writing $\Xi\qm$ for the equivalent of $\Xi$ defined using $\kkm-$
in place of $\kk$,
by the first part of the proof, $Z_k\qm/n\pto\P(\Xi\qm=k)$. Thus
\whp
\begin{equation}\label{ju2}
  |Z_k\qm/n-\P(\Xi\qm=k)|<\eps.
\end{equation}

Finally, we have $\E\Xi=\ints\gl(x)\dd\mu(x)=\iikxy$. Since
$\gl\qm(x)\le\gl(x)$, we can assume that $\Xi\qm\le\Xi$, and thus
\begin{multline}
  \label{ju3}
\P(\Xi\neq\Xi\qm)
=
\P(\Xi-\Xi\qm\ge1)
\le
\E(\Xi-\Xi\qm)
=\iint_{\sss^2}\kk
-\iint_{\sss^2}\kkm-
<\eps.
\end{multline}

Combining \eqref{ju1}, \eqref{ju2} and \eqref{ju3}, we see that
$ |Z_k/n-\P(\Xi=k)|<4\eps$ \whp.
\end{proof}

Let $\gL$ be the random variable $\gl(\xi)$, where $\xi$ is a random
point in $\sss$ with distribution $\mu$. Then we can also describe the
mixed Poisson distribution of $\Xi$ as $\Po(\gL)$.
Under mild conditions,
the tail probabilities $\P(\Xi>t)$ and $\P(\gL>t)$ are similar for
large $t$.
We state this for the case of power-law tails; the result generalizes to
regularly varying tails.
As above, let $D$ be the degree of a random vertex in $\gnkxn$.
Let $Z_{\ge k}$ be the number of vertices with degree $\ge k$.
\begin{corollary}
  \label{CD}
Let $(\kkn)$ be a graphical sequence of kernels on a vertex
space $\vxs$ with limit $\kk$.
Suppose that $\P(\gL>t)=\mu\set{x:\gl(x)>t}\sim a t^{-\ga}$ as \ttoo{}
for some
$a>0$ and $\ga>1$.
Then
\[
 Z_{\ge k}/n \pto  \P(\Xi\ge k)\sim a k^{-\ga},
\]
where the first limit is for $k$ fixed and $n\to\infty$,
and the second for \ktoo.
In particular,
$\lim_{\ntoo}\P(D \ge k) \sim a k^{-\ga}$
as \ktoo.
\end{corollary}
\begin{proof}
It suffices to show that $\P(\Xi\ge k)\sim a k^{-\ga}$;
the remaining conclusions then follow from \refT{TD}.
  For any $\eps>0$,
$\P\bigpar{\Po(\gL)>t\mid \gL>(1+\eps)t}\to1$
and
$\P\bigpar{\Po(\gL)>t\mid \gL<(1-\eps)t}=o(t^{-\ga})$
as \ttoo, for example by standard Chernoff estimates
\cite[Remark 2.6]{JLR}. It follows that
$\P(\Xi>t)=\P\bigpar{\Po(\gL)>t}\sim a t^{-\ga}$ as \ttoo{}.
\end{proof}

This result shows that our model does include natural cases with
power-law degree distributions. For example, taking $\sss=(0,1]$
with the Lebesgue measure, and
$\kkn(x,y)=\kk(x,y)=c/\sqrt{xy}$ for $c>0$ constant, we have
$\gl(x)=2c/\sqrt{x}$, so $\P(\gL>t)=4c^2/t^2$ for $t\ge 2c$.
Thus, by \refC{CD}, $\P(D\ge k)\sim 4c^2/k^2$ as $\ktoo$.
In fact, in this case $\gnkx$ is the `mean-field'
version of the Barab\'asi--Albert scale-free model; see
\refSS{SSrtij}. For other power laws, see \refSS{SSrank1}.

\section{Distances between vertices}\label{Sdist}

One of the properties of inhomogeneous graphs that has received
much attention is their `diameter'. For example, considering the
scale-free model of Barab\'asi and Albert~\cite{BAsc},
the diameter was determined heuristically and experimentally
to be $\Theta(\log n)$ in~\cite{N401,BAJp2,nswpp}; for a precise
version of this model, the LCD model, the value $(1+o(1))\log n/\log\log n$
was found rigorously in~\cite{diam}; later,
this value was also found heuristically
in~\cite{CH}.

Often, the diameter is taken to mean the average distance between
a random pair of vertices, or perhaps the `typical' distance,
although the usual graph theoretic definition (the maximum
distance between a pair of vertices) is also used. Here we shall 
consider both interpretations.

\subsection{Typical distances}
In this subsection we study the `typical' distance between vertices;
our aim is to prove \refT{Tdist},
giving upper and lower bounds on the distances
between almost all pairs of vertices, showing that almost all pairs
of vertices in the giant component are at distance roughly $\log n/\lntk$. 

Many related results have been published, concerning random graphs
with a fixed degree sequence, or random graphs with a given expected
degree sequence; we shall only describe a few here.
These models are similar to (and in some cases
special cases of) the rank 1 case of our model; see \refSS{SSrank1}.
For example, Chung and Lu~\cite{ChungLu:dist2002,ChungLu:dist2003}
studied distances in a `random graph with given expected degrees'.  For
power-law degrees with exponent $\beta>3$, where their model is a
special case of ours, they obtained an asymptotic diameter of $\log n$
over the log of the average of the squares of the degrees, a special case of \refT{Tdist};
see \refSS{SSrank1} for the connection to our model.

Van der Hofstad, Hooghiemstra and Van Mieghem~\cite{HHM_RSA}
(see also \cite{HHZ1,HHZ2,HHZ3})
studied a model where the vertex degrees are \iid{} with a certain
distribution, 
and the graph is chosen uniformly among all graphs with these degrees.
They analyze the growth of vertex neighbourhoods in this model by
using a branching process;
this process is single type, but the number of children of a particle
is not Poisson, 
so it is rather different from the one considered here. They obtain
very precise 
results on the distances between a random pair of vertices, showing that
it is $\log n/\log c+O_p(1)$, where $c$ is the
expectation of $d(v)(d(v)-1)$. (The $-1$ here comes from the degrees being
(conditionally) fixed, rather than essentially Poisson.)

There are many other papers in this area, both heuristic and
mathematical; we shall not 
attempt to list them. Let us mention only that Fernholz and
Ramachandran~\cite{FR:diam}, 
while mainly focussing on the diameter (see \refSS{SSdiam}) also treat
the typical 
distance between vertices. For further references we refer
the reader to the discussion of related work in~\cite{HHM_RSA}. 

\medskip
Let us now begin our preparation for the proof of \refT{Tdist}. 
Let $\kkn$ be a graphical sequence of kernels on a vertex space $\vxs$
with limit $\kk$, 
and let $G_n=\gnxx{\kkn}$. Note that we do not consider generalized
vertex spaces; arguing as in \refSS{SSgvs}, to prove \refT{Tdist} it suffices
to consider vertex spaces.
(For part (iv), we use also the fact that, by standard arguments, the
conclusion holds if and only if \eqref{pi} holds for $f(n)=\eta\log
n$, for every $\eta>0$.)
We shall write
$d_G(v,w)$ for the graph distance between two vertices $v$, $w$ of a graph
$G$, taking $d_G(v,w)=\infty$ if $v$ and $w$ are not in the same component
of $G$. When the graph $G$ is not specified, $G_n=\gnxx{\kkn}$ is to
be understood. 

\begin{lemma}
  \label{Ldfin}
If $\kk$ is \qir, then 
\begin{equation}\label{dfin}
 \bigl|\bigl\{\{v,w\} : d(v,w) <\infty\bigr\}\bigr|
= \frac{C_1(G_n)^2}2+o_p(n^2) = \frac{\rho(\kk)^2n^2}2 +o_p(n^2).
\end{equation}
\end{lemma}
\begin{proof}
As noted in \refSS{SSdist}, by \eqref{nconn} this is
immediate from Theorems \refand{T2}{T2b}.
\end{proof}

In particular, almost all pairs with either or both vertices outside the giant
component are not connected at all, so we shall study
the typical distance only in the supercritical case $\norm{\tk}>1$.

\begin{lemma}\label{Ldrf}
Let $\kk$ be a \rfin{} kernel on a vertex space \vxs\ with $\norm{\tk}>1$.
For any $\eps>0$,
\[
 \E\bigl|\bigl\{\{v,w\} : d(v,w)
  \le (1-\eps)\log n/\log\norm{\tk}\bigr\}\bigr|
= o(n^2).
\]
\end{lemma}
\begin{proof}
Changing only the notation, we may assume that the type space $\sss$ is finite,
say $\sss=\{1,2,\ldots,r\}$.
It turns out that, as usual, we may assume that $\mui>0$ for every $i$; however,
here we cannot simply ignore $o_p(n)$ edges, so an argument is needed.
Suppose that $\mui=0$ for some $i$. Taking
$\kk'(i,j)=\kk'(j,i)=\max\kk$ for all $j$, 
and $\kk'(j,k)=\kk(j,k)$ for $j,k\ne i$, we have $\kk\le \kk'$, so we may
couple $G_n\=\gnkx$ with $G_n'=\gnxx{\kk'}$ so that $G_n\subseteq G_n'$.
Note that $\norm{\tk}=\norm{T_{\kk'}}$, as $\kk=\kk'$ a.e.

For any $\eta>0$, define $\mu'$ by $\mu'\{i\}=\eta$ and
$\mu'\{j\}=(1-\eta)\muj$, $j\ne i$. 
Thus $\mu'$ is obtained from $\mu$ by shifting
some measure from types other than $i$ to type $i$. 
Changing the types of some vertices correspondingly, we obtain a vertex
space $\vxs'=(\sss,\mu',\xssd)$ such that whenever a vertex has type $j$
in $\vxs$, it has either type $j$ or type $i$ in $\vxs'$.
As $\kk'(j,k)$ is maximal when one or both of $j$ and $k$ is equal to $i$,
it follows that we can couple $G_n'$ and
$G_n''=\gnxxp{\kk'}$ so that $G_n'\subseteq G_n''$.
As $\eta\to 0$, the norm of $T_{\kk'}$ defined with respect to $\mu'$ tends
to the norm defined with respect to $\mu$. Since $G_n\subseteq G_n''$, to prove
\refL{Ldrf} for $G_n$, it thus suffices to prove the same result
for $G_n''$,
defined on a vertex space with $\mu'\{i\}>0$.
Iterating, it suffices
to prove \refL{Ldrf} in the case where $\mui>0$ for every $i$.

Let $\Ga_d(v)=\Ga_d(v,G_n)$ denote the {\em $d$-distance set} of $v$ in $G_n$,
i.e., the set of vertices of $G_n$ at graph
distance exactly $d$ from $v$, and let $\Gale d(v)=\Gale d(v,G_n)$ denote
the {$d$-neighbourhood} $\bigcup_{d'\le d}\Ga_{d'}(v)$ of $v$.

Let $0<\eps<1/10$ be arbitrary.
The proof of \eqref{uexact}
involved first showing that, for $n$ large enough,
the neighbourhood exploration process
starting at a given vertex $v$ of $G_n$ with type $i$
(chosen without inspecting $G_n$)
could be coupled with
the branching process $\bpx{(1+2\eps)\kk}(i)$ so that the branching process
dominates. In particular, the two processes can be coupled
so that
for every $d$, $|\Ga_d(v)|$ is at most the number $N_d$ of particles
in generation $d$ of $\bpx{(1+2\eps)\kk}(i)$.
Elementary properties of the branching process imply that
$\E N_d=O(\norm{T_{(1+2\eps)\gk}}^d)=O(((1+2\eps)\gl)^d)$,
where $\gl=\norm{\tk}>1$.

Set $D=(1-10\eps)\log n/\log\gl$. Then $D<(1-\eps)\log n/\log((1+2\eps)\gl)$
if $\eps$ is small enough, which we shall assume. Thus,
\[
 \E|\Gale D(v)| \le \E\sum_{d=0}^D N_d = O(((1+2\eps)\gl)^D) = O(n^{1-\eps}) = o(n).
\]
Summing over $v$, the expected number of pairs of vertices within distance
$D$ is $o(n^2)$, and the result follows.
\end{proof}

We now turn to the reverse bound, showing that most vertices in the giant
component are within distance roughly $\log n/\log\norm{\tk}$. First we
consider two random vertices.

\begin{lemma}\label{Ldrf2}
Let $\kk$ be a \qir\ \rfin{} kernel on a vertex space \vxs\ with
$\norm{\tk}>1$, 
and let $v$ and $w$ be two vertices of $\gnkx$ chosen independently
and uniformly. 
Then, for any $\eps>0$,
\[
 \P\bigpar{d(v,w) < (1+\eps)\log n/\log\norm{\tk}} \to \rho(\kk)^2
\]
as $n\to\infty$.
\end{lemma}

\begin{proof}
Note that an upper bound $\rho(\kk)^2+o(1)$ follows from \refL{Ldfin},
so it suffices to prove a corresponding lower bound. 
This time we may simply ignore types $i$ with $\mui=0$, working
entirely within the subgraph $G_n'$ of $G_n=\gnkx$
induced by vertices of the remaining types.
As there $o_p(n)$ vertices of types $i$ with $\mui=0$, changing $\eps$
slightly it suffices to prove the result for $G_n'$. Thus we shall
assume that $\mui>0$ for every $i$. Also, restricting to a suitable
subset of the types and renormalizing, we may and shall assume
that $\kk$ is irreducible.

Fix $0<\eta<1/10$. We shall assume that $\eta$ is small enough that
$(1-2\eta)\gl>1$, 
where $\gl=\norm{\tk}$.
In the argument leading to \eqref{n2} in proof of \refL{L:step1},
we showed that, given $\go(n)$ with $\go(n)=o(n)$
and a vertex $v$ of type $i$, the neighbourhood exploration process of
$v$ in $G_n$ 
could be coupled with the branching process $\bpx{(1-2\eta)\kk}(i)$
so that \whp\ the former dominates until it reaches size $\go(n)$.
More precisely, writing $N_{d,k}$ for the number of particles
of type $k$ in generation $d$ of $\bpx{(1-2\eta)\kk}(i)$,
and $\Ga_{d,k}(v)$ for the set of type-$k$ vertices
at graph distance $d$ from $v$, \whp\
\begin{equation}\label{coupl}
 |\Ga_{d,k}(v)|\ge N_{d,k},\, k=1,\ldots,r,
\text{ for all $d$ s.t. }|\Gale d(v)|<\go(n).
\end{equation}
The key point is that this coupling works because we have only `looked at'
$o(n)$ vertices at each step.

Let us call a kernel $\kk$ \emph{bipartite} if
\begin{equation}\label{bipdef_moved}
 \sss = L\cup R,\hbox{ with }\kk(i,j)=0\hbox{ whenever }i,j\in L\hbox{ or }
 i,j\in R,
\end{equation}
in which case the graph $G_n$ is bipartite.
For the moment, let us suppose that $\kk$ is not bipartite.
Let
$N_t(i)$ be the number of particles of type $i$ in the $t$th
generation of $\bpk$, and let $\vNt$ be the vector
$N_t(1),\ldots,N_t(r)$. 
Also, let $\vnu=(\nu_1,\ldots,\nu_r)$ be the eigenvector
of $\kk$ with eigenvalue $\gl$
(unique, up to normalization, as $\kk$ is irreducible).
From standard branching process results,
for example, \cite[Theorems V.6.1 and V.6.2]{AN},
we have
\begin{equation}\label{nt_moved}
 \vNt/\lambda^t\to X\vnu \quad\text{a.s.,}
\end{equation}
where $X\ge 0$ is a real-valued
random variable, $X$ is continuous except that it has some mass at
$0$, and $X=0$ if and only if the branching process
eventually dies out.

Let $D$ be the integer part of $\log(n^{1/2+2\eta})/\log((1-2\eta)\gl)$.
From \eqref{nt_moved},
\whp\ either $N_D=0$,
or $N_{D,k}\ge n^{1/2+\eta}$ for each $k$.
Furthermore, as
$\lim_{d\to\infty}\P(N_d\neq0)=\rho((1-2\eta)\kk)$ and $D\to\infty$,
we have $\P(N_D\neq0)\to \rho((1-2\eta)\kk)$.
Thus, if $n$ is large enough,
\[
 \P\left(\forall k: N_{D,k}\ge n^{1/2+\eta}\right) \ge \rho((1-2\eta)\kk)-\eta.
\]
By \refT{TappB}, the right-hand side tends to $\rho(\kk)$ as $\eta\to 0$.
Hence, given any fixed $\gamma>0$, if we choose $\eta$ small enough we have
\begin{equation}\label{bpl}
 \P\left(\forall k: N_{D,k}\ge n^{1/2+\eta}\right) \ge \rho(\kk)-\gamma
\end{equation}
for $n$ large enough.
It is easy to check that $\E(|\Gale D(v)|)=o(n^{2/3})$ if $\eta$ is
small enough; 
for example, we may argue as in the proof of \refL{Ldrf}. Hence,
\begin{equation}\label{GaleD}
 |\Gale D(v)|\le n^{2/3}\quad\text{\whp},
\end{equation}
and \whp\ the coupling described in \eqref{coupl} extends at least to the
$D$-neigh\-bour\-hood.

Now let $v$ and $w$ be two fixed vertices of $\gnkx$, of types $i$ and $j$
respectively. We explore both their neighbourhoods at the same time,
stopping either when we reach distance $D$ in both neighbourhoods,
or we find an edge from
one to the other, in which case $v$ and $w$ are within graph distance
$2D+1$. We consider two independent branching processes
$\bpx{(1-2\eta)\kk}(i)$, 
$\bpx{(1-2\eta)\kk}'(j)$, with $N_{d,k}$ and $N_{d,k}'$ vertices of type
$k$ in generation $d$ respectively. By \eqref{GaleD}, \whp\ we encounter
$o(n)$ vertices in the explorations so, by the argument leading to
\eqref{coupl}, 
\whp\ either the explorations meet, or
\begin{equation*}
 |\Ga_{D,k}(v)|\ge N_{D,k} \quad \text{and}\quad
 |\Ga_{D,k}(w)| \ge N_{D,k}',
 \quad k=1,\ldots,r.
\end{equation*}

Using \eqref{bpl} and the independence of the branching processes, it follows
that
\begin{multline}\label{eor}
 \P\left( d(v,w)\le 2D+1\hbox{ or } \forall k:
 |\Ga_{D,k}(v)|,|\Ga_{D,k}(w)| \ge n^{1/2+\eta}\right)
\\
  \ge (\rho(\kk)-\gamma)^2-o(1).
\end{multline}
Conditional on the second event in \eqref{eor} holding and not the first,
we have not examined any edges from $\Ga_D(v)$ to $\Ga_D(w)$, so these edges
are present independently with their original unconditioned probabilities.
For any $i'$, $j'$, the expected number of these edges is at least
$|\Ga_{D,i'}(v)||\Ga_{D,j'}(v)|\kk(i',j')/n$. Choosing $i'$, $j'$ such that $\kk(i',j')>0$,
this expectation is $\Omega((n^{1/2+\eta})^2/n)=\Omega(n^{2\eta})$.
It follows that at least one edge is present
with probability $1-\exp(-\Omega(n^{2\eta}))=1-o(1)$. If such an edge
is present, then $d(v,w)\le 2D+1$.
Thus, \eqref{eor} implies that
\[
 \P(d(v,w)\le 2D+1) \ge (\rho(\kk)-\gamma)^2-o(1) \ge \rho(\kk)^2-2\gamma-o(1).
\]
Choosing $\eta$ small enough, we have $2D+1\le (1+\eps)\log n/\log\gl$.
As $\gamma$ is arbitrary, we have
\[
 \P(d(v,w)\le (1+\eps)\log n/\log\gl) \ge \rho(\kk)^2-o(1),
\]
and the lemma follows.

The argument for the bipartite case is essentially the same, except that
if $v$ and $w$ are of types in the same class of the bipartition,
we should look for an edge between $\Ga_D(v)$ and $\Ga_{D-1}(w)$.
\end{proof}

Lemmas \refand{Ldfin}{Ldrf2} have the following immediate
consequence.

\begin{corollary}\label{Cdrf2}
Let $\kk$ be a \qir\ \rfin{} kernel on a vertex space \vxs\ with
$\norm{\tk}>1$.
For any $\eps>0$,
\begin{equation}\label{dub}
 \bigl|\bigl\{\{v,w\} : d(v,w) \le (1+\eps)\log n/\log\norm{\tk}\bigr\}\bigr| 
= \rho(\kk)^2n^2/2 +o_p(n^2).
\end{equation}
\end{corollary}

\begin{proof}
It follows from \refL{Ldfin} that the expected number
of vertex pairs \set{v,w} with $d(v,w)<\infty$ is $\rho(\kk)^2n^2/2+o(n^2)$.

Fix $\eps>0$.
{}From \refL{Ldrf2}, 
the expected number of pairs of vertices at distance less than
$d=(1+\eps)\log n/\log\norm{\tk}$ is $\rho(\kk)^2n^2/2+o(n^2)$.
Hence, the expected number of pairs with $d\le d(v,w)<\infty$ is $o(n^2)$,
so there are $o_p(n^2)$ such pairs. Using \eqref{dfin} again,
\eqref{dub} follows. 
\end{proof}

After this preparation it is easy to deduce \refT{Tdist}. As noted 
earlier, it suffices to consider vertex spaces, rather than
generalized vertex spaces.

\begin{proof}[Proof of \refT{Tdist}]
Let $\kkn$ be a graphical sequence of kernels on a vertex space $\vxs$ with limit $\kk$,
let $G_n=\gnxx{\kkn}$, and let $\eps>0$ be fixed.
We must prove four statements, which we recall separately below.

\pfitemx{\ref{Tdist0}}
The first part of \refT{Tdist} is exactly \refL{Ldfin} (but with
\qir\ replaced by irreducible), which we have already proved. 

\pfitemx{\ref{Tdistu}}
We must show that if $\supxyn<\infty$, then only $o_p(n^2)$
vertices of $G_n$ are within distance $(1-\eps)\log n/\log\norm{\tk}$.
As usual, we approximate with the \rfin\ case. Let $\kkm+$ be a sequence
of \rfin\ kernels on $\vxs$ with the properties guaranteed by \refL{L:uapprox}.
By \refL{L:uapprox}\ref{L:uapprox3}, $\norm{\Tx{\kkm+}}\to \norm{\tk}>1$,
so there is an $m$ such that
$(1-\eps/2)/\log\norm{\Tx{\kkm+}} \ge (1-\eps)/\log\norm{\tk}$.
Fixing such an $m$, we may couple $G_n$ and $\gnxx{\kkm+}$ so that
$G_n\subseteq \gnxx{\kkm+}$ for $n\ge m$, and the result follows by
applying \refL{Ldrf} to $\gnxx{\kkm+}$ with $\eps/2$ in place of $\eps$.

\pfitemx{\ref{Tdistl}}
This time we must show that if $\kk$ is irreducible and $\norm{\tk}<\infty$, then
$\rho(\kk)^2n^2/2+o_p(n^2)$ pairs of vertices of $G_n$ are within distance
$(1+\eps)\log n/\log\norm{\tk}$. By \eqref{dfin}, it suffices
to prove the lower bound. Again, we approximate with the \rfin\ case,
this time working with a graph $G_n'=\gnxx{\tkkm}\subseteq G_n$.
The argument is as above, but using the approximating
kernels $\tkkm$ given by \refL{L:approx} instead of $\kkm+$,
and applying \refC{Cdrf2} with $\eps/2$ in place of $\eps$:
by \refL{L:approx}\ref{L:approxa}, the \qirity\ condition of \refC{Cdrf2}
is satisfied, while \refL{L:approx}\ref{L:approxb} implies
$\norm{\Tx{\tkkm}}\upto \norm\tk$. 

\pfitemx{\ref{Tdistl2}}
This time we must show essentially that if $\kk$ is irreducible and
$\norm{\tk}=\infty$ 
then almost all pairs of vertices in the giant component
are within distance $o(\log n)$.
First, let $\eta>0$. By the same proof as for part \ref{Tdistl} above,
except that 
we have $\norm{\Tx{\tkkm}}\to\infty$, we see that
$\rho(\kk)^2n^2/2+o_p(n^2)$ pairs of vertices of $G_n$ are within distance
$\eta\log n$. Since $\eta$ is arbitrary, a standard argument shows
that we can replace $\eta\log n$ by some function $f(n)=o(\log n)$.
\end{proof}

\begin{remark} \label{Rsmalldist}
Part \ref{Tdistu} of \refT{Tdist} does not hold if we omit the condition
that $\supxyn<\infty$, even if $\kappa_n=\kappa$ for all $n$, 
with $\norm{\tk}<\infty$.
To see this,
let $\sss=[0,1]$ with $\mu$ the Lebesgue measure,
and let $\xss$ be disjoint deterministic 
sequences such that $\vxs=(\sss,\mu,\xss)$ is a vertex
space. We shall write $\xs$ as $(x^{(n)}_1,\ldots,x^{(n)}_n)$ to
emphasize the dependence of
the terms on $n$;
for example, we may take $x^{(n)}_i=(i-\sqrt2/2)/n$.
Taking $\kk(x,y)=2$ for all $x$, $y$, 
the graph $\gnxx{\kk}$ is a supercritical Erd\H os--R\'enyi random
graph.

Forming $\kk'$ by modifying $\kk$ on a set of 
measure zero, we can effectively add or delete $o(n)$
given edges to/from $\gnxx{\kk}$ whilst keeping
$\kk'$ graphical on $\vxs$ with  $\kk'=\kk=2$ a.e.
In particular, given $f(n)=o(n)$,
taking $\kk'(x^{(n)}_1,x^{(n)}_i) = n^2$, say, for $1\le i\le f(n)$,
we  may ensure that in the graph $G_n'=\gnxx{\kk'}$,
\whp\ the vertex $1$ is joined to all of the vertices $2,3,\ldots,f(n)$.
By \refT{T2}, the giant component still has $\rho(2)n+o_p(n)$ vertices.

For any $\omega(n)\to\infty$ it is easy to check that if we choose
$f(n)$ large enough,
all but $o_p(n)$ vertices in the giant component
are within distance $\omega(n)$ of one of the vertices
$1,2,\dots,f(n)$, and thus,
all but $o_p(n^2)$ pairs of
vertices in the giant component
are within distance $2\omega(n)+2$. Hence, even if $\norm{\tk}$ is bounded,
the typical distance between vertices may be smaller than any given
function tending to infinity.

Even if we allow $\norm\tk=\infty$,
the typical distance 
cannot be as small as a constant:
one can check that when
$\kk$ is irreducible, for any $C$ there are \whp{}
$\Theta(n^2)$ pairs of vertices in the giant component
at distance at least $C$. In fact, there
are $\Theta(n)$ vertices in the giant component
whose $C$-neighbourhood is a path. This can be proved using
a combination of the arguments leading to Theorems \refand{T2}{T3}.
\end{remark}

\begin{remark}\label{Rlargedist}
Using the same vertex space as in \refR{Rsmalldist},
for any $f(n)=o(n)$,
we can define a graphical sequence of kernels $\kkn''$ with (irreducible)
limit $\kk=2$ such that \whp\ $\gnxx{\kkn''}$ is obtained from
the Erd\H os-R\'enyi graph $G(n,2/n)$ by deleting $f(n)+1$ vertices
and replacing them with a path of length $f(n)$ not joined to the rest
of the graph. In this graph there are at least $(f(n)/3)^2$ pairs
of vertices at distance at least $f(n)/3$. Hence
the average distance between vertices
(counting only pairs at finite distance)
is at least $\Omega(f(n)^3)/n^2$, which may be much larger
than $O(\log n)$: in fact, it may be larger than any
given function that is $o(n)$.
This is the reason for considering the distances between
almost all pairs rather than the average distance
in \refT{Tdist}\ref{Tdistu}.
\end{remark}

\begin{remark}\label{ologn}
In certain cases, we know better bounds than $o(\log n)$
on the typical distances between vertices.
For example, the $(1+o(1))\log n/\log\log n$ formula
for the diameter of the $m=2$ LCD model proved by Bollob\'as and Riordan~\cite{diam}
certainly holds as a bound on the typical distances
in the much simpler `mean-field' case described in \refSS{SSrtij},
where $\kk(x,y)=1/\sqrt{xy}$ and $x_i=i/n$.

However, without further restrictions (which could be on $\kk$,
or on the distributions of the $x_i$),
we cannot strengthen the $o(\log n)$
bound in part \ref{Tdistl2} of \refT{Tdist}. Indeed, given a graphical kernel
$\kk$ on a vertex space $\vxs=(\sss,\mu,\xss)$ with $\norm{\tk}=\infty$,
and any function $g(n)\to\infty$, setting $\kkn=\kk\bmin g(n)$
we have 
$\E|\Ga_d(v)|\le g(n)^d$ for a given vertex $v$ and
every $d\ge0$. Indeed, we have $\gnkn\subseteq\gnx{g(n)/n}$.
The argument in the proof of \refL{Ldrf}
shows that the typical distance is at least $(1+o(1))\log n/\log g(n)$.
Hence the $o(\log n)$ bound in part \ref{Tdistl2} is best possible.
A similar example may be constructed with a fixed kernel $\kk$ by modifying the
sequences $\xss$ appropriately; 
take for example $\kk(x,y)=1/\sqrt{xy}$ and
$x_i=\max\bigset{i/n,1/g(n)}$.
\end{remark}

\subsection{The diameter}\label{SSdiam}

Let $\kk$ be a kernel on a (generalized) vertex space $\vxs$ 
in which the set of types is finite.
In this subsection we study the diameter of $G_n=\gnkx$,
measured in the usual
graph theoretical sense for disconnected graphs:
\[
 \diam(G_n) \= \max\{d(v,w)\: :\: v,w\in V(G),\:d(v,w)<\infty\},
\]
where $d(v,w)$ is the graph distance between $v$ and $w$ in $G_n$.

The following is a partial list of existing
work on the diameter of sparse random graphs: 
Bollob\'as and Fernandez de la Vega~\cite{BF} 
found the asymptotic diameter of random $r$-regular graphs,
\L uczak~\cite{Luczak:diam} obtained detailed results
for $G(n,c/n)$ with $c<1$, Chung and Lu~\cite{ChungLu:diam}
studied $G(n,c/n)$, $c>1$, and Fernholz and Ramachandran~\cite{FR:diam}
obtained a precise result for random graphs with \iid{} degrees (see below). 

When $G_n$ has finite-type, provided $\kk$
is not critical we can easily find the diameter of
$G_n$ in the form $(c+o(1))\log n$.
The constant $c=c(\kk)$ will be obtained from the branching
process $\bpk$ in a simple way, different in the sub- and super-critical cases;
see Theorems \refand{th_dbelow}{th_dabove} below. Together, these results,
which we shall prove separately, constitute \refT{th_diam}.

The results in this subsection correspond to those of 
Fernholz and Ramachandran~\cite{FR:diam} for a different model,
where the distribution of the vertex degrees is fixed, the vertex
degrees are sampled independently from this distribution and,
conditional on the degree sum being even,
the graph is then chosen uniformly at random from all graphs with the given degree sequence.
The special case of our model where $\kk$ has rank one is a special case
of this model; see \refSS{SSrank1}. In general, the two models are different.
The proofs in \cite{FR:diam} are much more complicated
than those we shall present here, because their model does not
have independence built in. Thus,
roughly speaking, Fernholz and Ramachandran  have to work to
get the branching process approximation
that we have here as our starting point.
Also, here we keep things simple by considering only the finite type 
case; as noted in \refSS{SSdist}, even the single type case is non-trivial.

Throughout this subsection, when exploring the neighbourhoods of a vertex $v$
in $G_n=\gnkx$, we fix in advance an arbitrary order on the vertices
of $G_n$. At each step in the exploration, among unexplored
vertices at minimal distance from $v$, we choose the first vertex $w$
in this order, and reveal all edges from $w$ to vertices
not yet reached by the exploration.
In this way we reveal the vertex sets of the
neighbourhoods $\Ga_t(v)$, $t=1,2,\ldots$ successively. Furthermore,
the graph we reveal is always a tree, rooted at $v$. We shall denote this
graph by $T(v)$, and call it the {\em reduced component} of $v$.
For $t\ge 0$ we write $T_t(v)$ for $T(v)\cap \Gale{t}(v)$, the
{\em reduced $t$-neighbourhood} of $v$.

Specifying which unexplored vertex to choose next does not affect
the coupling arguments leading to \eqref{n2}, for example, where
any unexplored vertex could be chosen at each step. The advantage
is that the tree $T(v)$ is uniquely specified even if the component
containing $v$ has cycles; below we shall sum the probability
that the {\em reduced} component of a vertex is a particular
tree over all trees. Using the reduced component guarantees
that the corresponding events are disjoint. Note that if $v$ and $w$
lie in the same component, then $d(v,w)$ is the same as the graph
distance in $T(v)$ between the root, $v$, and $w$.

We shall first prove the subcritical case of \refT{th_diam}, restated below.

\begin{theorem}\label{th_dbelow}
Let $\kk$ be a kernel on a (generalized) vertex space $\vxs=(\sss,\mu,(\xs))$, 
with $\sss=\{1,2,\ldots,r\}$ finite and
$\mui>0$ for each $i$.
If $0<\norm\tk<1$, then
\[
 \frac{\diam(G_n)}{\log n} \pto \frac{1}{\lntki}
\]
as $n\to\infty$, where $G_n=\gnkx$.
\end{theorem}

\begin{proof}
We may assume without loss of generality that $\kk$ is irreducible.
Also, by conditioning on the sequences $(\xs)$, we may assume that
$\vxs$ is a vertex space, and that the number $n_i$ of vertices
of type $i$ is deterministic (see \refSS{SSgvs}). 

Let $p_{d,i}$ be the probability that the branching process $\bpk(i)$
survives for at least $d$ generations. 
As the number of particles in the first generation that have descendants
in generation $d+1$ has a Poisson
distribution, we have
\[
 p_{d+1,i} = 1-\exp\Bigpar{-\sum_j \kk(i,j)\muj p_{d,j}}.
\]
Recalling that $p_{d,i}\downto \rho(\kk,i)=0$ as $d\to\infty$, and
using $1-\exp(-x)=x+O(x^2)$, 
it follows easily that, for each $i$,
\begin{equation}\label{pdi}
 p_{d,i} = (\norm\tk+o(1))^d \quad \text{as}\quad d\to\infty.
\end{equation}
Let $\go=\go(n)=A\log n$, where $A$ is a constant, chosen large enough
that the estimates below hold.
As $\bpk$ is subcritical, $\rhogek(\kk)$ decays exponentially with $k$; see
\eqref{ineq-GW}. Thus $\rho_{\ge\go}(\kk)=o(n^{-2})$, say.
Let $d=(1\pm\eps)\log n/\lntki$, where $\eps$ is a small positive constant;
we shall consider both choices of sign below.

Let $T$ be a rooted tree where each vertex has a type from $\sss=\{1,2,\ldots,r\}$.
We shall say that a tree $T$ is \emph{relevant} if it has height at least $d$
and contains at most $\go$ vertices. Let $\pi(T)$ be the probability
that $\bpk$ is isomorphic to $T$, in the natural sense.

The sum of $\pi(T)$ over relevant $T$ is 
the probability that $\bpk$ survives at least $d$ generations and
contains at most $\go$ particles in total, which is
$(\norm\tk+o(1))^d+o(n^{-2})=n^{-1\mp\eps+o(1)}$.

Let $p(T)$ be the probability that the reduced component of a random vertex $v$
of $G_n$ is isomorphic to $T$ in the natural sense.
{}From the step-by-step exploration,
one can check that
\[
 p(T) = (1+o(1))^{|T|}\pi(T) = n^{o(1)}\pi(T)
\]
for any relevant $T$: the proof is similar to that of \refL{l:fincouple},
but one shows that at each step the conditional probability of finding
the right number $a$ of new neighbours of a particular type in the graph 
is within a factor $(1+o(1))^{a+1}$ of the corresponding Poisson
probability, as long as 
both $a$ and the number of previously uncovered vertices are $o(n)$.
Let $\sigma$ be the sum of $p(T)$ over relevant $T$. Then it follows that
$\sigma = n^{-1\mp \eps+o(1)}$.
In particular, $n\sigma=o(1)$ if we take the plus sign in $d=(1\pm\eps)\log n/\lntki$,
and $n\sigma\to\infty$ if we take the minus sign.

Using \eqref{ineq-GW} again, the expected number of vertices in $G_n$ with
more than $\go(n)$ vertices in their (reduced or unreduced) component is
$o(1)$ (in fact, $o(n^{-100})$), so \whp\ there no such vertices.
Taking the plus sign in $d$, the expected number of vertices in $G_n$
whose reduced component is a relevant tree is $n\sigma=o(1)$.
Together, these bounds show that \whp\ every vertex $v$ is such
that all $w$ with $d(v,w)<\infty$ have $d(v,w)<d$. Thus, \whp\,
$\diam(G_n)<d=(1+\eps)\log n/\lntki$. As $\eps>0$ was arbitrary,
this proves the upper bound in \refT{th_dbelow}.

For the lower bound we take the minus sign in $d$, so $n\sigma=n^{\eps+o(1)}\to\infty$,
and use the second moment method. The key point is that if $T$, $T'$
are relevant trees, and $p(T,T')$ is the probability that independently
chosen random vertices $v$, $w$
have \emph{vertex-disjoint} reduced neighbourhoods isomorphic to $T$, $T'$ respectively,
then
\begin{equation}\label{pTT}
 p(T,T') = (1+O(\go/n))^\go p(T) p(T') \sim p(T)p(T').
\end{equation}
This again follows from the step-by-step exploration, as finding one tree uses up at
most $\go$ vertices. (Note that we do not have $p(T)\sim \pi(T)$,
because all we know about $n_i$ is that $n_i\sim n\mui$.)

Let $X_v$ be the indicator function of the event that
the reduced neighbourhood of $v$ is a relevant tree,
and let $N=\sum_{v=1}^n X_v$,
so $\E(N)=n\sigma$.
Expanding $\E(N^2)=\sum_v\sum_w\E(X_vX_w)$, the contribution from pairs $v$, $w$ in the same
component is at most $\sum_v\E(X_v\omega)=n\sigma \go$: if the component containing
$v$ is relevant, then by definition it contains at most $\go$ vertices $w$.
Using \eqref{pTT} above,
it follows that
\[
 \E(N^2) = (1+o(1))n^2\sigma^2 + O(n\sigma\go).
\]
Now $\E(N)=n\sigma=n^{\eps+o(1)}$, which
is much larger than $\go$, so $\E(N^2)\sim \E(N)^2$,
and \whp\ $N> 0$. So \whp\ there is a vertex $v$ whose
neighbourhood is a relevant tree, and thus includes a vertex
$w$ at distance at least $d=(1-\eps)\log n/\lntki$,
completing the proof of \refT{th_dbelow}.
\end{proof}

\begin{remark}\label{ReinT}
Recall that $p(T)$ was defined as the probability
that the {\em reduced} neighbourhoods of a random
vertex $v$ are isomorphic to $T$, allowing the possibility
that there are some edges within each $\Ga_d(v)$.
The reason was that, in proving
the upper bound on the diameter, we must
rule out components of large diameter that contain
cycles, as well as components that are trees.
If we redefine $p(T)$ to exclude edges within each $\Ga_d(v)$,
then $p(T)$ changes by a factor $(1+O(|T|/n))^{|T|}$. 
As we only ever consider $T$ with $|T|=O(\log n)$, this factor
is $1+o(1)$, and all our estimates go through. In particular,
\whp\ $G_n$ contains a tree component of diameter
at least $(1-\eps)\log n/\lntki$.
\end{remark}

\begin{remark}
One might expect \refT{th_dbelow} to generalize immediately from
$\sss$ finite to (at least) the case $\kk$ bounded. However,
$\diam(G)$ does {\em not} always decrease when an edge is added to
$G$, as the new edge might join two components. Thus
one cannot just sandwich $\gnkx$ between finite-type graphs
and apply \refT{th_dbelow}.  In the unbounded case, the construction
described in \refR{Rlargedist} shows that for any
$\go(n)=o(n)$ one can construct a graph $\gnkx$ with $\kk$
supercritical, such that $\diam(\gnkx)\ge \go(n)$ \whp.  Modifying
the construction by starting with a subcritical Erd\H os-R\'enyi graph
gives an example with $\kk$ subcritical.
\end{remark}

We now turn to the supercritical case of \refT{th_diam}, restated 
as \refT{th_dabove} below.
Recall that,
given a supercritical kernel $\kk$ on a ground space $(\sss,\mu)$, there
is a `dual' kernel $\kkd$ on a ground space $(\sss,\mud)$, defined as follows:
as a function on $\sss\times\sss$, $\kkd=(1-\rho(\kk))\kk$,
while $\dd\mud(x)=(1-\rho(k;x))/(1-\rho(\kk))\dd\mu(x)$.
In particular, in the finite-type case,
$\mud(\{i\})=(1-\rho(k;i))/(1-\rho(\kk))\mui$. 
Note that we have chosen to renormalize
the dual kernel defined in \refD{Ddual} so that $\mud$ is a
probability measure. As discussed after \refD{Ddual},
this makes essentially no difference;
however, it allows us to speak of the branching process $\bpx{\kkd}$ started
with a particle whose type is chosen according to $\mud$.
(In the remark after \refD{Ddual} we wrote $\kkd'$ and $\mud'$ for the
renormalized dual kernel and associated measure; here we write $\kkd$
and $\mud$ 
for notational convenience.)
When we write $\norm{T_\kkd}$, we mean the norm of $T_\kkd$ defined
with respect 
to $(\sss,\mud)$.

\begin{theorem}\label{th_dabove}
Let $\kk$ be an irreducible kernel on a (generalized) vertex space 
$\vxs=(\sss,\mu,(\xs))$, with $\sss=\{1,2,\ldots,r\}$ finite and
$\mui>0$ for each~$i$. 
If $\norm\tk>1$, then
\[
 \frac{\diam(G_n)}{\log n} \pto   \frac{2}{\lntkdi} +\frac{1}{\lntk},
\]
where $G_n=\gnkx$.
\end{theorem}

The relevance of the dual kernel is that it describes components other than
the giant component. In particular, it follows from \refT{Tdual} and
\refT{th_dbelow} that
the diameter of the largest `small' component of $G_n$ will be
$(1+o(1))\log n/\lntkdi$. As we shall see, the same quantity will give
the height of the tallest tree attached to the two-core. The diameter
will be given by two such trees attached to vertices of the two-core
at typical distance, $(1+o(1))\log n/\lntk$.

The idea of the proof is as follows: instead of considering the event
that the neighbourhoods of a vertex $v$ form a tree of height at least
$d$, we consider the event 
that the neighbourhoods are \emph{thin} for $d$ generations, meaning
that each generation 
has size at most $\go$, with $\go=A\log n$ as before. For the upper
bound, we will 
show that \whp\ no vertex has neighbourhoods that are thin for more than
$d=(1+\eps)\log n/\log \alpha^{-1}$ generations, where
$\alpha=\norm{T_\kkd}$. For the lower 
bound, we will find two trees of height roughly $d$ attached to
typical vertices 
of the two-core.

The reason for considering thin neighbourhoods is that, once $\Ga_t(v)$
is larger than $A\log n$ for some $t$, the neighbourhoods $\Ga_s(v)$, $s\ge t$,
grow reasonably rapidly.

From now on, we assume that $\kk$ is an irreducible kernel on a finite
ground 
space $(\sss,\mu)$, with $\mui>0$ for each $i$, and that 
the number $n_i$ of vertices of each type $i$ is deterministic, with
$n_i/n\to \mui$. 
As before, it suffices to prove \refT{th_dabove} under these assumptions.
Let 
\[
 t(v)\=\min\{r:|\Ga_r(v)|\ge \go\}
\]
denote the index of the first {\em thick}  neighbourhood of a vertex
$v$, when there is one.

\begin{lemma}\label{Ltt}
For any $\eps>0$,
\whp\ the graph $G_n$ does not contain two vertices $v$, $w$ with the
properties 
that $t(v)$, $t(w)$ are defined, $t(v)$, $t(w)\le n^{1/2}$,
and $d(v,w)\ge t(v)+t(w)+(1+\eps)\log n/\log\norm\tk$.
\end{lemma}

\begin{proof}
We show that the expected number of pairs is $o(1)$, by showing that
the probability 
that a random pair $v$, $w$ has the properties is $o(n^{-2})$.
Explore the neighbourhoods of $v$ and $w$ simultaneously, stopping at the first
thick neighbourhood of each, if there is one. Suppose, as we may, that
$t(v)$ and $t(w)$ are defined and at most $n^{1/2}$. Then with very high probability
 we have seen $o(n)$
vertices (the neighbourhoods can't have grown too much in the last step). If the neighbourhoods
have already joined, we are happy.
Otherwise, continue exploring. Simple Chernoff bounds show that for any $\eta>0$,
if $A$ is chosen large enough, with probability $1-o(n^{-100})$
the number of vertices of each type found at each subsequent step is within a factor
$1\pm \eta$ of its expectation. It follows that the neighbourhoods grow
by a factor of $(1\pm 2\eta)\norm\tk$ at each step, after a few steps to allow
the distribution of types 
to converge to the relevant eigenvector of $\kk$. Once both neighbourhoods
reach size $n^{1/2+\eta}$, they join at the next step with
very high probability.
\end{proof}

For the rest of this section,
let $\alpha=\ntkd$ denote the norm of the dual kernel $\kkd$.
As $\sss$ is finite, $\alpha<1$ by \refT{Tdualsub}.
Recall, from the discussion before \refL{Ldual},
that $\bpx{\kkd}$ has the same distribution as $\bpk$ conditioned
on extinction.
Let $t_{d,n}$ be the probability that $\bpk$ stays alive but thin for $d$ generations:
\[
 t_{d,n} \= \P\bigpar{1\le |X_t| \le \omega\: :\:1\le t\le d},
\]
where $X_t$ is generation $t$ of $\bpk$. Note that $t_{d,n}$ depends on $n$, via the definition
of $\omega$.

\begin{lemma}\label{Ltd}
For any $\eta>0$, if $n$ is large enough, then
\begin{equation}\label{tdn}
 (\alpha-\eta)^d \le t_{d,n} \le (\alpha+\eta)^d
\end{equation}
holds for all $d$ in the range $\frac{1}{2}\log n/\log\alpha^{-1} \le d \le 2\log n/\log\alpha^{-1}$.
\end{lemma}

\begin{proof}
We start with the lower bound.

Let $p_{r,i}'$ be the probability that the branching process $\bpkd(i)$,
started with a particle of type $i$,  survives for at least
$r$ generations. By \eqref{pdi}, we have $p_{r,i}'=(\alpha+o(1))^r$ as $r\to\infty$.
Let $p_{r,j,i}''$ be the probability that generation $r$ of $\bpkd(j)$
consists of a single particle of type $i$.
As the branching process is subcritical, one can check
that $p_{r,i,i}''=\Theta(p_{r,i}')$ for
$r$ even. (The restriction $r$ even is only needed if the kernel $\kkd$
is bipartite, in the sense of \eqref{bipdef_moved}.)
We shall need only the much weaker statement that
$p_{r,i,i}''=(\alpha+o(1))^r$ as $r\to\infty$ with $r$ even;
this can be proved along the same lines as \eqref{pdi}:
let $N_1$ be the number of particles $x$ in the first generation
of $\bpkd(j)$ with the property that the descendants of $x$
in generation $d+1$ consist of a single particle of type $i$.
Let $N_2$ be the number of particles $x$ in the first generation
of $\bpkd(j)$ with more than one descendant in generation $d+1$ of
type $i$, or any descendants of types other than $i$.
Note that $N_1$ and $N_2$ have independent Poisson distributions, with
\[
 \E(N_1) = \sum_k \kkd(j,k)\mudk p_{d,k,i}'', 
\]
and
\[
 \E(N_2) = \sum_k \kkd(j,k)\mudk (p_{d,k}'-p_{d,k,i}'').
\]
Since $p_{d,k,i}''\le p_{d,k}'\le \rho_{\ge d}(\kkd;k)\to 0$,
we have $\E(N_1)$, $\E(N_2)\to 0$.
By definition,
\[
 p_{d+1,j,i}'' = \P(N_1=1, N_2=0) = \E(N_1)\exp(-\E(N_1)-\E(N_2)) \sim \E(N_1).
\]
In other words,
\[
 p_{d+1,j,i}'' = (1+o(1)) \sum_k \kkd(j,k)\mudk p_{d,k,i}''.
\]
Recalling that $\alpha$ is the norm of $T_\kkd$ defined with respect to $\mud$,
it follows that
$p_{r,i,i}''=(\alpha+o(1))^r$ as $r\to\infty$ with $r$ even.

The probability that $\bpkd$ has any thick generation at all
is at most the expected total size of $\bpkd$ divided by $\omega$, which is
$O(1/\omega)=O(1/\log n)$. 
Taking $r\to\infty$ slowly enough, for example $r=2\ceil{\log\log\log n}$,
this probability is much smaller than $p_{r,i,i}''$.
Hence, with probability $p_{r,i}'''=(\alpha+o(1))^r$ the branching process
$\bpkd(i)$ remains thin for $r$ generations, and the $r$th generation
is a single particle of type $i$. Restarting, we see that for any $d$
we have $t_{d,n}\ge \mui(p_{r,i}''')^{\ceil{d/r}}$,
and the lower bound in \eqref{tdn} follows.

Let $c>0$ be a (small) constant. Simple Chernoff bounds show
that, given that a certain generation $t$ of $\bpk$ has size at least $c\log n$
(and given the numbers of particles of each type),
generation $t+c_1$ has size at least $\omega$ with probability at least
$1-n^{-c_2}$, for some constants $c_1$, $c_2$ depending on $c$ and $\kk$.
Hence, if $L=L(\kk,c)$ is chosen large enough,
the probability that $\bpk$ stays thin for $d$ generations and has size at least
$c\log n$ for the last $L$ of these is at most $n^{-100}$: given that generation
$d-L$ has size at least $c\log n$, the probability that generation $d-L+c_1$ is still
thin is at most $n^{-c_2}$. If this generation is thin but also has size
larger than $c\log n$, generation $d-L+2c_1$ is unlikely to be thin, and so on.

Hence, $t_{d,n}$ is within $n^{-100}$ of
the probability that $\bpk$ stays thin for $d$ generations and one of the
last $L$ of the first $d$ generations has size at most $c\log n$.
Let $q_r=q_{r,c}$ be the probability of the event that $\bpk$ survives for
at least $r$ generations, that the first $r$ generations are thin, and that generation
$r$ has size at most $c\log n$. We have shown that
\begin{equation}\label{tq}
 t_{d,n} \le \sum_{d-L\le r\le d} q_r + n^{-100}.
\end{equation}
As a generation of size $c\log n$ has probability $n^{-O(c)}$ of dying out immediately,
the probability that $\bpk$ survives for {\em exactly} $r$ generations
is at least $q_rn^{-O(c)}$.
Hence, $p_{r,i}'$, the probability that $\bpx{\kkd}$, which is just $\bpk$
conditioned on dying out, survives for $r$ generations
is also at least $q_rn^{-O(c)}$. Using $p_{r,i}'=(\alpha+o(1))^r$,
it follows that $q_r\le n^{O(c)}(\alpha+o(1))^r$.
In particular, choosing $c$ small enough (depending on $\kk$ and $\eta$),
for $r\ge \frac{1}{3}\log n/\log \alpha^{-1}$, say, we have
$q_r\le (\alpha+\eta/2)^r$.
Using \eqref{tq}, it follows that
\[
 t_{d,n} \le L(\alpha+\eta/2)^{d-L} +n^{-100} \le
 (\alpha+\eta)^d
\]
for $n$ large enough and $d$ in the range considered in \eqref{tdn},
completing the proof of \refL{Ltd}.
\end{proof}

Let $T$ be a rooted tree in which each vertex has a type from $\sss$.
Let $h(T)$ be the {\em height} of $T$, i.e., the maximal distance
of a vertex form the root.
We shall say that $\bpk$ is {\em consistent with} $T$ if the first
$h(T)$ generations of $\bpk$ are isomorphic to $T$,
and write $\pi'(T)$ for the probability of this event.
Analogously, we say that the reduced neighbourhoods of a
vertex $v$ of $G_n$ are {\em consistent with} $T$
if $T_{h(T)}(v)$, the tree formed by the first $h(T)$ such
neighbourhoods, is isomorphic to $T$.
We write $p'(T)$ for the probability that the reduced neighbourhoods
of a random vertex of $G_n$ are consistent with $T$.
The definitions are almost the same as those of $\pi(T)$, $p(T)$, except
that we do not care what happens after the first $h(T)$ generations.
(We could have used these definitions in the subcritical case -- there
it was not essential that we explored the whole component.)

We are now ready to prove \refT{th_dabove}.
\begin{proof}[Proof of \refT{th_dabove}.]
As before, we may assume that $\vxs$ is a vertex space, and that the number 
$n_i$ of vertices of type $i$ is deterministic.

Let $\eps>0$ be fixed, and let $d=(1\pm\eps)\log n/\log \alpha^{-1}$,
where $\alpha=\norm{T_\kkd}$, as above.
As before, let $\omega=A\log n$, where $A$ is a constant chosen large enough
for our bounds to hold. Having chosen $A$, let $A'$ be another constant
chosen large enough that the bounds below hold.

Let $t_{h,n}^+$ be the probability that the branching process $\bpk$ is
alive and thin for $h$ generations (as in the definition of
$t_{d,n}$), but that 
these first $h$ generations contain more than $A'\log n$ particles.
The proof of Lemma~\ref{Ltd} shows that if $A'$ is chosen large enough,
then $t_{h,n}^+\le n^{-99}$ for any $h$. Indeed,
the argument leading to \eqref{tq} gives a corresponding bound
on $t_{h,n}^+$ with $q_r$ replaced by the probability $q_r^+$ of
an appropriate event, defined as $q_r$, but with the extra condition
that there are
at least $A'\log n/2$ particles in the first $r$ generations. (We take
$A'$ large enough that $L\omega\le A'\log n/2$.)
Then, as before, $q_r^+$ can be bounded by $n^{O(1)}$ times
the probability that $\bpx{\kkd}$ survives for exactly $r$ generations
and contains at least $A'\log n/2$ particles.
Using only the exponential decay of $\rhogek(\kkd)$, this is at most
$n^{-100}$ if 
$A'$ is large enough, for any $r$.

Let us say that a tree $T$ is {\em relevant} if $T$ has height $d$,
is thin, and contains at most $A'\log n$ vertices.
As $T$ has $O(\log n)$ vertices,  we have
$p'(T)=n^{o(1)}\pi'(T)$, as before.
The sum of $\pi'(T)$ over all relevant $T$ is exactly
\[
 t_{d,n}-t_{d,n}^+ = (\alpha+o(1))^d - O(n^{-99}) = n^{-1\mp\eps+o(1)},
\]
by \refL{Ltd} and our bound on $t_{d,n}^+$ above. As before, the sign above
is the opposite of the sign we choose in $d=(1\pm\eps)\log n/\log \alpha^{-1}$.
It follows that the sum $\sigma$ of $p'(T)$
over relevant trees is $n^{-1\mp\eps+o(1)}$.

To prove the upper bound in \refT{th_dabove}, let $d=(1+\eps)\log n/\log \alpha^{-1}$.
Then we have $n\sigma=o(1)$, so \whp\ no vertex of $G_n$ has neighbourhoods
consistent with a relevant tree. It is easy to check that \whp\ no vertex
has neighbourhoods consistent with a thin tree $T'$ of height $d$ that
is not relevant (because it contains more than $A'\log n$ vertices): any such
$T'$ contains a subtree $T''$ given by the first $h$ generations
of $T'$ for some $h$, such that $T''$ is thin and has between $A'\log n$
and $A'\log n+\omega=O(\log n)$ vertices. But then $p'(T'')=n^{o(1)}\pi'(T'')$,
and $\sum_{T''} \pi'(T'')$ is at most $\sum_{h\le n}t_{h,n}^+\le n^{-98}$.
It follows that \whp\ no vertex of $v$ has neighbourhoods consistent
with any thin tree of height $d$. Thus, for every $v$, either
$\Ga_d(v)$ is empty, or $t(v)\le d$. Applying \refL{Ltt}, the upper
bound on $\diam(G_n)$ claimed in \refT{th_dabove} follows.

\medskip
For the lower bound, we aim to find many tall thin trees attached to the two-core.
As in \refR{ReinT}, for this part of the proof we modify our notion of the
consistency of the neighbourhoods
of a vertex $v$ of $G_n$ with a tree $T$ of height $d$, by disallowing
edges within each $\Ga_t(v)$, $t\le d$. We redefine $p'(T)$ correspondingly;
as in \refR{ReinT}, this changes $p'(T)$ by a factor $1+o(1)$ for trees
of the size we consider, so our estimate $p'(T)=n^{o(1)}\pi'(T)$ goes through.

A \emph{good} tree $T$ will be a relevant tree with
height $d=(1-\eps)\log n/\log\alpha^{-1}$ in which generation $d$
consists of a single vertex $a$. 
Changing $\eps$ slightly, we shall assume that $d$ is a multiple of
the quantity 
$r$ considered in the proof of \refL{Ltd}.
Then this proof shows that with probability at least
$(\alpha+o(1))^d=n^{-1+\eps+o(1)}$ 
the first $d$ generations of $\bpk$ are thin and generation $d$
consists of a single 
particle. With our bound on $t_{d,n}^+$, it follows that $\sum_T
\pi'(T)=n^{-1+\eps+o(1)}$, 
where the sum is over good trees. As $p'(T)=n^{o(1)}\pi'(T)$ for each good $T$,
we have $\sum_T p'(T)=n^{-1+\eps+o(1)}$.
As in the subcritical case, the second moment method gives us many good trees in the
graph, but this is not enough -- they might not be attached to the two-core.

For $v\in V(G)$ and $T$ a good tree, let $E_2(v,T)$ be the event that
the following all hold, where $d_1=C\log\log n$, and $C$, $C_1$ are constants to be chosen below:
the first $d$ neighbourhoods of $v$ form the tree $T$,
the single $a\in \Ga_d(v)$ has two neighbours $b_1$, $b_2$ in $\Ga_{d+1}(v)$
each of which has at least $\omega$ `descendants' in $\Ga_{d+d_1}(v)$,
and $|\Gale{d+d_1}(v)|\le (\log n)^{C_1}$. We claim that if the constants $C$
and $C_1$ are chosen large enough, 
and $v$ is a random vertex of $G_n$, then $\P(E_2(v,T))=\Omega(p'(T))$.
To see this, note that $p'(T)$ is exactly the probability that the first condition
is satisfied. Conditional on this happening, bounding the neighbourhood
exploration below by a supercritical branching process shows
that the existence of $b_1$, $b_2$ with the required properties has probability bounded
away from zero.
Finally, given that $|\Ga_d(v)|=1$, the expected size
of the next $d_1$ generations is at most $(2\sup \kk)^{d_1}=(\log n)^{O(1)}$,
so the claim follows.

Let $E_2(v)$ be the event that $E_2(v,T)$ holds for some good $T$. 
As $\sum_T p'(T)=n^{-1+\eps+o(1)}$, we have
\begin{equation}\label{sgs}
 \P(E_2(v)) = \sum_T \P(E_2(v,T)) = n^{-1+\eps+o(1)},
\end{equation}
where the sum is over good $T$.

We would like to show that \whp\ there are two (in fact, many) vertices
$v$, $w$ for which $E_2(v)$, $E_2(w)$ hold. We could use the second moment
method, but as we shall need the relevant neighbourhoods of $v$ and $w$
to be disjoint, it turns out to be easier to test vertices one by one.

Whether $E_2(v,T)$ holds can be determined by exploring the neighbourhoods
of $v$, stopping when at most $M=(\log n)^{C_1}$ vertices have been uncovered.
Let us construct a sequence of $t$ tests, $t=n/M^3$, as follows.
Each test starts from a vertex $v_i$, where the $v_i$ are chosen independently
and uniformly at random from $V(G_n)$. In the $i$th test, we explore
the neighbourhoods of $v_i$, uncovering at most $M$ vertices,
and attempting to verify that $E_2(v_i)$ holds. We abort the attempt
if we reach a vertex uncovered in a previous attempt.
As at most $tM=n/M^2$ vertices have previously been uncovered,
for each vertex we reach, the probability that it was previously
uncovered is at most $O(M^{-2})$. As the $i$th test involves
examining at most $M$ vertices, conditional on everything so far,
the $i$th test succeeds with probability
\[
 (1-O(M^{-2}))^M\P(E_2(v))\sim \P(E_2(v)) = n^{-1+\eps+o(1)}.
\]
The number of tests that succeed dominates a binomial random variable
with mean $tn^{-1+\eps+o(1)}=n^{\eps+o(1)}\to\infty$,
so \whp\ at least two tests succeed.
Hence, \whp\ there are vertices $v$, $w$ in $G_n$ for which $E_2(v)$, $E_2(w)$
hold, with the relevant neighbourhoods disjoint.

Now \whp\ $G_n$ has the property that whenever $E_2(v)$ holds, the corresponding
$a\in \Ga_d(v)$ is in the two-core of $G_n$. The argument is as for \refL{Ltt}:
we may continue expanding the large neighbourhoods of $b_1$, $b_2$ until they meet.
Hence, \whp\ $G_n$ contains two vertices $v_1$, $v_2$ belonging
to separate trees $T_1$, $T_2$ of height $d$, each attached to the two-core by the
single vertex at distance $d$ from the root.
We shall need only this last fact, basic properties
of the model, and \refT{Tdist}.

Recall that $G_n$ is a graph on labelled vertices $\{1,2,\ldots,n\}$,
each of which has a type in $\sss$. Given $G_n$ and the vertex types,
let us separate $G_n$ into the two-core $\Gtc$, a list of trees $T_i$
each attached to the two-core at some attachment vertex $a_i$,
and the rest of $G_n$. So far, we remember the label of each vertex. Now, however,
let us forget the labels of the attachment vertices $a_i$ of $T_i$, while
remembering their types. To reconstruct $G_n$, we should identify
each attachment vertex with a vertex of $\Gtc$ of the same type.
Moreover, we may pick these vertices independently and uniformly
at random from the allowed vertices of $\Gtc$; this is because
all possible (labelled) graphs formed in this way have the same
number of vertices of each type, and the same number of edges between
vertices of each pair of types, and hence the same probability in our model.

We have shown above that \whp\ our list $T_i$ contains two trees, say $T_1$ and $T_2$,
each of which has a vertex $v_i$ at distance $d$ from the corresponding attachment
vertex $a_i$. By \eqref{tclarge}, \whp\ $\Gtc$ contains $\Theta(n)$ vertices.
Hence, by
\refT{Tdist}\ref{Tdistu},
\whp\ almost all pairs of vertices of $\Gtc$ are at distance
at least $d'=(1-\eps)\log n/\lntk$. Hence, \whp\ the vertices
of $\Gtc$ at which we reattach the $a_i$ are at distance at least $d'$.
Thus,
\whp, $\diam(G_n)\ge d(v_1,v_2)\ge 2d+d'$, completing the proof of \refT{th_dabove}.
\end{proof}

\section{The phase transition}\label{Stransitionpf}

Our main aim in this section is to prove \refT{T5}, which claims that if 
a kernel $\kk$
on a ground space $(\sss,\mu)$ 
is irreducible and
satisfies \eqref{t5a}, i.e.,
\begin{equation}\label{t5a'}
  \sup_x \int_\sss \kk(x,y)^2\dd\mu(y) <\infty,
\end{equation}
then the function $c\mapsto
\rho(c)\=\rho(c\kk)$
is analytic except at $c=c_0\=\norm{\tk}^{-1}$, that $\tk$ has an
eigenfunction $\psi$
of eigenvalue $\norm{\tk}$, and that every such eigenfunction is
bounded and satisfies
\eqref{t5}.

\begin{proof}[Proof of \refT{T5}]
Note that our
assumption \eqref{t5a'} on $\kappa$  implies that $\tk$
is a Hilbert--Schmidt
operator, and thus compact in $L^2$; see \refL{Lcomp}.
It further implies, by the \CS{} inequality,
that $\tk$ is bounded
$L^2\to L^\infty$, and that
\eqref{bb1} holds for all $x$.

\pfitem{i}
It is trivial that the function $c \mapsto \rho(c)=\rho (c\kappa)$
  is analytic
  for $c<c_0$, so we shall assume that $c>c_0$.
We shall show that we can extend this function
to a suitable neighbourhood of $c$ in the complex plane,
and that this
extension is (complex) analytic; this implies that $\rho$ is a real
analytic function at $c$.
Actually,
we will show that there is an analytic map $z\mapsto\extrho_z$ into
$L^2(\mu)$, defined in a neighbourhood of $c$,
such that $\extrho_z=\rho_{z\kk}$ when $z$ is real. Here, as before,
$\rho_{z\kk}$
is the function defined by $\rho_{z\kk}(x)=\rho(z\kk;x)$, the survival
probability of the branching process $\bpxx{z\kk}$, which starts with a
particle of type
$x$.
We may then take
$\extrho(z):=\int\extrho_z(x)\dd\mu(x)$ as the extension of $\rho$.

To show that the claimed extension exists,
we will use the implicit function theorem for complex analytic functions
in the Banach space $L^2$.
Of course, in this proof we use the complex version of $L^2$.
Recall that a function $f$ mapping
an open subset $U$ of a complex Banach space $E$ into another Banach
space $F$ is {\em analytic} if
and only if it is differentiable (in the Fr\'echet sense) at every
point in $U$; the derivative $f'(x)$ at a point $x\in U$ then is a
continuous linear operator $E\to F$, see,
\eg, Herv\'e \cite[Section 3.1]{Herve}.
For background on differentiable functions in (real or complex)
Banach spaces, see Cartan \cite{Cartan}; in particular, as a special case of
\cite[Theor\`eme 4.7.1]{Cartan} (which holds in both the real and
complex cases), we have the following.

\begin{lemma}[The implicit function theorem]
  \label{Limplicit}
Let $B$ be a complex Banach space and let $f_0\in B$, $z_0\in\bbC$.
Let $\Psi:\bbC\times B\to B$ be an analytic function,
and denote by
$D_2\Psi(z,f)$ the partial derivative of $\Psi$ with
respect to the second variable,
\ie{}, let $D_2\Psi(z,f)$ be the derivative of $f\mapsto \Psi(z,f)$. Suppose that
$\Psi(z_0,f_0)=0$ and that $D_2\Psi(z_0,f_0)$ is invertible.
Then there exists a neighbourhood $U$ of $z_0$ and an analytic function
$z\mapsto f(z)$ defined in $U$ such that $f(z_0)=f_0$ and
$\Psi(z,f(z))=0$, $z\in U$.
\nopf
\end{lemma}

For convenience, note that by
replacing $\kk$ by $c\kk$, we may assume that $c=1$ (and thus
$c_0<1$).
We then apply \refL{Limplicit} with $B=L^2(\mu)$, $z_0=c=1$, $f_0=\rhokk$
and $\Psi(z,f)=\Phik(zf)-f$, where $\Phik(f)=1-e^{-\tk f}$ as above.
Since $\tk$ is a bounded linear map
$L^2\to L^\infty$,
and $g\mapsto e^g$ is analytic $L^\infty \to L^\infty$,
$\Phi$ is analytic $L^2\to L^\infty\subseteq L^2$ and
thus $\Psi$ is analytic.
It is also easily seen that $D_2 \Psi(z,f)=z\Phik'(zf)-I$.
It remains to show that that the partial derivative at $(1,\rhokk)$ is
invertible; we state this as another lemma.

\begin{lemma}
  \label{Linv}
Assume that $\kk$ is irreducible,
\eqref{t5a'} holds and $c_0=\norm{\tk}\qi<1$.
Let $\tmu$ be the measure $\dd\tmu=e^{\tk\rhokk}\dd\mu$ on $\sss$. Then
$\norm{\Phik'(\rhokk)}_{L^2(\tmu)} <1$, and hence
$D_2\Psi(1,\rhokk)=\Phik'(\rhokk)-I$ is invertible in $L^2(\tmu)$ and in
$L^2(\mu)$.
\end{lemma}

\begin{proof}
Since $\tk$ maps
$L^2$ into $L^\infty$, we have $\tk\rhokk\in L^\infty$. Hence
$L^2(\tmu)= L^2(\mu)$
with equivalent norms.

We have
\begin{equation}
  \label{tt5a}
\Phik'(f)(g)=e^{-\tk f} \tk g.
\end{equation}
Hence, for any $g,h\in L^2(\tmu)=L^2(\mu)$, writing $\bar h$ for the
complex conjugate of $h$, we see that
\begin{equation*}
  \innprod{ \Phik'(\rhokk)(g),h}_{L^2(\tmu)}
=\int_\sss e^{-\tk \rhokk} (\tk g) \bar h e^{\tk \rhokk } \dd\mu
=\int_\sss (\tk g) \bar h  \dd\mu
\end{equation*}
is a hermitian form in $g$ and $h$, because $\tk$ is a symmetric
operator in $L^2(\mu)$.
Hence $\Phik'(\rhokk)$ is a symmetric operator in $L^2(\tmu)$.
Furthermore, as remarked above, $\tk$ is compact in $L^2(\mu)$ and
thus also in
$L^2(\tmu)$.

If we had $\norm{\Phik'(\rhokk)}_{L^2(\tmu)} \ge1$, there
would thus be an eigenfunction $g$ with an eigenvalue $\gl$ with
$|\gl|\ge1$. Let $h\=|g|\ge0$, so $\int h>0$.
Since $\Phik'(\rhokk)$ has a non-negative kernel by \eqref{tt5a}, we obtain
\begin{equation}\label{tt5b}
  h \le |\gl g|
=|\Phik'(\rhokk)g|
\le\Phik'(\rhokk) h
= e^{-\tk\rhokk} \tk h
.
\end{equation}
Moreover, $\rhokk>0$ \aex{} and thus, by part (iii) of \refL{L1b},
$\tk\rhokk>0$ \aex{} and
\begin{equation}\label{tt5c}
  \tk\rhokk < e^{\tk\rhokk} -1 = e^{\tk\rhokk}\Phik \rhokk
= e^{\tk\rhokk} \rhokk
\qquad  \text{a.e.}
\end{equation}
Multiplying \eqref{tt5b} and \eqref{tt5c} and integrating, recalling
that $\int h>0$, we obtain
\begin{equation*}
  \int h \tk\rhokk\dd\mu < \int \rhokk \tk h\dd\mu,
\end{equation*}
which contradicts the symmetry of $\tk$.

This contradiction shows that
$\norm{\Phik'(\rhokk)}_{L^2(\tmu)} <1$. Hence
$I-\Phik'(\rhokk)$ is invertible in $L^2(\tmu)= L^2(\mu)$, completing
the proof of \refL{Linv}.
\end{proof}

Continuing the proof of \refT{T5},
we can now apply the implicit function theorem (\refL{Limplicit}) to
conclude the existence of an analytic function
$z\mapsto\extrho_z\in L^2(\mu)$  with
$\Phik(z\extrho_z)=\extrho_z$, defined in a complex neighbourhood $U$ of $1$.
We may further (by continuity) assume that $U$ is so small that
$\normll{\extrho_z}>0$.
For real
$z\in U$, we thus have $\Phix{z\kk}(\extrho_z)=\Phik(z\extrho_z)=\extrho_z$
(in $L^2$, \ie{}, a.e.),
so $\extrho_z=\rho_{z\kk}$ \aex{} by \refT{T:GW} and \refR{Rae},
which completes the proof of (i).

\pfitem{ii}
This time we scale so that $\norm{\tk}=1$, and thus $c_0=1$, and write
$\rhoe$ for $\rho_{(1+\eps)\kk}$; we assume
below that $0\le \eps<1$. Thus, by \refT{T:GW},
\begin{equation*} 
1-e^{-(1+\eps)\tk\rhoe} = \Phix{(1+\eps)\kk}(\rhoe) = \rhoe,
\end{equation*}
i.e.,
\begin{equation}
  \label{xb}
(1+\eps)\tk\rhoe = -\ln(1-\rhoe)=\rhoe+R(\rhoe)
\end{equation}
where
\begin{equation}
  \label{xc}
R(f)\=\frac{f^2}2 + \frac{f^3}3+\dots
\end{equation}

By \refT{T:GW},
$\rhoe>0$ \aex{} when $\eps>0$, but $\rhoo=0$.
By \refL{Lcomp}, there exists an eigenfunction $\psi\in L^2$ with
$\psi=\tk\psi$, which now implies $\psi\in L^\infty$. Furthermore,
$\psi$ is determined up to a constant factor, so the coefficient
$\int\psi\int\psi^2/\int\psi^3$ does not depend on the choice of
$\psi$.
We will for convenience assume that $\psi$ is chosen with $\psi\ge0$
and $\int\psi^2\dd\mu=\normll{\psi}^2=1$.

The operator $\tk$ maps the subspace $\psi^\perp\subset L^2(\mu)$ into itself;
let $\tkq$ denote the restriction of $\tk$ to this subspace.
Then $1$ is not an eigenvalue of $\tkq$, and thus (since
$\tk$ is compact), $1$ does not belong to the spectrum of $\tkq$, \ie{},
$I-\tkq$ is invertible. By continuity, $I-(1+\eps)\tkq$ is also
invertible for small $\eps$, and there exists $\gd>0$ and $C<\infty$
such that $\norm{\bigpar{I-(1+\eps)\tkq}\qi}_{\psi^\perp} \le C$ for
$0\le \eps<\gd$, \ie{},
\begin{equation}
  \label{xba}
\normll{f} \le C \normll{\bigpar{I-(1+\eps)\tk}f},
\qquad 0\le\eps<\gd,
\; f\in\psi^\perp.
\end{equation}

\refT{TappB} implies that $\rhoe\downto\rhoo$ \aex{} as $\eps\downto0$,
and thus,
by dominated convergence,
\begin{equation}
  \label{xd}
\normll{\rhoe}\to0.
\end{equation}
We also have, by \eqref{xb}, $\rhoe\le(1+\eps)\tk\rhoe$ and thus,
as $\tk$ is bounded from $L^2$ to $L^\infty$,
\begin{equation}
  \label{xe}
\normoo{\rhoe} \le (1+\eps) \normoo{\tk\rhoe}\le C_1\normll{\rhoe}.
\end{equation}
In particular, by \eqref{xd}, $\normoo{\rhoe}\to0$ as $\eps\to0$.

Assume that $\eps>0$ is small enough to ensure that $\normoo{\rhoe}<1/2$. Then, by
\eqref{xc}, $R(\rhoe)\le\rhoe^2$, and thus,
using \eqref{xe},
\begin{equation}
  \label{xf}
\normll{R(\rhoe)}
\le \normoo{\rhoe^2}
=\normoo{\rhoe}^2
\le C_2 \normll{\rhoe}^2.
\end{equation}

Let $Q$ be the orthogonal projection onto $\psi^\perp$ and let
$\rhoex\=Q\rhoe$.
We thus have the orthogonal decomposition
\begin{equation}
  \label{xg}
\rhoe=a_\eps\psi+\rhoex,
\end{equation}
where
\begin{equation}
  \label{xh}
a_\eps=\innprod{\rhoe,\psi}=\int_\sss\rhoe\psi\dd\mu.
\end{equation}
Hence $0\le a_\eps \le\normll{\rhoe}$.

Applying the projection $Q$ to \eqref{xb} we find,
since $Q\tk =\tk Q$,
\begin{equation*}
 (1+\eps)\tk \rhoex = \rhoex + Q(R(\rhoe))
\end{equation*}
and thus by \eqref{xba} and \eqref{xf}, for $\eps<\gd$,
\begin{equation}
  \label{xfa}
\normll{\rhoex}
\le C \normll{\bigpar{I-(1+\eps)\tk}\rhoex}
=C\normll{ Q(R(\rhoe))}
\le C_3 \normll{\rhoe}^2.
\end{equation}
Consequently, by \eqref{xg} and \eqref{xd}, as $\eps\to0$,
\begin{equation}
  \label{xi}
  \begin{split}
a_\eps
&
=\normll{a_\eps\psi}
=\normll{\rhoe-\rhoex}
\\&
=\normll{\rhoe}+O(\normll{\rhoex})
=\normll{\rhoe}+O(\normll{\rhoe}^2)
\sim\normll{\rhoe}.
  \end{split}
\end{equation}

Furthermore, recalling \eqref{xh} (twice),  $\psi=\tk\psi$, and
\eqref{xb},
\begin{equation*}
  \begin{split}
(1+\eps)a_\eps
&
= (1+\eps)\innprod{\tk\psi,\rhoe}
= \innprod{\psi,(1+\eps)\tk\rhoe}
\\&
= \innprod{\psi,\rhoe}  + \innprod{\psi,R(\rhoe)}
= a_\eps + \innprod{\psi,R(\rhoe)}.
  \end{split}
\end{equation*}
Therefore, appealing to \eqref{xc}, \eqref{xg}, \eqref{xe}, \eqref{xfa}
and \eqref{xi}, we find that
\begin{equation*}
  \begin{split}
\eps a_\eps
&
= \innprod{\psi,R(\rhoe)}
= \innprod{\psi,\tfrac12\rhoe^2} + O(\normoo{\rhoe^3})
\\&
= \innprod{\psi,\tfrac12a_\eps^2\psi^2)} +
  \innprod{\psi,a_\eps\psi\rhoex}
  +  \innprod{\psi,\tfrac12(\rhoex)^2} + O(\normoo{\rhoe}^3)
\\&
= \tfrac12a_\eps^2\int\psi^3\dd\mu
  + O(a_\eps\normll{\rhoex})
  +  O(\normll{\rhoex}^2) + O(\normoo{\rhoe}^3)
\\&
= \tfrac12a_\eps^2\int\psi^3\dd\mu
  + O(a_\eps\normll{\rhoe}^2)
  +  O(\normll{\rhoe}^4) + O(\normoo{\rhoe}^3)
\\&
= \tfrac12a_\eps^2\int\psi^3\dd\mu
  + O(a_\eps^3).
  \end{split}
\end{equation*}
If $\eps>0$ is small enough then, by \eqref{xi}, we have
$a_\eps>0$, and so we can conclude that
\begin{equation}
\label{xj}
\eps  = \tfrac12a_\eps \int\psi^3\dd\mu
  + O(a_\eps^2).
\end{equation}

Finally, let $\eps\downto0$. Then \eqref{xi} and \eqref{xd} imply that
$a_\eps\to0$, and so from \eqref{xj} we see that
$\eps\sim \tfrac12a_\eps\int\psi^3\dd\mu=\Theta(a_\eps)$, and,
more precisely,
\begin{equation}
\label{xk}
a_\eps  = \tfrac2{\int\psi^3\dd\mu} \eps+  O(\eps^2).
\end{equation}
Consequently, using \eqref{xg}, \eqref{xfa}, \eqref{xi} and \eqref{xk},
\begin{equation*}
\begin{split}
\rho(1+\eps)
&
=\int_\sss\rhoe\dd\mu
=a_\eps\int_\sss\psi\dd\mu + \int_\sss\rhoex\dd\mu
\\&
=a_\eps\int_\sss\psi\dd\mu + O(\normll{\rhoex})
= 2\tfrac{\int\psi\dd\mu}{\int\psi^3\dd\mu} \eps+  O(\eps^2),
\end{split}
\end{equation*}
proving \eqref{t5}, and so completing the proof of \refT{T5}.
\end{proof}

As noted in \refSS{SStrans}, \refT{T5} has a simple consequence,
\refC{C5}, showing 
that the rate $c_0\rho'_+(c_0)$ of emergence of the giant component
at the phase transition is maximal in the Erd\H os--R\'enyi case,
and, more generally, when \eqref{homo} holds; see \refE{Ehomo}.

\begin{proof}[Proof of \refC{C5}]
Our aim is to show that if $\kk$ is an irreducible kernel
on a ground space $(\sss,\mu)$ 
for which
\eqref{t5a'} holds, and $c_0\=\norm{\tk}\qi>0$, then $c_0\rho'_+(c_0)\le 2$.

By \refT{T5},
$c_0\rho'_+(c_0)=2\xfrac{\int_\sss\psi \int_\sss\psi^2}{\int_\sss\psi^3}$.
By \refL{Lcomp}, we may assume that $\psi\ge0$. Then,
by \Holder's inequality, $\ints\psi\le\bigpar{\ints\psi^3}^{1/3}$
and $\ints\psi^2\le\bigpar{\ints\psi^3}^{2/3}$, with equality if and
only if $\psi$ is \aex{} constant, \ie, if and only if the constant function
$1$ is an eigenfunction of $\tk$, which is equivalent to \eqref{homo}.
\end{proof}

Turning to the number of edges at the phase transition, we shall
next prove \refP{PC6}, which says that if $\kkn$ is
a graphical sequence of kernels on a vertex space $\vxs$ with limit $\kk$,
and $\norm{\tk}=1$, 
then $\frac1n e(\gnxx{\kkn})\pto\hiik\le 1/2$,
with equality
if and only if \eqref{homo} holds.

\begin{proof}[Proof of \refP{PC6}]
By \refP{PE} we have $e\bigpar{\gnxx{\kkn}}/n\pto \gamma$, where $\gam\=\hiik$.
If $\norm{\tk}=1$, then
\begin{equation*}
\hiik=  \tfrac12\innprod{1,\tk1}
\le \tfrac12 \norm{\tk}=\tfrac12,
\end{equation*}
with equality if and only if the constant function
$1$ is an eigenfunction, \ie{}, if
and only if \eqref{homo} holds.
\end{proof}

\begin{remark}\label{Rextreme}
\refP{PC6} says that the number of edges at the phase
transition is \emph{largest} in the classical Erd\H os--R\'enyi
case, and is \emph{strictly smaller} in all other cases except some
very homogeneous ones.
Together, \refC{C5} and \refP{PC6} say roughly that inhomogeneities
make the giant component appear sooner,
but grow more slowly (initially, at
least).
\end{remark}

\begin{remark}\label{Rproc2} 
Another way to study the number of edges when the giant component
is born is to consider the graph
process in \refR{Rproc}. Let
us stop the growth when the largest component first has at least
$\go(n)$ vertices, where $\go(n)$ is a function chosen
in advance, with $\go(n)=o(n)$, and
$\go(n)$ increasing sufficiently rapidly with $n$. (If $\kk$ is
bounded, we can take any $\go$ with $\log n\ll\go(n)\ll n$, see
\refT{T4}.)
Then, for any $\eps>0$, \whp\ we stop at a time between
$(c_0-\eps)/n$ and $(c_0+\eps)/n$, where $c_0\=\norm{\tk}\qi$, and
it follows easily from \refP{PE} that if $N$ is the number of
edges when we stop, then $N/n\pto\hiik$. Again
we see that the number of edges required for a giant component is
largest in the homogeneous case.
\end{remark}

\begin{remark}
  The proof of part (i) of \refT{T5}
 shows that the function $\rho(c\kk;\cdot)$ defined before \eqref{rho}
  depends
  analytically on $c\neq c_0$ as an element of $L^2$. Hence
  $\nuni$ in \refT{T1A} also depends analytically on $c\neq c_0$.
\end{remark}

We expect the following extensions of \refT{T5} to hold.

\begin{conjecture}
\refT{T5}(i) holds without the condition \eqref{t5a},
\ie{}, $\rho(c)$ is always analytic except at $c_0$.
\end{conjecture}

\begin{conjecture}
Let $\kk$ be an irreducible kernel on a ground
space $(\sss,\mu)$. Then equation \eqref{t5} holds
with the larger
error term $o(\eps)$
whenever $\tk$
has an eigenfunction $\psi$ of eigenvalue $\norm{\tk}$ with
$\ints\psi^3<\infty$; conversely, $\rho'_+(c_0)=0$ if $c_0>0$
but no such $\psi$ exists or $\psi$ exists with $\ints\psi^3=\infty$.
Cf.\
\refSS{SSrank1}. 
\end{conjecture}

\section{Applications and relationship to earlier results}\label{Sapp}

In this section we apply our general results to several specific models 
that have been studied in recent years, and describe the relationships
between our results and various earlier results.

\subsection{Dubins' model}\label{SSDubins} 
A common setting is the following: the vertex space
$\vxs$ is $(\sss,\mu,\xss)$, where $\sss=(0,1]$, $\mu$ is the Lebesgue measure,
and $\xs=(x_1,\ldots,x_n)$ with $x_i=i/n$.  In this case, \eqref{pij}
gives $\pij=\kk(i/n,j/n)/n \bmin1$ for the probability of an edge
between vertices $i$ and $j$.  We shall consider several choices of
$\gk$ in some detail.

Observe first that if $\gk$ is a positive function on  $(0,\infty)^2$ that is
homogeneous of
degree $-1$, then \eqref{pij} yields $\pij=\gk(i,j)\bmin1$.
Since this does not depend on $n$, in this case we can also consider
the infinite graph
$\gxx{\infty}{\kk}$, defined in the same way as $G_n=\gnkx$ but on the
vertex set
\set{1,2,\dots}.
Note that the graphs $\gnkx$ are induced subgraphs of
$\gxx{\infty}{\kk}$ and that we can construct them by successively
adding new vertices, and for each new vertex an appropriate random set
of edges to earlier vertices.

We first consider $\gk(x,y)=c/(x\bmax y)$ with $c>0$, so that if $j\ge c$ then
\begin{equation}
  \label{b3a}
\pij=c/j
\qquad
\text{for $i<j$}.
\end{equation}
In this case we can regard $\gnkx$ as a sequence of graphs grown by
adding new vertices
one at a time where,
when vertex $k$ is added, it gets $\Bi(k-1,c/k)$ edges, whose
other endpoints are chosen uniformly among the other vertices.
(We might instead take $\Po(c)\bmin (k-1)$
new edges, without any difference in
the asymptotic results below.)

This infinite graph $\gxx{\infty}{\kk}$ was considered
by Dubins in 1984, who asked when $\gxx{\infty}{\kk}$ is a.s.\
connected. Dubins' question
was answered partially by Kalikow and Weiss \cite{KW}. A little later
Shepp \cite{Shepp} 
proved that $\gxx{\infty}{\kk}$ is a.s.\ connected if and only if $c>1/4$.
This result was generalized to more general homogeneous kernels by Durrett
and Kesten \cite{DK}.

The finite random graph $\gnkx$ with this $\kk$, i.e.,
with edge probabilities given by \eqref{b3a},
has been studied by
Durrett \cite{Durrett}, who points out that it has the same critical
value $c=1/4$ for the emergence of a giant component as the infinite version
has for connectedness, and by Bollob\'as, Janson and
Riordan
\cite{BJR}
who rigorously show that this example has a phase transition with
infinite exponent. More precisely, denoting $\rho(\kk)$ by
$\rho(c)$, it was shown by Riordan \cite{Rsmall} that
\begin{equation}\label{rsmall}
\rho(1/4+\eps) = \exp\bigpar{-\tfrac{\pi}2\eps^{-1/2} + O(\log\eps)} .
\end{equation}

A similar formula for the closely related CHKNS model (see \refSS{SSchkns}),
introduced by Callaway, Hopcroft, Kleinberg, Newman and Strogatz~\cite{CHKNS},
had been given earlier by
Dorogovtsev, Mendes and Samukhin~\cite{DMS-anomalous}
using non-rigorous methods.

To find the critical value by our methods, we have to find the norm of
$\tk$ on $L^2(0,1)$.
Using the isometry $U:f\mapsto e^{-x/2}f(e^{-x})$ of $L^2(0,1)$ onto
$L^2(0,\infty)$, we may instead consider $\ttk\=U\tk U\qi$, which by a
simple calculation is the integral operator on $L^2(0,\infty)$ with
kernel
\begin{equation}\label{e2aoo}
\tilde\kk(x,y)
= e^{-x/2}\kk(e^{-x},e^{-y})e^{-y/2}
= c e^{-x/2-y/2+x\bmin y}
= c e^{-|x-y|/2}.
\end{equation}
Hence $\ttk$ is the restriction to $(0,\infty)$ of the convolution
with $h(x)\=c e^{-|x|/2}$.
Because of translation invariance, it is easily seen that $\ttk$ has
the same norm as convolution with $h$ on $L^2(-\infty,\infty)$, and
taking the Fourier transform we find
\begin{equation*}
\norm{\tk}
=
\norm{\ttk}
=
\norm{f\mapsto h * f}_{L^2(-\infty,\infty)}
=\sup_{\xi\in\bbR} |\hat h(\xi)|
=\intooo h(x)\dd x = 4c.
\end{equation*}

Thus, \refT{T2} shows that there is a giant component if and only if
$c>1/4$, as shown in Durrett \cite{Durrett} and \cite{BJR}.

To find the size of the giant component is more challenging, and
we refer to Riordan~\cite{Rsmall} for a proof of \eqref{rsmall}.
Note that the hypothesis \eqref{t5a} of \refT{T5} fails, as do the
conclusions in part (ii). Indeed, it is
easy to see that $\tk$ is a non-compact operator, and that it has
no eigenfunctions at all in $L^2$. We suspect that this is
connected to the fact that the phase transition has infinite
exponent.

\subsection{The mean-field scale-free model}\label{SSrtij} 

Another interesting case with a homogeneous kernel as in \refSS{SSDubins} is
$\gk(x,y)=c/\sqrt{xy}$ with $c>0$;
then, for $ij\ge c^2$, we have
\begin{equation}
  \label{b3b}
\pij=c/\sqrt{ij}.
\end{equation}

This model has been studied in detail by Riordan~\cite{Rsmall}.
Considering the sequence $\gnkx$ as a growing graph, in this case,
together with each new vertex
we add a number of edges that has approximately a Poisson
$\Po(2c)$ distribution;
the other endpoint of each edge is chosen with probability proportional
to $i\qhi$, which is approximately proportional to the degree of vertex $i$.
Hence, this random graph model resembles the growth with preferential
attachment model of
Barab\'asi and Albert \cite{BAsc}, which was made precise
as the \emph{LCD model} by
Bollob\'as and Riordan~\cite{diam};
see also \cite{Rsmall}.
In fact, up to a factor of $1+o(i^{-1})$
in the edge probabilities, the model defined by \eqref{b3b}
is the so called  `mean-field' version of the
Barab\'asi--Albert model, having the same individual edge probabilities,
but with edges present independently. This (by now common) use 
of `mean-field' is not the standard one in physics, where it normally
means that all vertices interact equally. (So the mean-field
random graph model is $G(n,p)$.)

In this case, $\tk$ is an unbounded operator, because
$x^{-1/2}\not\in L^2(0,1)$, and thus there is no threshold.
In other words, $\rho(c)\=\rho(\kk)>0$ for every $c>0$.

As shown by Riordan~\cite{Rsmall},
$\rho(c)$ grows very slowly at first
in this case too; more precisely, 
\begin{equation}\label{rsmall1}
\rho(c) \sim 2e^{1-\gam}\exp\bigpar{-1/(2c)}
\qquad \text{as } c\to0,
\end{equation}
where $\gam$ is Euler's constant;
see also \refSS{SSrank1} below.
The result in \cite{Rsmall} for the
Barab\'asi--Albert model is different, showing that in this model the dependence
between edges is important.

\begin{remark}
Random graphs related to the ones defined here and in \refSS{SSDubins}
but with some dependence between edges
(and thus not covered by the present paper)
can be obtained by adding at each
new vertex a number of edges with some other distribution, for example
$\Bi(m,p)$
for some fixed $m$ and $p$. Such random graphs have been considered in
\cite{robust,giant2,DevroyeMcR,Rsmall}, and
these papers show that not only the expected numbers of edges
added at each step
are important, but also the variances; the edge dependencies shift the
threshold.
\end{remark}

\subsection{The CHKNS model}\label{SSchkns}
We next consider the CHKNS model of
Callaway, Hopcroft, Kleinberg, Newman and Strogatz~\cite{CHKNS}.
Here, the graph grows from a single
  vertex; vertices are added one by one, and after each vertex is
  added, an edge is added with probability $\gd$; the endpoints are
  chosen uniformly among all existing vertices. (Multiple edges are
  allowed; this does not matter for the asymptotics.)

Following Durrett \cite{Durrett}, we consider a modification
(which is perhaps at least as natural):
after adding each vertex,
add a Poisson $\Po(\gd)$ number of edges to the graph, again choosing
the endpoints
of these edges uniformly at random.
Thus, when vertex $k$ is added, each existing pair of vertices acquires
$\Po\bigpar{\gd/\binom k2}$ new edges, and these numbers are
independent.
When we have reached $n$ vertices, the number of edges between
vertices $i$ and $j$, with $1\le i\le j\le n$, is thus Poisson with
mean
\begin{equation}\label{chk1}
 e_{ij}
\=\sum_{k=j}^n \frac{\gd}{\binom k2}
= 2\gd \sum_{k=j}^n \frac1{k(k-1)}
=2\gd\Bigpar{\frac1{j-1}-\frac1n},
\end{equation}
and the probability that there is one or more edges between $i$ and $j$ is
$p_{ij}\=1-\exp(-e_{ij})$.

Hence, ignoring multiple edges, we have a graph $G_n$ of our type,
with
$\sss=(0,1]$, $\mu$ Lebesgue measure, $x_i=i/n$ and
\begin{align}
\kk_n(x,y)
&\=
n\Bigpar{1-\exp\Bigpar{-2\gd\Bigpar{\frac1{n(x\bmax y)-1}-\frac1n}}}
\notag
\\&
\to
\kk(x,y)
\=2\gd\Bigpar{\frac1{x\bmax y}-1}.
\label{chkns}
\end{align}

The conditions of \refT{T2} are immediately verified, and thus\linebreak
$C_1(G_n)/n\pto\rho(\kk)$.

Instead of adding a Poisson number of edges at each step, 
we could add a
binomial number by adding, after vertex  $k$, each possible edge with
probability $\gd/\binom k2$. We obtain the same results with slightly
different $\kk_n$ but the same $\kk$.

The original CHKNS model, $\tgn$, say, can be treated by the argument
in \cite{BJR}. It follows that $C_1(\tgn)/n\pto\rho(\kk)$ holds for
the CHKNS model too.

In particular, the threshold for the CHKNS model, as
well as for Durrett's modification, is given by $\norm{\tk}=1$, or
$2\gd=\norm{T}\qi$, where $T$ is the integral operator with kernel
$1/(x\bmax y)-1$ on $L^2(0,1)$.
This kernel is strictly smaller that the kernel $1/(x\bmax y)$
considered in \refSS{SSDubins}.
However, changing variables as in \eqref{e2aoo}, we see that $T$ is
equivalent to the operator on $L^2(0,\infty)$ with kernel
$e^{-|x-y|/2}-e^{-(x+y)/2}$.
Using translational invariance of the
operator with kernel
$e^{-|x-y|/2}$ considered in \refSS{SSDubins}, considering
functions supported in $(R,\infty)$ and letting $R\to\infty$,
it is easily seen that $T$ has the same norm as this operator,
namely 4.

Thus the thresholds for the CHKNS model and Durrett's modification
are both given by $2\gd=1/4$, \ie{} $\gd=1/8$, as
was found by non-rigorous arguments by
Callaway, Hopcroft, Kleinberg, Newman and Strogatz
\cite{CHKNS} and
Dorogovtsev, Mendes and Samukhin~\cite{DMS-anomalous},
and first proved rigorously by
Durrett \cite{Durrett}; see also \cite{BJR}.

\medskip
To study the size of the giant component in these models, let us write
$\kk_0(x,y)\=1/(x\bmax y)$ and
$\kk_1(x,y)\=1/(x\bmax y)-1$.
Then $\kk_1< \kk_0$, and thus $\rho(c\kk_1)\le\rho(c\kk_0)$ for each
$c>0$.
(We have strict inequality for $c>1/4$, see \refR{Rmono};
note that,
as pointed out by Durrett \cite{Durrett}, we have the same threshold
$1/4$ for both kernels although we have twice as many edges in
$\gnx{c\kk_0}$ as in $\gnx{c\kk_1}$.)
On the other hand, let $\eta>0$ and
 consider only vertices $i\le j\le
\eta n$. Then
\begin{equation*}
  n\qi c\kk_1(i/n,j/n)
= c\Bigpar{\frac1j-\frac1n}
\ge \frac{(1-\eta)c}j.
\end{equation*}
Hence, \cf{} \eqref{b3a},
$\gnx{c\kk_1} \supseteq \gxx{\eta n}{(1-\eta)c\kk_0}$. (Note that
these graphs have different numbers of vertices.) Thus,
for every $\eta$ with $0<\eta<1$,
\begin{equation*}
  \rho(c\kk_0)
\ge \rho(c\kk_1)
\ge
\eta\rho\bigpar{(1-\eta)c\kk_0}.
\end{equation*}
Taking $c=1/4+\eps$ and $\eta=\eps^2$, relation
\eqref{rsmall} for $\rho(c\kk_0)$ implies the same estimate for $\rho(c\kk_1)$;
in other words, if $\gd=1/8+\eps$, then the size of the giant
component is given by
\begin{equation*}
\rho(\kk) = \rho(2\gd\kk_1)
= \exp\Bigpar{-\frac{\pi}{2\sqrt2}\eps^{-1/2} + O(\log\eps)} .
\end{equation*}

As noted in \refSS{SSDubins}, a similar formula (with no error term, and a
particular constant in front of the exponential) was given by
Dorogovtsev, Mendes and Samukhin~\cite{DMS-anomalous}, with a derivation part 
of which can be made rigorous; see Durrett~\cite{Durrett}.

\subsection{The rank 1 case}\label{SSrank1} 
In this subsection we consider a special case of our general model
that, while very restrictive, is also very natural, and includes
or is closely related to many random graph models considered by other
authors. This is the \emph{rank 1} case, where the kernel $\kk$
has the form
$\kk(x,y)=\psi(x)\psi(y)$ for some function $\psi>0$ on
$\sss$. We shall assume that 
the kernel is graphical; in particular we assume
$\int\psi\dd\mu<\infty$, but not necessarily
that $\int\psi^2\dd\mu<\infty$.

The function $\psi(x)$ can be interpreted as the ``activity'' of a
vertex at $x$,
with the probability of an edge between two vertices proportional to
the product of their activities.
In the rank 1 case,  $\tk f=\bigpar{\int f\psi} \psi$, so 
\begin{equation}\label{r1n}
 \norm{\tk}=\normll{\psi}^2=\int\psi^2\dd\mu\le\infty.
\end{equation}
Thus $\tk$ is bounded if and only if $\psi\in L^2$, in which case
$\tk$ has rank 1, so it is compact, and $\psi$ is the unique (up to
multiplication by constants) eigenfunction with non-zero eigenvalue.

By \refT{TD}, the distribution of vertex degrees is governed by the
distribution of the function $\gl(x)=(\int\psi\dd\mu)\psi(x)$ on
$(\sss,\mu)$. In particular, by \refC{CD},
the degree sequence will (asymptotically) have a power-law tail
if the distribution of $\gl(x)$ has; for example, if
$\sss=(0,1]$ with $\mu$ Lebesgue measure, and
$\psi(x)=cx^{-1/p}$.
(Another, perhaps more canonical, version is to take $\psi(x)=x$
on $\sss=[0,\infty)$, with a suitable finite Borel
measure $\mu$.
Note that every random graph considered in this example may be defined
in this way, since we may map $\sss$ to $[0,\infty)$ by
  $x\mapsto\psi(x)$.
Alternatively, we may map by   $x\mapsto\gl(x)$ and have
$\psi(x)=cx$ with $c>0$ and $\gl(x)=x$.)

Random graphs of this type have been studied in several
papers; we shall not attempt a complete list, mentioning only several 
examples.
Chung and Lu \cite{ChungLu} and
Norros and Reittu \cite{NorrosR} give results on the existence and
size of a giant 
component.
Britton, Deijfen and Martin-L\"of \cite{BrittonDML} use
\eqref{pij''} with $\kk(x,y)=\psi(x)\psi(y)$ to define a random graph,
and observe that conditioned on the vertex degrees, the
resulting graph is uniformly distributed over all graphs with the
given degree sequence; they further prove a version of \refT{TD} for
this case.

Actually, in
\cite{ChungLu} and \cite{NorrosR}
the edge probabilities $\pij$ are given by $\pij\=w_iw_j/\sum_{i=1}^n w_i$, with $w_i$
deterministic in \cite{ChungLu}
and random in \cite{NorrosR}. 
Under suitable conditions on the $w_i$, these examples are also 
special cases of our general model.
For suitable deterministic sequences
$(w_i)_1^n$, we use \refD{Dg2};  we omit the details.
For random \iid{} $w_i$, as in \cite{NorrosR}, if we further
assume $\E w_i=\omega<\infty$, we can, for example, let
$\sss=[0,\infty)$, $\mu=\cL(w_1)$,
$x_i=w_i\bigpar{\sum_jw_j/n\omega}\qhi$
and $\psi(x)=\omega\qhi x$. Then $\kk(x,y)=xy/\omega$, and
we have $\gl(x)=x$ in \refT{TD}, and thus $\gL=w_1$ in \refC{CD}
and $\Xi\sim\Po(w_1)$.
Furthermore, from \eqref{r1n} the norm of $\tk$ is (essentially)
the `second order average degree' ${\overline d}=\sum w_i^2/\sum w_i$.
Thus, for example,
the result of Chung and Lu~\cite{ChungLu:dist2002,ChungLu:dist2003}
that, under certain assumptions, the typical distance between
two vertices of the model $G({\bfw})$ studied in \cite{ChungLu} is
$\log n/\log({\overline d})$ 
corresponds to~\refT{Tdist}. (Chung and Lu also study sequences $w_i$ falling
outside the scope of our model.)

Chung and Lu~\cite{ChungLu:volume} give a result for the `volume'
$\sum_{i\in \cC_1} w_i$ of the giant component $\cC_1$ of $G({\bfw})$.
This result corresponds to \refT{T1A} with $f(x)=\gl(x)$; indeed,
under certain assumptions on the $w_i$, it is
implied by \refT{T1A}. 
Unfortunately, the statement of the result in \cite{ChungLu:volume} is
incomplete, 
as no conditions on the $w_i$ are given. It is not clear what the
right conditions 
are; certainly some restrictions are needed.

\medskip
The random graphs \gnk{} obtained from rank 1 kernels should be compared to
the random graphs with a given (suitably chosen) degree sequence $(d_i)_1^n$,
studied by, for example, Luczak \cite{Luczak:sparse},
Molloy and Reed \cite{MR1,MR2} and
(in the power-law case) Aiello, Chung and Lu \cite{AielloCL}.
Note that in this model, the probability of an edge between $i$ and
$j$ is roughly $d_id_j/n$, but there are dependencies between the edges.
It was shown by Molloy and Reed \cite{MR1} that, 
under some conditions, the threshold for the existence of a giant
component in this model is $\sum_i d_i(d_i-2)=0$. This fits well with
our result, although we see no strict implication: \refT{TD}
shows that, for our model, the degree of a random vertex
converges in distribution to a random variable $\Xi$
with the mixed Poisson distribution $\int_\sss\Po(\gl(x))\dd\mu(x)$.
If $X\sim \Po(\gl)$, then $\E(X(X-2))=\E(X(X-1)-X)=\gl^2-\gl$, so
\begin{equation*}
\E\bigpar{ \Xi(\Xi-2)}
=\ints\bigpar{\gl(x)^2-\gl(x)}\dd\mu(x)
=\Bigpar{\ints\psi\dd\mu}^2\Bigpar{\ints\psi^2\dd\mu-1},
\end{equation*}
which vanishes when $\ints\psi^2=1$. As $\norm{\tk}=\ints\psi^2$,
this is indeed the threshold for the emergence of a giant
component in our model.
The result of Molloy and Reed~\cite{MR1} that, in the supercritical case, 
the second largest component has size $O(\log n)$ corresponds to
\refT{T4}\ref{T4p2}; 
again there is no strict implication, but the kernels $\kk$ corresponding
to the graphs studied by Molloy and Reed satisfy $\inf \kk(x,y)>0$.
In a subsequent paper,
Molloy and Reed \cite{MR2} gave further results on the size of the
giant component and on the structure of the remainder of the graph,
corresponding to our Theorems \refand{T2}{Tdual}.

The following variant of this model has also been studied: 
the degrees are first chosen according to some
distribution, and then the graph is chosen uniformly among all graphs
with the resulting degree sequence; see, for example,
the results of Van der Hofstad, Hooghiemstra and Van Mieghem~\cite{HHM_RSA}
and of Fernholz and Ramachandran~\cite{FR:diam} on distances and diameter,
respectively,
mentioned in \refS{Sdist}.

Yet another variant of the rank 1 case of $\gnkx$ was studied 
rather earlier by Khokhlov and Kolchin~\cite{KK1,KK2}, 
who proved results about the number of cycles; see \refS{Scycles}.

\medskip
In the rank 1 case, the size of 
the giant component (if any) of $\gnkx$ 
can be found rather easily.
In order to study the phase transition,
let us consider the kernel $c\kk(x,y)=c\psi(x)\psi(y)$, with $c>0$
a parameter.  
By \refC{C1}, the threshold for $c$ is
$c_0=\norm{\tk}\qi=\Bigpar{\ints\psi^2}\qi$.
For $c\ge c_0$,
let
\begin{equation}
  \label{sea1}
\ga(c)\=c\int\psi\rhock\dd\mu
,
\end{equation}
where, as before,
$\rhock(x)=\rho(c\kk;x)$
is the survival probability of the branching process $\bpxx{c\kk}$.

We have $\Tx{c\kk}\rhock = c\tk\rhock = \ga(c)\psi$. Thus, by \refT{T:GW}
and \eqref{sea1},
\begin{equation}
  \label{ea}
\rhock = \Phix{c\kk}(\rhock)
= 1-e^{-\Tx{c\kk}\rhock}
= 1-e^{-\ga(c)\psi}.
\end{equation}
(The condition \eqref{bb1} holds for every $x$.)
Let
\begin{equation}
  \label{eb}
\gb(t)\=\int_{\sss}\Bigpar{1-e^{-t\psi(x)}}\psi(x)\dd\mu(x),
\qquad t\ge0.
\end{equation}
Then, by  \eqref{sea1} and \eqref{ea},
\begin{equation}\label{ebba}
  \ga(c)=c\int_\sss\rhock\psi\dd\mu = c\gb\bigpar{\ga(c)},
\end{equation}
so $c=\ga(c)/\gb\bigpar{\ga(c)}$, \ie{}, $\ga$ is the inverse function
to $t\mapsto \gam(t)\=t/\gb(t)$.
Since $\gb$ is explicitly given by \eqref{eb}, for any $c>c_0$ 
this gives (at least in principle)
$\ga(c)$, and hence, by \eqref{ea}, the function $\rhock$.
Then $\rho(c\kk)=\ints\rhock(x)\dd\mu(x)$ determines the asymptotic number of vertices in the giant component. 
Similarly, by \refT{Tedges},
the asymptotic number of edges in the giant
component is determined 
by 
$\edgeno(c\kk)$, which
by the definition \eqref{Zdef} and \eqref{ebba} is given by
\begin{equation}\label{gron}
  \begin{split}
\edgeno(c\kk)
&
=c\ints\psi\rhock\dd\mu \ints\psi\dd\mu
-\frac c2 \Bigpar{\ints\psi\rhock\dd\mu}^2
=\ga(c) \ints\psi\dd\mu
-\frac {\ga(c)^2}{2c}
.
  \end{split}
\end{equation}
Moreover, the asymptotic value of $\sum_{i\in\cC_1}f(x_i)/n$ is given
by \refT{T1A} for suitable functions $f$. 

We now turn to asymptotics as $c\downto c_0$, in order to study the
phase transition more closely. 
Note that $\rhock\downto0$ \aex{} as $c\downto c_0$ by \refT{TappB}
so, by dominated convergence,
\begin{equation}\label{eaa}
 \ga(c)/c\downto 0\quad\text{as}\quad c\downto c_0.
\end{equation}
Further, by \eqref{ea} and dominated convergence,
\begin{equation}\label{rovera}
  \frac{\rho(c)}{\ga(c)}
=\int_{\sss}   \frac{\rhock(x)}{\ga(c)} \dd\mu(x)
\to\int_\sss \psi(x)\dd\mu(x)>0
\qquad\text{as }c\downto c_0.
\end{equation}
Consequently, the behaviour of $\gb$ at 0 determines, through $\gam$
and $\ga=\gam\qi$, the behaviour of $\rho(c)$ as $c\downto c_0$.
Note that, by \eqref{eb}, $\gb$ is continuous with
$\gb(0)=0$ and
\begin{equation}
  \label{ebi}
\gb'(t)=\int_{\sss}e^{-t\psi(x)}\psi^2(x)\dd\mu(x),
\qquad t>0.
\end{equation}
Moreover, from \eqref{gron}
we have
\begin{equation*}
  \begin{split}
\edgeno(c\kk)
&
=\ga(c) \ints\psi\dd\mu
-\frac {\ga(c)^2}{2c}
\sim \ga(c) \ints\psi\dd\mu
\sim\rho(c)
  \end{split}
\end{equation*}
as $c\downto c_0$, where the first $\sim$ is from \eqref{eaa} and
the second from \eqref{rovera}.
Hence the asymptotics are the same as for the number of vertices;
see \refR{Retrans}. 

\medskip
Let us consider some concrete examples. Once again, we take
$\sss=(0,1]$ with $\mu$ Lebesgue measure,
and let $\psi(x)=x^{-1/p}$ where $1<p\le\infty$.
We shall use $C$, $C_1$, etc. to  denote various positive constants
that depend on $p$.

\step{Case 1: $1<p<2$}
In this case, $\normll{\psi}=\infty$, so $c_0=0$.
As $t\to0$, by \eqref{ebi},
\begin{equation}\label{ebi1}
\gb'(t)\=\intoi e^{-t x^{-1/p}}x^{-2/p}\dd x
= p\int_t^\infty e^{-y} t^{-2+p} y^{1-p} \dd y
\sim C t^{p-2},
\end{equation}
noting for the last step that the integral $\int_0^\infty
e^{-y}y^{1-p}\dd y$
is convergent.
Thus $\gb(t)\sim C_1 t^{p-1}$ and $\gam(t)=t/\gb(t) \sim C_2 t^{2-p}$.
Consequently, using \eqref{rovera},
\begin{equation*}
\rho(c)\sim C_3 \ga(c)=C_3\gam\qi(c) \sim C_4 c^{1/(2-p)}
\qquad \text{as } c\to0.
\end{equation*}
Note that this exponent $1/(2-p)$ may be any real number in $(1,\infty)$.

\step{Case 2: $p=2$}
This is the case  \eqref{b3b} studied in \refSS{SSrtij} and \cite{Rsmall}.
We still have $\normll{\psi}=\infty$ and thus $c_0=0$.
In analogy with \eqref{ebi1} we now find that  $\gb'(t)\sim 2
\ln(1/t)$ as $t\to0$.
This yields $\gb(t)\sim 2t \ln(1/t)$ and $\gam\sim1/(2\ln(1/t))$ as $t\to0$,
and thus $\ga(c)=\gam\qi(c)=e^{-(1+o(1))/2c}$ and
\begin{equation*}
\rho(c)=e^{-(1+o(1))/2c}
\qquad \text{as } c\to 0.
\end{equation*}
More refined estimates
can be obtained in the same way, see
\eqref{rsmall1} and \cite{Rsmall}.

\step{Case 3: $2<p<3$}
For $p>2$ we have $\int \psi^2\dd\mu<\infty$, and thus $c_0>0$, so we
have a phase transition. (In fact, $c_0=1-2/p$.)
By \eqref{ebi}, $\gb'(t)$ is continuous for $t\ge0$ with
$\gb'(0)=\int\psi^2\dd\mu=c_0\qi$. Differentiating once more we obtain
as $t\to0$
\begin{equation*} 
\gb''(t)\
=-\intoi e^{-t x^{-1/p}}x^{-3/p}\dd\mu(x)
= -p\int_t^\infty e^{-y} t^{-3+p} y^{2-p} \dd y
\sim -C t^{p-3}
\end{equation*}
and thus $\gb'(t)=c_0\qi -(C_1+o(1)) t^{p-2}$ and
$\gb(t)=c_0\qi t -(C_2+o(1)) t^{p-1}$.
Hence $\gam(t)=t/\gb(t) = c_0+ (C_3+o(1)) t^{p-2}$.
Consequently, using \eqref{rovera},
\begin{equation*}
\rho(c_0+\eps)\sim C_4 \ga(c_0+\eps)\sim C_5 \eps^{1/(p-2)}
\qquad \text{as } \eps\downto 0.
\end{equation*}
We thus have a phase transition at $c_0$ with exponent $1/(p-2)$. Note that
this exponent may be any real number in $(1,\infty)$.
(Taking instead \eg{} $\psi(x)=x^{-1/2}\ln\qi(e^3/x)$,
it is similarly seen that
there is a phase transition with infinite exponent.)

\step{Case 4: $p=3$}
Similar calculations show that,
as $t\to0$,
$\gb''(t)\sim 3\ln t$, $\gb'(t)=3-(3+o(1))t\ln(1/t)$,
$\gb(t)=3t-(3/2+o(1))t^2\ln(1/t)$, and
$\gam(t)=1/3+(1/6+o(1))t\ln 1/t$.
Consequently, with $c_0=1/3$,
\begin{equation*}
\rho(c_0+\eps)\sim C \ga(c_0+\eps) \sim C_1 \eps/\ln (1/\eps)
\qquad \text{as } \eps\downto 0,
\end{equation*}
so $\rho'(c_0)=0$.

\step{Case 5: $3<p\le \infty$}
In this case, $\int \psi^3\dd\mu<\infty$
and we find as $t\to0$,
$\gb''(t)\sim -C$, $\gb'(t)=c_0\qi-(C+o(1))t$,
$\gb(t)=c_0\qi t-(C+o(1))t^2/2$, and
$\gam(t)=c_0+(C_1+o(1))t$.
Consequently, $\rho(c_0+\eps)\sim C_2 \ga(c_0+\eps) \sim C_3 \eps$,
so we have a phase transition with exponent 1.
This is similar to
\refT{T5}, although \eqref{t5a} is not satisfied (except in the
classical case $p=\infty$). Indeed, it can be checked that
\eqref{t5} holds, except that the error term may be larger (it is
$\Theta(\eps^{p-2})$ for $3<p<4$).

More generally, the same argument shows that \eqref{t5} holds
for any rank 1 kernel $\psi(x)\psi(y)$ with
$\int\psi^3\dd\mu<\infty$, provided
the error term is weakened to $o(\eps)$.
(The error term is $O(\eps^2)$ if $\int\psi^4\dd\mu<\infty$.)

\subsection{Turova's model}\label{SSTurova}
Turova \cite{Turova,Turova2,Turova3,Turova4}
has studied a dynamical random graph
$G(t)$, $t\ge0$,
defined as
follows, using three parameters $\gam>0$, $\gl>0$ and $\mux\ge0$.
The graph starts with a single vertex at time $t=0$.
Each existing vertex produces new, initially isolated, vertices
according to a Poisson process with intensity $\gam$.
As soon as there are at least two vertices, each vertex sends out
edges according to another Poisson process with intensity $\gl$; the
other endpoint is chosen uniformly among all other existing vertices.
(Multiple edges are allowed, but this makes little difference.)
Vertices live for ever,  but edges die with intensity $\mux$, \ie{}, the
lifetime of an edge has an exponential distribution with mean
$1/\mux$. (All these
random processes and variables are independent. We use $\mux$ for
Turova's $\mu$ to avoid conflicts with our notation.)

By homogeneity we may assume $\gam=1$; the general case follows by
replacing $\gl$ and $\mux$ by $\gl/\gam$ and $\mux/\gam$ and changing
the time scale.

\medskip
Our analysis of the random graph $G(t)$ is very similar to that of
S\"oderberg \cite{Sod1}; our theorems enable us to add technical
rigour to his calculations.  The vertices proliferate according to a
Yule process (binary fission process): writing $N(t)$ for the number
of vertices at time $t$, the probability that a new vertex is added in
the infinitesimal time interval $[t,t+\dd t]$ is $N(t)\dd t$. It is
well-known (see,  \eg{}, Athreya and Ney \cite[Theorems III.7.1--2]{AN}) that
\begin{equation}\label{yulew}
  e^{-t} N(t) \asto W
\qquad \text{as \ttoo}
\end{equation}
for a random variable $W$ with $W>0$ a.s.
(In fact, $W\sim\Exp(1)$, 
but we do not need this.)

We condition on the vertex process, and assume, as we may by
\eqref{yulew}, that
\begin{equation}\label{yule}
  e^{-t} N(t) \to w
\qquad \text{as \ttoo}
\end{equation}
for some $w>0$.
We take $\sss=[0,\infty)$ and let $x_1,\dots,x_{N(t)}$ be the ages of
the particles existing at time $t$. For any fixed $s>0$, by
  \eqref{yule} we have
  \begin{equation*}
\nu_t[s,\infty)
\= \frac1{N(t)} \#\set{i:x_i\ge s}
= \frac1{N(t)} N(t-s) \to e^{-s},
  \end{equation*}
as $\ttoo$.
This means (see \refR{RR}) that
$\nu_t\to\mu$, where $\mu$ is the measure on $[0,\infty)$ given by
  $d\mu/dx=e^{-x}$ (the exponential distribution).

If $x_i\le x_j$, the number of edges at time $t$ between two vertices
of ages $x_i$ and
$x_j$ has a Poisson distribution with mean
\begin{equation*}
  e_{ij}
\= \int_{t-x_i}^t e^{-\mux(t-s)} \frac{2\gl}{N(s)-1}\dd s
= 2\gl\int_{0}^{x_i} e^{-\mux s} \frac{\dd s}{N(t-s)-1}.
\end{equation*}
Set
\begin{equation*}
\kktx(x,y)
\= 2\gl \int_{0}^{x\bmin y} e^{-\mux s} \frac{N(t)}{N(t-s)-1} \dd s,
\end{equation*}
and
\begin{equation*}
 \kkt(x,y)=N(t)\bigpar{1-\exp(-\kktx(x,y)/N(t))}.
\end{equation*}
Thus $e_{ij}=\kktx(x_i,x_j)/N(t)$, and the probability $p_{ij}$
that there is at least one edge between $i$ and $j$ is given by
$p_{ij}=1-e^{-e_{ij}}=\kkt(x_i,x_j)/N(t)$.

By \eqref{yule},
$N(t)/\bigpar{N(t-s)-1} \to e^{s}$ as \ttoo{} for every $s$, and
dominated convergence shows that if $\mux\ne 1$ and $x_t\to x$,
$y_t\to y$, then
\begin{equation}\label{kkt}
\kktx(x_t,y_t)
\to
\kk_\mux(x,y)
\= 2\gl \int_{0}^{x\bmin y} e^{-\mux s +s} \dd s
=\frac{2\gl}{1-\gd}\bigpar{e^{(1-\gd)(x\bmin y)}-1}
\end{equation}

For $\mux=1$, corresponding to $\mux=\gamma$ in the non-rescaled model,
let $\kk_1(x,y)\=2\gl(x\bmin y)$. Then
$\kktx(x_t,y_t)\to \kk_\mux(x,y)$ in this case also.

\refT{T2} thus applies to $G(t)$ conditioned on the process $(N(t))_{t\ge0}$,
and we find (conditioned on $(N(t))_{t\ge0}$, and thus also unconditionally) that
\begin{equation*}
  \frac{C_1\bigpar{G(t)}}{N(t)} \pto \rho(\kk_\mux),
\end{equation*}
with $\kk_\mux$ given by \eqref{kkt}.

To study $\rho(\kk_\mux)$ further, and in particular to investigate the
threshold as we vary $\gl$ keeping $\mu\ge0$ fixed, we thus have to
investigate the integral operator $T_{\kk_\mux}$ with kernel
$\kk_\mux$ given by
\eqref{kkt}.
The change of variables $x\to e^{-x}$ transforms $\sss$ and $\mu$ to
the standard setting $(0,1]$ with Lebesgue measure, and the kernel
\eqref{kkt} becomes
\begin{equation}
  \label{tur}
\tkkmux(x,y)\=\frac{2\gl}{1-\gd}\Bigpar{(x\bmax y)^{\gd-1}-1},
\end{equation}
with $\tkko\=2\gl\ln(1/(x\bmax y))$.

In the case $\gd=0$, this is the same as \eqref{chkns}; hence we have
the same critical value $1/8$ (for $\gl$) as for the CHKNS model and the same
$\rho(\kk)$ giving the size of the giant component; in particular, the
phase transition has infinite exponent.
(Indeed, with $\gd=0$ the model is very similar to
(Durrett's form of) the CHKNS model discussed in \refSS{SSchkns}; now
a geometric number 
of edges between random vertices is added at each step, rather than a
Poisson number.)
For $\gd>0$, the kernel $\tkkmux$ is in $L^2((0,1]^2)$, so $T_\tkkmux$ is
compact (see \refL{Lcomp}) and its norm can be found by finding its
eigenvalues.
By the discussion in \refSS{SSlund}
below, this is equivalent to solving \eqref{sl2} with the given
boundary values.  
In our case, denoting the eigenvalue by $\ga$, this means solving
$
\ga G''(x)= -2\lambda x^{\gd-2} G(x)
$
with boundary values $G(0)=G'(1)=0$ (since $\phi(1)=0$).

The general solution is easily written down as a linear combination of
two hypergeometric series, and $G(0)=0$ yields
(up to a constant factor)
\begin{equation*}
  \begin{split}
  g(x)
&
=G'(x)=\sum_{n=0}^\infty \frac{1}{n!\,\Gamma(n+1/\gd)}
\Bigpar{-\frac{2\gl}{\ga\gd^2}x^\gd}^n
\\&
=
\Bigpar{\frac{2\gl}{\ga\gd^2}x^\gd}^{-(1/\gd-1)/2}
J_{1/\gd-1}\Bigpar{2\Bigpar{\frac{2\gl}{\ga\gd^2}x^\gd}\qh },
  \end{split}
\end{equation*}
where $J_\nu$ is a Bessel function.

The condition $g(1)=G'(1)=0$
(which gives the formula in
Turova 
\cite[Corollary 4.1]{Turova} and \cite{Turova2})
thus leads to
$J_{1/\gd-1}\Bigpar{\Bigpar{\frac{8\gl}{\ga\gd^2}}\qh }=0$, so if
$z_\nu$ is the first positive zero of $J_\nu$, then
\begin{equation*}
  \norm{T_{\kk_\mux}}=\norm{T_\tkkmux}=\ga_1=\frac{8\gl}{\gd^2 z_{1/\gd-1}^2}.
\end{equation*}
In other words, the critical value of $\gl$ is
$\glcrit(\gd)=\gd^2 z_{1/\gd-1}^2/8$,
as given by a related argument by S\"oderberg~\cite{Sod1}.

\refT{T5} applies only when $\gd>1/2$, but the eigenfunctions are
continuous and bounded for every $\gd>0$, and we believe that
the phase transition has exponent 1, and that
\eqref{t5} holds, for every $\gd>0$.

We can easily find the asymptotics of $\glcrit(\gd)$ as $\gd\to0$ or
$\infty$;
see Turova \cite{Turova}.
If $\gl,\gd\to\infty$ with $\gl/\gd\to c>0$, then
$\tkkmux(x,y)\to 2c$, pointwise and in $L^2((0,1]^2)$, and thus
$\norm{T_\tkkmux-T_{2c}} \le \normHS{T_\tkkmux-T_{2c}} \to0$. It follows that
for large $\gd$, the graph is subcritical if $2c<1$ and
supercritical if $2c>1$. In other words, $\glcrit/\gd \to1/2$ as
$\gd\to\infty$.
Similarly, if $\gd\downto0$, then $\kk_\gd\upto \kk_0$ and it follows
easily, \eg{} by \refT{TappB}, that
$\norm{ T_{\kk_\gd}} \to \norm{ T_{\kk_0}} $, and thus
$\glcrit(\gd)\to\glcrit(0)=1/8$.
(The contrary assertion in \cite{Turova} is incorrect; see the
erratum.)

\subsection{Functions of $\max\{x,y\}$.} \label{SSlund} 
In several of the examples above
(see Subsections \ref{SSDubins}, \ref{SSchkns}, and \ref{SSTurova})
we have $\sss=(0,1]$, $\mu$ is the Lebesgue measure and
$\kk(x,y)=\phi(x\bmax y)$ for some function $\phi\ge0$
on $(0,1]$. The integral operators $\tk$ with such kernels have been
studied by Maz'ya and Verbitsky~\cite{MV} and Aleksandrov, Janson,
Peller and Rochberg~\cite{SJ139}.
In particular, these papers
prove that
$\tk$ is bounded if and only if
$\sup_{x>0} x\int_x^1 \phi(y)^2\dd y <\infty$,
and that
$\tk$ is compact if and only if
$ x\int_x^1 \phi(y)^2\dd y \to0$ as $x\to0$.

In the case when $\phi$ is decreasing (as in the examples above),
these criteria simplify to $\phi(x)=O(x\qi)$ and $\phi(x)=o(x\qi)$ as
$x\to0$, respectively.

Unfortunately, there is no general formula known for the norm of $\tk$.
(However, the criteria just given extend to estimates within constant
factors;
for example, if $\phi$ is decreasing, then
$\sup (x\phi(x)) \le\norm{\tk} \le 4\sup (x\phi(x))$.)
In the compact case,
at least if $\phi$ has a continuous derivative on $(0,1]$,
the eigenvalues, and thus the norm, can be found
by studying a Sturm--Liouville equation.
In fact, $g$ is an eigenfunction with eigenvalue $\gl$
if
\begin{equation}
  \label{sl1}
\gl g(x)=\phi(x)\int_0^x g(y)\dd y + \int_x^1 \phi(y)g(y)\dd y.
\end{equation}
If $\gl\neq0$, it is easily seen that then $g\in C^1(0,1]$ and, by
  differentiating, that \eqref{sl1} is equivalent to
  \begin{align}
\label{sl2}
g(x)=G'(x) ,
& & &
\gl G''(x)=\phi'(x) G(x),
  \end{align}
with the boundary conditions $G(0)=0$, $G'(1)=\gl\qi\phi(1)G(1)$;
see \cite[Section 9]{SJ139} and the example in
\refSS{SSTurova} above.

\section{Paths and cycles}\label{Scycles}

Let $P_k(G)$ and $Q_k(G)$ be the numbers of paths and cycles, respectively,
of length $k$ (\ie{}, with $k$ edges) in a graph $G$.
Note that $P_1(G)=e(G)$ is the number of edges, and that $Q_1=Q_2=0$
for simple graphs.
(If we allow multiple edges and loops as in \refR{Rmulti}, the results
below extend to $Q_2$ and, under an additional continuity assumption
on $\kk$, to $Q_1$.)

In this section we briefly study the numbers $P_k=P_k(G_n)$
and $Q_k=Q_k(G_n)$, where $G_n=\gnkx$.
The results are easily extended to a sequence $\kk_n$ as in \refD{Dg2}
under appropriate conditions, but we leave the details to
the reader.

For $k\ge1$ let
\begin{align*}
\ga_k(\kk)
&\=
\frac{1}{2}
\int_{\sss^{k+1}} \kk(x_0,x_1) \kk(x_1,x_2) \dotsm
\kk(x_{k-1},x_k)\dd \mu(x_0)\dotsm \dd \mu(x_k),
\\
\gb_k(\kk)
&\=
\frac{1}{2k}
\int_{\sss^{k}}  \kk(x_1,x_2) \dotsm \kk(x_{k-1},x_k) \kk(x_k,x_1)
\dd \mu(x_1)\dotsm \dd \mu(x_k).
\end{align*}

Note that $\ga_k(\kk)=\frac{1}{2}\innprod{1,\tk^k1}$. Clearly, $\ga_k(\kk)$
and $\gb_k(\kk)$  may be infinite. In this case, the limiting
statements
in the result below have their natural interpretations.
\begin{theorem}
   \label{TPC}
Let $\kk$ be an \aex{} continuous
kernel on a (generalized)
vertex space $\vxs$, and let $G_n=\gnkx$.
\begin{thmxenumerate}
\item
For $k$ fixed,
\begin{align*}
\liminf_\ntoo \E P_k(G_n)/n
&\ge \ga_k(\kk), && k\ge1,
\\
\liminf_\ntoo  \E Q_k(G_n)
&\ge \gb_k(\kk), && k\ge3.
\end{align*}
\item
Suppose further that $\vxs$ is a vertex space.
If\/ $\kk$ is bounded, or if
$x_1,\dots,x_n$ are \iid{} random points with the distribution $\mu$,
then
\begin{align}
\E P_k(G_n)/n
&\to \ga_k(\kk), && k\ge1,
\label{Pke}
\\
\E Q_k(G_n)
&\to \gb_k(\kk), && k\ge3.
\label{bk}
\end{align}
\end{thmxenumerate}
Moreover, whenever \eqref{Pke} holds and $\ga_k(\kk)$ is finite,
\begin{equation}\label{pp}
   P_k(G_n)/n\pto \ga_k(\kk).
\end{equation}
Similarly, whenever \eqref{bk} holds and the $\gb_k(\kk)$ are finite,
\begin{equation}\label{qd}
   Q_k(G_n)\dto \Po\bigpar{\gb_k(\kk)}, \qquad k\ge3,
\end{equation}
jointly for all $k\ge3$ with independent limits.
\end{theorem}

\begin{proof}
The argument for parts (i) and (ii) is as in the proof of \refL{LE} (a
special case), considering first the \rfin{} case and then
approximating with $\kkm-$ and $\kkm+$
defined in \eqref{k-} and \eqref{k+};
we omit the details. The case
of \iid{} $x_i$ with the distribution $\mu$ is immediate.

The convergence \eqref{pp} and the
asymptotic (joint) Poisson distribution \eqref{qd} of $Q_k(G_n)$ follow
easily in the \rfin{} case (first conditioning on $\xs$ as in \refR{Rcond}
if $\vxs$ is a generalized vertex space),
for example by the method of moments as for
\gnp{}, \cf{} \cite{BB,JLR}.
The general cases then follow by appealing to Billingsley 
\cite[Theorem 4.2]{Bill},
noting that if $\kkm-$ is defined by \eqref{k-}, then
$\ga_k(\kkm-)\to\ga_k(\kk)$ and $\gb_k(\kkm-)\to\gb_k(\kk)$
by the monotone convergence theorem, while,
from the assumption \eqref{bk} and part (i) (applied to $\kkm-$),
\begin{multline*}
   \limsup_\ntoo\E|Q_k\bigpar{G(n,\kk)} - Q_k\bigpar{G(n,\kkm-)}|
\\
=\lim_\ntoo\E Q_k\bigpar{G(n,\kk)}
-
\liminf_\ntoo\E Q_k\bigpar{G(n,\kkm-)}
\le \gb_k(\kk) - \gb_k(\kkm-)
\to0
\end{multline*}
as \mtoo; an analogous bound holds for $P_k/n$.
\end{proof}

Part (ii) holds for many generalized vertex spaces too, but not for 
all. Indeed, \refR{Rgenexp} extends easily to the present situation,
although, writing $v_n$ for the number of vertices of $\gnkx$,
the variance condition $\Var(v_n/n)\to 0$ should be replaced by a higher
moment condition $\E(v_n^\ell/n^{\ell})\to \mu(\sss)^\ell$, with $\ell=k+1$
for \eqref{Pke} and $\ell=k$ for \eqref{bk}.
It is easily seen that \eqref{Pke} and \eqref{bk} hold
in the situation in \refE{E42} too.
However, these relations may fail for the counterexample in \refE{Ebad}.

Khokhlov and Kolchin~\cite{KK1,KK2} studied a model closely related to
the rank 
1 case of $\gnkx$: each vertex has an activity $a_i$, and edges are added
one by one, with the endpoints of the edge chosen independently, and
the probability 
that a vertex is chosen proportional to its activity. They proved results about
the distribution of the numbers of short cycles in this model corresponding
to the last part of \refT{TPC}.

\begin{proposition}\label{Ppaths}
   \begin{thmenumerate}
\item
If\/ $\norm{\tk}\le1$, then $\ga_k(\kk)\le1/2$ for every $k$.
\item
If\/ $\norm{\tk}>1$, then $\ga_k(\kk)\to\infty$ as $\ktoo$.
   \end{thmenumerate}

Consequently, if $\kk$ is a
graphical kernel on a vertex space $\vxs$, then
\gnkx{} has a giant component if and only if
$\sup_k\ga_k(\kk)=\infty$.
\end{proposition}
\begin{proof}
The first statement is immediate, as $\innprod{1,\tk^k1}\le \norm{\tk}^k$.
For the second statement, we argue as in the proof of \refL{L6}: there
is a bounded
kernel $\kk_N\le \kk$ with $\norm{T_{\kk_N}}>1$, and $T_{\kk_N}$ has a bounded
eigenfunction $\psi$ with eigenvalue $\lambda>1$. Taking $\normoo{\psi}=1$
we have
$\innprod{1,\tk^k1}\ge\innprod{1,T_{\kk_N}^k1}\ge\innprod{1,T_{\kk_N}^k\psi}
=\lambda^k\innprod{1,\psi}\to \infty$.
The final statement follows by \refT{T2}(i),
which states that \gnkx\ has
a giant component
if and only if $\norm{\tk}>1$.
\end{proof}

Similarly, at least when $\iikk<\infty$, we have the following
consequence of well-known properties of Hilbert--Schmidt operators
(\cf{} \refL{Lcomp}).

\begin{proposition}
   \label{PC}
If\/ $\iikk<\infty$, then $\tk$ is compact and self-adjoint, and if
$\gl_i$ are its (real) eigenvalues (counted with multiplicities), then
\begin{equation*}
   \gb_k(\kk) = \frac1{2k} \Tr(\tk^k) = \frac1{2k}\sum_i \gl_i^k
   <\infty,
\qquad k\ge2.
\end{equation*}
\end{proposition}

\begin{example}
   \label{EChomo}
As in \refE{Ehomo}, let $\sss=(0,1]$ (regarded as a circle) with $\mu$
Lebesgue measure, and
   $\kk(x,y)=h(x-y)$, with $h\ge0$ an even periodic function that is integrable
   over $(0,1]$.
Then $\tk$ is the convolution operator $f\mapsto h*f$ with eigenvalues
$\hath(j)=\intoi e^{-2\pi i j x} h(x)\dd x$, $j\in\bbZ$, so
   $\gb_k=\sum_{-\infty}^\infty\hath(j)^k$.

Considering functions $h$ with small support, we can obtain
arbitrarily large $\gb_k$ with $\kk$ bounded and $\iik=1$.
Alternatively, we can take $\hath(j) = 1/\ln(2+|j|)$, for example;
this defines an integrable function $h>0$ which is continuous except
at 0 \cite[Theorems V.(1.5) and V.(1.8)]{Zygmund}, and thus a kernel
$\kk$ with $\norm{\tk}=\intoi h<\infty$ but $\gb_k(c\kk)=\infty$ for
every $k\ge2$ and $c>0$.
\end{example}

\begin{example}
   \label{ECT}
Let $\kk(x,y)=c/(x\bmax y)$ on $\sss=(0,1]$ as in \refSS{SSDubins}.
Then $\gb_k(\kk)=\infty$ for every $k$ and every $c>0$; indeed,
if, say, $x_i=i/n$,
the expected number of $k$-cycles with vertices in $(2^{-m-1}n,2^{-m}n)$
tends to a positive constant independent of $m\ge0$, and thus
$\E Q_k\to\infty$.

The same holds for $\kk(x,y)=c(1/(x\bmax y) -1)$ as in \eqref{chkns}
and the $\gd=0$ case of \eqref{tur}.
\end{example}

\begin{example}
   Let $\kk(x,y)=\phi(x\bmax y)$ with $\sss=(0,1]$, as in \refSS{SSlund},
   and assume that $\phi\ge0$ is non-increasing with $x\phi(x)$ bounded.
Assume that $x_1,\dots,x_n$ are \iid{} and uniformly distributed on
   $(0,1]$.
(It can be checked that the same conclusions hold for $x_i=i/n$.)
Then, results of Aleksandrov, Janson, Peller and
Rochberg~\cite[Theorems 4.1 and 4.6]{SJ139} imply that $\tk$ is a bounded
   positive operator, and it is compact with eigenvalues $\gl_i$
   satisfying $\sum_i \gl_i^k<\infty$
(which means that $\tk$ belongs to the Schatten--von Neumann class $S_k$)
if and only if
$\intoi\bigpar{x\phi(x)}^k/x\dd x<\infty$.

By \refT{TPC}, we have  $\E Q_k\to\gb_k(\kk)\le\infty$ for $k\ge3$.
\refP{PC} assumes $\iikk<\infty$, but it can be shown (using
truncations of $\phi$) that the result
extends to the present situation;
hence $\gb_k(\kk)<\infty$ if and
only if
$\intoi\bigpar{x\phi(x)}^k/x\dd x<\infty$.

Consequently, we may for any given $\ell$ choose $\phi$ such that
$\gb_k=\infty$ for $3\le k\le \ell$ but $\gb_k<\infty$ for $k>\ell$.
\end{example}

Under suitable conditions, the expected total number of cycles
converges to $\sum_{k=3}^\infty \gb_k$; we omit the details.
By \refP{PC}
this sum is given by the following formula.

\begin{corollary}
   \label{CPC}
If\/ $\iikk<\infty$, then
\begin{equation*}
\sum_{k=3}^\infty  \gb_k(\kk)
=
\begin{cases}
\sum_i\bigpar{-\tfrac12\ln(1-\gl_i)-\tfrac12\gl_i-\tfrac14\gl_i^2}
   <\infty,
&
\text{if\/ } \norm{\tk}<1,
\\
\infty,
&
\text{if\/ } \norm{\tk}\ge1.
\end{cases}
\end{equation*}
\end{corollary}

The sum on the right-hand side can be written as
$-\tfrac12\ln\bigpar{\det_3(I-\tk)}$, where $\det_3$ is a
renormalized Fredholm determinant \cite[\S9]{Simon}.

\begin{remark}
   \label{RTurova}
Turova \cite{Turova2} studies the number of cycles in the random graph
discussed in \refSS{SSTurova}, including a formula for $\lim\E Q_k=\gb_k$.
She conjectures 
that the threshold for the existence of a giant component
is the same as the threshold for $\sum_k \gb_k=\infty$. (This
conjecture inspired the present section.)
We now see from \refC{CPC} that this is true in great
generality; for example, if $\kk$ is bounded, then \refT{TPC} and
\refP{PC} imply that the threshold $c_0=\norm{\tk}\qi$ in \refC{C1}
may be written as
\begin{equation}
   \label{turc}
c_0=\sup\Bigset{c:\sum_{k=3}^\infty \gb_k(c\kk)<\infty}.
\end{equation}
Note, however, that exactly at the threshold, \ie{}, for $\gnk$ with
$\norm{\tk}=1$, there is no giant component although
$\sum_{k=3}^\infty \gb_k(\kk)=\infty$.
Moreover, the relation \eqref{turc} may fail for unbounded $\kk$, see
the examples above. In Turova's case \eqref{tur}, the relation
\eqref{turc} holds for $\gd>0$ (when $\tk$ is Hilbert--Schmidt), but
not for $\gd=0$ (when $\tk$ is not compact), see \refE{ECT}.
\end{remark}

\section{Further remarks}\label{Sconcl}

Random graphs defined via kernels appear in various other contexts.
One natural example is the `dense' case: let $\kk$ be a symmetric
function from $[0,1]^2$ to $[0,1]$ with some suitable `smoothness'
property, and form a graph on $[n]$ by taking the probability $\pij$
of the edge $ij$ to be $\kk(x_i,x_j)$, where
$x_i$ is the type of vertex $i$ (e.g., $x_i=i/n$),
and different edges are present independently. 
Thus, when $\kk=p$ is constant,
one recovers the dense Erd\H os--R\'enyi graph $G(n,p)$.
The study of this inhomogeneous dense model is as far from the
concerns of the present
paper as the study of $G(n,1/2)$, say, is from the study of $G(n,c/n)$.

Another case, dense but not so dense, is obtained from our model if we 
omit the restriction that
$  \kk\in L^1(\sss\times\sss,\mu\times\mu)$.
One particular example that might have interesting
behaviour is $\kk(x,y)=1/|x-y|$, 
with $x_i=i/n$, say, so $\pij=1/|i-j|$, for $i\ne j$.
A similar model (in the rank 1 case of \refSS{SSrank1}) has been studied by
Norros and Reittu \cite{NorrosR}. 
Newman and Schulman~\cite{NS86} 
studied percolation in a closely related infinite random graph:
two vertices $i,j\in \bbZ$ are joined with probability
$1-\exp(-\beta |i-j|^{-s})$, where $\beta$ and $s$ are parameters of the model.

Models with intermediate density (a number of edges that is more than 
linear but less than quadratic in the number $n$ of vertices)
could be obtained by defining the edge probabilities $\pij$
in terms of a kernel $\kk$ but with different scaling to that in
\eqref{pij}.
For example, we could take
$p_{ij}=\min \{\kappa (x_i, x_j)/n^{\alpha}, 1\}$, where $0\le \alpha \le
1$ is a fixed number, or $p_{ij}=n^{-\kappa(x_i, x_j)}$, say.
Although these definitions bear a formal resemblance to the one we have used,
they lead to very different models. Nevertheless, these models may also
repay close attention. In some cases such models might correspond
to, or resemble, graphs growing in time by the addition of vertices,
with the addition of an increasing number of edges at each step:
the case $\pij=j^{-\alpha}$ for $i<j$ is one particular example.

A different connection between graphs and symmetric functions $W$ from
$[0,1]^2$ to $[0,1]$ arises in the work of Lov\'asz and Szegedy~\cite{LS},
where the limit of a sequence of dense graphs is defined by considering
the number of subgraphs isomorphic to each fixed graph.

Another natural model is the following: take the type space
as $[0,1]^2$, say, with the Lebesgue measure, and take the types
of the vertices to be independent. (Or, more naturally,
generate the vertex types by a Poisson process
of intensity $n$, so the total number of vertices is random.)
Join two vertices with a probability $p(d,n)$ that is a function
of $d$, the Euclidean distance between the (types of the) vertices,
and $n$. Since the typical distances are order $n^{-1/2}$,
the natural normalization is $p=f(dn^{1/2})$, for example,
$p(d,n)=c_1 \exp(-c_2d^2n)$.

In many ways, a model defined in this way is similar to that
considered in this paper: if $f$ decays sufficiently fast, the
expected degrees are of order $1$, and the degree distribution will be
asymptotically Poisson.  However, in other ways this graph is very
different from the ones we have been studying: in particular, it has
many small cycles.  Determining the threshold for the emergence of the
giant component in this model is likely to be as hard as finding the
critical probability for a planar percolation model (indeed, it is
essentially the same task), and is thus likely to be impossible except
perhaps in very special cases.

\medskip
Another interesting property of a graph is the behaviour of
the contact process on the graph. Suppose that each vertex is either
{\em susceptible}, or {\em infected}: infected vertices
infect their susceptible neighbours with rate $\lambda$, and recover
with rate $1$, returning to the susceptible state.
The process starts with a single randomly chosen infected vertex.
When the average degree is of order $1$, one might expect that there
is a critical value $\lambda_c$ such that for $\lambda<\lambda_c$ constant,
the expected number of vertices ever infected is $O(1)$,
while for $\lambda>\lambda_c$ constant, with probability bounded
away from zero almost all (perhaps even all) vertices
in the giant component become infected at some point,
and the infection lasts an exponentially long time. This is
the case for the 4-regular grid graph on the torus, say; see Liggett~\cite{Liggett2}.

The behaviour of the contact process on the scale-free LCD graph
has been studied by Berger, Borgs, Chayes and Saberi~\cite{BBCS},
who gave detailed results showing in particular that there is
no threshold (i.e., $\lambda_c=0$).

For $\gnkx$, one might expect a positive threshold if and only if
$\norm{\tk}<\infty$. (Perhaps an additional condition would be
needed, such as $\kk$ bounded.) In fact, one can make a
more detailed prediction based on the contact process on infinite trees:
the threshold $\lambda_c$ should be the same as the threshold
for the process to continue forever on an infinite tree
generated by the branching process $\bpkx$. Note that it is likely
that there are {\em two} distinct thresholds for the behaviour of the contact
process on such trees (this is known only for certain classes
of trees, including regular trees; see \cite{Pemantle,Liggett1,Stacey}):
a lower threshold $\lambda_1$
above which the process has positive probability of never dying out,
and an upper threshold $\lambda_2$
above which a given vertex has a positive probability of becoming reinfected
infinitely often. (In both cases, we start with a single infected vertex.)
In the graph, $\lambda_1$ should be relevant: if the process survives
with drift in the infinite tree, it will eventually revisit a given vertex of $\gnkx$,
as the neighbourhoods of a vertex are only locally treelike.

Related results have been proved by Durrett and Jung \cite{DurrettJung} 
for a $d$-dimensional version of the small-world model of Bollob\'as and
Chung: the vertex set is a discrete torus, each vertex is connected
to all vertices within a fixed distance, and then all pairs in a random
matching of the vertices are added as `long-range' edges.
Durrett and Jung prove separation of $\gl_1$ and $\gl_2$ for an infinite
version of this graph. Also, they show that for $\gl>\gl_1$,
a modified contact process on the finite graph survives for an exponential time;
the modification is to allow an infected vertex to infect a randomly
chosen other vertex, at an arbitrarily small but positive rate $\gamma$.
The result is likely to hold with $\gamma=0$, since the `long-range' edges
already provide sufficient global randomness.

If the definition of the contact process on a graph $G$ is modified so
that when a vertex recovers it cannot be reinfected, 
one might expect that, starting with a single infected vertex $v$, the
set of vertices 
eventually infected have a `similar' distribution to the component of $G[p]$
containing $v$, where $G[p]$ is formed from $G$ by keeping each edge
independently with probability $p=\gl/(1+\gl)$. Roughly speaking, 
for each edge $ww'$ of $G$,
we may declare the edge $ww'$ to be open if whichever of $w$ and $w'$ is
first infected will try to infect the other before it recovers, an event with probability $\gl/(1+\gl)$.
The set of infected vertices is the component of $v$ in the graph formed by the open edges.
Unfortunately, since the probabilities of infection from $w$ to $w'$ and from $w$ to $w''$
both depend on the random time that $w$ remains infected,
the events that different edges are open are not independent.
This fact is missed by Newman~\cite{Newman:contact}, who states that this modified
contact process is equivalent to percolation on $G$; we should like to thank an
anonymous referee for drawing this paper to our attention.
Nevertheless, it may still be true that
the threshold in this modified contact process is close to
the percolation threshold on $G$, at least under certain conditions.

\appendix
\section{Probabilistic lemmas} \label{Smeasure}

In this appendix we prove three simple technical results concerning
sequences of random variables. The first and third are
used in the main body of the paper; the second is needed to prove
the third. The first
concerns random Borel measures.

Let $\sss$ be a separable metric space,
let $M(\sss)$ be the space
of all finite (positive) Borel measures on $\sss$,
and let $P(\sss)$ be the subspace
of all Borel probability measures on $\sss$. We equip $M(\sss)$ and
$P(\sss)$ with the usual (weak) topology:
$\mu_n\to\mu$ if and only if
$\int f\dd\mu_n\to\int f\dd\mu$ for every function $f$ in the space $C_b(\sss)$
of bounded continuous functions on $\sss$.
Alternatively, as is well known, $\mu_n\to\mu$ if and only if
$\mu_n(A)\to\mu(A)$ for every $\mu$-continuity set $A$,
\ie{}, every measurable set $A$ with $\mu(\ddd A)=0$. 

\begin{remark}
The case of probability measures is perhaps better known, 
and is treated in detail in, for example, 
Billingsley  \cite{Bill}.
Many 
results extend immediately to $M(\sss)$, either by inspecting
the proof, or because $\mu_n\to\mu$ in $M(\sss)$ if and only if 
$\mu_n(\sss)\to\mu(\sss)$ and either $\mu(\sss)=0$ or
$\mu_n/\mu_n(\sss)\to\mu/\mu(\sss)$ in $P(\sss)$.
\end{remark}

The spaces $P(\sss)$ and $M(\sss)$ are themselves
separable metric spaces; for $P(\sss)$, see \cite[Appendix III]{Bill}.

The characterizations above of convergence in $P(\sss)$  and $M(\sss)$
extend to random measures and convergence in
probability as follows; see Kallenberg
\cite[Theorem 16.16]{Kall} for a similar
(but stronger) theorem under a stronger hypothesis on $\sss$.
Note that both \ref{lma} and \ref{lmc} are special cases of \ref{lmac}. 

\begin{lemma}
  \label{Lmeas}
Let $\sss$ be a separable metric space, and suppose that $\nu_n$,
$n\ge 1$, are random measures in $M(\sss)$. 
Then the following assertions are equivalent:
\begin{romenumerate}
  \item \label{lm}
$\nu_n\pto\mu$;
\item \label{lma}
$\nu_n(A)\pto\mu(A)$ for every $\mu$-continuity set $A$;
\item \label{lmc}
$\int f\dd\nu_n \pto \int f\dd\mu$ for every bounded continuous
function $f:\sss\to\bbR$.
\item \label{lmac}
$\int f\dd\nu_n \pto \int f\dd\mu$ for every bounded $\mu$-\aex{}
continuous function $f:\sss\to\bbR$. 
\end{romenumerate}
\end{lemma}

\begin{proof}
\pfitemx{\ref{lm}$\implies$\ref{lma}}
If $A$ is a \mucs\ then $\nu\mapsto\nu(A)$ defines a measurable
function $M(\sss)\to\bbR$
which is continuous at $\mu$.
\pfitemx{\ref{lma}$\implies$\ref{lmac}}
We may suppose that $0\le f\le 1$; the general case follows by linearity. 
Let $N$ be the $\mu$-null
set consisting of points at which $f$ is discontinuous,
and let
$A_t\=f\qi(t,\infty)=\set{x:f(x)>t}$. If $x\in\overline{A_t}\setminus A_t$
and $x\notin N$, then $f(x)=t$ by continuity. Thus $\partial A_t\subseteq
N\cup  f\qi\set{t}$. Hence, the sets $\partial A_t\setminus N$ are
disjoint, and 
$\mu(\partial A_t)=0$ except for at most countably many $t$.
When $\mu(\partial A_t)=0$, we have $\nu_n(A_t)\pto \mu(A_t)$ by \ref{lma}, so
$\E|\nu_n(A_t)-\mu(A_t)|\to0$ by dominated convergence.
By dominated convergence again,
\begin{equation*}
  \begin{split}
\E\,&\biggl|\int_\sss f\dd\nu_n  - \int_\sss f\dd\mu \biggr|
=
\E\,\biggl|\int_0^1\nu_n(A_t)\dd t - \int_0^1\mu(A_t)\dd t \biggr|
\\&
\le
\E\int_0^1\bigl|\nu_n(A_t) - \mu(A_t) \bigr| \dd t
=
\int_0^1\E\bigl|\nu_n(A_t) - \mu(A_t) \bigr| \dd t
\\&
\to0.
  \end{split}
\end{equation*}
\pfitemx{\ref{lmac}$\implies$\ref{lmc}}
Trivial.
\pfitemx{\ref{lmc}$\implies$\ref{lm}}
The topological space $M(\sss)$ is metrizable, but the topology is
also defined by
the functionals $\mu\mapsto\int f\dd\mu$, $f\in C_b(\sss)$.
Hence, if $U$ is a neighbourhood of $\mu$ in $M(\sss)$, there is a
finite set of
functions $f_1,\dots,f_N\in C_b(\sss)$ and $\eps>0$ such that if
$\bigl|\int f_i\dd\nu - \int f_i\dd\mu \bigr| <\eps$, $i=1,\dots,N$,
then $\nu\in U$.
Consequently,
\begin{equation*}
  \P(\nu_n\notin U) \le\sum_{i=1}^N
\P\biggpar{\biggl|\int f_i\dd\nu - \int f_i\dd\mu \biggr| \ge\eps}
\to0.
\end{equation*}
\vskip-\baselineskip
\end{proof}

\begin{remark}
  \label{RR}
To verify condition \ref{lma}, it often suffices to consider $A$
in a suitably selected family
of subsets. For example, it is well-known that on $\bbR$, it suffices
to consider \mucs{s} of the form $(-\infty,x]$ 
and, for $M(\sss)$, $\bbR$ itself; see
\cite[Section 3]{Bill}.
\end{remark}

\medskip
Recall that if $X_n$ is a sequence of random variables and $a_n$ a sequence 
of positive real numbers, then $X_n=O(a_n)$ \whp{} means that there is
a constant 
$C$ such that $|X_n|\le Ca_n$ \whp. 
Our final technical result (Lemma~\ref{Loc2} below) is simple, but perhaps
a little surprising: we shall show that under suitable assumptions, if
$X_n=O(a_n)$ 
holds conditionally (after conditioning on the sequences $\xs$ in our
model), then it holds 
unconditionally, i.e., that the implicit constant may be assumed to be
deterministic. 
We start with a preparatory lemma.

\begin{lemma}
  \label{Loc}
Let $\cA_1,\cA_2,\ldots$, be non-empty families of random variables such 
that for any sequence $X_n \in \cA_n$ we have $X_n = O(a_n)$ \whp. Then
there is a 
constant $C$ such that $\sup_{X\in \cA_n} \P(|X|>Ca_n) \to 0$ as
$n\to\infty$. In other 
words, the 
implicit constant in $X_n=O(a_n)$ \whp{}
can be chosen uniformly for $X_n\in \cA_n$.
\end{lemma}

\begin{proof}
Replacing $X_n$ by $X_n/a_n$, we may assume that $a_n=1$.
Suppose the conclusion fails. 
Then, for every $m$ there is an $\eps_m>0$ such that there are
arbitrarily large $n$ for which there is an $X_n \in \cA_n$ with 
$\P(|X_n|>m) >\eps_m$. Let $(m_k)$ be a sequence of integers where
each positive integer appears infinitely many times. Select
inductively an increasing sequence $(n_k)$ and $X_{n_k} \in \cA_{n_k}$
such that $\P(|X_{n_k}|>m_k) > \eps_{m_k}$. For $n\notin \set{n_k}$,
choose $X_n$ from $\cA_n$ arbitrarily.

For any positive integer $m$, there are infinitely many $k$ such that
$m_k=m$, and thus infinitely many $n$ 
such that $\P(|X_n|>m) > \eps_m$. Hence $(X_n)$ is not $O(1)$ \whp,
which contradicts our assumption.
\end{proof}

The next lemma can be stated in terms of families
$\mu_n(y)$ of probability distributions (on $\bbR$) and mixtures 
$\E \mu_n(Y_n)$
of them, but we prefer a statement in terms of random variables
$X_n(y)\sim\mu_n(y)$; we consider a sequence of families $X_n(y)$ of random
variables defined for $y$ in a subset $\cM_n$ of a certain space $\cM$
as this is convenient when we apply the lemma to $\gnkx$.

\begin{lemma}  \label{Loc2}
Let $\cM$ be a metric space, and, for each $n\ge1$, let $X_n(y)$, $y\in \cM_n\subseteq \cM$,
be a (measurable) family of real-valued random variables.
Let $y_0\in\cM$, and suppose that for every sequence $(y_n)$ with $y_n\in\cM_n$ and $y_n\to y_0$
we have $X_n(y_n)=O(a_n)$ \whp. 
Then, if $(Y_n)$ is a sequence of $\cM_n$-valued random variables, independent of
all $X_n(y)$,
with $Y_n \pto y_0$,
we have $X_n(Y_n) = O(a_n)$ \whp. 
\end{lemma}

\begin{proof}
Since $Y_n \pto y_0$, there is a sequence $\delta_n\to 0$
such that $\P(d(Y_n,y_0)<\delta_n)\to 1$.
Set $U_n = \set{y\in \cM_n: d(y,y_0) < \delta_n}$, so
$\P( Y_n \in U_n) \to 1$.
Note that $y_n\in U_n$ implies $y_n\to y_0$, and thus $X_n(y_n)=O(a_n)$ \whp.
Let $\cA_n = \set{X_n(y) : y \in U_n}$. By \refL{Loc}, there exists $C$
such that $\eps_n:= \sup_{y\in U_n} \P(|X_n(y)|>Ca_n) \to 0$.
Finally, $\P(|X_n(Y_n)|>Ca_n) \le \eps_n + \P(Y_n \notin U_n)\to0$.
\end{proof}

\medskip
\begin{ack}
Part of this research was done during visits of S.J.\ to Cambridge,
supported by the Swedish Royal Academy of Sciences and the London
Mathematical Society. The paper was revised during a visit of
all three authors to the Institute for Mathematical Sciences, 
National University of Singapore.
\end{ack}

\newcommand\AAP{\emph{Adv. Appl. Probab.} }
\newcommand\JAP{\emph{J. Appl. Probab.} }
\newcommand\AMS{Amer. Math. Soc.}
\newcommand\JAMS{\emph{J. \AMS} }
\newcommand\MAMS{\emph{Memoirs \AMS} }
\newcommand\PAMS{\emph{Proc. \AMS} }
\newcommand\TAMS{\emph{Trans. \AMS} }
\newcommand\AnnMS{\emph{Ann. Math. Statist.} }
\newcommand\AnnPr{\emph{Ann. Probab.} }
\newcommand\AnnSt{\emph{Ann. Statist.} }
\newcommand\AnnAP{\emph{Ann. Appl. Probab.} }
\newcommand\CPC{\emph{Combin. Probab. Comput.} }
\newcommand\JMAA{\emph{J. Math. Anal. Appl.} }
\newcommand\JASA{\emph{J. Amer. Statist. Assoc.} }
\newcommand\RSA{\emph{Random Struct. Alg.} }
\newcommand\ZW{\emph{Z. Wahrsch. Verw. Gebiete} }
\newcommand\LNCS[2]{Lecture Notes Comp. Sci. #1, Springer, #2}
\newcommand\Springer{Springer}
\newcommand\Wiley{John Wiley \& Sons}
\newcommand\CUP{Cambridge Univ. Press}
\newcommand\Poznan{Pozna\'n}
\newcommand\jour{\emph}
\newcommand\book{\emph}
\newcommand\inbook{\emph}
\newcommand\pp{\relax}
\newcommand\vol{\textbf}
\newcommand\toappear{\unskip, to appear}
\newcommand\webcitesvante{\webcite{http://www.math.uu.se/\~{}svante/papers}}

\end{document}